\documentclass[10pt]{amsart}
\usepackage[margin=1.5in]{geometry}
\usepackage{xcolor}
\usepackage{xcolor}
\usepackage{ulem}
\newtheorem{theorem}{Theorem}[section]
\newtheorem{lemma}[theorem]{Lemma}
\newtheorem{proposition}[theorem]{Proposition}
\newtheorem{corollary}[theorem]{Corollary}
\theoremstyle{definition}
\newtheorem{definition}[theorem]{Definition}

\let\hide\iffalse

\usepackage{graphicx}
\usepackage{amssymb,amsmath,amsthm,}
\usepackage{a4wide}
\usepackage{amscd}
\usepackage{amsfonts}
\usepackage{amssymb}
\usepackage{latexsym}
\usepackage{color}
\usepackage{esint}
\usepackage{graphicx}
\usepackage{caption}
\usepackage{subcaption}
\usepackage{amsthm}
\usepackage[colorlinks]{hyperref}
\usepackage{ulem}
\usepackage{kotex}
\usepackage[makeroom]{cancel}
\usepackage{mathtools}

\newtheorem{remark}[theorem]{Remark}
\let\a=\alpha

\let\p=\partial

\let\g=\gamma

\let\b=\beta
\let\k=\kappa

\newcommand{\avv}{\langle v\rangle}
\newcommand{\avvv}{\langle v_*\rangle}

\newcommand{\R}{\mathbb{R}}

\renewcommand{\S}{\mathbb{S}}
\newcommand{\T}{\mathbb{T}}

\newcommand{\be}{\begin{equation}}
\newcommand{\bm}{\begin{multline}}
\newcommand{\ee}{\end{equation}}

\newcommand{\half}{\frac{1}{2}}

\newcommand{\Bes}{\begin{eqnarray*}}
	\newcommand{\Ees}{\end{eqnarray*}}
\newcommand{\Be}{\begin{equation}}
\newcommand{\Ee}{\end{equation}}

\def\p{\partial}

\def\R{\mathbb{R}}

\def\B{\begin{equation}}
\def\E{\end{equation}}
\def\BN{\begin{eqnarray*}}
\def\EN{\end{eqnarray*}}

\usepackage{color}

\theoremstyle{remark}

\numberwithin{equation}{section}

\allowdisplaybreaks

\begin{document}
\title[Large-Data Boltzmann Dynamics]{A Kinetic Concentration Criterion for Large-Data Boltzmann Dynamics}

\author{Lingbing He}
\address{Department of Mathematical Sciences, Tsinghua University, Beijing 100084, P.R. China}
\email{hlb@tsinghua.edu.cn}

\author{Jin-cheng Jiang}
\address{Department of Mathematics, National Tsing Hua University, Hsinchu, Taiwan 300040, R.O.C.
\\, and  National Center for Theoretical Sciences, National Taiwan University, Taipei, Taiwan 100046, R.O.C.}
\email{jcjiang@math.nthu.edu.tw}

\author{Jong-in Kim}
\address{Department of Mathematics, Pohang University of Science and Technology, Pohang, South Korea}
\email{kimjim@postech.ac.kr}

\author{Donghyun Lee}
\address{Department of Mathematics, Pohang University of Science and Technology, Pohang, South Korea}
\email{donlgee@postech.ac.kr}

\author{Sungbin Park}
\address{Basic Science Research Institute(BSRI), Pohang University of Science and Technology, Pohang, South Korea}
\email{parksb2942@postech.ac.kr}




\begin{abstract}
We establish a conditional theory of genuinely non-perturbative global dynamics for the cutoff Boltzmann equation with hard potentials on the three-dimensional torus, in weighted \(L^\infty\) spaces. In the stretched-exponential regime, we show that any global mild solution satisfying
\[
\sup_{t\ge0}\int_{\mathbb R^3} W(v)\sup_{x\in\mathbb T^3}F(t,x,v)\,dv<\infty,
\]
for a suitable stretched-exponential weight \(W\), necessarily enjoys uniform weighted \(L^\infty\) control.
Consequently, the solution is unique in the relevant dissipative class, relaxes qualitatively in relative entropy, and converges exponentially fast to the global Maxwellian equilibrium. In particular, any arbitrary-amplitude global mild solution admitting a uniform stretched-exponential velocity upper bound is unique in the dissipative class and converges exponentially fast to the global Maxwellian.

The key new mechanism is a \textit{tail-to-coercivity principle} : the weighted $L^1_vL^\infty_x$ bound forces quantitative positive lower bounds for the local density and, consequently, for the nonlinear collision frequency. Thus macroscopic coercivity is generated from an upper-tail condition rather than imposed as an independent a priori assumption. A second ingredient is a short-time positivity mechanism based on a double-Duhamel expansion, which generates a quantitative Gaussian lower bound immediately after the initial time. This removes the need for any pointwise lower bound on the initial distribution and, in particular, allows initial data with local vacuum. These principles reveal a rigidity mechanism by which a single additional macroscopic-type control determines the global behavior of arbitrarily large dissipative solutions, yielding weighted \(L^\infty\) bounds, uniqueness, and exponential relaxation.

We further develop a polynomially weighted counterpart. Assuming the initial datum admits a spatially uniform lower velocity profile, we show that the corresponding polynomially weighted velocity-integrability condition on the spatial supremum leads to the analogous weighted \(L^\infty\) control, uniqueness, and exponential convergence to equilibrium.

This framework provides a low-regularity counterpart to conditional global theories developed for smoother Boltzmann solutions. Rather than assuming a collection of pointwise macroscopic bounds or high regularity, the long-time dynamics are governed by a single weighted velocity-tail control. The mechanism is formulated within a weighted \(L^\infty\) framework and suggests a possible route toward conditional large-data dynamics for Boltzmann equations with physical boundary conditions.

\end{abstract}

\maketitle
\tableofcontents

\section{Introduction}
We consider the nonlinear Boltzmann equation which is a classical mathematical model for rarefied collisional gas on the three-dimensional torus $\mathbb{T}^3_x$. The non-negative distribution function $F(t,x,v)$ solves
\Be \label{BE}
	\p_t F + v\cdot\nabla_x F = Q(F,F)
\Ee
where the quadratic collision operator $Q(F,F)$ is defined by   
\Be \notag
Q(F,F) = \int_{\R^3}\int_{\mathbb{S}^2} B(v-v_*,\sigma)[F(t,x,v'_*)F(t,x,v') - F(t,x,v_*)F(t,x,v)] d\sigma dv_*. 
\Ee
We consider hard potentials with an angular cutoff, and the collision kernel $B$ is given by 
\Be \notag
	B(v-v_*,\sigma) = |v-v_*|^\g b(\cos\theta),\quad \cos\theta = \frac{(v-v_*)\cdot\sigma}{|v-v_*|},\quad 0 < \gamma \leq 1 
\Ee
For cross section $b(\cos\theta)$, we assume
\Be \label{kernel}
	c_b \leq b(\cos\theta) \leq C_b
\Ee
for some uniform positive constants $c_b$ and $C_b$. The post-collisional velocities $v'$ and $v'_*$ are given by 
\begin{align*}
     v' = \frac{v+v_*}{2}+\frac{|v-v_*|}{2}\sigma, \quad v'_* = \frac{v+v_*}{2}-\frac{|v-v_*|}{2}\sigma,
\end{align*} 
and they satisfy the following conservations :
\Be \label{loc conserv}
	v + v_* = v' + v'_* \quad\text{and} \quad |v|^2 +|v_*|^2 = |v'|^2 + |v'_*|^2,
\Ee

The above micro-level conservations give collision invariance relation
\Be \label{invariance}
	\int_{\R^3} Q(F,F) 
	\phi(v)
	dv = 0
	\quad \text{if and only if}\quad
	\phi(v) \in \text{span}\{1, v, |v|^2\}
\Ee
formally, and we obtain (formal) macroscopic conservations : 
\Be \label{macro_conserv}
	\iint_{\mathbb{T}^3\times\R^3} F(t,x,v)\begin{pmatrix}
		1 \\ v \\ |v|^2
	\end{pmatrix} dv dx
	=
	\iint_{\mathbb{T}^3\times\R^3} F_0(x,v)\begin{pmatrix}
		1 \\ v \\ |v|^2
	\end{pmatrix} dv dx
\Ee  
where $F_0(x,v)$ is initial distribution. 
	
A fundamental property of the Boltzmann equation is the monotonicity of its entropy functional. The Boltzmann entropy functional $H(F)$ is defined by 
\Be \notag
	H(F) := \iint_{\mathbb{T}^3\times\R^3} F\log F dvdx.
\Ee 
When $F$ is sufficiently smooth and decays sufficiently rapidly in high velocity, it is well-known that the entropy functional is non-increasing in time :
\Be \label{H-theorem}
	\frac{d}{dt}H(F)(t) + \mathcal{D}(F)(t) = 0.
\Ee
where entropy dissipation functional $\mathcal{D}(F)$ is defined by 
\begin{equation*}
	\begin{split}
	\mathcal{D}(F)(t) &:= \int_{\mathbb{T}^3}{D}(F)(t) dx 
	 := -\iint_{\mathbb{T}^3\times\R^3} Q(F,F)\log F \\
	 &= \frac{1}{4}\iiiint_{\mathbb{T}^3\times\R^3\times\R^3\times \mathbb{S}^2} B(v-v_*,\sigma)(F'_* F' - F_*F)\log\frac{F'_* F'}{F_*F} \geq 0.
	\end{split}
\end{equation*}

This is called the H-theorem and this is why we expect a Boltzmann solution converges to a state with the lowest entropy. For more explanation, let us suppose a global (smooth) solution $F(t,x,v)$ of \eqref{BE}. We recall the standard definition of macroscopic quantities $(\rho, U, T)$ :
\begin{equation} \label{def macro}
	\begin{split}
		\rho(t,x) &:= \int_{\R^3} F(t,x,v) dv,\\
		\rho(t,x)U(t,x) &:= \int_{\R^3} vF(t,x,v) dv, \\
		3\rho(t,x)T(t,x) &:= \int_{\R^3}|v-U(t,x)|^2 F(t,x,v) dv.		
	\end{split}
\end{equation}

Let us consider a Gaussian function whose macroscopic quantities are same as $F$. We use $M_{F}(t,x,v)$ to denote the Gaussian function and call the local Maxwellian of $F$:
\Be \notag
	M_{F}(t,x,v) := \frac{\rho(t,x)}{\sqrt{2\pi T(t,x)}^3} e^{-\frac{|v-U(t,x)|^2}{2T(t,x)}}
\Ee

To explain entropy minimizing property, we first define relative entropy 
\Be \label{def rel entropy}
	H(F  \ \vert \ G) = \iint_{\mathbb{T}^3\times\R^3} \Big( F\log\frac{F}{G} - (F - G)\Big) dvdx.
\Ee
Note that among all probability density functions with same macroscopic quantities $(\rho(t,x), U(t,x), T(t,x))$, $M_{F}$ minimizes entropy. This is obtained by 
\[
	H(F \ \vert \ M_{F}) \geq 0\quad\text{and}\quad H(F \ \vert \ M_{F}) = H(F) - H(M_F).
\]
See Lemma \ref{lem_min entropy} for detail. From \eqref{H-theorem}, it is natural to expect that a global solution $F$ converges to a Maxwellian $M_F$, and then from local conservation \eqref{loc conserv}, we obtain
\Be
	Q(M_F, M_F) = 0.
\Ee
Thus we have
\[
	\p_t M_F + v\cdot \nabla_x M_F = 0,
\]
and a local Maxwellian which solves the free transport equation should be a uniform Maxwellian (See \cite{Desvillettes1990, Desvillettes}). Combining with \eqref{macro_conserv}, the expected global equilibrium is 
\[
	\mu(v) := \frac{1}{\sqrt{2\pi}^3}e^{-\frac{|v|^2}{2}}
\]
where we normalize the initial datum so that its total mass, momentum, and energy agree with those of $\mu$. \\

In sharp contrast with the spatially homogeneous theory, the global large-data dynamics of the spatially inhomogeneous Boltzmann equation remain far from being understood. More precisely, although global renormalized solutions for arbitrary physical large data are available (\cite{DL}), global well-posedness in a strong or mild class, uniqueness, and unconditional convergence to equilibrium for general large data are still among the central open problems in the theory. The spatially homogeneous equation, in which the transport mechanism and the associated spatial concentration effects are absent, has a substantially more complete theory. Its mathematical study goes back to the pioneering work of Carleman \cite{Carleman1933Acta.Math.}, who established global existence and uniqueness for hard spheres, initially under symmetry assumptions, and also investigated the large-time approach to equilibrium. This theory was subsequently developed substantially by Arkeryd \cite{Arkeryd1972,Arkeryd1972a} who constructed global \(L^1\)-type solutions for hard potentials and initiated a systematic large-data theory far from equilibrium. Further developments by Povzner \cite{Povzner1962}, Elmroth \cite{Elmroth1983}, Wennberg \cite{Wennberg1994}, Desvillettes \cite{Desvillettes1993}, and others provided refined moment estimates, regularity properties, and asymptotic information. A particularly definitive well-posedness result was obtained by Mischler and Wennberg \cite{Mischer1999AIHP}, who proved that, for cutoff hard potentials, every nonnegative initial datum with finite mass and energy generates a unique energy-conserving solution. The long-time theory was developed in parallel: compactness and entropy methods yielded qualitative \(L^1\) convergence to the Maxwellian, while Arkeryd \cite{Arkeryd1988ARMA} obtained nonconstructive exponential convergence in \(L^1\), later extended to broader \(L^p\) settings \cite{Wennberg1993}. For hard potentials with Grad's cutoff, Mouhot \cite{Mouhot2006CMP} subsequently established a quantitative exponential relaxation theory, connecting the large-data nonlinear dynamics to the spectral gap of the linearized collision operator. 

The situation changes drastically once the transport term \(v\cdot\nabla_x F\) is present. Near equilibrium, global classical and mild solutions have been constructed through perturbative theories originating in the works of Ukai \cite{Ukai1974} and subsequent developments by Nishida \cite{Nishida1976}, Guo \cite{Guo2003}, and many others \cite{Shizuta1983,Ukai2006,Ukai1982,Liu2004b}; a different perturbative regime near vacuum was treated by Illner and Shinbrot \cite{Illner1984CMP}, see also~\cite{He2023,He2024} and reference therein. For arbitrary large data, the landmark work of DiPerna and Lions \cite{DL} established global renormalized solutions under the natural finite mass, energy, and entropy assumptions, but this theory does not provide uniqueness or the strong pointwise control required to determine the global dynamics. The convergence problem is equally delicate. Desvillettes and Villani \cite{Desvillettes} proved a powerful conditional theory of convergence to global equilibrium for spatially inhomogeneous Boltzmann solutions, assuming strong a priori bounds on smoothness, velocity decay, and positivity, while more recent work of Imbert and Silvestre \cite{Imbert2021} developed a conditional regularity theory in the non-cutoff setting based on pointwise control of macroscopic quantities. Thus, between the perturbative theories and the large-data renormalized framework lies a substantial gap: for genuinely large spatially inhomogeneous solutions, it remains largely open to derive from weak or macroscopic information the coercivity and weighted pointwise bounds needed for uniqueness and quantitative relaxation to equilibrium. The present work is aimed precisely at this intermediate regime. 
 
 The work of Desvillettes and Villani \cite{Desvillettes} shows, in particular, that sufficiently smooth solutions satisfying their a priori assumptions approach the global Maxwellian at an almost exponential rate. This was complemented by the perturbative theory of \cite{Gualdani2017}, developed in weighted \(L^1_vL^\infty_x\)-type spaces, which upgrades the eventual relaxation to exponential convergence once the solution enters a sufficiently small neighborhood of equilibrium. Such spaces are natural for the cutoff Boltzmann equation: the collision operator involves a single integration in the velocity variable, while the free transport preserves the spatial \(L^\infty\) structure. This is also consistent with the classical spatially homogeneous theory, whose natural large-data framework is based on \(L^1_v\) control.
 
 There is, however, an important limitation of this viewpoint in the low-regularity regime considered below. The monotonicity of the Boltzmann entropy provides a natural mechanism for measuring proximity to equilibrium, but small relative entropy is fundamentally an integral statement in phase space. When sufficient spatial regularity is available, this information can be combined with regularity and velocity decay to rule out concentration and thereby yield smallness in a stronger norm such as a weighted \(L^1_vL^\infty_x\) norm. Without such regularity, this implication fails in general: a solution may have very small entropy while retaining narrow spatial spikes of large amplitude, which are invisible to a global integral quantity. To exclude this scenario one would need, at a minimum, some uniform spatial non-concentration mechanism, for instance a global modulus of continuity. This is one reason why a direct \(L^1_vL^\infty_x\) perturbative approach is less suited to the rough mild-solution framework pursued here.

 The situation is particularly compelling for the Boltzmann equation in bounded domains with physical boundary conditions, where sufficiently smooth solutions cannot in general be expected, even for cutoff collision kernels. The interaction of kinetic transport with grazing trajectories and repeated reflections at the boundary produces singular behavior that is absent in the periodic problem and is strongly sensitive to the geometry of the domain. For inflow, diffuse, and bounce-back boundary conditions, Kim \cite{Kim2011} showed that discontinuities can be created at non-convex grazing points and subsequently propagate along characteristics, while Guo, Kim, Tonon, and Trescases \cite{Guo2015} established the corresponding optimal \(BV\) regularity in general non-convex domains. In strictly convex domains, Guo, Kim, Tonon, and Trescases \cite{Guo2017Invent.Math.} developed a weighted \(C^1\) and \(W^{1,p}\) regularity theory away from the grazing set and, at the same time, exhibited the obstruction to second-order regularity up to the boundary. For specular reflection, where the geometric structure is different, Kim and Lee \cite{Lee2024CPAM}  proved \(C^{0,\frac12-}\) Hölder regularity past a convex obstacle; this line of work was further developed by An and Lee \cite{An2026Ann.PDE.} through the study of mixed convex-concave geometries and, more recently, the optimal \(C^{0,\frac12}\) regularity in non-convex domains. These results illustrate that, in the presence of physical boundaries, a low-regularity framework is not merely convenient but is often intrinsic to the dynamics.

These regularity results naturally connect with another line of research that develops global theories directly in weighted \(L^\infty_{x,v}\) spaces. A fundamental step was Guo's work \cite{Guo2010ARMA} on the Boltzmann equation in bounded domains, which established global mild solutions near a global Maxwellian and exponential \(L^\infty\) decay for several basic physical boundary conditions, including inflow, bounce-back, specular, and diffuse reflection. Since then, weighted \(L^\infty\) methods, based on characteristic representations together with \(L^2\)-\(L^\infty\) estimates and repeated Duhamel expansions, have become a principal tool for studying global dynamics of Boltzmann boundary problems in the perturbative regime. 

More recently, this framework has been pushed beyond the small-\(L^\infty\) regime. Duan, Huang, Wang, and Zhang \cite{Duan2019} developed a global stability theory for cutoff soft potentials with non-isothermal diffuse reflection, including a class of large perturbations, and obtained long-time decay in weighted \(L^\infty\). For hard potentials with diffuse reflection, Duan and Wang \cite{Duan2019b} constructed global solutions whose weighted \(L^\infty\) amplitude may be arbitrarily large, provided that the initial perturbation remains sufficiently small in an \(L^2\)-type sense; such solutions eventually enter the perturbative regime and converge exponentially to equilibrium. Analogous large-amplitude dynamics for specular reflection were subsequently obtained by Duan, Ko, and Lee \cite{Ko2023JSP} under a small relative-entropy condition. These works reveal a characteristic mechanism of the boundary theory: large pointwise amplitude can be compatible with global dynamics once it is coupled to a weaker notion of proximity to equilibrium. 

These developments suggest that weighted \(L^\infty\) spaces are a natural setting in which to seek a theory beyond the perturbative regime, especially with future applications to physical boundary problems in mind. The challenge is to bridge the gap between the large-data existence theory, which offers only weak compactness and entropy information, and the perturbative \(L^\infty\) theory, which relies on strong coercivity and pointwise control. 

It is this gap that the present work addresses. We seek a conditional theory of global uniqueness and relaxation to equilibrium in a low-regularity weighted \(L^\infty\) framework, rather than in a smooth solution class. A fundamental obstruction in the genuinely non-perturbative mild formulation is the lack of an a priori lower bound for the nonlinear collision frequency. 
Without such a bound, the loss term does not provide the coercivity needed to close global weighted \(L^\infty\) estimates, and the standard perturbative mechanisms are no longer directly available.



%

To overcome this difficulty, we introduce what we call the tail-to-coercivity principle. Roughly speaking, the uniform upper-tail condition
\Be \label{rough M}
    \sup_{t\ge0}\int_{\mathbb R^3} W(v)\sup_{x\in\mathbb T^3}F(t,x,v)\,dv<\infty,
\Ee
is converted into quantitative lower bounds for the local density and, consequently, for the nonlinear collision frequency after a Gaussian lower bound is generated. The constants in this compactness step are solution-dependent and non-constructive. Thus the macroscopic coercivity required by the weighted \(L^\infty\) theory is generated from upper-tail information rather than imposed as an independent hypothesis.	
This gives a different conditional framework: the global mild solution is not assumed smooth, and the coercive macroscopic bounds are deduced from \eqref{rough M} under the initial hypotheses stated below. The criterion is a single weighted $L^1_vL^\infty_x$ bound; placing the spatial supremum inside the velocity integral makes it stronger than a bound on the local density or local velocity moments alone.

In particular, neither a positive lower bound on the local density nor a pointwise upper bound on the local entropy density is imposed independently. The tail-to-coercivity principle converts the rough velocity-tail information contained in \eqref{rough M} into quantitative lower bounds for the local density and the nonlinear collision frequency, thereby supplying the coercivity needed for the weighted \(L^\infty\) theory. Thus the present framework replaces the combination of smoothness and several pointwise macroscopic bounds by a single rough upper-tail condition on a mild solution, and from this derives conditional weighted \(L^\infty\) control, uniqueness, and relaxation to equilibrium in a genuinely non-perturbative regime. 

In the polynomially weighted setting, a spatially uniform initial lower profile is still required to initiate the argument. However, if a global mild solution is additionally assumed to belong to a polynomially weighted $L^\infty$ space at each time, this initial lower profile assumption can be removed. Moreover, it suffices to take $W(v)=\avv^2$ to establish uniform-in-time weighted $L^\infty$ bounds and convergence to equilibrium.


Overall, our results suggest a different organizing principle for the cutoff Boltzmann dynamics: the decisive structural input is not smoothness itself, but uniform control of the velocity tails. Once \eqref{rough M} is available, weighted \(L^\infty\) bounds, macroscopic control, non-concentration, uniqueness, and relaxation to equilibrium all emerge from the same hypothesis. 

More precisely, once a global mild solution satisfies the tail condition \eqref{rough M}, we prove directly that the solution remains uniformly bounded in the relevant weighted \(L^\infty\) space. Along the way, the same tail control generates the macroscopic information that is usually imposed separately: quantitative lower bounds for the density and collision frequency, upper control of the macroscopic fields, and the non-concentration properties needed for the long-time analysis. In this sense, \eqref{rough M} acts as a single structural condition that organizes the global dynamics, replacing a collection of smoothness and macroscopic assumptions by one rough velocity-tail hypothesis.

A further feature of the stretched-exponential theory is that no pointwise lower bound on the initial datum is required; in particular, the initial distribution is allowed to contain local vacuum. The missing positivity is generated dynamically. Combining the Lebesgue density theorem with a double-Duhamel expansion, we show that, after an arbitrarily short positive time, the solution acquires a uniform Gaussian lower bound in velocity. This short-time positivity mechanism is a crucial input into the tail-to-coercivity principle: once a quantitative amount of mass has been created uniformly in space and velocity, the assumed upper-tail control prevents that mass from escaping to large velocities and yields a uniform positive lower bound for the local density and hence for the nonlinear collision frequency. Thus, in the stretched-exponential setting, the coercivity needed for the global weighted \(L^\infty\) theory is generated dynamically even from initial data with vacuum. \\

\subsection{Main results}

We first recall the two notions of solutions that will be used throughout the paper. We solve the cutoff Boltzmann equation
$$
\partial_t F+v\cdot\nabla_x F=Q(F,F),
\qquad (t,x,v)\in (0,\infty)\times \mathbb T^3\times\mathbb R^3,
$$
with initial datum \(F_0\). We write the collision operator in the gain--loss form
$$
Q(F,F)=Q_+(F,F)-L(F)F,
$$
where $L(F)$ is defined
\Be \label{def LF}
L(F)(t,x,v)
:=
\int_{\mathbb R^3}\int_{\mathbb S^2}
B(v-v_*,\sigma)F(t,x,v_*)\,d\sigma\,dv_* .
\Ee

The first notion is that of a renormalized solution in the sense of DiPerna--Lions \cite{DL}. When $G = \b(F) := \log(1+F)$, a nonnegative distribution \(F\) is called a renormalized solution if 
\[
	\frac{1}{1+F}Q_{\pm}(F,F) \in L^1_{loc}((0,\infty)\times \mathbb{T}^3 \times \mathbb{R}^3)
\]
and
\[
	\p_t G + v\cdot\nabla_x G = \frac{1}{1+F}Q(F,F)
\]
holds in the sense of distributions. Actually, if $F$ is a renormalized solution, then for every admissible renormalization
$$
\beta\in C^1([0,\infty)),
\qquad
|\beta'(z)|\lesssim \frac{1}{1+z},
$$
the renormalized equation
$$
\partial_t\beta(F)+v\cdot\nabla_x\beta(F)
=
\beta'(F)Q(F,F)
$$
holds in the sense of distributions. This formulation is designed to make sense under the natural finite mass, energy, and entropy bounds, even when \(Q(F,F)\) itself lacks sufficient integrability to be interpreted directly as a distribution. The global solutions constructed by the DiPerna--Lions approximation scheme \cite{DL} belong to this class and satisfy the corresponding entropy and energy inequalities. 

\begin{definition} 
	We say $F \in \mathcal{D}_{DL}(F_0)$ if $F$ is a renormalized solution of \eqref{BE} with initial distribution $F_0$ which is obtained by approximation scheme as in \cite{DL}. 
\end{definition}

In particular, the renormalized framework provides a global large-data existence theory without requiring pointwise control or high regularity of the solution.

The second notion, which is more closely adapted to the weighted \(L^\infty\) analysis of the present paper, is that of a mild solution. 

\begin{definition} 
A nonnegative function \(F\) is called a mild solution on \([0,T]\) if it satisfies the Boltzmann equation along the free-transport characteristics. More precisely, 
\Be \label{def mild}
\begin{aligned}
	F(t,x,v)
	={}&
	F_0(x-tv,v)
	+
	\int_0^t Q(F,F)\bigl(s,x-(t-s)v,v\bigr)\,ds .
\end{aligned}
\Ee
If $L(F) \in L^1((0,T)\times B_{R} \times B_R)$ for all $R, T < \infty$, $F$ satisfies exponential multiplier form
\Be \label{def mild exp}
\begin{aligned}
	F(t,x,v)
	={}&
	F_0(x-tv,v)
	\exp\!\left[
	-\int_0^t
	L(F)\bigl(s,x-(t-s)v,v\bigr)\,ds
	\right]
	\\
	&+
	\int_0^t
	\exp\!\left[
	-\int_s^t
	L(F)\bigl(\tau,x-(t-\tau)v,v\bigr)\,d\tau
	\right]
	\\
	&\hspace{2.5cm}\times
	Q_+(F,F)\bigl(s,x-(t-s)v,v\bigr)\,ds .
\end{aligned}
\Ee
In this paper, the mild formulations \eqref{def mild} and \eqref{def mild exp} are equivalent, and we use them interchangeably, since our a priori assumption on $F$ satisfies $L(F) \in L^1((0,T)\times B_{R} \times B_R)$ for all $R, T < \infty$.
\end{definition}

The distinction between these two solution concepts is important for our purposes. Renormalized solutions are available globally for arbitrary large physical data, but in general do not come with the pointwise information required by a weighted \(L^\infty\) theory. Mild solutions, on the other hand, expose directly the competition between the gain term and the nonlinear damping generated by \(L(F)\). The main results below identify conditions under which this mild formulation becomes globally coercive, unique, and asymptotically stable, while remaining compatible with the large-data renormalized framework. 

\begin{definition} 
\label{def_dissipative soln}
	For \(T\in(0,\infty]\), let \(\mathcal{U}_T(F_0)\) denote the class of all nonnegative mild solutions $G$ emanating from \(F_0\), such that
	\[
	G \in L^\infty_{\mathrm{loc}}
	([0,T);L^1_{x,v}(\langle v\rangle^2))
	\]
	and
	\[
	\iint\langle v\rangle^2 G(t,x,v)\,dx\,dv
	\le
	\iint\langle v\rangle^2 F_0(x,v)\,dx\,dv
	\]
	for almost every \(t\in[0,T)\). Note that we do not impose entropy condition on \(G\). 
\end{definition}

\begin{definition} 
	We say $G$ is a physically admissible solution if $G$ satisfies mass, momentum, and energy conservation and entropy inequality
	\Be \label{entropy ineq}
		H(G)(T) - H(G_0)+ \int_{0}^{T} \mathcal{D}(G)(t)dt \leq 0,\quad \forall T\geq 0.
	\Ee
	
\end{definition}

For stretched exponential theory, we impose a priori assumption
\Be \label{apriori}
	\sup_{t\geq 0}\int_{v}\sup_{x}F(t,x,v) e^{{\alpha}|v|^{{\kappa}}} \avv^k dv \leq A_0 < \infty,\quad {\alpha}>0, \quad k\geq 1,\quad 1< {\kappa} < 2 .
\Ee

\begin{remark} 
	Consider the spatially homogeneous Boltzmann equation. If $F_0$ is bounded by a Gaussian $e^{-a_0|v|^2+c_0}$, then $F \leq e^{-a|v|^2+c}$ for some $a>0$ and $c$ by \cite{GPV}. Therefore,  \eqref{apriori} holds whenever $F_0$ has a Gaussian upper bound. When $\kappa=2$, however, $\alpha$ in \eqref{apriori} should be very small so that $\alpha  < a$, because $a$ is very small in general and hard to specify.
\end{remark}


We consider initial data \(F_0\in L^\infty(\avv^\b e^{\a |v|^\k})\), with \(\mu\) denoting the associated Maxwellian, i.e.,
\Be \label{initial}
	\|\avv^\b e^{\a |v|^\k} F_0\|_{L^\infty_{x,v}} < \infty ,\quad  \b > 2,  \ \ \a > 0, \ \ 1 < \k < 2
\Ee
and 
\Be \label{initial normal}
\iint_{\mathbb{T}^3\times\R^3} F_0(x,v)\begin{pmatrix}
	1 \\ v \\ |v|^2
\end{pmatrix} dv dx
=
\iint_{\mathbb{T}^3\times\R^3} \mu(v)\begin{pmatrix}
	1 \\ v \\ |v|^2
\end{pmatrix} dv dx,
	\Ee
where $\alpha$ and $\kappa$ are given as in \eqref{apriori} and $\b> 2$ is to be determined. Note that \eqref{initial} and \eqref{initial normal}  imply finite initial entropy $H(F_0) <\infty$.

We begin with the stretched-exponential case. All global existence assertions below are conditional: the theorem starts from a global mild solution satisfying the mixed-norm bound. Neither the existence of such a solution nor propagation of that bound follows from the initial conditions alone. Supremum norms of measurable solutions are understood as essential supremum norms. 

\begin{theorem} \label{thm:cond} 
	Let $F_0\geq 0$ satisfy \eqref{initial} and \eqref{initial normal}. Suppose that $F$ is a global mild solution  to \eqref{BE} emanating from $F_0$ and satisfying \eqref{apriori}. Then $F$ is the unique global dissipative mild solution with this datum and is physically admissible. Moreover, $F$ converges exponentially to the associated Maxwellian $\mu$:
	\Be \label{decay F}
		\|w_1(F(t)-\mu)\|_{L^\infty_{x,v}}
	\leq C_{F}e^{-\lambda_\infty t},
	\Ee
	where $w_1(v):=\langle v\rangle^{\beta_1}e^{\a|v|^{\k}}$ and $\b_1 \geq 0$ is specified in Proposition~\ref{prop uniform bound}. The exponential decay rate $\lambda_\infty$ is constructive and independent of $F$, whereas the constant $C_F$ depends on $F$ and is non-constructive.
\end{theorem}

If $F(t,x,v)$ has a (large) stretched exponential velocity upper bound, then we obtain the following Corollary.
\begin{corollary} \label{coro:cond}
	Let $F_0\geq 0$ satisfy \eqref{initial normal}. If a global mild solution $F(t,x,v)\geq 0$ of \eqref{BE} emanating from $F_0$ satisfies 
	\Be \label{cond Coro}
		\sup_{[0,\infty)\times \T^3\times \R^3} \avv^{\b_1} e^{\a|v|^\k} F(t,x,v)  < \infty, \quad \b_1\geq 0, \ \ \a > 0, \ \ 1 < \k < 2,
	\Ee
	then $F$ is the unique global dissipative mild solution emanating from $F_0$ and is physically admissible. Moreover, $F$ satisfies
	\Be \label{decay F coro}
		\|w_1(F(t)-\mu)\|_{L^\infty_{x,v}}
	\leq C_{F}e^{-\lambda_\infty t},\quad w_1(v) = \avv^{\b_1} e^{\a|v|^\k},
	\Ee
	with no loss in the velocity weight. The exponential decay rate $\lambda_\infty$ is constructive and independent of $F$, whereas the constant $C_F$ depends on $F$ and is non-constructive.
\end{corollary}

	\begin{remark} 
		Theorem~\ref{thm:cond} is conditional only through the global tail
		criterion \eqref{apriori}. No smallness assumption is imposed on the
		weighted $L^\infty$ norm, the relative entropy, or the distance to
		equilibrium. Thus the theorem is genuinely non-perturbative. Its content
		is a rigidity statement: once a global mild solution satisfies
		\eqref{apriori}, its global dynamics are completely determined.
	\end{remark}
	
	\begin{remark} 
		Let $\mathcal M(F_0)$ denote the set of global mild solutions emanating
		from $F_0$ and satisfying \eqref{apriori}. Then Theorem~\ref{thm:cond}
		and the weak--strong uniqueness argument imply that
		\[
		\mathcal M(F_0)\neq\varnothing
		\quad\Longrightarrow\quad
		\mathcal M(F_0)
		=
		\mathcal D_{\rm DL}(F_0)
		=
		\{F\}.
		\]
		Thus the tail criterion not only yields uniqueness within the mild class,
		but also selects the DiPerna--Lions approximation limit. In particular, $F$ is physically admissible and satisfies the entropy inequality together with the conservation of mass, momentum, and energy. 
	\end{remark}
	
	\begin{remark} 
        In particular, $F_0$ may vanish on nontrivial subsets of phase space. The positivity needed for the subsequent coercivity argument is generated dynamically through the mechanism developed in Section~3. More precisely, we choose a time $T_{\rm glb}>0$ within the local existence interval at which a Gaussian lower bound has been established, and then proceed with the lower-bound argument starting from $T_{\rm glb}$. Thus, neither a strictly positive density nor a Gaussian lower bound is assumed at the initial time $t=0$, and the initial datum may contain vacuum. Finally, the short-time upper estimate is used to extend the subsequent uniform bounds back to the interval $[0,t_s]$.
	\end{remark}
	
	\begin{remark}  
		\label{rem:structural-consequences}
		The conclusion of Theorem~\ref{thm:cond} is obtained through several
		additional uniform estimates which are of independent interest.
		More precisely, the tail condition \eqref{apriori} implies the existence
		of constants $\rho_->0$ and $\lambda_F>0$ such that
		\[
		\rho(t,x)
		:=
		\int_{\mathbb R^3}F(t,x,v)\,dv
		\geq \rho_-,
		\]
		and
		\[
		L(F)(t,x,v)
		\geq
		\lambda_F\langle v\rangle^\gamma
		\]
		for almost every $(t,x,v)\in
		[0,\infty)\times\mathbb T^3\times\mathbb R^3$ where we assumed $F_0$ has a Gaussian lower bound without loss of generality by previous remark.
		Moreover, the exponent $\beta_1\geq0$ determined in Proposition~\ref{prop uniform bound},
		\[
		\sup_{t\geq0}
		\left\|
		\langle v\rangle^{\beta_1} 
		e^{\alpha|v|^\kappa}F(t)
		\right\|_{L^\infty_{x,v}}
		<\infty.
		\]
		Thus the upper-tail condition \eqref{apriori} dynamically generates
		both macroscopic coercivity and global pointwise control.
	\end{remark}
	
	\begin{remark}  
		\label{rem:homogeneous-compatibility}
		The tail condition \eqref{apriori} is consistent with the classical
		spatially homogeneous theory. Indeed, for the homogeneous Boltzmann
		equation with hard potentials, the Gaussian-weighted $L^1$ estimate of
		Bobylev, later extended to the relevant class of hard potentials and
		used by Gamba--Panferov--Villani in their pointwise Maxwellian theory,
		shows that
		\[
		\int_{\mathbb R^3} e^{a_0|v|^2}F_0(v)\,dv<\infty
		\]
		implies
		\[
		\sup_{t\geq0}
		\int_{\mathbb R^3} e^{a|v|^2}F(t,v)\,dv<\infty
		\]
		for some $a>0$; see, for instance,
		\cite{Bobylev,GPV}. The propagated Gaussian rate $a$ need not coincide
		with the initial rate $a_0$ and may be smaller.
		
		This loss of Gaussian rate causes no difficulty for the
		stretched-exponential criterion used here. Indeed, whenever
		$1<\kappa<2$, for every $\alpha>0$ and $k\geq0$,
		\[
		\langle v\rangle^k e^{\alpha|v|^\kappa}
		\leq C e^{a|v|^2}
		\]
		for any $a>0$. Hence the classical homogeneous Gaussian-moment
		propagation immediately implies \eqref{apriori}, independently of the
		size of the propagated Gaussian exponent. In this sense,
		\eqref{apriori} is a natural inhomogeneous counterpart of the classical
		homogeneous high-energy tail estimates.
	\end{remark}
	
	\begin{remark} 
		The exponent $\b_1$ in Theorem~\ref{thm:cond} is determined by the
		weighted gain estimates of Section~5. In particular, when the tail
		order $k$ in \eqref{apriori} is chosen sufficiently large relative to
		the continuation weight, the argument closes with $\b_1=\beta$, so that
		there is no loss in the polynomial part of the initial
		stretched-exponential weight.
	\end{remark}

	We next turn to the polynomially weighted setting, where an additional lower-bound assumption on the initial datum is required. For initial data we consider 
	\Be \label{initial-poly}
		\|\avv^\b F_0\|_{L^\infty_{x,v}} < \infty, \quad \b  > 3
	\Ee
	and assume
	\Be \label{apriori-poly}
		\sup_{t\geq 0}\int_{v}\sup_{x}F(t,x,v)\avv^k dv \leq A_0 < \infty
	\Ee
    for some $k\geq 2$. The range of $k$ will be specified in the following theorems. We additionally assume that the initial datum admits a spatially uniform lower profile in velocity:
	\Be \label{ini profile-poly}
		\exists \varphi(v) \geq 0 \ \ \text{such that} \ \ F_0(x,v) \geq \varphi(v) \ \ \text{and}  \ \ \int_v \varphi(v) dv  > 0.
	\Ee
	Under the additional assumption \eqref{ini profile-poly}, we obtain an analogue of Theorem~\ref{thm:cond}. Note that $H(F_0) < \infty$ if $\b  > 3$ in \eqref{initial-poly}.

    \begin{theorem} \label{thm:cond-p}
		Let $F_0\geq0$ satisfy \eqref{initial-poly}, \eqref{ini profile-poly}, and \eqref{initial normal}. Let $F$ be a global mild solution of \eqref{BE} with this datum and satisfying \eqref{apriori-poly}. Choose $\beta$ and $k$ sufficiently large that
		\[
		k\geq2+\gamma,\qquad
		\beta_1(\beta, k)>\max\{11,\beta_\infty\},
		\]
		where $\beta_1$ is the exponent in Proposition~\ref{prop uniform bound}, applied as in Remark~\ref{rmk unif-p}, and $\beta_\infty$ is defined in Theorem~\ref{Mainresultlargeamp-p}. Then $F$ is the unique global dissipative mild solution with this datum and is physically admissible. Moreover, $F$ converges exponentially to the associated Maxwellian $\mu$:
		\Be \label{decay F-p}
		\|w_1(F(t)-\mu)\|_{L^\infty_{x,v}}
		\leq C_{F}e^{-\lambda_\infty t},
		\Ee
		where $w_1(v):=\langle v\rangle^{\beta_1}$ and $\b_1 \geq 0$ is specified in Proposition~\ref{prop uniform bound}. The exponential decay rate $\lambda_\infty$ is constructive and independent of $F$, whereas the constant $C_F$ depends on $F$ and is non-constructive.
	\end{theorem}


    The threshold $11$ in Theorem~\(\ref{thm:cond-p}\) is not claimed to be optimal; it results from choosing $p=\frac{3}{2}$ in Section~6.
    
			If $F(t,x,v)$ has a polynomially decaying velocity upper bound, then we obtain the following Corollary.
		\begin{corollary} \label{coro:cond-p}
			There exists a positive constant $\b_*$ such that if $\b_1 \geq \b_*$ and a global mild solution $F(t,x,v)\geq 0$ of \eqref{BE} satisfies  \eqref{initial normal} and 
			\Be \label{cond Coro-p}
			\sup_{[0,\infty)\times \T^3\times \R^3} \avv^{\b_1} F(t,x,v)  < \infty,
			\Ee
			then $F$ is the unique global dissipative mild solution emanating from $F_0$ and physically admissible. Moreover, $F$ satisfies
			\Be \label{decay F coro-p}
			\|w_1(F(t)-\mu)\|_{L^\infty_{x,v}}
			\leq C_{F}e^{-\lambda_\infty t},\quad w_1(v) = \avv^{\b_1}
			\Ee
		with no loss in the velocity weight. The exponential decay rate $\lambda_\infty$ is constructive and independent of $F$, whereas the constant $C_F$ depends on $F$ and is non-constructive. 
		\end{corollary}

	\begin{remark}[Relation with the stretched-exponential case]
		\label{rem:poly-relation}
		Apart from the initial positivity assumption, the remarks following
		Theorem~\ref{thm:cond} apply equally to the polynomially weighted
		Theorem~\ref{thm:cond-p}. In particular, the polynomial tail criterion
		leads to the same structural consequences: uniqueness in the dissipative
		mild class, macroscopic coercivity, uniform weighted $L^\infty$ control,
		and relaxation to equilibrium. The essential difference is that, in the
		polynomial setting, we impose the additional lower profile assumption
		\eqref{ini profile-poly}.
		
		This additional assumption reflects a limitation of the argument developed in the present paper rather than a fundamental obstruction to the positivity mechanism. Indeed, the double-Duhamel and spreading argument of Section~3 can still generate a Gaussian lower bound after an arbitrarily short positive time, provided that the necessary local pointwise bounds and moment and loss-frequency controls are available. However, in the purely polynomially weighted setting, we do not establish the large-data local propagation estimate required to apply this argument before obtaining the global weighted $L^\infty$ bound.
        Consequently, even if a Gaussian lower bound is produced at some
		$T_{\mathrm{glb}}>0$, we cannot justify that
		\[
		F(T_{\mathrm{glb}})
		\in L^\infty_{x,v}(\langle v\rangle^{\beta'})
		\]
		for a sufficiently large polynomial exponent $\beta'$ needed in the
		subsequent weighted estimates. Thus, unlike in the stretched-exponential
		setting, we cannot simply take $T_{\mathrm{glb}}$ as a new initial time
		while retaining the required solution class.
		
		For this reason, Theorem~\ref{thm:cond-p} assumes
		\eqref{ini profile-poly} directly. This hypothesis supplies the finite-time
		nondegeneracy needed to initiate the tail-to-coercivity argument without
		appealing to a restart of the solution at a positive time. If an
		appropriate large-data local theory in polynomially weighted
		$L^\infty$ spaces were available, the lower profile assumption could in
		principle be replaced by the dynamically generated short-time
		positivity mechanism used in Theorem~\ref{thm:cond}.
	\end{remark}
	
	\begin{remark} 
		\label{rem:homogeneous-poly}
		In the spatially homogeneous setting, the polynomial tail condition
		\eqref{apriori-poly} reduces simply to the uniform polynomial moment bound
		\[
		\sup_{t\geq0}
		\int_{\mathbb R^3}
		\langle v\rangle^k F(t,v)\,dv
		<\infty .
		\]
		Such bounds are classical for the homogeneous Boltzmann equation with
		hard potentials. Starting from the Povzner moment inequalities, global
		propagation and generation of polynomial moments were developed in the
		works of Elmroth, Desvillettes, Wennberg, and subsequent authors; see,
		for instance,
		\cite{Elmroth1983,Desvillettes1993,Wennberg1994,LuMouhot2012}
		and the references therein.
		Thus \eqref{apriori-poly} is fully consistent with the classical
		homogeneous large-data theory. The additional difficulty in the
		spatially inhomogeneous problem is the possible concentration in the
		spatial variable, which is encoded here through the stronger quantity
		involving $\sup_x F$.
	\end{remark}

   If a global mild solution $F$ is additionally assumed to satisfy $F\in L^\infty([0,T);L^\infty_{\b_1})$ for any finite $T <\infty$, then the initial lower profile assumption \eqref{ini profile-poly} can be removed, as discussed in Remark~\ref{rem:poly-relation}. Moreover, it suffices to assume \eqref{apriori-poly} with $k=2$.
	
	\begin{theorem} \label{thm:cond-low 2}
    Let $F$ be a global nonnegative mild solution of \eqref{BE} satisfying the normalization condition \eqref{initial normal}, the second-order moment condition \Be \label{M2}
			\sup_{t\geq 0}\int_{v}\sup_{x}F(t,x,v)\avv^2 dv \leq A_0 < \infty.
		\Ee
      and
    \Be \notag
    \sup_{0\le t\le T}
    \|\langle v\rangle^{\beta_1}F(t)\|_{L^\infty_{x,v}}<\infty
    \qquad\text{for every } \ T<\infty.
    \Ee
    Then, there exists a positive constant $\b_*$ such that if $\b_1 > \b_*$, then $F$ satisfies
		\Be \label{decay F M2}
		\|w_1(F(t)-\mu)\|_{L^\infty_{x,v}}
		\leq C_{F}e^{-\lambda_\infty t}.
		\Ee
		The exponential decay rate $\lambda_\infty$ is constructive and independent of $F$, whereas the constant $C_F$ depends on $F$ and is non-constructive.
	\end{theorem}
	If we want uniqueness in Theorem \ref{thm:cond-low 2} in addition, it is enough to add an extra condition
	\[
		\sup_x \int_{v} \avv^{2+\g} F(t,x,v) dv \in L^1_{loc}[0,\infty)
	\]
	by Lemma \ref{lem_unique}.  
	
	Theorem \ref{thm:cond-low 2} can be viewed as a counterpart to \cite{Imbert2021} and \cite{Desvillettes}. More precisely, assuming only that $F(t)$ is locally uniformly bounded in time in a polynomially weighted $L^\infty_{\beta_1}$, we establish convergence to equilibrium in the same space under the second-order moment condition \eqref{apriori-poly} with $k=2$. This moment condition is analogous to the uniform upper bound on the local energy assumed in \cite{Imbert2021} and \cite{CPDE}, although it is stronger in the sense that the spatial supremum is taken inside the velocity integral. On the other hand, unlike those works, our result does not require a uniform positive lower bound on the local density or a uniform upper bound on the local entropy.

	\subsection{Proof Strategy and Main Ideas}
	
	In this subsection, we describe the main ideas in the proof of
	Theorem~\ref{thm:cond}. We first explain the uniqueness mechanism and
	the generation of a Gaussian lower bound near the initial time. 
	
	\medskip
	
	\subsubsection{Uniqueness in the dissipative mild class.}
	The first step is to establish uniqueness in a class substantially larger
	than the weighted $L^\infty$ class used in the later analysis.
	This part relies on a weak--strong uniqueness argument of Lu type and is
	independent of the positivity mechanism developed below.
	
	Let $F$ be the global mild solution in Theorem~\ref{thm:cond},
	satisfying the tail condition \eqref{apriori}, and let
	$G\in \mathcal U_\infty(F_0)$ be any dissipative mild solution with the
	same initial datum. The assumption \eqref{apriori} implies, in particular,
	that
	\[
	\Lambda_F(t)
	:=
	\sup_{x\in\mathbb T^3}
	\int_{\mathbb R^3}
	\langle v\rangle^{2+\gamma}
	F(t,x,v)\,dv
	<\infty ,
	\]
	uniformly in time. It also provides sufficient velocity integrability to
	justify the collision invariants for $F$, and hence $F$ conserves mass and
	energy. Since $G$ is dissipative,
	\[
	\iint_{\mathbb T^3\times\mathbb R^3}
	\langle v\rangle^2 G(t,x,v)\,dx\,dv
	\leq
	\iint_{\mathbb T^3\times\mathbb R^3}
	\langle v\rangle^2 F(t,x,v)\,dx\,dv .
	\]
	
	Following the argument of Lu~\cite{Lu}, we consider the one-sided
	difference
	\[
	Y(t)
	:=
	\iint_{\mathbb T^3\times\mathbb R^3}
	\langle v\rangle^2
	(F-G)_+(t,x,v)\,dx\,dv .
	\]
	The elementary decomposition
	\[
	|F-G|
	=
	(G-F)+2(F-G)_+
	\]
	together with the preceding mass--energy inequality allows the full
	weighted $L^1$ distance to be controlled by $Y(t)$. The one-sided
	algebraic inequality for the collision operator then yields
	\[
	Y(t)
	\leq
	C\int_0^t \Lambda_F(s)Y(s)\,ds .
	\]
	Since $Y(0)=0$, Gr\"onwall's inequality gives $Y\equiv0$, and hence
	$F=G$. Thus $F$ is the unique solution in the whole dissipative mild
	class $\mathcal U_\infty(F_0)$. This is the mechanism developed in
	Lemma~\ref{lem_unique}.
	
	This uniqueness result also identifies the physically relevant solution.
	Indeed, every DiPerna--Lions approximation limit considered here belongs
	to the dissipative mild class. Therefore every
	\[
	F_{\mathrm{DL}}\in\mathcal D_{\mathrm{DL}}(F_0)
	\]
	must coincide with $F$. Consequently, $F$ inherits the entropy inequality
	and the other admissibility properties of the DiPerna--Lions solution.
	In particular, the entropy inequality is not assumed in the uniqueness
	argument, but is instead obtained as a consequence of the identification
	\[
	F=F_{\mathrm{DL}}.
	\]
	This observation will play an important role in the compactness and
	long-time arguments below.
	
	\medskip
	\subsubsection{Generation of a Gaussian lower bound near the initial time.}
	A new difficulty appears at the beginning of the positivity argument.
	We impose no pointwise lower bound on $F_0$; in particular, the initial
	distribution may contain local vacuum. Therefore, the classical spreading
	argument for $Q^+$ cannot be applied directly.
	
	We first recover the short-time regularity needed for the positivity
	construction. Proposition~\ref{prop_loc} gives a unique mild solution
	$\widetilde F$ on some interval $[0,\widehat t_0]$ satisfying
	\[
	\sup_{0\leq t\leq \widehat t_0}
	\left\|
	\langle v\rangle^\beta
	e^{\alpha |v|^\kappa}
	\widetilde F(t)
	\right\|_{L^\infty_{x,v}}
	\leq
	2
	\left\|
	\langle v\rangle^\beta
	e^{\alpha |v|^\kappa}
	F_0
	\right\|_{L^\infty_{x,v}} .
	\]
	By the uniqueness result established above, the global solution $F$ in
	Theorem~\ref{thm:cond} must coincide with $\widetilde F$ on
	$[0,\widehat t_0]$. Thus
	\[
	F=\widetilde F
	\qquad\text{on }[0,\widehat t_0],
	\]
	and the same short-time weighted $L^\infty$ bound is available for the
	global solution under consideration. This allows us to perform the
	positivity construction near $t=0$ without imposing any additional
	regularity assumption on $F$.
	
	The construction proceeds in three stages.
	
	First, conservation of the total mass and energy, together with the
	local $L^\infty$ bound, implies that a positive amount of the initial
	mass lies above a fixed positive level inside a bounded velocity region.
	Consequently, the set
	\[
	\left\{
	(x,v)\in\mathbb T^3\times B_R:
	F_0(x,v)\geq \mu_0
	\right\}
	\]
	has positive measure for suitable $R,\mu_0>0$. By the Lebesgue density
	theorem, one can therefore find a phase-space point $(x_*,v_*)$ and an
	anisotropic cylinder
	\[
	B_x(x_*,Ar)\times B_v(v_*,r)
	\]
	such that
	\[
	F_0(x,v)\geq \mu_0
	\]
	except on a subset of arbitrarily small relative measure. Thus the
	argument starts only from measure-theoretic positivity and does not use
	any pointwise lower bound on the initial datum.
	
	Second, this almost-everywhere positivity is converted into a genuine
	pointwise lower bound by a double-Duhamel argument. Along free-transport
	characteristics, the positive portion of the initial cylinder survives
	with a controlled loss factor. An averaging argument for the transported
	bad set shows that, for most spatial points in a suitable moving tube,
	only a small portion of the velocity ball is missing. The small-hole
	gain estimate then gives a pointwise lower bound for $Q^+(F,F)$.
	
	Inserting this gain estimate once more into the Duhamel formula and
	performing a second averaging argument removes the remaining exceptional
	spatial set. We consequently obtain, at some arbitrarily small positive
	time $t_0$,
	\[
	F(t_0,x,v)\geq \kappa_0>0
	\]
	on a genuine product set
	\[
	(x,v)\in
	B_x(x_0,r_x)\times B_v(v_0,r_0).
	\]
	This is the short-time positivity mechanism developed in
	Proposition~\ref{prop:2.8}.
	
	The third step is to propagate this local product-ball lower bound.
	Repeated applications of the velocity-spreading property of the gain
	operator enlarge the velocity ball:
	\[
	B_v(v_0,r_n)
	\longrightarrow
	B_v(v_0,r_{n+1}),
	\qquad r_{n+1}>r_n .
	\]
	The associated time increments are chosen sufficiently small so that
	all relevant backward characteristics remain inside the spatial region
	where the previous lower bound is available.
	
	Once the velocity radius has become sufficiently large, free transport
	provides the spatial mixing needed to cover the torus. Indeed, since
	$\mathbb T^3$ has finite diameter, sufficiently large velocities can
	connect the original positive spatial region to an arbitrary target
	point in a fixed positive time. We therefore obtain, at some
	$T_1\in(0,\widehat t_0)$, a uniform local velocity profile at every spatial
	point:
	\[
	\forall\,x\in\mathbb T^3,
	\qquad
	F(T_1,y,v)
	\geq \beta_0>0
	\]
	for $y$ in a fixed neighborhood of $x$ and for $v$ in a velocity ball
	whose parameters are uniformly controlled.
	
	In particular, this yields a strictly positive spatially uniform density
	lower bound at time $T_1$:
	\[
	\rho(T_1,x)
	=
	\int_{\mathbb R^3}F(T_1,x,v)\,dv
	\geq \rho_1>0,
	\qquad x\in\mathbb T^3.
	\]
	The short-time weighted $L^\infty$ estimate also gives the required
	uniform upper bounds on the local macroscopic quantities. In particular,
	the weighted $L^\infty$ bound implies the local entropy-density bound
	needed in the standard quantitative positivity theory.
	
	We may therefore apply the Gaussian lower-bound theorem for the
	spatially inhomogeneous cutoff Boltzmann equation
	(see, for instance,~\cite{CPDE}) and obtain some
	\[
	0<T_{\mathrm{glb}}<\widehat t_0
	\]
	such that
	\[
	F(T_{\mathrm{glb}},x,v)
	\geq
	a_0 e^{-c_0|v|^2},
	\qquad
	(x,v)\in\mathbb T^3\times\mathbb R^3,
	\]
	for some constants $a_0,c_0>0$. The constants may depend on the initial
	datum, the a priori tail bound, and $T_{\mathrm{glb}}$.
	
	Since $T_{\mathrm{glb}}$ lies inside the local existence interval,
	the weighted $L^\infty$ estimate is still available:
	\[
	\left\|
	\langle v\rangle^\beta
	e^{\alpha |v|^\kappa}
	F(T_{\mathrm{glb}})
	\right\|_{L^\infty_{x,v}}
	\leq
	2
	\left\|
	\langle v\rangle^\beta
	e^{\alpha |v|^\kappa}
	F_0
	\right\|_{L^\infty_{x,v}} .
	\]
	We may therefore restart the subsequent argument at
	$t=T_{\mathrm{glb}}$. After this harmless time shift, we may assume,
	without loss of generality, that the initial datum satisfies a uniform
	Gaussian lower bound.

	\medskip
	\subsubsection{Tail-to-coercivity: uniform macroscopic lower bounds.}
	Having generated a Gaussian lower bound at a positive time, we next show
	that the tail condition \eqref{apriori} prevents the solution from
	developing macroscopic degeneracy at later times. More precisely, we prove
	that there exist constants $\varepsilon_{\rho} > 0$ and $\lambda_{F} > 0$ such that
	\[
	\rho(t,x)
	:=
	\int_{\mathbb R^3}F(t,x,v)\,dv
	\geq \varepsilon_{\rho},
	\]
	and
	\[
	L(F)(t,x,v)
	\geq
	\lambda_{F} \langle v\rangle^\gamma
	\]
	for almost every $(t,x,v)\in [0,\infty)\times\mathbb T^3\times\mathbb R^3$.
	The constants obtained in this step are non-constructive.
	
	The argument is based on a contradiction and compactness mechanism.
	For fixed $R>0$ and $a\in\mathbb R^3$, introduce the truncated collision
	frequency
	\[
	u^{(a)}(t,x)
	:=
	\int_{|v|\leq R}
	|a-v|^\gamma F(t,x,v)\,dv .
	\]
	Since a Gaussian lower bound is available at the new initial time, the
	mild formulation and the upper bound on the loss coefficient give, on
	every finite time interval,
	\[
	u^{(a)}(t,x)>0
	\]
	uniformly in $x$ and locally uniformly in $a$. In particular, a failure
	of a global uniform lower bound cannot occur at a bounded sequence of
	times; any degenerating sequence must escape to infinity.
	
	Suppose, toward a contradiction, that no uniform lower bound exists.
	Then one can find $t_n\to\infty$ such that
	\[
	\inf_{|a|\leq R_1}u^{(a)}(t,x)
	\]
	becomes arbitrarily small on a set of positive $(t,x)$-measure in the
	time slabs $[t_n,t_n+1]$. We introduce the time-translated solutions
	\[
	F_n(t,x,v):=F(t_n+t,x,v),
	\qquad
	(t,x,v)\in (0,1)\times\mathbb T^3\times\mathbb R^3 .
	\]
	
	The identification with the DiPerna--Lions solution established in the
	uniqueness step is crucial here. It supplies the entropy and energy
	inequalities, while the tail condition \eqref{apriori} gives uniform
	velocity tightness. Standard large-time compactness results for the
	Boltzmann equation then imply, after extraction of a subsequence, that
	$F_n$ converges weakly in $L^1$ to a local Maxwellian solving the free
	transport equation. By the conservation laws and our normalization, this
	limit is necessarily the global Maxwellian $\mu$.
	
	To turn this weak convergence into pointwise information on the
	collision frequency, we use the $L^1$ velocity averaging lemma. For
	\[
	u_n^{(a)}(t,x)
	:=
	\int_{|v|\leq R}
	|a-v|^\gamma F_n(t,x,v)\,dv ,
	\]
	velocity averaging gives relative compactness in $L^1_{t,x}$ for every
	fixed $a$. Moreover,
	\[
	\|u_n^{(a)}-u_n^{(b)}\|_{L^1_{t,x}}
	\leq C_R |a-b|^\gamma ,
	\]
	uniformly in $n$. A diagonal extraction over a countable dense subset
	of $B_{R_1}$, followed by this parameter continuity, therefore produces
	one common subsequence such that
	\[
	\sup_{|a|\leq R_1}
	\left|
	u_n^{(a)}(t,x)-c_{R,a}
	\right|
	\longrightarrow 0
	\]
	for almost every $(t,x)$, where
	\[
	c_{R,a}
	:=
	\int_{|v|\leq R}|a-v|^\gamma\mu(v)\,dv .
	\]
	Since
	\[
	\inf_{|a|\leq R_1}c_{R,a}>0,
	\]
	this already gives a uniform positive lower bound away from an
	arbitrarily small exceptional set in $(t,x)$.
	
	The main difficulty is to control the solution on this exceptional set.
	For this purpose, we extend the translated solutions to a slightly larger
	time interval and use Egoroff's theorem to isolate a bad space--time set
	of arbitrarily small measure. Outside this set, two pieces of information
	are available simultaneously: the truncated collision frequency is
	uniformly positive, and the local entropy density
	\[
	\int_{|v|\leq R}F_n(t,x,v)|\log F_n(t,x,v)|\,dv
	\]
	is bounded. We refer to this complement as the good set.
	
	We then return to the mild formulation and perform a double-Duhamel
	expansion. The first Duhamel step expresses $u_n^{(a)}(t,x)$ in terms of
	$Q^+(F_n,F_n)$ at earlier points along characteristics, while a second
	Duhamel expansion is applied to one of the factors entering the gain
	term. After the Carleman representation, the resulting expression
	contains planar averages of the form
	\[
	\int_{E_{v,w}\cap B_R}
	F_n(s,y,v')\,dE(v').
	\]
	The local entropy bound on the good set prevents the mass of $F_n$ from
	concentrating in arbitrarily thin regions. A directional averaging
	argument, combined with an Orlicz--H\"older estimate, then shows that
	there exists a set of directions, and hence a set of velocities $w$ of
	uniformly positive measure, on which these planar averages admit a
	strictly positive lower bound. This is the non-concentration mechanism
	that replaces any assumed pointwise lower bound on the macroscopic
	fields.
	
	At the second gain level, we write
	\[
	Q^+(F_n,F_n)
	=
	Q^+(\mu,\mu)
	+
	\bigl(Q^+(F_n,F_n)-Q^+(\mu,\mu)\bigr).
	\]
	The Maxwellian contribution is strictly positive on bounded velocity
	regions. On the other hand, the compactness theory for the gain operator
	gives convergence in measure
	\[
	Q^+(F_n,F_n)\longrightarrow Q^+(\mu,\mu)
	\]
	locally in $(t,x,v)$. Since the family $Q^+(F_n,F_n)$ is uniformly
	integrable, Vitali's theorem upgrades this convergence to the
	$L^1_{\mathrm{loc}}$ convergence needed in the double-Duhamel estimate.
	Consequently, for all sufficiently large $n$, the positive Maxwellian
	contribution dominates the error term.
	
	The smallness of the bad space--time set, together with the transport
	changes of variables appearing in the two Duhamel integrations, ensures
	that a definite portion of the integration domain lies in the good set.
	Combining this geometric observation with the preceding Carleman
	lower bound yields
	\[
	\inf_{|a|\leq R_1}
	u_n^{(a)}(t,x)
	\geq \varepsilon_{R_1}>0
	\]
	also on the exceptional set. This contradicts the assumed degenerating
	sequence. We conclude that
	\[
	\inf_{\substack{t\geq0,\;x\in\mathbb T^3\\|a|\leq R_1}}
	\int_{|v|\leq R}
	|a-v|^\gamma F(t,x,v)\,dv
	\geq
	\varepsilon_{R_1}>0 .
	\]
	
	The remaining velocities are treated by elementary geometry. Choosing
	$R_1$ sufficiently larger than $R$, for $|a|\geq R_1$ and $|v|\leq R$
	we have
	\[
	|a-v|\geq |a|-R\gtrsim \langle a\rangle .
	\]
	The lower bound already obtained on the mass contained in $B_R$ therefore
	gives
	\[
	L(F)(t,x,a)
	=
	\int_{\mathbb R^3}|a-v|^\gamma F(t,x,v)\,dv
	\geq
	\lambda_{F} \langle a\rangle^\gamma .
	\]
	Likewise, taking a fixed parameter value in the truncated estimate gives
	a uniform lower bound for the density,
	\[
	\rho(t,x)\geq \varepsilon_{\rho} > 0.
	\]
	
	Thus the single upper-tail condition \eqref{apriori}, together with the
	short-time positivity generated in the preceding step, rules out
	macroscopic degeneration at arbitrarily large times. In particular, the
	nonlinear loss term acquires the global coercivity
	\[
	L(F)(t,x,v)\geq \lambda_{F} \langle v\rangle^\gamma ,
	\]
	which is the key input for the global weighted $L^\infty$ estimate in
	the next step.

	\medskip
	
	\subsubsection{Uniform weighted $L^\infty$ estimate.}
	Once the uniform coercive estimate
	\[
	L(F)(t,x,v)\geq \lambda_{F} \langle v\rangle^\gamma
	\]
	has been established, we close the global weighted $L^\infty$ estimate
	by an Arkeryd-type argument. Set
	\[
	f(t,x,v):=e^{\alpha |v|^\kappa}F(t,x,v).
	\]
	Using the Carleman representation of the gain operator and the
	submultiplicative property of the stretched-exponential weight under
	collisions, the gain term is reduced to the product of a singular
	three-dimensional integral and a planar velocity integral. The former
	is controlled by the assumed tail bound \eqref{apriori}, while the
	latter is estimated by combining the mild formulation with the
	coercive damping generated by $L(F)$.
	
	Here an important feature of the a priori condition \eqref{apriori}
	enters the argument. When the mild formulation is inserted into the
	Carleman representation, a pre--post collisional change of variables
	does not preserve the transported spatial point
	\[
	x-v(t-s).
	\]
	After the change of variables, the spatial argument depends on the
	corresponding pre- or post-collisional velocity. Consequently, the resulting velocity integrals naturally require control of
	\[
	\int_{\mathbb R^3}
	e^{\alpha|v|^\kappa}\langle v\rangle^k
	\sup_{x\in\mathbb T^3}F(t,x,v)\,dv,
	\]
	which is precisely the quantity controlled by \eqref{apriori}. This choice is also consistent with the functional framework of Gualdani, Mischler, and Mouhot \cite{Gualdani2017}, who employed weighted \(L^1_vL^\infty_x\) spaces in their perturbative Boltzmann theory, providing further support for the naturality of placing the spatial supremum inside the velocity integral.
	Thus the placement of the spatial supremum inside the velocity
	integral is essential for the collision change-of-variables argument,
	rather than merely a technical strengthening of the hypothesis.
	
	With this observation, a first Carleman estimate gives a uniform
	$L^\infty$ bound with a lower polynomial weight. Feeding this estimate
	back into the gain estimate and exploiting the additional velocity
	decay supplied by \eqref{apriori}, we recover the desired polynomial
	weight and obtain
	\[
	\sup_{t\geq0}
	\left\|
	\langle v\rangle^{\b_1}
	e^{\alpha |v|^\kappa}F(t)
	\right\|_{L^\infty_{x,v}}
	\leq C_{\infty,F},
	\]
    where $C_{\infty,F}$ depends on $\alpha,\kappa,\beta,k,\gamma, F_0, A_0,\varepsilon_L,\lambda_F$ and $\b_1>0$ is specified in Proposition~\ref{prop uniform bound}.
	This completes the global pointwise estimate required in the
	subsequent long-time analysis. 
    
	\medskip
	\subsubsection{Entropy decay.}
	We next prove the qualitative decay of the relative entropy.
	The uniform weighted $L^\infty$ estimate obtained above provides
	velocity non-concentration, while the previous sections give a uniform Gaussian lower bound. These estimates yield uniform upper and lower bounds for the density and temperature, together with a uniform bound on the bulk velocity. Together with the entropy-production estimate,
	these bounds imply that, for every sequence $t_n\to\infty$,
	\[
	F_n(t,x,v):=F(t_n+t,x,v),
	\qquad 0<t<1,
	\]
	satisfies, after extraction of a subsequence,
	\[
	\|F_n-M_{F_n}\|_{L^1((0,1)\times\mathbb T^3\times\mathbb R^3)}
	\longrightarrow0,
	\]
	where $M_{F_n}$ denotes the corresponding local Maxwellian.
	
	The $L^1$ velocity averaging lemma, combined with the uniform velocity
	tails, gives strong compactness of the macroscopic fields and hence
	strong $L^1$ compactness of $\{M_{F_n}\}$. On the other hand,
    time-shift identification proved in Proposition~\ref{Desvillettes}, using the direct measurable-Maxwellian rigidity Lemma~\ref{lem:r1-Maxwellian-rigidity}, shows that every such weak limit is the global Maxwellian $\mu$. 
    Since $F_n-M_{F_n}\to0$ strongly in $L^1$, the weak
	limit of $F_n$ and the strong limit of $M_{F_n}$ coincide. Consequently,
	\[
	F_n\longrightarrow\mu
	\qquad\text{in }
	L^1((0,1)\times\mathbb T^3\times\mathbb R^3)
	\]
	along the extracted subsequence. A reverse entropy estimate, together
	with the monotonicity of the relative entropy, then yields
	\[
	H(F(t)\mid\mu)\longrightarrow0,
	\qquad t\to\infty.
	\]
	
	In fact, the same macroscopic nondegeneracy and weighted \(L^\infty\) control also allow one to implement a constructive entropy-decay scheme. Thus, the above qualitative argument can be strengthened to yield a quantitative decay rate. We note, however, that the macroscopic lower bounds obtained in the preceding section are themselves non-constructive; consequently, this refinement does not improve the quantitative statement of Theorem~\ref{thm:cond}. Therefore, we refer to the forthcoming work~\cite{DH-entropy} for the constructive entropy-decay argument.
    
	\medskip
	\subsubsection{Exponential decay.}
	Finally, we combine the qualitative entropy relaxation with the
	perturbative theory near the global Maxwellian. At this stage, the
	solution is uniformly bounded in the stretched-exponential weighted
	$L^\infty$ space, and the previous step gives
	\[
	H(F(t)\mid\mu)\longrightarrow0
	\qquad\text{as } t\to\infty .
	\]
    The entropy becomes smaller than a large-amplitude threshold at some finite time $T_{ent}$. Small entropy alone does not imply smallness in weighted $L^\infty$. The gain estimate of Lemma~\ref{pointwiseQ+est} and the Duhamel argument of Section~7 supply the additional bridge: after a further finite waiting time, the weighted norm itself enters the small-perturbation ball. The small-perturbation theory then gives the shifted exponential estimate, and hence
	\[
	\|F(t)-\mu\|_{L^\infty_{x,v}
		(\langle v\rangle^{\b_1} e^{\alpha|v|^\kappa})}
	\leq C_{F}e^{-\lambda t},
	\qquad t\geq0,
	\]
	{after enlarging the constant $C_{F}$ to cover the finite interval
	$[0,T_{ent}]$. Thus the qualitative convergence obtained above is upgraded
	to exponential relaxation in the same stretched-exponential weighted
	$L^\infty$ framework. 
    
	\medskip
	\subsubsection{The polynomially weighted case.}
	The proof of Theorem~\ref{thm:cond-p} follows essentially the same
	scheme as that of Theorem~\ref{thm:cond}. As explained in the preceding
	remark, the main difference at the initial stage is that we assume the
	uniform lower profile \eqref{ini profile-poly} instead of generating a
	near-initial Gaussian lower bound. Once this initial nondegeneracy is
	available, the weak--strong uniqueness argument, the tail-to-coercivity
	principle, the uniform polynomially weighted $L^\infty$ estimate, and
	the entropy-decay argument carry over with only minor modifications.
	
	A further difference appears in the final perturbative large-amplitude
	regime. In the polynomially weighted setting, estimates for the
	linearized Boltzmann operator naturally produce an additional
	collision-frequency factor of order $\langle v\rangle^\gamma$.
	Consequently, the same stretched-exponential Duhamel estimate cannot be
	used verbatim, and we employ a slightly different mild-formulation
	scheme adapted to the polynomial weight. This modification nevertheless
	leads to the same conclusion: once the solution has entered a sufficiently
	small neighborhood of the global Maxwellian, it converges exponentially
	fast to $\mu$ in the corresponding polynomially weighted $L^\infty$
	space.

    \medskip
    
    \subsubsection{The second-order moment criterion.} The main feature of Theorem 1.18 is that, if a global mild solution belongs to a sufficiently high-order polynomially weighted $L^\infty_{\beta_1}$ space at each time, without requiring its norm to be uniformly bounded in time, then the higher-order moment condition in Theorem 1.14 can be replaced by the second-order moment condition (1.27). Under this additional assumption on the solution class, the initial lower profile assumption (1.23) can also be removed. The compactness and tail-to-coercivity arguments of Section 4 remain applicable under the second-order moment condition, although the weak limit can initially be identified only as a global Maxwellian $m(v)$, rather than necessarily $\mu(v)$. This is sufficient for the coercivity estimates.

The key additional step is to establish a uniform-in-time bound in the original weighted $L^\infty_{\beta_1}$ space using only the second-order moment condition. We first obtain a uniform $L^\infty_5$ bound by the argument of Section 5 and then employ a sharper form of the Carleman representation to upgrade it to a uniform $L^\infty_{\beta_1}$ bound for sufficiently large $\beta_1$. This improved estimate allows us to identify $m$ with $\mu$ and establish the decay of the relative entropy. The polynomially weighted large-amplitude theory then yields exponential convergence to $\mu$ in the same weighted $L^\infty_{\beta_1}$ space.
	
	\subsection{Organization of the paper}
	
	The rest of the paper is organized as follows.
	Section~2 establishes uniqueness in the dissipative mild class and
	identifies the solution with the DiPerna--Lions approximation limit.
	Section~3 develops the short-time positivity mechanism and generates a
	Gaussian lower bound from general initial data, allowing local vacuum.
	Section~4 proves the uniform macroscopic lower bounds and establishes the
	tail-to-coercivity principle. Section~5 derives the global weighted
	$L^\infty$ estimate. Section~6 proves the decay of the relative entropy,
	and Section~7 upgrades this qualitative relaxation to exponential
	convergence toward the global Maxwellian. Appendix~C treats the
	polynomially weighted counterpart.  Finally, in Section 8, we prove the main theorems. The appendices collect the local
	well-posedness result and several auxiliary lemmas used throughout the
	paper. The convention for constants is as follows. Throughout the paper, constants denoted by $C$ may change from line to line. Structural constants depend on the fixed kernel, torus, equilibrium and the weight parameters in the estimate. After the compactness argument of Section~4, uniform-in-time constants may additionally depend on the solution-specific positive lower bounds $\varepsilon_\rho,\varepsilon_L,\lambda_F$ and the time $T_{glb}$ for a Gaussian lower bound.

\bigskip

\section{Uniqueness}
%

We start with a weak-strong uniqueness type lemma.
\begin{lemma} \label{lem_unique} \cite{Lu}
	Let $T\in(0,\infty]$, and let $f,g\in\mathcal U_T(f_0)$ be dissipative mild solutions with the same initial datum $f_0=g_0\in L^1_2$. Assume
    \Be \notag
	\iint_{\mathbb{T}^3\times \R^3} \avv^2 g(t) \leq \iint_{\mathbb{T}^3\times \R^3} \avv^2 f(t),\quad t\in[0,T),
	\Ee
	and
	\Be \label{cond_varphi}
	\sup_x \int_{\R^3}A(v-v_*)\avv^2 f(t,x,v) dv \leq \avvv^{2}\varphi(t),\quad \text{for some $\varphi \in L^1_{loc}[0,\infty)$} 
	\Ee
	Then, $f=g$ a.e. on $[0,T)\times\mathbb T^3\times\mathbb R^3$. Here, $A(v) := \int_{\mathbb{S}^2} B(v, \sigma) d\sigma$.
\end{lemma}
\begin{proof}
    We use only the stated ordering of the mass--energy integrals; conservation for the weaker solution $g$ is not required.
    We split (as in X. Lu~\cite{Lu})
	\[
	|f-g| = (g-f) + 2(f-g)^{+}
	\]
	and then
	\begin{equation*}
		\begin{split}
	\iint_{\mathbb{T}^3\times \R^3}(1+|v|^2)|f-g| &= \iint_{\mathbb{T}^3\times \R^3}(1+|v|^2)(g-f) + 2\iint_{\mathbb{T}^3\times \R^3}(1+|v|^2)(f-g)^{+} \\
	&\leq 2\iint_{\mathbb{T}^3\times \R^3}(1+|v|^2)(f-g)^{+}.
\end{split}
\end{equation*}
	
	We recall argument in \cite{Lu}.  Let $f'_* = f(v'_*)$, $f'=f(v')$, $g'_* = g(v'_*)$, and $g'=g(v')$. 
	\begin{equation*}
		\begin{split}
			[(f'f'_* - ff_*) - (g'g'_* - gg_*)]\mathbf{1}_{f\geq g} &\leq (f'f'_* - g'g'_*)\mathbf{1}_{f\geq g} - (ff_* - gg_*)^{+} \\
			&\quad + \underbrace{(ff_* - gg_*)^{+} - (ff_* - gg_*)\mathbf{1}_{f\geq g}}_{\Delta}. 
		\end{split}
	\end{equation*}
	For $\Delta$, we split into four cases :   \\
	If $f \geq g$ and $ff_* - gg_* \geq 0$, then $\Delta = (ff_* - gg_*) - (ff_* - gg_*) = 0$. \\
	If $f < g$ and $ff_* - gg_* \leq 0$, then $\Delta = 0 - 0 = 0$.  \\
	If $f \geq g$ and $ff_* - gg_* \leq 0$, then $\Delta = - ff_* + gg_* \leq - ff_* + fg_* \leq f|f_* - g_*|$. \\
	If $f < g$ and $ff_* - gg_* \geq 0$, then $\Delta = ff_* - gg_* \leq ff_* - fg_* \leq f|f_* - g_*|$. \\
	
	\noindent Hence
	\Be \label{lem Lu}
	[(f'f'_* - ff_*) - (g'g'_* - gg_*)]\mathbf{1}_{f\geq g} \leq (f'f'_* - g'g'_*)^{+} - (ff_* - gg_*)^{+} + |f||g_*- f_*|.
	\Ee
    
	Now we perform Gronwall type estimate : 
	\begin{equation*}
		\begin{split}
			&\iint_{x,v}\avv^2|f - g|(t,x,v) = \iint_{x,v}\avv^2|f - g|(t,x+vt,v) \\
			&\leq 2\iint_{x,v} \avv^2 (f-g)\mathbf{1}_{f\geq g}(t,x+vt,v) \\
			&\leq \cancel{2\iint_{x,v} \avv^2 (f_0 - g_0)}(x,v)  \\
			&\quad + 2\int_{0}^{t}\iiiint_{x,v,v_*, \sigma} \avv^2 B(v-v_*,\sigma)[(f'f'_* - ff_*) - (g'g'_* - gg_*)]\mathbf{1}_{f\geq g}(s,x+sv,v) \\
			&\leq 2\int_{0}^{t}\iiiint_{x,v,v_*,\sigma} \avv^{2} B(v-v_*,\sigma) \big[ (f'f'_* - g'g'_*)^{+} - (ff_* - gg_*)^{+} + |f||g_*- f_*| \big](s,x,v)\quad \text{using \eqref{lem Lu}}  \\
			&\leq 2\int_{0}^{t} \iiint_{x,v,v_*} \avv^{2} A(v-v_*)f(s,x,v)|f_* - g_*|(s,x,v_*)   \\
			&\leq 2\int_{0}^{t} \varphi(s) \int_{x,v_*} \avvv^{2}|f_* - g_*|(s,x,v_*) \quad \text{using \eqref{cond_varphi}.}  \\
		\end{split}
	\end{equation*}
	Using standard Gronwall inequality, we finish the proof. 
\end{proof}

\begin{proof} [Proof of uniqueness in Theorems \ref{thm:cond} and \ref{thm:cond-p}]
	If $F$ is a mild solution emanating from $F_0$ with assumption \eqref{apriori}, then $F$ is a global dissipative mild solution. In particular, it satisfies mass, momentum, and energy conservations by \eqref{apriori}. In Lemma \ref{lem_unique}, we put $f=F$ and $g=G$ for a global dissipative mild solution $G$. Since \eqref{apriori} implies
    \[  
        \sup_{t,x}\int\langle v\rangle^{2+\gamma}F(t,x,v)\,dv<\infty,
    \]
    \eqref{cond_varphi} also holds and we obtain uniqueness in the dissipative mild solution class $\mathcal{U}_{\infty}(F_0)$. Since $\emptyset\neq\mathcal{D}_{DL}(F_0)\subset \mathcal{U}_{\infty}(F_0)$ (e.g. \cite{DL, Golse02}), the unique solution $F\in \mathcal{D}_{DL}(F_0)$. So $F$ also satisfies the entropy inequality by \cite{DL2}.
\end{proof}

From now on, we focus on proving \eqref{decay F}. As written in above proof, $F$ is a unique global dissipative mild solution that satisfies entropy inequality and mass, momentum, and energy conservations. 

\bigskip


\section{Generation of Gaussian lower bound near initial time}
In this section, we prove the generation of a Gaussian lower bound near the initial time, with constants depending on the initial datum \(F_0\). Whenever possible, we formulate the argument through separate, self-contained lemmas in order to identify the minimal assumptions required for the Gaussian lower bound and to track explicitly how the resulting constants depend on the relevant physical quantities.

Let us assume the following conservation laws.
	\begin{equation}\label{assum:mass}
		\iint_{\mathbb T^3\times\mathbb R^3}
		F(t,x,v)\,dv\,dx
		=
		M,
		\qquad
		0<M<\infty.
	\end{equation}
	Second, the total momentum is finite and conserved:
	\begin{equation}\label{assum:moment}
		\iint_{\mathbb T^3\times\mathbb R^3}
		vF(t,x,v)\,dv\,dx
		=
		U,
		\qquad
	|U|<\infty.
	\end{equation}
	Third, the total energy is finite and conserved:
	\begin{equation}\label{assum:energy}
		\iint_{\mathbb T^3\times\mathbb R^3}
		|v|^2 F(t,x,v)\,dv\,dx
		=
		E,
		\qquad
		0<E<\infty.
	\end{equation}

We work on the three-dimensional flat torus \(\mathbb T^3\) and denote its
diameter by
\[
C_{\mathbb T}:=\operatorname{diam}(\mathbb T^3).
\]
All spatial variables are understood modulo \(\mathbb T^3\). For \(0\le s\le t\),
we write
\[
X(s;t,x,v):=x-v(t-s)
\]
for the backward in time trajectory of free transport at time $s$ which was at $(t,x,v)$.

We shall use the following notation for spatial balls, velocity balls, and
phase-space cylinders. For \(x_0\in\mathbb T^3\),
\(v_0\in\mathbb R^3\), and \(r>0\), we write
\[
B_x(x_0,r):=\{x\in\mathbb T^3: |x-x_0|<r\},
\qquad
B_v(v_0,r):=\{v\in\mathbb R^3: |v-v_0|<r\}.
\]
For \(z_0=(x_0,v_0)\), we define the phase-space cylinders
\[
C_r(z_0):=B_x(x_0,r)\times B_v(v_0,r),
\qquad
C_{A,r}(z_0):=B_x(x_0,Ar)\times B_v(v_0,r),
\]
where \(A>1\) in the second one.  \\

Throughout this section, we also assume there exists a unique local solution on $[0,T_{loc})$ such that
\begin{equation}\label{assum:linfty}
	\|F\|_{L^\infty([0,T_{loc})\times\mathbb T^3\times\mathbb R^3)}
	\le
	K
	<\infty.
\end{equation}
(This is true for Theorem \ref{thm:cond} from Proposition \ref{prop_loc}.) And, for collision frequency,
we have
\Be \label{assum:LF}
	L(F) := \iint_{\R^3\times \mathbb{S}^2} B(v-v_*,\sigma) F(t,x,v_*) d\sigma dv_* \leq C_{L}\avv^\g
\Ee
by \eqref{apriori} obviously.

Throughout this section, we work with a pointwise representative satisfying this
bound everywhere.
These are the only a priori bounds used in the construction below.

\subsection{Basic Inequalities}
The first estimate is the basic velocity-spreading mechanism of the hard-potential
gain operator. It says that a positive lower bound on a velocity ball produces,
through \(Q^+\), a positive lower bound on a larger velocity ball.


\begin{lemma} \label{lem:2.2}
	For any \(v_0\in\mathbb R^3\), \(r>0\), and \(\xi\in(0,1)\),
	\[
	Q^+\!\bigl(\mathbf 1_{B_v(v_0,r)},\mathbf 1_{B_v(v_0,r)}\bigr)(v)
	\ge
	c_\xi r^{3+\g}
	\mathbf 1_{B_v(v_0,r\sqrt2(1-\xi))}(v),
	\]
	where \(c_\xi:=c_0\xi^{5/2}\).
	Here \(c_0>0\) depends only on the hard-potential collision kernel.
\end{lemma}

\begin{proof}
	This is the quadratic hard-potential spreading estimate of
	Pulvirenti--Wennberg~\cite[Lemma~3.2]{PulvirentiWennberg1997} and Mouhot \cite{CPDE}, specialized to
	dimension three and to the hard-potential kernel. The factor \(r^{3+\gamma}\) follows from
	the hard-potential scaling of \(Q^+\).
\end{proof}

We will also need a stability version of the gain lower bound when a small subset
of the velocity ball is removed. The first step is to control the mixed term in
which one input is supported on the full ball and the other on a small bad set.
This is where the Carleman representation enters.


\begin{lemma}\label{lem:2.3}
	There exists a constant \(C_{C_b}>0\), depending only on the hard-potential collision
	kernel, such that the following holds. Let \(B:=B_v(v_0,r)\), and let
	\(H\subset B\) be measurable. If
	\(
	|H|\le \eta |B|
	\)
	for some \(\eta\in(0,1)\), then
	\[
	Q^+(\mathbf 1_B,\mathbf 1_H)(u)
	\le
	C_{C_b}\eta^{\frac{1+\g}{3}}r^{3+\g}.
	\]
	The estimate holds for every \(u\in\mathbb R^3\).
\end{lemma}

\begin{proof}
	We use the Carleman representation. It gives
	\[
	Q^+(\mathbf 1_B,\mathbf 1_H)(u)
	\leq
	C_b
	\int_{\mathbb R^3}
	\frac{\mathbf 1_H(z)}{|z-u|}
	\left(
	\int_{\Pi_{u,z}}  
	\frac{1}{\sqrt{|z-u|^2 + |y-u|^2}^{1-\g}}
	\mathbf 1_B(y)\,d\pi_{u,z}(y)
	\right)
	dz,
	\]
	where \(\Pi_{u,z}\) is the plane passing through \(u\) and orthogonal to \(z-u\).
	
	Every planar section of \(B_v(v_0,r)\) has area bounded by \(Cr^2\). Hence
	\[
	\int_{\Pi_{u,z}}\mathbf 1_B(y)\,d\pi_{u,z}(y)
	\le
	Cr^2.
	\]
	Therefore
	\[
	Q^+(\mathbf 1_B,\mathbf 1_H)(u)
	\le
	Cr^2\int_H\frac{dz}{|z-u|^{2-\g}}.
	\]
	
	By the Hardy-Littlewood inequality, among all measurable sets of measure \(|H|\), the integral of \(|z-u|^{-2+\gamma}\) is maximized by the ball
	centered at \(u\) having the same measure. Choose \(\rho\ge0\) so that
	\(
	|B_v(u,\rho)|=|H|.
	\)
	Then
	\[
	\int_H\frac{dz}{|z-u|^{2-\g}}
	\le
	\int_{B_v(u,\rho)}\frac{dz}{|z-u|^{2-\g}}
	=
	4\pi\int_0^\rho s^{\g} \,ds
	=
	\frac{4\pi}{1+\g} \rho^{1+\g}.
	\]
	Since \(\rho^3\lesssim |H|\), while
	\(|H|\le\eta|B|\) and \(|B|\sim r^3\), we obtain
	\[
	\int_H\frac{dz}{|z-u|^{2-\g}}
	\le
	C|H|^{\frac{1+\g}{3}}
	\le
	C\eta^{\frac{1+\g}{3}}r^{1+\g}.
	\]
	Consequently,
	\[
	Q^+(\mathbf 1_B,\mathbf 1_H)(u)
	\le
	C_{C_b}\eta^{\frac{1+\g}{3}}r^{3+\g}.
	\]
	This proves the claim.
\end{proof}

We now combine the full-ball lower bound with the mixed-term estimate to obtain
the small-hole gain estimate. This result will be applied twice in the
double-Duhamel construction.


\begin{lemma} \label{lem:2.4}
	There exist constants \(\eta_*\in(0,1)\) and \(c_B>0\), depending only on the
	hard-potential collision kernel, such that the following holds.
	
	Let \(v_0\in\mathbb R^3\), \(r>0\), and set \(B:=B_v(v_0,r)\). If \(H\subset B\)
	is measurable and
	\(
	|H|\le \eta_*|B|,
	\)
	then
	\[
	Q^+\!\left(\mathbf 1_{B\setminus H},\mathbf 1_{B\setminus H}\right)(u)
	\ge
	c_Br^{3+\g}.
	\]
	The estimate holds for every \(u\in B\).
\end{lemma}

\begin{proof}
	Choose a fixed number
	\(\xi\in(0,1-1/\sqrt2)\).
	Then \(B\subset B_v(v_0,r\sqrt2(1-\xi))\). By Lemma~\ref{lem:2.2},
	\[
	Q^+(\mathbf 1_B,\mathbf 1_B)(u)
	\ge
	c_\xi r^{3+\g}
	\]
	for every \(u\in B\).
	
	Since \(\mathbf 1_{B\setminus H}=\mathbf 1_B-\mathbf 1_H\), the bilinearity
	and symmetry of \(Q^+\), together with nonnegativity, imply
	\[
	\begin{aligned}
		Q^+(\mathbf 1_{B\setminus H},\mathbf 1_{B\setminus H})
		&=
		Q^+(\mathbf 1_B,\mathbf 1_B)
		-Q^+(\mathbf 1_B,\mathbf 1_H) - Q^+(\mathbf 1_H,\mathbf 1_B)
		+Q^+(\mathbf 1_H,\mathbf 1_H) \\
		&\ge
		Q^+(\mathbf 1_B,\mathbf 1_B)
		-Q^+(\mathbf 1_B,\mathbf 1_H) - Q^+(\mathbf 1_H,\mathbf 1_B).
	\end{aligned}
	\]
	Choose \(\eta_*\in(0,1)\) so small that
	\[
	   \big( C_{C_b} + \frac{C_b}{c_b}C_{C_b} \big) \eta_*^{\frac{1+\g}{3}}\le\frac{c_\xi}{2}.
	\]
	If \(|H|\le\eta_*|B|\), then Lemma~\ref{lem:2.3} gives
	\[
	Q^+(\mathbf 1_B,\mathbf 1_H)(u)
	\le
	C_{C_b} \eta_*^{\frac{1+\g}{3}}r^{3+\g}
	\]
	for every \(u\in B\). Also taking $\sigma\leftrightarrow -\sigma$, 
    \[
        Q^+(\mathbf 1_H,\mathbf 1_B)(u)
	\le
	   \frac{C_b}{c_b}C_{C_b} \eta_*^{\frac{1+\g}{3}}r^{3+\g}.
    \]
    Therefore
	\[
	Q^+\!\left(\mathbf 1_{B\setminus H},\mathbf 1_{B\setminus H}\right)(u)
	\ge
	\frac{c_\xi}{2}r^{3+\g}
	\]
	for every \(u\in B\). The claim follows by setting
	\(c_B:=c_\xi/2\).
\end{proof}
\begin{remark}\label{rem:smallhole}
	It is very important to note that the threshold \(\eta_*\) in Lemma~\ref{lem:2.4} may be chosen independently
	of the radius \(r\). Indeed, the hard-potential gain operator has the exact scaling
	property
	\begin{equation}\label{eq:Qscaling}
		Q^+\!\left(\mathbf 1_G,\mathbf 1_F\right)(v_0+r\hat u)
		=
		r^{3+\g}
		Q^+\!\left(
		\mathbf 1_{(G-v_0)/r},
		\mathbf 1_{(F-v_0)/r}
		\right)(\hat u)
	\end{equation}
	for measurable sets \(G,F\subset\mathbb R^3\), \(v_0\in\mathbb R^3\),
	\(r>0\), and \(\hat u\in\mathbb R^3\). This follows from the change of
	variables \(v_*=v_0+r\hat v_*\): the velocity measure contributes a factor
	\(r^3\), while the hard-potential kernel contributes an additional factor \(r^\gamma\).
	
	Consequently, both the full-ball lower bound from Lemma~\ref{lem:2.2} and the
	mixed-term estimate from Lemma~\ref{lem:2.3} scale like \(r^{3+\g}\). Hence the
	smallness requirement on the removed set depends on its relative measure
	\(|H|/|B|\), but not on the radius \(r\).
\end{remark}

\subsection{Formation of a Local Profile}
\label{sec:2.2}

In this section we construct, in arbitrarily short time, a local product-ball lower
bound. The starting point is only measure-theoretic positivity of the initial datum:
from the global bounds we find a positive-measure level set of \(f_0\), and then a
density point of that set. The goal is to turn this weak information into a
pointwise lower bound on a genuine product set in phase space.

The argument proceeds as follows. We first choose an anisotropic density cylinder
around the density point, so that the complement of the level set occupies only a
small relative measure. We then follow this bad set along free-transport
characteristics and show, by a moving-tube averaging argument, that most velocity
fibers remain almost full. The small-hole gain estimate from Lemma~\ref{lem:2.4}
then gives a pointwise gain lower bound. A second use of the same mechanism,
together with the Duhamel formula, produces the desired local product-ball profile.

For \(A>1\), \(r>0\), and \(z_*=(x_*,v_*)\), we recall the notation
\(C_{A,r}(z_*):=B_x(x_*,Ar)\times B_v(v_*,r)\).



\begin{lemma}\label{lem:2.6}
	Let \(F\) satisfy the a priori bounds
	\eqref{assum:mass}--\eqref{assum:linfty}.
	Fix \(A>1\) and \(\varepsilon>0\). Define
	\[
	R_V:=\sqrt{\frac{4E}{M}},
	\qquad
	V:=|\mathbb T^3|\,|B_v(0,R_V)|,
	\qquad
	\mu_0:=\frac{M}{2V}.
	\]
	Then there exist \(x_*\in\mathbb T^3\), \(v_*\in B_v(0,R_V)\),
	\(r\in(0,1]\), and a measurable set
	\(E\subset C_{A,r}(z_*)\), with \(z_*=(x_*,v_*)\), such that
	\begin{equation}\label{eq:local-bad-set-small}
		|E|
		\le
		\varepsilon |B_x(x_*,Ar)|\,|B_v(v_*,r)|
	\end{equation}
	and
	\begin{equation}\label{eq:initial-cylinder-lower}
		F_0(x,v)
		\ge
		\mu_0\,\mathbf 1_{C_{A,r}(z_*)\setminus E}(x,v)
	\end{equation}
	for every \((x,v)\).
	Moreover, \(B_x(x_*,Ar)\) may be chosen to lie in a single coordinate chart of
	\(\mathbb T^3\).
\end{lemma}

\begin{proof}
	With \(R_V\), \(V\), and \(\mu_0\) as in the statement, the energy identity
	gives
	\[
	\iint_{\{|v|\ge R_V\}} F_0(x,v)\,dv\,dx
	\le
	\frac{1}{R_V^2}
	\iint_{\mathbb T^3\times\mathbb R^3}
	|v|^2 F_0(x,v)\,dv\,dx
	=
	\frac{E}{R_V^2}
	=
	\frac{M}{4}.
	\]
	Hence the total mass identity implies
	\[
	\iint_{\mathbb T^3\times B_v(0,R_V)}
	F_0(x,v)\,dv\,dx
	\ge
	\frac{3M}{4}.
	\]
	Define
	\[
	\Omega_{\mu_0}
	:=
	\left\{
	(x,v)\in\mathbb T^3\times B_v(0,R_V):
	F_0(x,v)\ge \mu_0
	\right\}.
	\]
	Using the \(L^\infty\) bound, we obtain
	\[
	\begin{aligned}
		\frac{3M}{4}
		&\le
		\iint_{\mathbb T^3\times B_v(0,R_V)}
		F_0(x,v)\,dv\,dx \\
		&=
		\iint_{(\mathbb T^3\times B_v(0,R_V))\setminus\Omega_{\mu_0}}
		F_0(x,v)\,dv\,dx
		+
		\iint_{\Omega_{\mu_0}} F_0(x,v)\,dv\,dx \\
		&\le
		\mu_0V+K|\Omega_{\mu_0}| \\
		&=
		\frac{M}{2}+K|\Omega_{\mu_0}|.
	\end{aligned}
	\]
	Hence
	\[
	|\Omega_{\mu_0}|
	\ge
	\frac{M}{4K}
	>
	0.
	\]
	
	Since \(\Omega_{\mu_0}\) has positive measure, the Lebesgue density theorem gives a
	density point \(z_*=(x_*,v_*)\in\Omega_{\mu_0}\). In particular,
	\(v_*\in B_v(0,R_V)\).
	
	We now take the density limit along anisotropic cylinders. Since \(A>1\) is
	fixed, the cylinders \(C_{A,r}(z_*)\) shrink to \(z_*\) with uniformly bounded
	eccentricity as \(r\downarrow0\). Equivalently, they are comparable to
	Euclidean balls at scale \(r\), with constants depending only on \(A\).
	Hence the density-point property may be tested on these cylinders. Since
	\(z_*\) is a density point of \(\Omega_{\mu_0}\), we have
	\[
	\lim_{r\downarrow0}
	\frac{|\Omega_{\mu_0}\cap C_{A,r}(z_*)|}
	{|C_{A,r}(z_*)|}
	=
	1.
	\]
	Consequently, for \(r>0\) sufficiently small,
	\[
	|\Omega_{\mu_0}\cap C_{A,r}(z_*)|
	\ge
	(1-\varepsilon)|C_{A,r}(z_*)|.
	\]
	We may decrease \(r\), if necessary, so that \(r\le1\) and
	\(B_x(x_*,Ar)\) lies in a single coordinate chart of \(\mathbb T^3\).
	These additional restrictions do not affect the density estimate.
	
	Set \(E:=C_{A,r}(z_*)\setminus\Omega_{\mu_0}\). Then
	\[
	|E|
	\le
	\varepsilon |C_{A,r}(z_*)|
	=
	\varepsilon |B_x(x_*,Ar)|\,|B_v(v_*,r)|,
	\]
	which proves \eqref{eq:local-bad-set-small}. Moreover,
	\(C_{A,r}(z_*)\setminus E=C_{A,r}(z_*)\cap\Omega_{\mu_0}\). Hence
	\(F_0(x,v)\ge\mu_0\) on \(C_{A,r}(z_*)\setminus E\), while the
	right-hand side of \eqref{eq:initial-cylinder-lower} vanishes outside this
	set. Thus \eqref{eq:initial-cylinder-lower} holds for every \((x,v)\).
\end{proof}

The previous lemma gives an initial lower bound outside a small bad set in an
anisotropic cylinder. The next proposition transports this information to a
short time interval and converts it into a pointwise lower bound for the gain
term. The only geometric point is that, along a suitable moving spatial tube,
the free-transport trajectory remains inside the density cylinder.


\begin{proposition}\label{prop:2.7}
	Let \(t\in(0,T_{loc})\), \(\varepsilon>0\), and let
	\(a>0\), \(A>1\) satisfy \(A>a+t\).
	Assume that $\eta_*$ is given in Lemma~\ref{lem:2.4}, and \(x_*,v_*,r,\mu_0\), and \(E\) are given by
	Lemma~\ref{lem:2.6} with these choices of \(A\) and \(\varepsilon\).
	For \(s\in[0,t]\), set \(D_s:=B_x(x_*+sv_*,ar)\), and let
	\(\nu:=C_{L}(2+R_V)^{\g}\).
	
	Then, for each \(s\in[0,t]\), there exists a measurable set
	\(X_{\mathrm{bad}}(s)\subset D_s\) such that
	\begin{equation}\label{eq:Xbad-est}
		|X_{\mathrm{bad}}(s)|
		\le
		\frac{\varepsilon}{\eta_*}|B_x(x_*,Ar)|
	\end{equation}
	and, whenever \(y\in D_s\setminus X_{\mathrm{bad}}(s)\),
	\begin{equation}\label{eq:first-gain-lower}
		Q^+(F,F)(s,y,u)
		\ge
		\mu_1r^{3+\g}
	\end{equation}
	for every \(u\in B_v(v_*,r)\), where
	\(\mu_1:=e^{-2\nu t}\mu_0^2c_B\).
\end{proposition}

\begin{proof}
	We first verify the moving-tube inclusion. If \(s\in[0,t]\), \(y\in D_s\),
	and \(w\in B_v(v_*,r)\), then
	\[
	\begin{aligned}
		|y-sw-x_*|
		&\le
		|y-(x_*+sv_*)|+s|v_*-w| \\
		&<
		ar+tr \\
		&<
		Ar.
	\end{aligned}
	\]
	Thus
	\begin{equation}\label{eq:footpoint-in-cylinder}
		y-sw\in B_x(x_*,Ar).
	\end{equation}
	
	For \(s\in[0,t]\) and \(y\in D_s\), define
	\[
	E_{s,y}
	:=
	\left\{
	w\in B_v(v_*,r):
	(y-sw,w)\in E
	\right\},
	\]
	and set
	\[
	X_{\mathrm{bad}}(s)
	:=
	\left\{
	y\in D_s:
	|E_{s,y}|>\eta_*|B_v(v_*,r)|
	\right\}.
	\]
	
	The average size of these fibers is controlled by the size of \(E\). Indeed,
	Tonelli's theorem gives
	\[
	\begin{aligned}
		\int_{D_s}|E_{s,y}|\,dy
		&=
		\int_{B_v(v_*,r)}
		\int_{D_s}
		\mathbf 1_E(y-sw,w)\,dy\,dw.
	\end{aligned}
	\]
	For fixed \(s\) and \(w\), the change of variables \(x'=y-sw\) gives
	\[
	\int_{D_s}\mathbf 1_E(y-sw,w)\,dy
	=
	\int_{D_s-sw}\mathbf 1_E(x',w)\,dx'.
	\]
	By \eqref{eq:footpoint-in-cylinder},
	\(D_s-sw\subset B_x(x_*,Ar)\). Hence
	\[
	\int_{D_s}\mathbf 1_E(y-sw,w)\,dy
	\le
	\int_{B_x(x_*,Ar)}\mathbf 1_E(x',w)\,dx'.
	\]
	Since \(E\subset B_x(x_*,Ar)\times B_v(v_*,r)\), integration in \(w\)
	yields
	\begin{equation}\label{eq:fiber-average-bound}
		\int_{D_s}|E_{s,y}|\,dy
		\le
		|E|.
	\end{equation}
	
	By Markov's inequality, \eqref{eq:fiber-average-bound}, and
	\eqref{eq:local-bad-set-small},
	\[
	|X_{\mathrm{bad}}(s)|
	\le
	\frac{|E|}{\eta_*|B_v(v_*,r)|}
	\le
	\frac{\varepsilon}{\eta_*}|B_x(x_*,Ar)|.
	\]
	This proves \eqref{eq:Xbad-est}.
	
	Let \(y\in D_s\setminus X_{\mathrm{bad}}(s)\). Then
	\[
	|E_{s,y}|
	\le
	\eta_*|B_v(v_*,r)|.
	\]
	For every \(w\in B_v(v_*,r)\setminus E_{s,y}\), we have
	\((y-sw,w)\in C_{A,r}(z_*)\setminus E\).
	Moreover, since \(v_*\in B_v(0,R_V)\) and \(r\le1\), we have
	\(
	|w|\le R_V+1.
	\)
	Hence \eqref{assum:LF} bounds the loss coefficient by
	\(
	\nu=C_{L}(2+R_V).
	\)
	Using the mild formulation and \eqref{eq:initial-cylinder-lower}, we obtain
	\[
	F(s,y,w)
	\ge
	e^{-\nu s}F_0(y-sw,w)
	\ge
	e^{-\nu t}\mu_0.
	\]
	Equivalently,
	\begin{equation}\label{eq:tube-fiber-lower}
		F(s,y,\cdot)
		\ge
		e^{-\nu t}\mu_0
		\mathbf 1_{B_v(v_*,r)\setminus E_{s,y}}.
	\end{equation}
	
	Applying Lemma~\ref{lem:2.4} with \(B=B_v(v_*,r)\) and \(H=E_{s,y}\),
	and using \eqref{eq:tube-fiber-lower} together with the bilinearity and
	monotonicity of \(Q^+\), we obtain
	\[
	\begin{aligned}
		Q^+(F,F)(s,y,u)
		&\ge
		e^{-2\nu t}\mu_0^2
		Q^+\!\left(
		\mathbf 1_{B_v(v_*,r)\setminus E_{s,y}},
		\mathbf 1_{B_v(v_*,r)\setminus E_{s,y}}
		\right)(u) \\
		&\ge
		e^{-2\nu t}\mu_0^2c_Br^{3+\g} \\
		&=
		\mu_1r^{3+\g}.
	\end{aligned}
	\]
	This holds for every \(u\in B_v(v_*,r)\), and proves
	\eqref{eq:first-gain-lower}.
\end{proof}

We now insert the preceding gain lower bound into a second Duhamel step. The
main point is again an averaging argument: for most velocities \(w\), the
backward tube entering the first gain lower bound avoids the bad spatial set
\(X_{\mathrm{bad}}(s)\) for a positive amount of time.


\begin{proposition}\label{prop:2.8}
	Let \(F(t,x,v)\) be a nonnegative pointwise mild solution to the Boltzmann
	equation satisfying the a priori bounds
	\eqref{assum:mass}--\eqref{assum:LF} on \([0,T_{loc})\).
	Fix \(t>0\) with \(2t < T_{loc}\).
	Then there exist \(x_*\in\mathbb T^3\), \(v_*\in\mathbb R^3\),
	\(r>0\), \(c_x,c_v\in(0,1)\), and \(\mu_{\mathrm{sol}}>0\) such that
	\begin{equation}\label{eq:local-product-ball-profile}
		F(2t,x,v)\ge \mu_{\mathrm{sol}}
	\end{equation}
	for every
	\[
	(x,v)
	\in
	B_x(x_*+2tv_*,c_xr)\times B_v(v_*,c_vr).
	\]
	Moreover, one may take
	\begin{equation}\label{eq:musol}
		\mu_{\mathrm{sol}}
		=
		\frac{t^3}{8}
		e^{-\frac{17}{2}\nu t}
		\mu_0^4c_B^3r^{9+3\g},
	\end{equation}
	where \(\nu=C_{L}(2+R_V)^{\g}\), \(R_V\) and \(\mu_0\) are as in
	Lemma~\ref{lem:2.6}, and \(c_B\) is the constant from
	Lemma~\ref{lem:2.4}.
\end{proposition}

\begin{proof}
	The parameters are chosen successively as follows. We first fix the small-hole
	threshold \(\eta_*\) from Lemma~\ref{lem:2.4}. By
	Remark~\ref{rem:smallhole}, this threshold is independent of the radius \(r\).
	We then choose \(c_x,c_v\), followed by \(a,A\), \(\varepsilon\), and finally
	the density cylinder. This order ensures that \(A\) and \(\varepsilon\) are
	already fixed when Lemma~\ref{lem:2.6} is applied.
	
	Choose \(c_x,c_v\in(0,1)\). Then choose \(a>0\) and \(A>1\) such that
	\begin{equation}\label{eq:a-A-choice}
		a>c_x+\frac{c_vt}{2}+2t,
		\qquad
		A>a+t.
	\end{equation}
	Finally, choose \(\varepsilon>0\) sufficiently small that
	\begin{equation}\label{eq:eps-choice}
		\frac{16A^3}{t^3}\frac{\varepsilon}{\eta_*}
		\le
		\eta_*.
	\end{equation}
	
	Having fixed \(A\) and \(\varepsilon\), apply Lemma~\ref{lem:2.6}.
	This gives \(z_*=(x_*,v_*)\), \(r>0\), and
	\(E\subset C_{A,r}(z_*)\) satisfying
	\eqref{eq:local-bad-set-small} and \eqref{eq:initial-cylinder-lower} for $F_0$.
	Apply Proposition~\ref{prop:2.7} with the parameter \(t,\varepsilon,a,A\) chosen above. We obtain
	measurable sets
	\(X_{\mathrm{bad}}(s)\subset D_s\), where
	\(D_s:=B_x(x_*+sv_*,ar)\), such that
	\begin{equation}\label{eq:Xbad-bound}
		|X_{\mathrm{bad}}(s)|
		\le
		\frac{\varepsilon}{\eta_*}|B_x(x_*,Ar)|
	\end{equation}
	for every \(s\in[0,t]\), and
	\begin{equation}\label{eq:first-gain-lower-used}
		Q^+(F,F)(s,y,u)\ge \mu_1 r^{3+\g}
	\end{equation}
	whenever \(y\in D_s\setminus X_{\mathrm{bad}}(s)\) and
	\(u\in B_v(v_*,r)\), where
	\(\mu_1:=e^{-2\nu t}\mu_0^2c_B\).
	
	Fix \(x\in B_x(x_*+2tv_*,c_xr)\),
	\(v\in B_v(v_*,c_vr)\), and \(\tau\in[3t/2,2t]\).
	For \(s\in[0,t]\) and \(w\in B_v(v_*,r)\), define
	\[
	y(\tau):=x-v(2t-\tau),
	\qquad
	z(s;w):=y(\tau)-w(\tau-s).
	\]
	Since \(v_*\in B_v(0,R_V)\), \(r\le1\), and \(c_v<1\), we have
	\(|w|\le R_V+1\) and \(|v|\le R_V+1\). Hence
	\eqref{assum:LF} bounds all loss coefficients used below by
	\(\nu=C_{L}(2+R_V)^{\g}\).
	
	We first check that \(z(s;w)\in D_s\).
	Since \(\tau\in[3t/2,2t]\) and \(s\in[0,t]\), we have
	\(2t-\tau\le t/2\) and \(\tau-s\le2t\). Therefore
	\[
	\begin{aligned}
		|z(s;w)-(x_*+sv_*)|
		&\le
		|x-(x_*+2tv_*)|
		+|v-v_*|(2t-\tau)
		+|w-v_*|(\tau-s) \\
		&<
		c_xr+\frac{c_vt}{2}r+2tr \\
		&=
		\left(c_x+\frac{c_vt}{2}+2t\right)r \\
		&<
		ar,
	\end{aligned}
	\]
	where the last inequality follows from \eqref{eq:a-A-choice}. Thus
	\(z(s;w)\in D_s\).

	For the fixed \(x,v,\tau\), define
	\[
	E_{\mathrm{vel}}(\tau)
	:=
	\left\{
	w\in B_v(v_*,r):
	\int_0^t
	\mathbf 1_{X_{\mathrm{bad}}(s)}(z(s;w))\,ds
	>
	\frac{t}{2}
	\right\}.
	\]
	By Markov's inequality and Tonelli's theorem,
	\[
	\begin{aligned}
		\frac{t}{2}|E_{\mathrm{vel}}(\tau)|
		&\le
		\int_0^t
		\int_{B_v(v_*,r)}
		\mathbf 1_{X_{\mathrm{bad}}(s)}
		\bigl(y(\tau)-w(\tau-s)\bigr)\,dw\,ds.
	\end{aligned}
	\]
	For fixed \(s\), the change of variables
	\(z=y(\tau)-w(\tau-s)\) gives \(dw=(\tau-s)^{-3}dz\). Since
	\(\tau-s\ge t/2\),
	\[
	\int_{B_v(v_*,r)}
	\mathbf 1_{X_{\mathrm{bad}}(s)}
	\bigl(y(\tau)-w(\tau-s)\bigr)\,dw
	\le
	\frac{8}{t^3}|X_{\mathrm{bad}}(s)|.
	\]
	Using \eqref{eq:Xbad-bound} and \eqref{eq:eps-choice}, we obtain
	\begin{equation}\label{eq:Evel-small}
		|E_{\mathrm{vel}}(\tau)|
		\le
		\frac{16A^3}{t^3}
		\frac{\varepsilon}{\eta_*}
		|B_v(v_*,r)|
		\le
		\eta_*|B_v(v_*,r)|.
	\end{equation}
	
	Let \(w\in B_v(v_*,r)\setminus E_{\mathrm{vel}}(\tau)\). Then
	\[
	I_{\mathrm{good}}(w)
	:=
	\{s\in[0,t]:z(s;w)\notin X_{\mathrm{bad}}(s)\}
	\]
	satisfies \(|I_{\mathrm{good}}(w)|\ge t/2\).
	For every \(s\in I_{\mathrm{good}}(w)\), we have
	\(z(s;w)\in D_s\setminus X_{\mathrm{bad}}(s)\), and therefore
	\[
	Q^+(F,F)(s,z(s;w),w)\ge \mu_1 r^{3+\g}.
	\]
	Using the mild formulation at \((\tau,y(\tau),w)\), we obtain
	\[
	\begin{aligned}
		F(\tau,y(\tau),w)
		&\ge
		\int_{I_{\mathrm{good}}(w)}
		e^{-\nu(\tau-s)}
		Q^+(F,F)(s,z(s;w),w)\,ds \\
		&\ge
		\frac{t}{2}e^{-2\nu t}\mu_1 r^{3+\g}.
	\end{aligned}
	\]
	Define
	\begin{equation}\label{eq:mu2-definition}
		\mu_2
		:=
		\frac{t}{2}e^{-2\nu t}\mu_1
		=
		\frac{t}{2}e^{-4\nu t}\mu_0^2c_B.
	\end{equation}
	Then
	\begin{equation}\label{eq:local-profile-intermediate}
		F(\tau,y(\tau),w)\ge \mu_2r^{3+\g}
	\end{equation}
	for every \(w\in B_v(v_*,r)\setminus E_{\mathrm{vel}}(\tau)\).
	
	We now apply Lemma~\ref{lem:2.4} with
	\(B=B_v(v_*,r)\) and \(H=E_{\mathrm{vel}}(\tau)\). By
	\eqref{eq:Evel-small}, the small-hole condition is satisfied. Using
	\eqref{eq:local-profile-intermediate}, together with the bilinearity and
	monotonicity of \(Q^+\), we obtain
	\begin{equation}\label{eq:Qfinal}
		Q^+(F,F)(\tau,y(\tau),u)
		\ge
		\mu_2^2c_Br^{9+3\g}
	\end{equation}
	for every \(u\in B_v(v_*,r)\). Since
	\(v\in B_v(v_*,c_vr)\subset B_v(v_*,r)\), we may take \(u=v\).
	
	Finally, using \eqref{eq:Qfinal} in the mild formulation at \((2t,x,v)\), we
	obtain
	\[
	\begin{aligned}
		F(2t,x,v)
		&\ge
		\int_{3t/2}^{2t}
		e^{-\nu(2t-\tau)}
		Q^+(F,F)(\tau,y(\tau),v)\,d\tau \\
		&\ge
		e^{-\nu t/2}\frac{t}{2}\mu_2^2c_Br^{9+3\g}.
	\end{aligned}
	\]
	Set
	\(
	\mu_{\mathrm{sol}}
	:=
	e^{-\nu t/2}\frac{t}{2}\mu_2^2c_Br^{9+3\g}.
	\)
	By \eqref{eq:mu2-definition}, this is precisely the quantity stated in
	\eqref{eq:musol}. Since the terminal point \((x,v)\) was arbitrary in the
	target product ball, this proves \eqref{eq:local-product-ball-profile}.
\end{proof}

\begin{remark}[Measurability of the bad sets]
	Since \(E\) is measurable, Fubini--Tonelli implies that
	\[
	(s,y)
	\longmapsto
	|E_{s,y}|
	=
	\int_{B_v(v_*,r)}
	\mathbf 1_E(y-sw,w)\,dw
	\]
	is measurable. Moreover, the set
	\(
	\{(s,y)\in[0,t]\times\mathbb T^3:y\in D_s\}
	\)
	is measurable, since \(D_s=B_x(x_*+sv_*,ar)\). Hence the space-time set
	\(
	\{(s,y):s\in[0,t],\ y\in X_{\mathrm{bad}}(s)\}
	\)
	is measurable. Since \((s,w)\mapsto (s,z(s;w))\) is continuous, it follows that
	\(E_{\mathrm{vel}}(\tau)\) is measurable.%
\end{remark}


Proposition~\ref{prop:2.8} gives a pointwise product-ball lower bound. Setting
\[
t_0:=2t,
\qquad
x_0:=x_*+2tv_*,
\qquad
v_0:=v_*,
\]
and
\[
r_x:=c_xr,
\qquad
r_0:=c_vr,
\qquad
\kappa_0:=\mu_{\mathrm{sol}},
\]
we may rewrite its conclusion as
\[
F(t_0,x,v)\ge \kappa_0
\]
for every \((x,v)\in B_x(x_0,r_x)\times B_v(v_0,r_0)\). This gives a local profile which will be used in the next subsection.


\subsection{Spreading argument and spatial mixing}

\paragraph{\textbf{Local velocity spreading}}

We begin with a solid local bump. Assume that, for some time \(t_0\in[0,T_{loc})\),
there exist
\[
x_0\in\mathbb T^3,\qquad
v_0\in\mathbb R^3,\qquad
r_x>0,\qquad
r_0>0,\qquad
\kappa_0>0
\]
such that
\begin{equation}\label{P1-initial-seed}
	F(t_0,x,v)\ge \kappa_0
\end{equation}
for every \(x\in B_x(x_0,r_x)\) and every \(v\in B_v(v_0,r_0)\).

The purpose of the first step is to enlarge the velocity ball by repeated
applications of the spreading estimate in Lemma~\ref{lem:2.2}. At each step the
velocity radius increases by a fixed factor larger than one. The price is that the
spatial ball must shrink, because backward characteristics with larger velocities
can travel farther during the short time interval used in the Duhamel formula.
The time steps below are chosen precisely so that these backward characteristics
remain inside the spatial region where the previous lower bound is available.

Fix
\(\xi\in(0,1-1/\sqrt2),\)
and define
\begin{equation}\label{P1-lambda}
	\lambda:=\sqrt2(1-\xi)>1,
	\qquad
	r_n:=\lambda^n r_0,
	\qquad n\ge0.
\end{equation}
Let \((\Delta_n)_{n\ge0}\) be a sequence of positive numbers satisfying
\begin{equation}\label{P1-time-choice}
	(|v_0|+r_{n+1})\,\Delta_n
	\le
	\frac{r_x}{2^n}-\frac{r_x}{2^{n+1}}
	=
	\frac{r_x}{2^{n+1}}
\end{equation}
for every \(n\ge0\). We set
\[
t_n:=t_0+\sum_{k=0}^{n-1}\Delta_k,
\qquad n\ge1.
\]


\begin{proposition}\label{prop:2.10}
	Assume \eqref{P1-initial-seed}. Then, for every integer \(n\ge0\) such that
	\(t_n < T_{loc}\), there exists a constant \(\kappa_n>0\) such that
	\[
		F(t_n,x,v)\ge \kappa_n
	\]
	for every
	\(x\in B_x(x_0,r_x/2^n)\) and \(v\in B_v(v_0,r_n)\).
	Moreover, whenever \(t_{n+1} < T_{loc}\), the constants may be chosen recursively by
	\[
	\kappa_{n+1}
	=
	c_\xi r_n^{3+\g} \Delta_n
	e^{-3\nu_n\Delta_n}
	\kappa_n^2,
	\qquad n\ge0,
	\]
	where
	\(
	\nu_n
	:=  C_{K,A_0}\bigl(1+|v_0|+r_{n+1}\bigr),
	\)
	and \(c_\xi>0\) is the spreading constant from Lemma~\ref{lem:2.2}.
\end{proposition}

\begin{proof}
	We argue by induction on \(n\). For \(n=0\), the conclusion is exactly
	\eqref{P1-initial-seed}. Assume that, for some \(n\ge0\) with
	\(t_{n+1}\le T\),
	\begin{equation}\label{eq:local-spread-induction}
		F(t_n,x,v)\ge \kappa_n
	\end{equation}
	whenever
	\(x\in B_x(x_0,r_x/2^n)\) and \(v\in B_v(v_0,r_n)\).
	
	Fix
	\(x\in B_x(x_0,r_x/2^{n+1})\), \(v\in B_v(v_0,r_{n+1})\).
	For \(s\in[t_n,t_{n+1}]\), set
	\[
	y:=X(s;t_{n+1},x,v)
	=
	x-v(t_{n+1}-s).
	\]
	
	We first transfer the induction hypothesis from time \(t_n\) to the
	intermediate time \(s\). Let \(w\in B_v(v_0,r_n)\). Since
	\(r_n\le r_{n+1}\), we have
	\[
	\begin{aligned}
		|y-w(s-t_n)-x_0|
		&\le
		|x-x_0|
		+|v|(t_{n+1}-s)
		+|w|(s-t_n) \\
		&\le
		\frac{r_x}{2^{n+1}}
		+\bigl(|v_0|+r_{n+1}\bigr)
		\bigl((t_{n+1}-s)+(s-t_n)\bigr) \\
		&=
		\frac{r_x}{2^{n+1}}
		+\bigl(|v_0|+r_{n+1}\bigr)\Delta_n \\
		&\le
		\frac{r_x}{2^n},
	\end{aligned}
	\]
	where the last inequality follows from \eqref{P1-time-choice}. Hence
	\(y-w(s-t_n)\in B_x(x_0,r_x/2^n)\).
	Moreover,
	\(|w|\le |v_0|+r_n\le |v_0|+r_{n+1}\).
	By \eqref{assum:LF}, the loss coefficient along the corresponding
	characteristic is therefore bounded by \(\nu_n\). Using the mild formulation
	between \(t_n\) and \(s\), discarding the nonnegative gain contribution, and
	applying the induction hypothesis, we obtain
	\[
	\begin{aligned}
		F(s,y,w)
		&\ge
		\exp\!\left(
		-\int_{t_n}^s
		L[F](\tau,y-w(s-\tau),w)\,d\tau
		\right)
		F(t_n,y-w(s-t_n),w) \\
		&\ge
		e^{-\nu_n(s-t_n)}\kappa_n \\
		&\ge
		e^{-\nu_n\Delta_n}\kappa_n.
	\end{aligned}
	\]
	Equivalently,
	\begin{equation}\label{eq:local-spread-intermediate}
		F(s,y,\cdot)
		\ge
		e^{-\nu_n\Delta_n}\kappa_n\,
		\mathbf 1_{B_v(v_0,r_n)}.
	\end{equation}
	
	We now use the gain part of the mild formulation between \(t_n\) and
	\(t_{n+1}\). Since
	\(
	|v|
	\le
	|v_0|+r_{n+1},
	\)
	\eqref{assum:LF} also bounds the outer loss factor by \(\nu_n\). Therefore
	\[
		F(t_{n+1},x,v)
	\ge
	e^{-\nu_n\Delta_n}
	\int_{t_n}^{t_{n+1}}
	Q^+(F,F)\bigl(s,y,v\bigr)\,ds.
	\]
	By \eqref{eq:local-spread-intermediate}, together with the bilinearity and
	monotonicity of \(Q^+\),
	\[
	\begin{aligned}
		Q^+(F,F)\bigl(s,y,v\bigr)
		&\ge
		e^{-2\nu_n\Delta_n}\kappa_n^2
		Q^+\!\left(
		\mathbf 1_{B_v(v_0,r_n)},
		\mathbf 1_{B_v(v_0,r_n)}
		\right)(v).
	\end{aligned}
	\]
	Since \(v\in B_v(v_0,r_{n+1})\) and \(r_{n+1}=\lambda r_n\),
	Lemma~\ref{lem:2.2} gives
	\[
	Q^+\!\left(
	\mathbf 1_{B_v(v_0,r_n)},
	\mathbf 1_{B_v(v_0,r_n)}
	\right)(v)
	\ge
	c_\xi r_n^{3+\g}.
	\]
	Consequently,
	\[
		F(t_{n+1},x,v)
	\ge
	c_\xi r_n^{3+\g}\Delta_n
	e^{-3\nu_n\Delta_n}
	\kappa_n^2.
	\]
	Thus the claim holds at step \(n+1\) with
	\(
	\kappa_{n+1}
	:=
	c_\xi r_n^{3+\g}\Delta_n
	e^{-3\nu_n\Delta_n}
	\kappa_n^2.
	\)
	This closes the induction.
\end{proof}

\begin{remark}
	A convenient choice of the time steps is
	\[
	\Delta_n
	:=
	\vartheta\frac{r_x}{2^{n+1}(|v_0|+r_{n+1})},
	\qquad 0<\vartheta\le1.
	\]
	Then \eqref{P1-time-choice} is automatically satisfied. Moreover,
	\[
	\sum_{n=0}^{\infty}\Delta_n
	\lesssim
	\vartheta r_x
	\sum_{n=0}^{\infty}
	\frac{2^{-n}}{|v_0|+\lambda^{n+1}r_0}
	<\infty.
	\]
	In particular, the total time required for the local velocity spreading can be
	made arbitrarily small by choosing \(\vartheta>0\) sufficiently small.
	
\end{remark}


The preceding spreading step gives a lower bound on a large velocity ball, but it
is still localized near the original spatial point \(x_0\). The purpose of the
spatial mixing step is to use free transport to move this information throughout
the torus. Since the torus has finite diameter, a sufficiently large velocity ball
contains velocities that connect the original spatial region to any prescribed
target point in a fixed time.

Fix \(T_{\mathrm{mix}}>0\) so that \(t_0+T_{\mathrm{mix}} < T_{loc}\), and choose \(N\)
large enough that
\begin{equation}\label{eq:mix-radius-choice}
	r_N>\frac{C_{\mathbb T}}{T_{\mathrm{mix}}}.
\end{equation}
Since \(N\) is now fixed, by choosing the parameter \(\vartheta\) in the admissible
time steps \(\Delta_n\) sufficiently small, we may further arrange that
\(t_N+T_{\mathrm{mix}} < T_{loc}\). Then Proposition~\ref{prop:2.10} gives
\begin{equation}\label{eq:mix-input}
	F(t_N,x,v)\ge \kappa_N
\end{equation}
for every
\(x\in B_x(x_0,r_x/2^N)\) and \(v\in B_v(v_0,r_N)\).

The next proposition transports this lower bound to a neighborhood of each
spatial point. The velocity center required to reach a point \(\zeta\) may depend
on \(\zeta\), and this dependence will be removed in the recentering step.


\begin{proposition}\label{prop:2.12}
	Assume that \eqref{eq:mix-input} holds for every
	\(x\in B_x(x_0,r_x/2^N)\) and every \(v\in B_v(v_0,r_N)\), that
	\eqref{eq:mix-radius-choice} holds, and that
	\(t_N+T_{\mathrm{mix}} < T_{loc}\). Then there
	exist constants \(\delta_0>0\), \(\rho_0>0\), and \(\beta_0>0\) such that, for
	every \(\zeta\in\mathbb T^3\), there exists
	\(v_\zeta\in B_v(v_0,r_N)\) for which
	\begin{equation}\label{eq:mix-output}
		F(t_N+T_{\mathrm{mix}},x,v)\ge \beta_0
	\end{equation}
	whenever \(x\in B_x(\zeta,\delta_0)\) and \(v\in B_v(v_\zeta,\rho_0)\). Hence, we have uniform density lower bound at time $t_N +T_{\mathrm{mix}}$ : 
	\[
		\rho(t_N +T_{\mathrm{mix}}, x) \geq \frac{4}{3}\pi \rho_0^3 \beta_0 > 0.
	\]
\end{proposition}

\begin{proof}
	Fix \(\zeta\in\mathbb T^3\). We choose \(v_\zeta\) by perturbing \(v_0\)
	so that free transport carries \(x_0\) to \(\zeta\) in time
	\(T_{\mathrm{mix}}\).
	By the definition of \(C_{\mathbb T}\), the element
	\(\zeta-x_0-T_{\mathrm{mix}}v_0\in\mathbb T^3\) admits a representative
	\(d_\zeta\in\mathbb R^3\) with \(|d_\zeta|\le C_{\mathbb T}\). Set
	\(
	v_\zeta:=v_0+\frac{d_\zeta}{T_{\mathrm{mix}}}.
	\)
	Then, on \(\mathbb T^3\),
	\begin{equation}\label{eq:mix-center-relation}
		\zeta=x_0+T_{\mathrm{mix}}v_\zeta,
	\end{equation}
	and
	\(|v_\zeta-v_0|\le C_{\mathbb T}/T_{\mathrm{mix}}\).
	
	Choose
	\[
	0<\rho_0<
	\min\left\{
	r_N-\frac{C_{\mathbb T}}{T_{\mathrm{mix}}},
	\frac{r_x}{2^NT_{\mathrm{mix}}}
	\right\},
	\]
	which is possible by \eqref{eq:mix-radius-choice}. Then choose
	\(
	0<\delta_0<\frac{r_x}{2^N}-T_{\mathrm{mix}}\rho_0.
	\)
	These choices imply
	\begin{equation}\label{eq:mix-velocity-margin}
		\frac{C_{\mathbb T}}{T_{\mathrm{mix}}}+\rho_0\le r_N
	\end{equation}
	and
	\begin{equation}\label{eq:mix-spatial-margin}
		\delta_0+T_{\mathrm{mix}}\rho_0\le \frac{r_x}{2^N}.
	\end{equation}
	The first inequality keeps the velocities inside \(B_v(v_0,r_N)\), while the
	second keeps the backward in time position inside
	\(B_x(x_0,r_x/2^N)\).
	
	Fix \(x\in B_x(\zeta,\delta_0)\) and \(v\in B_v(v_\zeta,\rho_0)\). By
	\eqref{eq:mix-velocity-margin},
	\[
	|v-v_0|
	\le
	|v-v_\zeta|+|v_\zeta-v_0|
	<
	\rho_0+\frac{C_{\mathbb T}}{T_{\mathrm{mix}}}
	\le r_N.
	\]
	Thus \(v\in B_v(v_0,r_N)\). Moreover, using \eqref{eq:mix-center-relation} and
	\eqref{eq:mix-spatial-margin},
	\[
	\begin{aligned}
		|x-T_{\mathrm{mix}}v-x_0|
		&\le
		|x-\zeta|+T_{\mathrm{mix}}|v-v_\zeta|  \\
		&<
		\delta_0+T_{\mathrm{mix}}\rho_0
		\le
		\frac{r_x}{2^N}.
	\end{aligned}
	\]
	Therefore \(x-T_{\mathrm{mix}}v\in B_x(x_0,r_x/2^N)\), and
	\eqref{eq:mix-input} gives
	\[
		F(t_N,x-T_{\mathrm{mix}}v,v)\ge\kappa_N.
	\]
	
	Set
	\(
	\nu_{\mathrm{mix}}
	:=
	C_{K,A_0}\bigl(1+|v_0|+r_N\bigr)^{\g}.
	\)
	By \eqref{assum:LF}, this is a uniform upper bound for the loss
	coefficient on
	\(
	|v|\le |v_0|+r_N.
	\)
	Using the mild formulation between \(t_N\) and
	\(t_N+T_{\mathrm{mix}}\), we obtain
	\[
	\begin{aligned}
		F(t_N+T_{\mathrm{mix}},x,v)
		&\ge
		e^{-\nu_{\mathrm{mix}}T_{\mathrm{mix}}}
		F(t_N,x-T_{\mathrm{mix}}v,v) \\
		&\ge
		e^{-\nu_{\mathrm{mix}}T_{\mathrm{mix}}}\kappa_N.
	\end{aligned}
	\]
	The proof is complete by setting
	\(
	\beta_0:=e^{-\nu_{\mathrm{mix}}T_{\mathrm{mix}}}\kappa_N.
	\)
\end{proof}


Proposition~\ref{prop:2.12}, together with \eqref{assum:linfty}, yields a spatially uniform positive lower bound on the density at $t=t_N+T_{\mathrm{mix}}<T_{\mathrm{loc}}$, where $T_{\mathrm{loc}}$ denotes the local existence time given by Proposition~\ref{prop_loc}. Proposition~\ref{prop_loc} further ensures the existence of a possibly different spatially uniform positive lower bound on the density over a sufficiently short time interval around $t_N+T_{\mathrm{mix}}$.

Moreover, the a priori assumption \eqref{apriori} provides the uniform upper bounds on the local density and the local energy required in Theorem 1.1 of \cite{CPDE}. The bound \eqref{assum:linfty} also prevents concentration on $[0,T_{\mathrm{loc}})$ and allows us to control the local entropy density, thereby ensuring the remaining assumptions of Theorem 1.1 of \cite{CPDE}. Alternatively, the required local entropy bound can be established directly from Proposition~\ref{prop_loc}. Consequently, we obtain a Gaussian lower bound at some time $T_{\mathrm{glb}}\in(0,T_{\mathrm{loc}})$.


\begin{proposition} \label{prop:glb} 
    Assume a mild solution $F$ that satisfies \eqref{assum:mass} -- \eqref{assum:LF}.
    There exists $0 < T_{glb} < T_{loc}$ such that 
	\[
		F(T_{glb}, x, v) \geq e^{-a_0|v|^2+c_0} \quad \text{and}\quad \|\avv^{\b}e^{\alpha|v|^\k}F(T_{glb}, x, v)\| \leq 2\|\avv^{\b}e^{\alpha|v|^\k}F_0(x, v)\| ,
	\]
	where constants $a_0, c_0$ depend on $F$, $A_0$, and $T_{glb}$. 
\end{proposition}

Thus, in the setting of Theorem \ref{thm:cond}, a Gaussian lower bound is established after a sufficiently short positive time, since \eqref{assum:linfty} and \eqref{assum:LF} follow directly from \eqref{initial}, Proposition~\ref{prop_loc}, and \eqref{apriori}. 

By Proposition \ref{prop:glb}, the entropy inequality \eqref{entropy ineq}, and Proposition \ref{prop_loc}, we may shift the initial time to \(T_{\rm glb}\). Hence, throughout the subsequent sections, we assume without loss of generality that \(F_0(x,v)\) admits a spatially uniform Gaussian lower bound.

\bigskip

\section{Macroscopic lower bounds}

\subsection{Velocity Averaging, Rigidity, and Gain Compactness}

Recall that the solution $F$ selected in Section~2 is a global mild solution of the Boltzmann equation which satisfies the physical conservation laws and the entropy inequality. Although these properties were established under \eqref{apriori} or \eqref{apriori-poly} with $k\geq 2$, if they are assumed independently, the arguments in this section require only the weaker a priori bound
\begin{equation}\label{cond L1 vel avg}
\sup_{t,x}\int_{\mathbb R^3}
\langle v\rangle^\gamma F(t,x,v)\,dv\leq C_1.
\end{equation}


Set
\begin{equation} \label{set Fn}
    F_n(t,x,v):=F(t_n+t,x,v),\quad (t,x,v)\in (0,1)\times \T^3\times \R^3.
\end{equation}
Under assumption of \eqref{cond L1 vel avg}, we apply velocity averaging lemma which is known in \cite{DL, GLPS}.
\begin{lemma}\label{lem L1 vel avg}
	Assume $F$ is a global mild solution of \eqref{BE} which satisfies the physical conservation laws, the entropy inequality, and \eqref{cond L1 vel avg}. Then, for any sequence $t_n \rightarrow \infty$ and every fixed
	$\psi\in L^\infty((0,1)\times\mathbb T^3\times\mathbb R^3)$,
	the sequence
	\[
	\left\{
	\int_{\mathbb R^3}F_n(t,x,v)\psi(t,x,v)\,dv
	\right\}_{n\geq1}
	\]
	is relatively compact in $L^1((0,1)\times\mathbb T^3)$.
\end{lemma}

\begin{proof}
    We apply the
	$L^1$ velocity-averaging lemma of \cite{DL}.
	To this end, we first establish the relative weak compactness
	of $\{F_n\}$ in $L^1((0,1)\times\mathbb T^3\times\mathbb R^3)$ and of
	$\{Q_\pm(F_n)\}$ in $L^1((0,1)\times\mathbb T^3\times B_{\bar R})$
	for every $\bar R>0$.
	
	The mass and energy bounds, together with the entropy inequality,
	imply that
	\[
	\mathcal E:=
	\sup_n \sup_{0<t<1}
	\int_{\mathbb T^3\times\mathbb R^3}
	\langle v\rangle^2F_n(t,x,v)\,dx\,dv<\infty
	\]
	and
	\[
	\mathcal H:=
	\sup_n\int_{(0,1)\times\mathbb T^3\times\R^3}
	F_n\log^+F_n\,dv\,dt\,dx<\infty.
	\]
	Hence, $\{F_n\}$ is uniformly integrable. Moreover,
	\[
	\sup_n\int_{(0,1)\times\mathbb T^3}\int_{|v|>R}F_n\,dv\,dt\,dx
	\leq \frac{\mathcal E}{1+R^2},
	\]
	which gives uniform tightness in velocity.
	The Dunford--Pettis theorem therefore yields relative weak
	compactness in $L^1((0,1)\times\mathbb T^3\times\mathbb R^3)$.
	
	By \eqref{cond L1 vel avg},
	\[
	m_{\gamma,n}(t,x):=
	\int_{\mathbb R^3}\langle v\rangle^\gamma
	F_n(t,x,v)\,dv\leq C_1
	\]
	uniformly in $n$, $t$, and $x$. In particular,
	$\rho_n(t,x)\leq C_1$.
	Using the kernel bound and pre-/post-collisional symmetry, we obtain
	\[
	\begin{aligned}
		\|Q_\pm(F_n)\|_{L^1(D\times\mathbb R^3)}
		&\leq C\int_{(0,1)\times\mathbb T^3}
		\iint_{\mathbb R^3\times\mathbb R^3}
		\langle v\rangle^\gamma\langle v_*\rangle^\gamma
		F_nF_{n*}\,dv_*\,dv\,dt\,dx\\
		&\leq C C_1^2.
	\end{aligned}
	\]
	It remains to establish uniform integrability on bounded
	velocity sets.
	
	For the loss term, if $|v|\leq\bar R$, then
	\[
	L(F_n)(t,x,v)
	\leq C\langle v\rangle^\gamma
	m_{\gamma,n}(t,x)
	\leq C_{\bar R}C_1.
	\]
	Consequently,
	\[
	0\leq Q_-(F_n)\leq C_{\bar R}C_1F_n
	\qquad\text{on } (0,1)\times\mathbb T^3\times B_{\bar R},
	\]
	and the uniform integrability of $\{Q_-(F_n)\}$ follows
	from that of $\{F_n\}$.
	
	We next consider the gain term. Write
	\[
	P_n:=F_nF_{n*},
	\qquad
	P_n':=F_n'F_{n*}'.
	\]
	Since
	\[
	P_n\log^+P_n
	\leq F_nF_{n*}
	\bigl(\log^+F_n+\log^+F_{n*}\bigr),
	\]
	we have
	\begin{equation}\label{UI}
		\begin{aligned}
			&\int_{(0,1)\times\mathbb T^3}
			\iiint_{\mathbb R^3\times\mathbb R^3\times\mathbb S^2}
			P_n\log^+P_n\,d\sigma\,dv_*\,dv\,dt\,dx\\
			&\qquad\leq
			2|\mathbb S^2|\int_{(0,1)\times\mathbb T^3}
			\rho_n(t,x)
			\int_{\mathbb R^3}F_n\log^+F_n\,dv\,dt\,dx\\
			&\qquad\leq 2|\mathbb S^2|C_1\mathcal H.
		\end{aligned}
	\end{equation}
	By pre-/post-collisional symmetry in the full collision variables,
	the same bound holds with $P_n$ replaced by $P_n'$.
	
	For $M\geq1$, decompose
	\[
	Q_+(F_n)=Q^{\mathrm{tail}}_{+,n,M}
	+Q^{\mathrm{core}}_{+,n,M},
	\]
	where the two terms correspond to the regions $|v_*|>M$
	and $|v_*|\leq M$, respectively.
	
	For the tail term, conservation of kinetic energy gives
	\[
	|v'|^2+|v_*'|^2
	=|v|^2+|v_*|^2>M^2
	\qquad\text{when }|v_*|>M.
	\]
	Thus, at least one of $|v'|$ and $|v_*'|$ exceeds $M/\sqrt2$.
	Using pre-/post-collisional symmetry and
	\eqref{cond L1 vel avg}, we obtain
	\begin{equation}\label{UI tail}
		\begin{aligned}
			\|Q^{\mathrm{tail}}_{+,n,M}\|_{L^1((0,1)\times\mathbb T^3\times B_{\bar R})}
			&\leq C
			\int_{(0,1)\times\mathbb T^3} \iint_{\mathbb R^3\times\mathbb R^3}
			\langle v\rangle^\gamma\langle v_*\rangle^\gamma
			F_nF_{n*}\\
			&\qquad\qquad\times
			\left(
			\mathbf 1_{\{|v|>M/\sqrt2\}}
			+\mathbf 1_{\{|v_*|>M/\sqrt2\}}
			\right)
			\,dv_*\,dv\,dt\,dx\\
			&\leq C C_1
			\int_{(0,1)\times\mathbb T^3} \int_{|v|>M/\sqrt2}
			\langle v\rangle^\gamma F_n\,dv\,dt\,dx\\
			&\leq C C_1\mathcal E
			\langle M/\sqrt2 \rangle^{\gamma-2}.
		\end{aligned}
	\end{equation}
	Since $\gamma<2$, this bound tends to zero as $M\to\infty$,
	uniformly in $n$.
	
	For the core term, let $E\subset (0,1)\times\mathbb T^3\times B_{\bar R}$
	be measurable. On the region $|v|\leq\bar R$ and $|v_*|\leq M$,
	\[
	B(v-v_*,\sigma)\leq C_b(\bar R+M)^\gamma.
	\]
	For any $K>1$, splitting into the regions $P_n'\leq K$
	and $P_n'>K$ gives
	\begin{equation}\label{UI core}
		\begin{aligned}
			\int_E Q^{\mathrm{core}}_{+,n,M}\,dt\,dx\,dv
			&\leq C_b(\bar R+M)^\gamma
			\int_{E\times B_M\times\mathbb S^2}P_n'\\
			&\leq C_b(\bar R+M)^\gamma
			\left(
			K|B_M||\mathbb S^2||E|
			+\frac{C}{\log K}
			\right),
		\end{aligned}
	\end{equation}
	where $C$ is independent of $n$, $E$, and $K$.
	Here, the last term follows from \eqref{UI} and
	pre-/post-collisional symmetry.
	
	Given $\varepsilon>0$, first choose $M$ sufficiently large
	so that \eqref{UI tail} is smaller than $\varepsilon/2$.
	With this $M$ fixed, choose $K$ sufficiently large and then
	$|E|$ sufficiently small in \eqref{UI core}.
	This proves the uniform integrability of $\{Q_+(F_n)\}$
	in $L^1((0,1)\times\mathbb T^3\times B_{\bar R})$.
	Thus, both $\{Q_+(F_n)\}$ and $\{Q_-(F_n)\}$
	are relatively weakly compact in this space.
	
	Since
	\[
	(\partial_t+v\cdot\nabla_x)F_n
	=Q_+(F_n)-Q_-(F_n),
	\]
	the $L^1$ velocity-averaging lemma, applied after smooth
	localization, yields relative compactness of
	\[
	\left\{
	\int_{\mathbb R^3}F_n\varphi\,dv
	\right\}_{n\geq1}
	\]
	in $L^1_{\mathrm{loc}}((0,1)\times\mathbb T^3)$ for every
	$\varphi\in C_c^\infty((0,1)\times\mathbb T^3\times\mathbb R^3)$.
	
	For each fixed $\delta\in(0,1/2)$ and $R>0$,
	approximation in measure by uniformly bounded smooth functions,
	together with the uniform integrability of $\{F_n\}$,
	extends this compactness to bounded measurable weights
	supported in
	$[\delta,1-\delta]\times\mathbb T^3\times B_R$.
	The corresponding averages are relatively compact in $L^1((0,1)\times\mathbb T^3)$.
	
	Finally, the mass and energy bounds imply
	\[
	\begin{aligned}
		\sup_n\bigg[
		&\int_{((0,\delta)\cup(1-\delta,1))
			\times\mathbb T^3\times\mathbb R^3}F_n
			+
		\int_{(0,1)\times\mathbb T^3}\int_{|v|>R}F_n\,dv\,dt\,dx
		\bigg]
		\leq
		2\delta\mathcal E+\frac{\mathcal E}{1+R^2}.
	\end{aligned}
	\]
	Multiplying this estimate by $\|\psi\|_{L^\infty}$ controls
	the error incurred by truncating $\psi$ in time and velocity.
	Letting $\delta\downarrow0$ and $R\to\infty$ proves the
	claimed relative compactness in $L^1((0,1)\times\mathbb T^3)$.
\end{proof}

\begin{lemma} 
	\label{Qn vanish}
    Assume that the hypotheses of Lemma~\ref{lem L1 vel avg} hold. $\{F_n\}$ of Lemma \ref{lem L1 vel avg} satisfies
	\begin{equation}\label{eq:r1-Qshift-zero}
		\|Q(F_n,F_n)\|_{L^1((0,1)\times \T^3\times\mathbb R^3)}\longrightarrow0.
	\end{equation}
\end{lemma}
\begin{proof}
	
    The total entropy dissipation is finite, so
	\[
	d_n:=\int_{t_n-1}^{t_n+1}\mathcal D(F)(s)\,ds\longrightarrow0.
	\]
	Put $A_n=F_n'F_{n*}'$, $C_n=F_nF_{n*}$ and
	$d\Lambda=B\,dt\,dx\,dv\,dv_*\,d\sigma$.  The pre--post collisional
	change of variables and \eqref{cond L1 vel avg} give
	\[
	\int(A_n+C_n)\,d\Lambda\le C_{B,\gamma} C_1^2.
	\]
	Using $(a-c)\log(a/c)\ge4(\sqrt a-\sqrt c)^2$ and
	Cauchy--Schwarz, with the usual extended-value convention at zero, gives
	\begin{equation}\label{eq:r1-collision-defect}
		\begin{split}
			\int|A_n-C_n|\,d\Lambda
			&\le \left(\int(\sqrt{A_n}-\sqrt{C_n})^2\,d\Lambda\right)^{1/2}
			\left(\int(\sqrt{A_n}+\sqrt{C_n})^2\,d\Lambda\right)^{1/2}\\
			&\le C_{B,C_1}\,d_n^{1/2}\longrightarrow0.
		\end{split}
	\end{equation}
	Integration in $(v_*,\sigma)$ proves \eqref{eq:r1-Qshift-zero}.
\end{proof}

We explain rigidity of a local Maxwellian which solves free-transport equation. We also refer \cite{Desvillettes1990}.
\begin{lemma} 
	\label{lem:r1-Maxwellian-rigidity}
	Let $I$ be a nonempty open interval, and let $m\ge0$ be locally integrable on
	$I\times\mathbb T^3\times\mathbb R^3$, with finite positive total mass,
	and solve $(\partial_t+v\cdot\nabla_x)m=0$ in distributions.  Suppose
	that, for almost every $(t,x)$, $m(t,x,\cdot)$ is either a Maxwellian
	with positive finite temperature or the zero function.  Then $m$ is
	one stationary global Maxwellian.  No regularity of its macroscopic
	coefficients is required.
\end{lemma}
\begin{proof}
	The characteristic change of variables gives
	$m(t,x,v)=g(x-tv,v)$ almost everywhere.  Choose two distinct good
	times $t_0,t_1$, put $\delta=t_1-t_0\ne0$, and write $m_j=m(t_j)$.
	Then
	\begin{equation}\label{eq:r1-two-slices}
		m_1(x,v)=m_0(x-\delta v,v)\quad\hbox{for a.e. }(x,v).
	\end{equation}
	Let $E_j$ be the set of spatial points where $m_j$ is not zero.
	A nonzero Maxwellian is positive at every velocity; hence
	$\mathbf1_{E_1}(x)=\mathbf1_{E_0}(x-\delta v)$ almost everywhere.
	Average in $v$ over $\{v:\delta v\in[0,1)^3\}$.  This gives
	$\mathbf1_{E_1}(x)=|E_0|$ almost everywhere.  Positive total mass
	excludes $|E_0|=0$, so both $E_j$ have full measure.
	
	We may therefore write
	\[
	\log m_j(x,v)=a_j(x)+b_j(x)\cdot v+c_j(x)|v|^2,
	\qquad c_j(x)<0,
	\]
	with finite measurable coefficients.  Intersect the full-measure
	sets on which \eqref{eq:r1-two-slices} holds at the seven velocities
	$v$ and $v\pm e_i/\delta$, $i=1,2,3$.  All seven velocities give the
	same point $y=x-\delta v$ on the torus.  Thus the polynomial
	\[
	P(w)=a_1(x)-a_0(y)+(b_1(x)-b_0(y))\cdot w
	+(c_1(x)-c_0(y))|w|^2
	\]
	vanishes at those seven points.  Its second differences give
	$c_1(x)=c_0(y)$, its first differences then give $b_1(x)=b_0(y)$,
	and $P(v)=0$ gives $a_1(x)=a_0(y)$.
	Restricting $v$ to the same scaled fundamental cube shows that these
	coefficient vectors agree for almost every pair $(x,y)$.
	Fubini's theorem makes them a single constant vector $(a,b,c)$,
	with $c<0$.  Consequently $m_0(v)=\exp(a+b\cdot v+c|v|^2)$;
	the characteristic formula propagates this same Maxwellian to all
	times.  No integrability or differentiation of $\log m$ was used.
\end{proof}

From Lemma \ref{lem L1 vel avg}, we have a weak limit $m(t,x,v)$ such that along that same
subsequence,
\begin{equation}\label{eq:r1-average-limit}
	\int F_n\psi\,dv\longrightarrow\int m\psi\,dv
	\quad\hbox{strongly in }L^1((0,1)\times \T^3).
\end{equation}

We also introduce tightness condition for local energy.
	\begin{equation}\label{eq:r1-energy-tightness}
		\lim_{R\to\infty}\sup_{t,x}
		\int_{|v|>R}|v|^2F(t,x,v)\,dv=0.
	\end{equation}

\begin{proposition} 
\label{Desvillettes}
	Assume that the hypotheses of Lemma~\ref{lem L1 vel avg} hold. For every sequence $t_n\to\infty$ and every fixed $T>0$, there is a
	subsequence and stationary global Maxwellian $m$ such that
	\[
	F_n\rightharpoonup m
	\quad\hbox{weakly in }L^1((0,1)\times \T^3\times\mathbb R^3).
	\]
	If we also assume \eqref{eq:r1-energy-tightness}, \eqref{eq:r1-average-limit} holds with $m=\mu$.
\end{proposition}
\begin{proof}
	Take a weak limit $m$ using Lemma~\ref{lem L1 vel avg}.
	Its equation passes to free transport by \eqref{eq:r1-Qshift-zero} in Lemma \ref{Qn vanish}.
	We now identify its collision structure without assuming lower
	semicontinuity of entropy dissipation under bare weak convergence.
	
	Write $z=(t,x)$ and take a smooth compactly supported test
	$\Phi(z,v,v_*,\sigma)$ whose support is away from $v=v_*$.
	Set $V=(v+v_*)/2$, $r=|v-v_*|$, $k=(v-v_*)/r$.
	Angular integration gives the loss-product kernel
	\[
	K_C(z,v,v_*)=r^\gamma\int_{\mathbb S^2}
	b(k\cdot\sigma)\Phi(z,v,v_*,\sigma)\,d\sigma.
	\]
	For the gain product, write the collision measure in
	$(V,r,k,\sigma)$ and exchange $k$ and $\sigma$.  Both the measure
	and $b(k\cdot\sigma)$ are preserved, and the resulting kernel is
	\[
	K_A(z,v,v_*)=r^\gamma\int_{\mathbb S^2}b(k\cdot\sigma)
	\Phi\left(z,V+\frac r2\sigma,V-\frac r2\sigma,k\right)\,d\sigma.
	\]
	Both kernels are continuous and compactly supported.  Indeed,
	approximate the angular factor by smooth zonal functions in
	$L^1(\mathbb S^2)$.  The kernel error is bounded uniformly by
	$C_\Phi\|b-b_j\|_{L^1(\mathbb S^2)}$.  The support condition removes
	the direction singularity at $r=0$, and energy conservation preserves
	compact velocity support in the gain change of variables.
	
	Each kernel is uniformly approximable by finite sums
	$\eta_\ell(z)\psi_\ell(v)\widetilde\psi_\ell(v_*)$, with smooth
	compactly supported velocity factors.  By
	\eqref{eq:r1-average-limit}, the corresponding averages converge
	strongly in $L^1((0,1)\times \T^3)$ and are uniformly bounded in
	$L^\infty((0,1)\times \T^3)$ by \eqref{cond L1 vel avg}.  Their products
	therefore converge strongly in $L^1((0,1)\times \T^3)$.  The error in uniform
	kernel approximation is controlled by $C_1^2$.  We obtain
	\[
	\int\Phi A_n\,d\Lambda\longrightarrow\int\Phi m' m_*'\,d\Lambda,
	\qquad
	\int\Phi C_n\,d\Lambda\longrightarrow\int\Phi m m_*\,d\Lambda.
	\]
	Here $A_n,C_n,d\Lambda$ have the meanings used in
	\eqref{eq:r1-collision-defect}.  That estimate implies
	$B(m'm_*'-mm_*)=0$ almost everywhere.  Since $b\ge c_b>0$, detailed
	balance holds outside the null diagonal $v=v_*$.  The measurable
	equality case of the Boltzmann $H$-theorem implies that
	$m(z,\cdot)$ is a Maxwellian or zero for almost every $z$; see
	\cite{Cercignani1994}.  Its finite local mass, energy and entropy follow
	by velocity truncation, convex lower semicontinuity, and Fubini's
	theorem.  Thus this equality characterization does not presume
	strict positivity or smoothness of $m$. The total mass remains positive in the weak limit.  Apply
	Lemma~\ref{lem:r1-Maxwellian-rigidity} to conclude that $m$ is a
	stationary global Maxwellian.  
	
	If we also assume \eqref{eq:r1-energy-tightness}, to identify its parameters, first test
	weak convergence with bounded velocity truncations of $1,v,|v|^2$,
	multiplied by arbitrary time tests.  Condition
	\eqref{eq:r1-energy-tightness} makes the omitted tails tend to zero
	uniformly in $n$.  The limit therefore has the conserved mass,
	momentum and energy in \eqref{macro_conserv}, and hence equals $\mu$.
\end{proof}

Convergence in measure of $Q_{+}(F_n, F_n)$ in $(0,1)\times \T^3\times B_{\bar R}$ are wel-known in \cite{lion1}. Using  Vitali convergence lemma, it can be improved to strong $L^1$ convergence. We provide an alternate proof for this fact. 
\begin{lemma} 
\label{lem:r1-gain-strong}
	Assume that the hypotheses of Lemma~\ref{lem L1 vel avg} and \eqref{eq:r1-energy-tightness} hold. Along the subsequence in Proposition~\ref{Desvillettes}, for every
	$\bar R<\infty$,
	\[
	Q_+(F_n,F_n)\longrightarrow Q_+(\mu,\mu)
	\quad\hbox{strongly in }L^1((0,1)\times \T^3\times B_{\bar R}).
	\]
	If we do not assume \eqref{eq:r1-energy-tightness}, we obtain simlar result replacing $\mu$ to a stationary global Maxwellian $m$.
\end{lemma}
\begin{proof}
	We prove only for $\mu$. Proof for $m$ is also same. \\
	Choose $0\le\chi_M\le1$ smooth, equal to one on $B_M$ and supported
	in $B_{2M}$.  Put $G_n^M=\chi_MF_n$ and $q_n^M=Q_+(G_n^M,G_n^M)$;
	use the subscript $\mu$ for the analogous fixed functions.
	The gain weak form places the cutoff error on the input product as
	$F_n(v)F_n(v_*)[1-\chi_M(v)\chi_M(v_*)]$.
	Consequently
	\begin{equation}\label{eq:r1-gain-tail}
		\sup_n\|Q_+(F_n,F_n)-q_n^M\|_{L^1((0,1)\times \T^3\times\mathbb R^3)}
		\le C_{b,C_1}\sup_n\int_{(0,1)\times \T^3}\int_{|v|>M}
		\langle v\rangle^\gamma F_n\,dv\,dz
		\le\omega_M\longrightarrow0.
	\end{equation}
	For the uniform estimate on the right, one may use the energy bound
	and $0<\gamma\le1<2$, obtaining $\omega_M\le C_T M^{\gamma-2}$ for
	$M\ge1$.  The same estimate holds for $\mu$.
	
	For $N>e$, set $G_n^{M,N}=\min\{G_n^M,N\}$ and
	$q_n^{M,N}=Q_+(G_n^{M,N},G_n^{M,N})$.  The entropy bound gives
	\[
	\sup_n\|G_n^M-G_n^{M,N}\|_{L^1_{z,v}}
	\le \frac{C_T}{\log N}.
	\]
	For inputs supported in $B_{2M}$ the gain weak form gives
	$\|Q_+(g,h)\|_1\le C_{M,b}\|g\|_1\|h\|_1$.
	Bilinearity and the uniform local mass cap therefore imply
	\begin{equation}\label{eq:r1-gain-amplitude}
		\sup_n\|q_n^M-q_n^{M,N}\|_{L^1_{z,v}}
		\le \frac{C_{M,T,b,C_1}}{\log N}.
	\end{equation}
	
	At fixed $M,N$,
	$\|G_n^{M,N}(z)\|_2^2\le NC_1$ for almost every $z$.
	Since $b$ is bounded, it belongs to $L^2(\mathbb S^2)$.
	The gain estimate of \cite[Theorem~3.5]{MouhotVillani2004},
	in dimension three with $s=\eta=0$, gives for inputs supported in
	$B_{2M}$
	\begin{equation}\label{eq:r1-velocity-gain-regularity}
		\|Q_+(g,g)\|_{H^1(\mathbb R^3)}
		\le C_{M,b,\gamma}\bigl(\|g\|_2^2+\|g\|_1^2\bigr).
	\end{equation}
	Its kinetic-factor assumption is satisfied by $\Phi(r)=r^\gamma$.
	All $q_n^{M,N}(z,\cdot)$ thus lie in one relatively compact subset
	of $L^1(B_{\bar R})$.
	
	For completeness, take finite families of nonnegative smooth
	functions $\theta_{J,j}\in C_c^\infty(B_{\bar R})$, with support
	diameters tending uniformly to zero, such that
	$0\le\sum_j\theta_{J,j}\le1$ and this sum tends to one locally
	uniformly in the ball.  Such families follow from a fine smooth
	partition of unity and an interior cutoff.  Omitting zero functions,
	define the finite-rank operators
	\[
	P_Jh(v)=\sum_j\frac{\theta_{J,j}(v)}{\int\theta_{J,j}}
	\int\theta_{J,j}(w)h(w)\,dw.
	\]
	Then $\|P_J\|_{L^1\to L^1}\le1$.  Uniform continuity proves
	$P_Jh\to h$ first for continuous compactly supported $h$, and density
	proves it for every $h\in L^1(B_{\bar R})$.  The convergence is uniform
	on compact subsets, and hence
	\begin{equation}\label{eq:r1-gain-projection}
		\sup_n\|(I-P_J)q_n^{M,N}\|_{L^1((0,1)\times \T^3\times B_{\bar R})}
		\longrightarrow0\qquad(J\to\infty).
	\end{equation}
	
	For $\phi\in C_c^\infty(B_{\bar R})$, the weak collision form is
	\[
	\int q_n^M(z,v)\phi(v)\,dv
	=\iint G_n^M(z,v)G_n^M(z,v_*)K_\phi(v,v_*)\,dv\,dv_*,
	\]
	where
	\[
	K_\phi(v,v_*)=|v-v_*|^\gamma\int_{\mathbb S^2}
	b(k\cdot\sigma)\phi\left(\frac{v+v_*}{2}
	+\frac{|v-v_*|}{2}\sigma\right)\,d\sigma.
	\]
	Here $k=(v-v_*)/|v-v_*|$ away from $v=v_*$.
	This kernel is continuous on compact velocity sets.  Away from
	$v=v_*$, approximate $b$ by smooth zonal factors in $L^1$; the induced
	kernel error is uniform.  At the diagonal,
	$|K_\phi|\le4\pi C_b\|\phi\|_\infty|v-v_*|^\gamma\to0$.
	Uniform approximation by finite sums of separated smooth velocity
	factors, followed by \eqref{eq:r1-average-limit} with $m=\mu$, shows
	\begin{equation}\label{eq:r1-gain-coefficients}
		\int q_n^M\phi\,dv\longrightarrow\int q_\mu^M\phi\,dv
		\quad\hbox{strongly in }L^1((0,1)\times \T^3).
	\end{equation}
	Here products of averages converge because both factors are uniformly
	bounded in $L^\infty((0,1)\times \T^3)$.  Importantly, averaging is applied to
	$F_n$, not to its nonlinear amplitude truncation.
	
	With all norms on $(0,1)\times \T^3\times B_{\bar R}$,
	\begin{align*}
		\|q_n^M-q_\mu^M\|_1
		&\le2\|q_n^M-q_n^{M,N}\|_1
		+\|(I-P_J)q_n^{M,N}\|_1\\
		&\quad+\|P_J(q_n^M-q_\mu^M)\|_1
		+\|(I-P_J)q_\mu^M\|_1.
	\end{align*}
	For fixed $M$, first choose $N$ using
	\eqref{eq:r1-gain-amplitude}, then $J$ using
	\eqref{eq:r1-gain-projection} and the fixed function $q_\mu^M$,
	and finally let $n\to\infty$ using the finitely many coefficients in
	\eqref{eq:r1-gain-coefficients}.  Thus $q_n^M\to q_\mu^M$ strongly.
	Letting $M\to\infty$ in \eqref{eq:r1-gain-tail} proves the claim.
\end{proof}

\subsection{Diagonal Extraction}

We define 
\[
	u^{(a)}_{n} := \int_{|v|\leq R} |a - v|^{\g} F_n(t,x,v) dv
\]
and
\[
c_{R,a}:=\int_{|v|\le R}|a-v|^{\g} \mu(v)\,dv
\]
for some $R > 0$ and $a\in \R^3$. Using Proposition~\ref{Desvillettes} and Lemma~\ref{lem L1 vel avg}, for each fixed $a\in\mathbb R^3$, there exists a subsequence $\{n_k\}$ such that
\Be \label{un weak conv}
	u_{n_k}^{(a)} \ \rightarrow\ c_{R,a}\qquad \text{in } \quad L^1((0,1)\times\mathbb T^3_x),
\Ee
where we have used fact that $F$ is physically admissible and $m(t,x,v)=\mu$ by Proposition \ref{Desvillettes}. 


\begin{lemma} 
\label{lem:Lip}
	For all $a,b\in\mathbb R^3$,
	\[
	\|u_n^{(a)}-u_n^{(b)}\|_{L^1_{t,x}}
	\le\ |a-b|^{\g}\,M_R
	\]
	for some finite $M_R$. Consequently, the constant limits satisfy
	\[
	|c_{R,a}-c_{R,b}|
	\le |a-b|^{\g} \int_{|v|\le R}\mu(v)\,dv .
	\]
\end{lemma}
\begin{proof}
	We just check
	\Be
    \begin{split}
        \|u_n^{(a)}-u_n^{(b)}\|_{L^1_{t,x}}&\le \; \iiint_{|v|\le R} \left||a-v|^{\g} - |b-v|^{\g}\right|F_n(t,x,v)\,dv\,dx\,dt\\
		&\le |a-b|^{\g}\;\iiint_{|v|\le R} F_n(t,x,v)\,dv\,dx\,dt
		\ \le\ |a-b|^{\g}\,M_R,
    \end{split}
    \Ee
	and 
	\[
		|c_{R,a}-c_{R,b}|
	\leq \int_{|v|\le R}\left||a-v|^{\g} - |b-v|^{\g}\right|\,\mu(v)\,dv
	\le |a-b|^{\g} \int_{|v|\le R}\mu(v)\,dv
	\]
	since 
	\[
	\left||a-v|^\g - |b-v|^\g\right|\leq |a-b|^\g
	\]
	for $\gamma \leq 1$.	
\end{proof}

We now show how to construct one common subsequence (independent of $a$) along which $u_n^{(a)}$ converges strongly in $L^1_{t,x}$ to $c_{R,a}$ for all $a\in B_{R_1}$.   

\begin{lemma} 
\label{lem:diag}
	Consider $\{F_n\}$ with $t_n \rightarrow \infty$. There exists a subsequence $\{n_j\}$, independent of $a$, such that
	\[
		\Big\| \sup_{a\in B_{R_1}} |u_{n_j}^{(a)} - c_{R,a}| \Big\|_{L^1_{t,x}} \xrightarrow[j\to\infty]{} 0.
	\]
	And there exists a further subsequence $\{n_{j_k}\}$, independent of $a$, such that
		\[
		\sup_{a\in B_{R_1}}\,|u_{n_{j_k}}^{(a)}(t,x)-c_{R,a}|\ \xrightarrow[k\to\infty]{}\ 0
		\quad\text{for a.e.\ }(t,x),
		\]
\end{lemma}

\begin{proof}
	Let $\mathcal{D} :=\mathbb Q^3 \cap B_{R_1}$ Enumerate $\mathcal D=\{q_m\}_{m\ge1}$.
	By relative compactness Lemma \ref{lem L1 vel avg}, for $q_1$ pick a subsequence (still denoted $n$) so that $u_n^{(q_1)}\to c_{R,q_1}$ in $L^1_{t,x}$. For $q_2$, extract a further subsequence so that $u_n^{(q_1)}\to c_{R,q_1}$ and $u_n^{(q_2)}\to c_{R,q_2}$ in $L^1_{t,x}$. Continue inductively and define the diagonal subsequence $n_j$ so that
	\[
	u_{n_j}^{(q_m)} \ \to\ c_{R,q_m}\quad\text{in }L^1_{t,x}\quad\text{for every }m\in\mathbb N.
	\]
    
		From compactness of $B_{R_1}$, there exists a finite $\delta$-cover
		\[ 
			a_1,\dots,a_{N(\delta)}\in \mathcal{D}\cap B_{R_1}.
		\]
		Then, for any $a\in B_{R_1}$, there exists some $a_\ell$ such that $|a-a_\ell|<\delta$ and
		\begin{equation}
		\begin{split}
		 |u_{n_j}^{(a)}-c_{R,a}| &\le |u_{n_j}^{(a)}-u_{n_j}^{(a_\ell)}| + |u_{n_j}^{(a_\ell)}-c_{R,a_\ell}| + |c_{R,a_\ell}-c_{R,a}|\\ 
        &\le |a-a_{\ell}|^{\g} \int_{|v|\leq R} F_{n_j}(t,x,v)dv + |u_{n_j}^{(a_\ell)}-c_{R,a_\ell}| + |a-a_{\ell}|^{\g} \int_{|v|\leq R} \mu(v) dv  \\
		  &\le |u_{n_j}^{(a_\ell)}-c_{R,a_\ell}| + (A_0 + C_{\mu,R})\delta^\gamma. 
		\end{split}
		\end{equation}
		Therefore, pointwise,
		\[ 
			\sup_{|a|\le R_1} |u_{n_j}^{(a)}(t,x)-c_{R,a}| \le \max_{1\le\ell\le N} |u_{n_j}^{(a_\ell)}(t,x)-c_{R,a_\ell}| + C_R\delta^\gamma. 
		\]
		Now we integrate over $(0,1)\times\T^3$ to obtain
		\[ 
			\Big\| \sup_{|a|\le R_1} |u_{n_j}^{(a)}-c_{R,a}| \Big\|_{L^1((0,1)\times\T^3)} \le \sum_{\ell=1}^{N} \|u_{n_j}^{(a_\ell)}-c_{R,a_\ell}\|_{L^1((0,1)\times\T^3)} + (A_0 + C_{\mu,R})\delta^\gamma.
		\]
		Since the first term RHS is of finite sum, 
        \[
            \limsup_{j\rightarrow \infty} \Big\| \sup_{|a|\le R_1} |u_{n_j}^{(a)}-c_{R,a}| \Big\|_{L^1((0,1)\times\T^3)} \leq (A_0 + C_{\mu,R})\delta^\gamma.
        \]
        This holds for any $\delta>0$, and we send $\delta\downarrow0$ to obtain
        \[
		      \Big\| \sup_{a\in B_{R_1}} |u_{n_j}^{(a)} - c_{R,a}| \Big\|_{L^1_{t,x}} \xrightarrow[j\to\infty]{} 0.
	    \]
\end{proof}

  \subsection{Uniform macroscopic lower bounds}
We combine the compactness results established above with either the positive initial lower profile specified below or the initial Gaussian lower bound assumed at the end of Section~3 to establish uniform positive lower bounds for the loss frequency $L(F)$ and the density $\rho$.  

\begin{proposition} \label{prop: macro lb}
Assume $F(0,x,v)\geq\varphi(v)\geq0$ with $\int\varphi\,dv>0$. If $F$ is a global dissipative mild solution of \eqref{BE} which satisfies entropy inequality, \eqref{cond L1 vel avg} and \eqref{eq:r1-energy-tightness}. Then
there exist uniform constants $\varepsilon_{L}, \varepsilon_{\rho} > 0$ such that
		\Be \label{lb_Lf}
			L(F)(t,x,a) \geq \varepsilon_{L}\quad  \text{a.e.} \  (t,x)\quad \forall a\in \R^3 
		\Ee
		and
		\Be \label{lb_rho}
			\rho(t,x) \geq \varepsilon_{\rho}\quad \text{a.e.} \  (t,x). 
		\Ee
		Therefore,
		\[
			L(F)(t,x,v) \geq \lambda_{F} \avv^\g \quad  \text{a.e.}\  (t,x). 
		\]
		for some constant $\lambda_{F} > 0$ which depends on $F$ in non-constructive manner.
\end{proposition}
\begin{proof}
	For a unit vector $e$, define
	\[
	\beta_b:=\int_{\mathbb S^2}b(e\cdot\sigma)\,d\sigma
	=2\pi\int_{-1}^1 b(z)\,dz>0,
	\qquad C_L:=2^{\g/2}\beta_b.
	\]
	The value of $\beta_b$ is independent of $e$, and \eqref{def LF} gives
	\[
	L(F)(t,x,a)=\beta_b\int_{\mathbb R^3}|a-v|^\g F(t,x,v)\,dv
	\leq C_LA_0\langle a\rangle^\g.
	\]
	Choose $R$ large enough that $\int_{B_R}\varphi(v)\,dv>0$, in addition to the geometric requirements imposed below.
    Define
	\[
		u^{(a)}(t,x) := \int_{|v_*|\leq R} |a-v_*|^{\g} F(t,x,v_*) dv_*,\quad |a|\leq R_1
	\]
	 for $R_1 > R$. 
     We also set
	\[
	c_{R,*}:=\min_{|a|\leq R_1}c_{R,a}>0.
	\]
	Indeed, $a\mapsto c_{R,a}$ is continuous on the compact ball $\{|a|\leq R_1\}$ and is strictly positive at every $a$, since $\mu>0$ and $|a-v|^\g$ vanishes only at one point. \\
	\textit{Step 1.} 
	

The preceding loss bound yields
	\[
	\begin{aligned}
	&\int_{s}^{t}L(F)(\tau, x-v(t-\tau), v)\,d\tau\\
	&\quad=\beta_b\int_{s}^{t}\int_{v_*}|v-v_*|^{\g}
	 F(\tau, x-v(t-\tau), v_*)\,dv_*\,d\tau\\
	&\quad\leq C_LA_0\langle v\rangle^{\g}(t-s).
	\end{aligned}
	\]
	
	Now, for $|v_*|\leq R$, 
	\begin{equation}
		\rho(t,x) \geq \int_{|v_*|\leq R} F dv_* \geq e^{-C_LA_0\langle R \rangle^{\g} t}\int_{|v|\leq R}\varphi(v) dv,
	\end{equation}
    and
	\begin{equation} \label{lower bound in time}
		u^{(a)}(t,x) \geq e^{-C_LA_0\langle R\rangle^\g t}
		\int_{|v_*|\leq R}|a-v_*|^\g\varphi(v_*)\,dv_*
		\geq e^{-C_LA_0\langle R\rangle^\g t}m_{\varphi,R,R_1},
	\end{equation}
	where
	\[
	m_{\varphi,R,R_1}:=\min_{|a|\leq R_1}
	\int_{B_R}|a-v|^\g\varphi(v)\,dv>0.
	\]
	This positivity follows from continuity in $a$, compactness, and the positive mass of $\varphi$ on $B_R$; no radial symmetry of $\varphi$ is required. In particular, the infimum of $u^{(a)}$ over $|a|\leq R_1$ has a positive lower bound on each bounded time interval.
     We prove the stronger statement that $\inf_{a\in B_{R_1}}u^{(a)}$ has a positive uniform lower bound. This implies \eqref{lb_Lf} for $a\in B_{R_1}$ because $L(F)(t,x,a)\geq\beta_bu^{(a)}(t,x)$; equivalently,
		\[
			\{(t,x) \ \vert \ L(F)(t,x,a) \leq \varepsilon\} \subset  \{(t,x) \ \vert \ u^{(a)}(t,x) \leq \varepsilon/\beta_b\}.
		\]
		Suppose, to the contrary, that this stronger uniform bound for $u^{(a)}$ fails. Then for any $\varepsilon>0$ we have
    \Be \label{if not}
		\mathfrak{m}\{(t,x) : \inf_{a\in B_{R_1}} u^{(a)}(t,x) < \varepsilon \} > 0,
	\Ee
	where $\mathfrak{m}$ is 4-dimensional Lebesgue measure in $(t,x)$. Now let us consider a positive sequence $\varepsilon_n \downarrow 0$. Note that \eqref{lower bound in time} implies that there exists a sequence $T_n\rightarrow\infty$ such that
	\[
	t \leq T_n \implies (t,x)\notin \{(t,x) : \inf_{a\in B_{R_1}} u^{(a)}(t,x) < \varepsilon_n \},
	\]
	which implies
	\Be \label{rho est}
	\mathfrak{m}\{(t,x) : \inf_{a\in B_{R_1}} u^{(a)}(t,x) < \varepsilon_n,\ t \geq T_{n}\} > 0.
	\Ee
	Therefore, for each $n$, there exist at least one $t_{n}$ such that $T_n\leq t_n$ and
	\Be \label{positive meas}
	\mathfrak{m}\{(t,x) : \inf_{a\in B_{R_1}} u^{(a)}(t,x) < \varepsilon_n, \ t_n\leq t \leq t_n+1\} > 0.
	\Ee
	\\
	Now, we consider the sequence $\{t_n\}$, $t_n\rightarrow\infty$. 
   By Lemma~\ref{lem L1 vel avg}, for each $a\in B_{R_1}$, $\{u_{n}^{(a)}\}$ is relatively compact in $L^1_{t,x}$, where
	\[
		u_{n}^{(a)}(t,x):= \int_{|v_*|\leq R} |a-v_*|^{\g} F_n(t,x,v_*) dv_*,\quad |a|\leq R_1,
	\]
	where
	\[
		F_{n}(t,x,v) := F(t_{n}+t,x, v),\quad (t,x,v)\in [0,1]\times \mathbb{T}^3_x \times \R^3_v.
	\]
	\\
	Now, using Lemma \ref{lem:diag}, we can choose a subsequence 
	\[
		u_{n_{j_{k}}}^{(a)} \longrightarrow c_{R,a} \quad \text{a.e.}\quad [0,1]\times \mathbb{T}^3 \ \ \text{uniform in $a\in B_{R_1}$}
	\]
	since $u^{(a)}_n \rightarrow c_{R,a}$ in $L^1$. Using Egoroff Theorem, for any $\eta > 0$, there exists $A\subset [0,1]\times \mathbb{T}_x^3$ with $|A| \leq \eta$ such that 
	\Be \label{01domain}
			u_{n_{j_{k}}}^{(a)} \longrightarrow c_{R,a} \quad \text{uniformly in $A^c := \{[0,1]\times \mathbb{T}_x^3\}\backslash A$}, \quad \text{uniform in $a\in B_{R_1}$}.
	\Ee
	To obtain a global uniform lower bound of $\inf_{a\in B_{R_1}} u^{(a)}_n$, we have to estimate the quantity at $(t,x)\in A$, since convergence to $c_{R,a}$ in $A$ can be arbitrarily slow.   \\
	
	\noindent \textit{Step 2.}  For given subsequence $n_{j_k}$ we extend the domain of $u^{(a)}_{n_{j_k}}$ to $[-1,1]\times \mathbb{T}_x^3$. And for the sequence $n_{j_k}$, we repeat \textit{Step 1} and Lemma \ref{lem:diag}. Then, for any $\eta>0$, there exists $\bar{A}\subset [-1,1]\times \mathbb{T}_x^3$ with $|\bar{A}|\leq 2\eta$ such that
	\Be \label{extdomain}
		u_{n_{j_{k_{p}}}}^{(a)} \longrightarrow c_{R,a} \quad \text{uniformly in $\bar{A}^c := \{[-1,1]\times \mathbb{T}_x^3\}\backslash \bar{A}$}, \quad \text{uniform in $a\in B_{R_1}$}.
	\Ee
	WLOG, we relabel the sequence as $u_n^{(a)}$. $u_n^{(a)}$ satisfies \eqref{extdomain} in extended domain $[-1,1]\times \mathbb{T}_x^3$ and \eqref{01domain} in $[0,1]\times \mathbb{T}_x^3$. We can choose $\bar{A}$ so that $A \subset \bar{A}$. \\ 
	
	Now, let us measure the size of bad sets. For given $n$, we define 
	\Be \label{def B}
	\mathcal{B}^{(n)} := \{ (t,x) \in [-1,1]\times \mathbb{T}_x^3 \ \vert \ \inf_{a\in B_{R_1}} u^{(a)}_n(t,x) \leq \frac{c_{R,*}}{2} \ \text{or}  \  \int_{B_R} F_n |\log F_n| \geq H_*  \}
	\Ee
	and
	\Be \label{def G}
	\mathcal{G}^{(n)} := \{ (t,x) \in [-1,1]\times \mathbb{T}_x^3 \ \vert \ \inf_{a\in B_{R_1}} u^{(a)}_n(t,x) \geq \frac{c_{R,*}}{2}\ \text{and}  \  \int_{B_R} F_n |\log F_n| \leq H_*  \}.
	\Ee
	
	By Step 1, for any given $\eta >0$, we can choose a large $N_\eta$ and $H_* = H_*(f_0,R,\eta)\geq 1$ such that if $n > N_\eta$, then 
	\Be \label{small B1}
		\mathfrak{m}\{ (t,x) \in [-1,1]\times \mathbb{T}_x^3 \ \vert \ \inf_{a\in B_{R_1}} u^{(a)}_n(t,x) \leq \frac{c_{R,*}}{2} \} \leq |\bar{A}| \leq 2\eta
	\Ee 
	and 
	\Be \label{small B2}
		\mathfrak{m}\{ (t,x) \in [-1,1]\times \mathbb{T}_x^3 \ \vert \  \int_{B_R} F_n |\log F_n| \geq H_*  \} \leq \frac{2 H(F_0\vert\mu) + C|B_{R}|}{H_*}  \leq \eta.
	\Ee	
	 because 
		\[
		\int_{t,x,v} F|\log F| = \int_{t,x,v} F\log F - 2\int_{\substack{t,x,v \\ F<1 }} F\log F \leq 2H(F_0\vert\mu) + C|B_{R}| \quad \text{for some generic $C>0$}
		\] 
		with abbreviation $\int_{t,x,v} = \int_{-1}^1 dt \int_{\mathbb{T}_x^3} dx \int_{B_R} dv$. \\
	From \eqref{def B}, \eqref{small B1}, and \eqref{small B2}, we have
	\Be \label{small B}
		|\mathcal{B}^{(n)}| \leq 3\eta, \quad  \text{for $n\geq N_\eta$ and $H_* = H_*(f_0, R, \eta)$}.
	\Ee
	
	
	\noindent \textit{Step 3.} In this step, we will perform double Duhamel iteration for $u_n^{(a)}$.  Choose $(t,x)\in A\subset [0,1]\times \mathbb{T}_x^3$. For fixed $a\in B_{R_1}$, 
	\Be \label{double Q+}
		\begin{split}
			&\int_{|v|\leq R} |a-v|^{\g} F_n(t,x,v) dv  \\
			&\geq \delta^{\g}\int_{\substack{v\in B_R\\|v-a|\geq\delta}}dv
			 \int_{t-\Delta_1}^{t-\delta_1}e^{-C_LA_0\langle v\rangle^{\g}(t-s)}\,ds\\
			&\qquad\times Q_+(F_n)(s,x-v(t-s),v)\\
			&\geq c_b\delta^{\g}e^{-2C_LA_0\langle R\rangle^{\g}}
			 \int_{\substack{v\in B_R\\|v-a|\geq\delta}}dv
			 \int_{t-\Delta_1}^{t-\delta_1}ds\\
			&\qquad\times\int_{v_*}\int_{\sigma}|v-v_*|^{\g}
			 F_n(s,y_1(s),v'_*)F_n(s,y_1(s),v')\\
			&\geq c_b\delta^{\g}e^{-2C_LA_0\langle R\rangle^{\g}}
			 \int_{\substack{v\in B_R\\|v-a|\geq\delta}}dv
			 \int_{t-\Delta_1}^{t-\delta_1}ds\\
			&\qquad\times\int_{v_*}\int_{\sigma}|v-v_*|^{\g}F_n(s,y_1(s),v')\\
			&\qquad\qquad\times\int_{s-\Delta_2}^{s-\delta_2}
			 e^{-2C_LA_0\langle v'_*\rangle^{\g}}
			 Q_+(F_n)(s_1,y_2(s_1),v'_*)\,ds_1,
 		\end{split}
	\Ee
	where $y_1(s)=x-v(t-s)$, $y_2(s_1)=y_1(s)-v'_*(s-s_1)$ and  small positive constants $\delta, \delta_1, \Delta_1, \delta_2, \Delta_2$ will be determined later. Now, we perform Carleman estimate to obtain
	\Be \label{Carleman}
	\begin{split}
		& \int_{v_*}\int_{\sigma}|v-v_*|^{\g} F_n(s, y_1(s), v') Q_+(F_n)(s_1, y_2(s_1), v'_*) \\
		&=4 \int_{w\in \R^3} dw Q_+(F_n)(s_1, y_2(s_1), w) \frac{1}{|v-w|} \int_{v'\in E_{v,w}} \frac{1}{|v'-w|^{1-\g}}F_n(s ,y_1(s), v') dE(v') \\
	\end{split}
	\Ee
	where $E_{v,w} = \{ v' \ \vert \ (v'-v)\cdot (w-v) = 0 \}$ and $dE(v')$ is 2D plane integral on $v'\in E_{v,w}$. When $(s, y_1(s))\in \mathcal{G}^{(n)}$, we define
	\Be \label{def Omega}
		\Omega_K^{(n)}(s, y, v) := \{ w\in B_{\bar{R}} \ \vert \ \int_{E_{v,w}\cap B_{R}} F_n(s, y, v') dE(v')\geq c_0  \},
	\Ee
	where $c_0$ is to be determined later so that independent of $(n,s,y,v)$ and $\bar{R} > R$. 
	Applying to \eqref{double Q+},  and discarding the harmless factor $4$ in the lower bound,  we have
	\Be 
	\begin{split}
		\eqref{double Q+} &\geq \underbrace{ \frac{c_b c_0 \delta^{\g}}{ (R+\bar{R})^2 } e^{-2 C_LA_0 \langle R \rangle^{\g}} e^{-2 C_LA_0 \langle \bar{R} \rangle^{\g}} }_{:= c_*} \int_{t-\Delta_1}^{t-\delta_1} ds \int_{s-\Delta_2}^{s-\delta_2} ds_1 \int_{\substack{v\in B_{R} \\ |v-a|\geq \delta}}dv \int_{w\in B_{\bar{R}}} dw \\ &\quad \times Q_+(F_n)(s_1, y_2(s_1), w)  \mathbf{1}_{\mathcal{G}^{(n)}(s, y_1(s))} \mathbf{1}_{\Omega_{K}^{(n)}}(w)  \\
		&= c_* \underbrace{\int_{s, s_1, v, w} Q_+(\mu)  \mathbf{1}_{\mathcal{G}^{(n)}(s, y_1(s))} \mathbf{1}_{\Omega_{K}^{(n)}}(w)}_{I_n(t,x,a)}  + c_* \underbrace{\int_{s, s_1, v, w} (Q_+(F_n) - Q_+(\mu))  \mathbf{1}_{\mathcal{G}^{(n)}(s, y_1(s))} \mathbf{1}_{\Omega_{K}^{(n)}}(w)}_{II_n(t,x,a)}.
	\end{split}
	\Ee
	
	Let us estimate $I_n(t,x,a)$ first. For $Q_+(\mu)(w)$, we set 
	\[
		q_{\bar{R}} := \inf_{|w|\leq \bar{R}} Q_+(\mu)(w) > 0.
	\]
	\Be \label{I lower bound}
		I_n(t,x,a)
		\ge
		q_{\bar R} (\Delta_2 - \delta_2)
		\int_{t-\Delta_1}^{t-\delta_1}ds
		\int_{\substack{v\in B_R\\|v-a|\ge\delta}}
		\mathbf 1_{\mathcal G^{(n)}}(s,y_1(s))
		|\Omega_K^{(n)}(s,y_1(s),v)|\,dv.
	\Ee
		Let us define
		\[
		\mathcal P_n(t,x,a)
		:=
		\int_{t-\Delta_1}^{t-\delta_1}ds
		\int_{\substack{v\in B_R\\|v-a|\ge\delta}}
		\mathbf 1_{\mathcal G^{(n)}}(s,y_1(s)).
		\]
		We perform change of variable between $v$ and $y_1(s)=x- v(t-s)$. In $\T^3$, we have
		\[
		\int_{B_R}h(x-\tau v)\,dv
		=
		\frac1{\tau^3}
		\int_{\mathbb T^3}
		N_{\tau,x}(y)h(y)\,dy,
		\]
		for non-negative periodic function $h$ and $\tau \geq 0$. Here, $N_{\tau,x}(y)$ means number of lifts of $y$ which belong to $B_{\tau R}(x)$ If we choose $\delta_1 R$ large so that it is bigger than covering radius of $\T^3$, then $N_{\tau,x}(y)\ge1 $ for
		\[
			\tau=t-s\in[\delta_1,\Delta_1].
		\]
		For example, we can choose $R\ge \frac{10}{\delta_1}$. Therefore,
		\[
		\int_{B_R}
		\mathbf1_{\mathcal G^{(n)}}(s,x-v(t-s))\,dv
		\ge
		\frac1{\Delta_1^3}
		\int_{\mathbb T^3}\mathbf1_{\mathcal G^{(n)}}(s,y)\,dy.
		\]
		Now, we have 
		\[
		\begin{aligned}
			\mathcal P_n
			&\ge
			\frac1{\Delta_1^3}
			\int_{t-\Delta_1}^{t-\delta_1}\int_{\mathbb T^3}
			\mathbf1_{\mathcal G^{(n)}}(s,y)\,dy\,ds
			-
			\frac{4\pi}{3}\delta^3(\Delta_1-\delta_1).
		\end{aligned}
		\]
		By \eqref{small B},
		\[
		\int_{t-\Delta_1}^{t-\delta_1}\int_{\mathbb T^3}
		\mathbf1_{\mathcal G^{(n)}}(s,y)\,dy\,ds
		\ge
		(\Delta_1-\delta_1)-3\eta.
		\]
		and we obtain
		\Be \label{volume Gn}
		\begin{aligned}
			\mathcal P_n
			&\ge
			\frac{(\Delta_1-\delta_1)-3\eta}{\Delta_1^3}
			-
			\frac{4\pi}{3}\delta^3(\Delta_1-\delta_1)\\
			&\ge
			\frac1{\Delta_1^3}
			\left[
			(\Delta_1-\delta_1)
			-3\eta
			-\frac{4\pi}{3}\delta^3\Delta_1^4
			\right].
		\end{aligned}
		\Ee
	Once we prove
	\Be \label{U}
		|\Omega_K^{(n)}(s,y,v)|\ge\omega_0>0
		\quad
		\text{for }(s,y)\in\mathcal G^{(n)},\ v\in B_R,
	\Ee
	then from \eqref{I lower bound} and \eqref{volume Gn}
	\Be \label{I est}
		I_n(t,x,a)
		\ge
		q_{\bar R}\omega_0
		\frac{\Delta_2-\delta_2}{\Delta_1^3}\,
		\left[
		(\Delta_1-\delta_1)
		-3\eta
		-\frac{4\pi}{3}\delta^3\Delta_1^4
		\right] > 0.
	\Ee

	 \noindent \textit{Step 4.} 
    We estimate the absolute value of the signed remainder $II_n$.
For fixed $s,s_1,w$, the map
$v\mapsto y_2=x-w(s-s_1)-v(t-s)$ has Jacobian $(t-s)^3$.
On the torus its covering multiplicity is bounded by a constant
$M_2=M_2(R,\Delta_1)$, independent of $t,x,s,s_1,w,n,a$; for example
$M_2=(2\Delta_1R+3)^3$ is sufficient.  Since $t-s\ge\delta_1$,
\begin{align*}
 |II_n(t,x,a)|
 &\le\int_{t-\Delta_1}^{t-\delta_1}\!ds
       \int_{s-\Delta_2}^{s-\delta_2}\!ds_1
       \int_{B_R}\!dv\int_{B_{\bar R}}\!dw\,
       |Q_+(F_n)-Q_+(\mu)|(s_1,y_2,w)\\
 &\le\frac{M_2}{\delta_1^3}
       \int_{t-\Delta_1}^{t-\delta_1}\!ds
       \int_{s-\Delta_2}^{s-\delta_2}\!ds_1
       \int_{\mathbb T^3}\!dy\int_{B_{\bar R}}\!dw\,
       |Q_+(F_n)-Q_+(\mu)|(s_1,y,w).
\end{align*}
All times $s_1$ lie in $[-1,1]$ for the parameters chosen below.
Lemma~\ref{lem:r1-gain-strong}, therefore gives
\begin{equation}\label{II est}
 \sup_{\substack{t\in[0,1],\ x\in\mathbb T^3\\a\in B_{R_1}}}
 |II_n(t,x,a)|
 \le \frac{M_2(\Delta_1-\delta_1)}{\delta_1^3}
 \|Q_+(F_n)-Q_+(\mu)\|_{L^1((-1,1)\times\mathbb T^3\times B_{\bar R})}
 \longrightarrow0.
\end{equation}
The same bound applies for each fixed $a$, with a right-hand side
independent of $a$.  Thus the strong-gain limit supplies exactly the
uniform error control required in the contradiction argument.

	\noindent \textit{Step 5.} In this step, we claim \eqref{U} holds for some sufficiently large $n$. Inspired by the fact that \eqref{def Omega} depends on the direction of $w-v$ not $w$ itself, we perform directional averaging. Choose $n=\frac{v-w}{|v-w|}$ and we also define 
	\Be \label{def A}
		\mathcal{A}_v := B_{R}\backslash B_r(v),\quad \text{where $v\in B_{R}$}.
	\Ee
	Here, $r$ is to be chosen and $B_r(v)$ is a ball with radius $r$ centered at $v$.
	\Be \label{direc avg}
		\begin{split}
			&\int_{n\in \mathbb{S}^2} dn \int_{v'\in E_{v,w}\cap \mathcal{A}_v} F_n(s, y_1, v') dE(v') \\
			&=\int_{n\in \mathbb{S}^2} dn \int_{v' \in \mathcal{A}_v} dv' F_n(s, y_1, v') \delta((v'-v)\cdot n) \\
			&= \int_{v' \in\mathcal{A}_v} dv' F_n(s, y_1, v') \int_{n} dn \delta((v'-v)\cdot n) \\
			&= 2\pi \int_{v' \in\mathcal{A}_v}\frac{1}{|v-v'|} F_n(s, y_1, v') dv'. \\
		\end{split}
	\Ee
	Note that $(s, y_1) \in \mathcal{G}^{(n)}$, so using definition of \eqref{def G},
	\Be \label{density tilde R}
		\int_{|v'|\leq R} F_n(s,y_1, v') dv' \geq \frac{1}{(R+R_1)^\gamma} \int_{|v'|\leq R} |a-v'|^{\g} F_n(s, y_1, v') dv' \geq \frac{1}{(R + R_1)^\g} \frac{c_{R,0}}{2}.
	\Ee
	Also, using \eqref{def G}, 
    with all local entropy integrals restricted to $B_R$,
	\Be \label{non concen r}
		\begin{split}
			&\int_{v'\in B_r(v)\cap B_R}F_n(s,y_1,v')\,dv'\\
			&\leq\int_{\substack{v'\in B_r(v)\cap B_R\\F_n<j}}
			 F_n(s,y_1,v')\,dv'\\
			&\qquad+\frac{1}{\log j}
			 \int_{\substack{v'\in B_r(v)\cap B_R\\F_n\geq j}}
			 F_n(s,y_1,v')|\log F_n(s,y_1,v')|\,dv'\\
			&\leq\frac{4\pi r^3}{3}j+\frac{H_*}{\log j}
			 \leq\frac{1}{(R+R_1)^\g}\frac{c_{R,*}}{4}.
		\end{split}
	\Ee
	by choosing sufficiently large $j>1$ and then sufficiently small $r$. We note that the choice of $r$ depends only on $H_*$ and $(R + R_1)^{-\g}c_{R,*}/2$. 
     From \eqref{direc avg}, \eqref{density tilde R} and \eqref{non concen r}, we have 
	\Be \label{avg lb est}
		\begin{split}
		\int_{n\in \mathbb{S}^2} dn \int_{v'\in E_{v,w}\cap \mathcal{A}_v} F_n(s, y_1, v') dE(v') &\geq \frac{\pi}{R} \int_{v'\in\mathcal{A}_v} F_n(s, y_1, v') dv' \\
		&\geq  \frac{\pi}{R(R + R_1)^\g} \frac{c_{R,*}}{4} := C_1.
		\end{split}
	\Ee
	This gives a lower bound of directional average \eqref{direc avg}.  \\
 	
	Let us move on to upper bound estimate of \eqref{direc avg}. We define
	\Be \label{def Sigma}
		\Sigma_K^n := \{ n= \frac{w-v}{|w-v|} \ \vert \ \int_{v'\in E_{v,w}\cap B_{R}} F_n(s, y_1, v') dE(v') \geq c_0  \}.
	\Ee
	We have
	\Be \label{decomp AB}
		\begin{split}
			&\int_{n\in \mathbb{S}^2} dn \int_{v'\in E_{v,w}\cap \mathcal{A}_v} F_n(s, y_1, v') dE(v')  \\
			&\leq \int_{n\in \Sigma_K^n} dn \int_{v'\in E_{v,w}\cap \mathcal{A}_v} F_n(s, y_1, v') dE(v')  +  \int_{n\in \mathbb{S}^2\backslash \Sigma_K^n} dn \int_{v'\in E_{v,w}\cap \mathcal{A}_v} F_n(s, y_1, v') dE(v')  \\
			&:= (A) + (B).
		\end{split}
	\Ee
	
	For $(B)$, since
	\[
		\int_{v'\in E_{v,w}\cap \mathcal{A}_v} F_n(s, y_1, v') dE(v') \leq \int_{v'\in E_{v,w}\cap B_{R}} F_n(s, y_1, v') dE(v') \leq c_0,
	\]
	we choose $c_0 = \frac{C_1}{8\pi}$ so that
	\Be \label{B est}
		(B) \leq 4\pi c_0 \leq \frac{C_1}{2}.
 	\Ee
 	For $(A)$, we estimate as follows : 
 	\Be \label{A est}
 		\begin{split}
	 			(A) &= \int_{n\in \mathbb{S}^2} \mathbf{1}_{\Sigma_K^n} (n) dn \int_{v'\in \mathcal{A}_v} F_n(s, y_1, v') \delta((v'-v)\cdot n) dv' \\
	 			&= \int_{v'\in\mathcal{A}_v} F_n(s, y_1, v') \frac{1}{|v'-v|} \int_{n\in \Sigma_K^n \cap (v'-v)^{\perp}} dn dv'  \\
	 			&= \int_{v'\in\mathcal{A}_v} F_n(s, y_1, v') \Phi_v^n (v') dv',  \\
 		\end{split}
 	\Ee
 	where we used $\Phi_v^n (v')$ to denote $\frac{1}{|v'-v|} \int_{n\in \Sigma_K^n \cap (v'-v)^{\perp}} dn$. To apply Orlicz-H\"older inequality, we first observe Young's inequality via Legendre-Fenchel transform of $t\log(1+t)$ : 
 	\[
 		st \leq \Phi(t) + \Psi(s) := t\log(1+t) + e^{s}-1.
 	\]
 	We define Luxemburg norms
 	\Be \label{Luxemburg}
 		\begin{split}
	 		\|f\|_{L^{\Phi}} &:= \inf \Big\{ \lambda > 0 \ \vert \ \int_{\mathcal{A}_v} \Phi(\frac{|f(v')|}{\lambda})dv' \leq 1 \Big\} .\\
 		\end{split}
 	\Ee
 	By definition and the Young's inequality, we obtain Orlicz-H\"older inequality : 
 	\Be \label{Orlicz}
 		 \eqref{A est} \leq 2\|F_n(s, y_1,\cdot)\|_{L^\Phi} \|\Phi^n_v\|_{L^\Psi}. 
 	\Ee
 	 For $\|F_n(s, y_1,\cdot)\|_{L^\Phi}$, we use the following simple inequality (with $(s, y_1)\in \mathcal{G}^{(n)}$)
 	 \Be \notag
 	 \begin{split}
 	 	\int_{\mathcal{A}_v} \Phi(f) &:= \int_{\mathcal{A}_v} \big( f\log(1+f) \big) \leq \int_{\mathcal{A}_v} \big( f\log^+ f + f \big) \\
 	 	& \leq \int_{\mathcal{A}_v} \big( f|\log f| + f \big) \leq H_* + A_0.
 	\end{split}
 	\Ee
 	 (i) If $H_* + A_0 \geq 1$, we use convexity of $\Phi(t)=t\log(1+t)$ to get
 	 \[
 	 	\Phi\Big(\frac{f}{\lambda}\Big) \leq \frac{1}{\lambda}\Phi(f), \quad \forall \lambda > 1,
 	 \]
 	 and then
 	 \[
 	 	\int_{\mathcal{A}_v} \Phi\Big( \frac{f}{H_* + A_0} \Big) \leq \frac{1}{H_* + A_0} \int_{\mathcal{A}_v} \Phi(f) \leq 1.
 	 \]
 	 (ii) If $H_* + A_0 \leq 1$, $\lambda=1$ satisfies the condition of \eqref{Luxemburg}. \\
 	 Therefore, by above two cases, 
 	 \Be \label{L Phi est}
 	 	\|F_n(s, y_1, \cdot)\|_{L^\Phi} \leq \max\{1, H_* + A_0\} := \mathfrak{H}.
 	 \Ee
 	\\
 	We move on to the estimate of $\|\Phi^n_v\|_{L^\Psi}$.	 
 	By the coarea formula,
 	\Be \label{Phi order 1}
 	\begin{split}
 		\int_{A_v}\Phi_v^n(z)\,dz
 		&=\int_{\Sigma_K^n}
 		\mathcal H^2\bigl(A_v\cap(v+\omega^\perp)\bigr)\,d\omega\\
 		&\le \pi R^2|\Sigma_K^n|.
 	\end{split}
 	\Ee
 	Moreover, since \(A_v=B_R\setminus B_r(v)\),
 	\[
 	0\le\Phi_v^n(z)\le \frac{2\pi}{r}.
 	\]
 	Hence, for every integer \(k\ge1\),
 		\Be \label{Phi order k}
 	\begin{split}
 	\int_{A_v}(\Phi_v^n)^k
 	\le
 	\left(\frac{2\pi}{r}\right)^{k-1}
 	\pi R^2|\Sigma_K^n|.
 	\end{split}
 	\Ee
 	
 	
 	For \(\Psi(s)=e^s-1\),
 		\[
 		\|\Phi_v^n\|_{L^\Psi(A_v)}
 		=
 		\inf\left\{
 		\lambda>0:
 		\int_{A_v}
 		\left(
 		e^{\Phi_v^n/\lambda}-1
 		\right)dz
 		\le1
 		\right\}.
 		\]
 		From \eqref{Phi order k},
 		\Be \label{avg ub est}
 		\begin{split}
 			\int_{A_v}
 			\left(e^{\Phi_v^n/\lambda}-1\right)dz
 			&=
 			\sum_{k=1}^\infty
 			\frac1{k!\lambda^k}
 			\int_{A_v}[\Phi_v^n]^kdz\\
 			&\le
 			\pi R^2|\Sigma_K^n|
 			\sum_{k=1}^\infty
 			\frac1{k!\lambda^k}
 			\left(\frac{2\pi}{r}\right)^{k-1}\\
 			&=
 			\frac{R^2r}{2}|\Sigma_K^n|
 			\left(
 			e^{\frac{2\pi}{r\lambda}}-1
 			\right).
 		\end{split}
 		\Ee
 		Therefore, once we choose
 		\[
 		\lambda
 		=
 		\frac{2\pi/r}{
 			\log\left(
 			1+\dfrac{2}{R^2r|\Sigma_K^n|}
 			\right)},
 		\]
 		RHS of \eqref{avg ub est} becomes $1$. Hence,
 			\Be \label{L Psi est}
 	 			\|\Phi_v^n\|_{L^\Psi(A_v)}
 			\le
 			\frac{2\pi/r}{
 				\log\left(
 				1+\dfrac{2}{R^2r|\Sigma_K^n|}
 				\right)}.
 		\Ee
 	
	
	We combine \eqref{avg lb est}, \eqref{decomp AB}, \eqref{B est}, \eqref{A est}, \eqref{Orlicz}, \eqref{L Phi est}, and \eqref{L Psi est} to obtain 
%
	\Be
		\begin{split}
			C_1 &\leq \int_{n\in \mathbb{S}^2} dn \int_{v'\in E_{v,w}\cap \mathcal{A}_v} F_n(s, y_1, v') dE(v') \leq 2\mathfrak{H} \frac{2\pi / r}{\log\Big( 1 + \frac{2}{R^2 r |\Sigma_K^n| } \Big)} + \frac{C_1}{2}.
		\end{split}
	\Ee
	From 
	\Be
		\frac{C_1}{2} \leq 2\mathfrak{H} \frac{2\pi / r}{\log\Big( 1 +\frac{2}{R^2 r |\Sigma_K^n| } \Big)}, 
	\Ee
	we get a lower bound estimate of $|\Sigma_{K}^n|$ : 
	\Be
		|\Sigma_K^n| \geq \frac{2}{
			R^2r
			\left[
			\exp\!\left(\frac{8\pi \mathfrak{H}}{C_1r}\right)-1
			\right]
		} :=  \vartheta_0 > 0.
	\Ee

	Let us recover $|\Omega_{K}^{(n)}|$. 
    The plane-integral condition defining $\Omega_K^{(n)}$ depends only on the direction of $w-v$; the additional restriction $w\in B_{\bar R}$ truncates each ray. We use spherical coordinates centered at $v$ to obtain
	\Be \notag
		|\Omega_{K}^{(n)}| = \int_{n\in \Sigma_K^{(n)}} dn \int_{0}^{\ell(n)} \rho^2 d\rho,
	\Ee
	where $\ell(n)$ upper bound of radius in the direction of $n=\frac{w-v}{|w-v|}\in\Sigma_K^n$. Since $v\in B_{R}$ and $w\in B_{\bar{R}}$, we choose $\bar{R} > R$ so that $\ell(n) \geq \bar{R}-R  > 0$ for all $n\in \mathbb{S}^2$. Then,
	\Be \notag
		|\Omega_{K}^{(n)}| \geq |\Sigma_K^{(n)}| \frac{1}{3} (\bar{R}-R)^3 \geq \frac{1}{3} \vartheta_0 (\bar{R}-R)^3 > 0.
	\Ee
	Hence, we justified \eqref{U} : 
	\Be \label{Omega K lb}
		\inf_{\substack{n\ge N_\eta\\
		(s,y)\in\mathcal G^{(n)},\,v\in B_R}}
		|\Omega_K^{(n)}(s,y,v)|
		\ge \omega_0 := \frac{1}{3} \vartheta_0 (\bar{R}-R)^3 > 0
	\Ee
	and \eqref{I est} holds. \\
	
		We choose $N_1$ sufficiently large and then from \eqref{I est} and \eqref{II est},
		\Be \label{Est Q+ double}
		\begin{split}
			\eqref{double Q+} &\geq c_* (I_n(t,x,a) - |II_n(t,x,a)|)
            \\ 
			&\ge \frac{c_*}{2} q_{\bar R}\omega_0
			\frac{\Delta_2-\delta_2}{\Delta_1^3}\,
			\left[
			(\Delta_1-\delta_1)
			-3\eta
			-\frac{4\pi}{3}\delta^3\Delta_1^4
			\right] > 0.
		\end{split}	
		\Ee
    Before we finish this step, let us check the order in which all the parameters can be chosen. First we choose 
	\[
		\delta_1=\delta_2=\frac{1}{4},\quad \Delta_1=\Delta_2=\frac{1}{2} 
	\]
	so that $s, s_1 \in [-1, 1]$ since we extended the domain in \textit{Step 2}. Then, as mentioned before, we choose 
	\[
		R \geq \frac{10}{\delta_1}=40,\quad \int_{B_R}\varphi(v)\,dv>0, \quad R_1 > R, \quad \text{and}\quad \bar{R}>R.
	\] 
	Then, we choose  sufficiently small $\delta$ and $\eta$ so that satisfy
	\[
		\Delta_1 - \delta_1 - \frac{4}{3}\pi\delta^3\Delta_1^4 - 3\eta > 0
	\]
	as in \eqref{volume Gn}. $\eta$ also determines $H_*(f_0, R, \eta)$ and then $r$ is determined so that satisfies \eqref{non concen r}. The above procedure shows that all the parameters introduced above can indeed be chosen consistently. \\

	\noindent\textit{Step 6.} Now, we combine with \eqref{Est Q+ double} to obtain
	\Be \label{lb after N1}
		\begin{split}
		&\int_{|v|\leq R} |a-v|^{\g} F_n(t,x,v) dv \\
        &\geq \frac{c_*}{2}q_{\bar R}\omega_0
		\frac{\Delta_2-\delta_2}{\Delta_1^3}\,
		\left[
		(\Delta_1-\delta_1)
		-3\eta
		-\frac{4\pi}{3}\delta^3\Delta_1^4
		\right]  \\
		&:= \varepsilon_{R_1}, \quad  n\geq N_1 \quad \forall a\in B_{R_1}, \quad \forall (t,x)\in A.
		\end{split}  
	\Ee

	Now, we derive contradiction. From \eqref{extdomain}, there exists $N_0$ such that
	\Be \label{lb after N0}
		u^{(a)}_n \geq \frac{c_{R,a}}{2} \geq  \frac{c_{R,*}}{2},\quad n\geq N_0\quad \quad \forall (t,x)\in A^c.
	\Ee
	Since $\varepsilon_n \downarrow 0$, there exists sufficiently large $N_* \geq \max\{N_0, N_1\}$ such that $\varepsilon_{N_*} <  \min\{\frac{c_{R,*}}{2}, \varepsilon_{R_1}\}$. \\

	Then \eqref{lb after N0} implies that
	\Be 
	\begin{split}
		&\mathfrak{m}\{(t,x)\in A^c : \inf_{a\in B_{R_1}} u_{N_*}^{(a)}(t,x) < \varepsilon_{N_*} \}  = 0.
	\end{split}
	\Ee
	Similarly, \eqref{lb after N1} yields
	\Be 
	\begin{split}
		&\mathfrak{m}\{(t,x)\in A : \inf_{a\in B_{R_1}} u_{N_*}^{(a)}(t,x) < \varepsilon_{N_*} \}  = 0.
	\end{split}
	\Ee
	Therefore, we obtain
	\Be 
	\begin{split}
		&\mathfrak{m}\{(t,x)\in [0,1]\times\mathbb{T}_x^3 : \inf_{a\in B_{R_1}} u_{N_*}^{(a)}(t,x) < \varepsilon_{N_*} \}  = 0,
	\end{split}
	\Ee
	which is contradiction to \eqref{positive meas}. Hence there exists a global uniform lower bound $\varepsilon_{L_1} > 0$ such that
	\[
		\inf_{a\in B_{R_1}}u^{(a)}(t,x) \geq \varepsilon_{L_1} > 0,\quad \text{a.e.} \  (t,x).
	\]
%
	If $a\notin B_{R_1}$, since $R_1 > R$,
	\Be \label{L lower bound2}
	\begin{split}
	&\int_{|v|\leq R}|a-v|^{\g}F(t,x,v)\,dv\\
	&\quad\geq\frac{(R_1-R)^{\g}}{R^\g}
	 \int_{|v|\leq R}|v|^\g F(t,x,v)\,dv\\
	&\quad\geq\frac{(R_1-R)^{\g}}{R^\g}\varepsilon_{L_1}>0.
	\end{split}
	\Ee
	Since $L(F)(t,x,a)\geq\beta_b\int_{B_R}|a-v|^\g F(t,x,v)\,dv$, this proves \eqref{lb_Lf} by putting
	\[
		\varepsilon_L = \beta_b\min\{ \varepsilon_{L_1}, \frac{(R_1-R)^{\g}}{R^\g}\varepsilon_{L_1} \}.
	\]

	Density lower bound \eqref{lb_rho} is easier to get  because 
	\begin{equation*}
		\begin{split}
			\varepsilon_{L_1} &\leq \int_{|v| \leq R} |v|^{\g} F(t,x,v) dv  \leq R^{\g} \int_{|v| \leq R} F(t,x,v) dv  \\
			&\leq R^{\g} \rho(t,x).
		\end{split}
	\end{equation*}
	This proves \eqref{lb_rho} by putting $\varepsilon_{\rho} := \frac{\varepsilon_{L_1}}{R^{\g}}$.  
 In fact, the same calculation gives
	\[
	\int_{B_R}F(t,x,v)\,dv\geq\varepsilon_\rho.
	\]
	For $|a|\geq2R$, we have $|a-v|\geq|a|/2$ on $B_R$, and therefore
	\[
	L(F)(t,x,a)\geq\beta_b\varepsilon_\rho\,2^{-\g}|a|^\g
	\geq\beta_b\varepsilon_\rho\,2^{-3\g/2}\langle a\rangle^\g,
	\]
	where $R\geq40$ ensures $|a|\geq1$. For $|a|\leq2R$, \eqref{lb_Lf} gives
	\[
	L(F)(t,x,a)\geq\varepsilon_L
	\geq\frac{\varepsilon_L}{\langle2R\rangle^\g}\langle a\rangle^\g.
	\]
	Combining these two ranges proves the final assertion with
	\[
	\lambda_F:=\min\left\{\frac{\varepsilon_L}{\langle2R\rangle^\g},
	\,2^{-3\g/2}\beta_b\varepsilon_\rho\right\}>0.
	\]

\end{proof}

\bigskip
  
\section{Uniform $L^\infty$  estimate}
In this section, we prove a uniform $L^\infty$ estimate of the solution. To simplify notations, we define
\Be \label{def f}
	f := e^{\a|v|^\k}F	
\Ee
for $\alpha\geq 0$. From \eqref{apriori} or \eqref{apriori-poly} with $k\geq 2$, we have
\Be \label{aprioriforf}
	\sup_{t\geq 0}\int_{\R^3}\sup_{x}f(t,x,v) \avv^k\,dv \leq A_0 < \infty, \quad k\geq 1.
\Ee
Also, for the nonlinear collisional operator, we define
\[
	\Gamma_{\pm}(f,f)(v) :=  e^{\a|v|^\k}Q_{\pm}(e^{-\a|v|^\k}f, e^{-\a|v|^\k}f).
\]
Using this notation, we write the mild formulation of the solution by
\begin{align*}
    f(t,x,v) &= f_0(x-vt, v) + \int_0^t e^{\a|v|^\k}Q(e^{-\a|v|^\k}f, e^{-\a|v|^\k}f)(\tau, x-v(t-\tau),v)\,d\tau\\
    &=f_0(x-vt, v) + \int_0^t \left(\Gamma_+(f,f) - L(F)\right)(\tau, x-v(t-\tau),v)\,d\tau,
\end{align*}
or equivalently
\begin{align*}
    f(t,x,v) &= f_0(x-vt,v)\exp\!\left[-\int_0^t L(F)\bigl(s,x-(t-s)v,v\bigr)\,ds\right] \\
    &\quad + \int_0^t \exp\!\left[-\int_s^t L(F)\bigl(\tau,x-(t-\tau)v,v\bigr)\,d\tau\right] \Gamma_+(f, f)(s,x-v(t-s),v)\,ds.
\end{align*}
Our goal is to prove the following proposition. The proof will be given at the end of this section. As mentioned in the end of Section 3, we assume $F_0$ has a Gaussian lower bound.
\begin{proposition} \label{prop uniform bound}
    Suppose $\beta>2$. Under the hypotheses of Proposition~\ref{prop: macro lb}, we have the uniform bound
	\Be \label{uniform est}
		 \operatorname*{ess\,sup}_{t\geq 0} \| \avv^{\b_1}f(t) \|_{L^\infty_{x,v}} \leq C_{\infty,F},
	\Ee
    where
	\[
	\begin{aligned}
	C_{\infty,F}&:=C_{\alpha,\kappa,\beta,k,\gamma}
	\big(F_0, A_0,\varepsilon_L,\lambda_F\big)<\infty,
	\end{aligned}
	\]
	where $\varepsilon_L$ and $\lambda_F$ are the positive lower-bound constants supplied by Proposition~\ref{prop: macro lb}, and  
    \begin{equation*}
		\begin{split}
			\b_1 &:= \min\{ \b, \mathfrak{b}+\g+\b-2, \mathfrak{b}+\mathfrak{b}_1 + k + 2\g \} \\
			&= \begin{cases}
				\min\{ \b, \mathfrak{b}+\g+\b-2\},\quad \b \leq \mathfrak{b}_1+k+\g + 2 \\
				\min\{ \b, \mathfrak{b}+\mathfrak{b}_1 + k +2\g \},\quad \b \geq \mathfrak{b}_1+k+\g + 2 \\
			\end{cases}
		\end{split}
	\end{equation*}
	with 
    \[
    \mathfrak{b} := \min\{2-\g, \frac{2-\g}{3}\mathfrak{a} + \frac{k(1+\g)}{3}\},\quad \mathfrak{b}_1 := \min\Big\{1-\g, \frac{1-\g}{3}\mathfrak{a} + \frac{k(2+\g)}{3} \Big\}.
    \]
	Here, $\mathfrak{a}:= \min\{\b-2+\g, k+ 2\g\}$. Note that $\b_1 \geq 0$ if $\b \geq 2$.
\end{proposition}
\begin{remark} \label{rmk unif-p}
	We do not crucially use exponential weight $e^{\a|v|^\k}$ in the kernel estimate and just take cancellation as in \eqref{weight interac}. Therefore, Proposition \ref{prop uniform bound} holds even for $F$, i.e., the $\alpha=0$ case, and can be used to study the polynomial weight case. 
\end{remark}

\subsection{Decay estimates - I}
We first prove that $\langle v\rangle^{\mathfrak{a}}f(t,x,v)\in L^\infty_{t,x,v}$ for $\mathfrak{a}:= \min\{\b-2+\g, k+ 2\g\}$. We start by listing some basic properties of $\Gamma_+$.
\begin{lemma} 
\label{lem_Carleman}
	We have Carleman type estimate 
	\[
		\Gamma_{+}(f,f)(v) \leq C \int_{\R^{3}_{\eta}} \frac{f(\eta)}{|v-\eta|^{2-\g}} \int_{w\in E_{v,\eta}}  
		 f(w) \mathbf{1}_{|w| \geq \frac{|v|}{\sqrt{2}}}\,dE_{w} d\eta
	\]
    for some generic constant $C>0$, where $E_{v,\eta} = \{w\in\R^3 : (w-v)\cdot(\eta-v)=0\}$ and $dE_w$ is the standard two-dimensional Lebesgue measure on the plane $E_{v,\eta}$, for $\eta\ne v$.
\end{lemma}
\begin{proof}
	From $|v|^2 + |u|^2 = |v'|^2 + |u'|^2$, $|v'|$ and $|u'|$ cannot be located in $B_{\frac{|v|}{\sqrt{2}}}(0)$ at the same time. Therefore, for $\sigma\in \mathbb{S}^2$, using symmetry $\sigma\leftrightarrow -\sigma$,
	\begin{equation}\label{Carleman_eq1}
		\begin{split}
            \Gamma_{+}(f,f) &\leq \iint_{\mathbb{R}^3\times\mathbb{S}^2} |v - v_*|^{\g} 
            e^{\a|v|^{\k}} e^{-\a|v^{\prime}|^{\k}} e^{-\a|v^{\prime}_*|^{\k}}  
            f(v'_*)f(v') \Big( \mathbf{1}_{|v'_*|\geq \frac{|v|}{\sqrt{2}}} + \mathbf{1}_{|v'|\geq \frac{|v|}{\sqrt{2}}} \Big)\,dv_*d\sigma \\
            &\leq 2 \iint_{\mathbb{R}^3\times\mathbb{S}^2} |v-v_*|^{\g}
            e^{\a|v|^{\k}} e^{-\a|v^{\prime}|^{\k}} e^{-\a|v^{\prime}_*|^{\k}} 
            f(v'_*)f(v')  \mathbf{1}_{|v'_*|\geq \frac{|v|}{\sqrt{2}}}\,dv_*d\sigma.
		\end{split}
	\end{equation}
    We treat
    \[
    e^{\a(|v|^2)^{\k/2}} e^{-\a(|v'|^2)^{\k/2}} e^{-\a(|v'_*|^2)^{\k/2}}.
    \]
    Using concavity of $(\cdot)^{\k/2}$ with $\k/2 <1$, 
    \[
    (|v|^2)^{\k/2} - [(|v'|^2)^{\k/2} + (|v'_*|^2)^{\k/2}] \leq (|v|^2)^{\k/2} - [(|v'|^2+|v'_*|^2)^{\k/2}] = (|v|^2)^{\k/2} - [(|v|^2+|v_*|^2)^{\k/2}] \leq 0,
    \]
    which means
    \Be \label{weight interac}
        e^{\a(|v|^2)^{\k/2}} e^{-\a(|v'|^2)^{\k/2}} e^{-\a(|v'_*|^2)^{\k/2}} \leq 1.
    \Ee
    Therefore, using the Carleman representation, we have
    \begin{equation*}
        \begin{split}
            \eqref{Carleman_eq1}&\leq 2 \iint_{\mathbb{R}^3\times\mathbb{S}^2} |v-v_*|^{\g} f(v'_*)f(v')  \mathbf{1}_{|v'_*|\geq \frac{|v|}{\sqrt{2}}}\,dv_*d\sigma\\
            &\leq C \int_{\R^{3}_{\eta}} \frac{f(\eta)}{|v-\eta|} \int_{w\in E_{v,\eta}} \frac{1}{\sqrt{|v-\eta|^2+|w-v|^2}^{1-\g}} f(w) \mathbf{1}_{|w| \geq \frac{|v|}{\sqrt{2}}}\,dE_{w} d\eta\\		
            &\leq C \int_{\R^{3}_{\eta}} \frac{f(\eta)}{|v-\eta|^{2-\g}} \int_{w\in E_{v,\eta}} f(w) \mathbf{1}_{|w| \geq \frac{|v|}{\sqrt{2}}}\,dE_{w} d\eta.
        \end{split}
    \end{equation*}
\end{proof}

\begin{lemma} \label{lem_int Gamma}
	For $0\leq e<2$, we have the following estimates.
    \begin{equation*}
	   \int_{\R^3}\frac{1}{|v_1- v|^e}\Gamma_{+}(f,f)(v_1)\,dv_1 \leq C \int_{\R^3} f(v_1)\,dv_1 \times \sup_{v_1}\int_{\R^3} {\color{black}\frac{1}{|v_1-u_1|^{e-\g}} } f(u_1)\,du_1.
    \end{equation*}
    Especially when $e=\gamma$, we have
	\begin{equation*}
        \int_{\R^3}\frac{1}{|v_1- v|^\g}\Gamma_{+}(f,f)(v_1)\,dv_1 \leq C A_0^2.
	\end{equation*}
\end{lemma}
\begin{proof}
	Using change of variable $(v_1,u_1)\leftrightarrow (v'_1, u'_1)$ and use \eqref{weight interac}
	\begin{equation} \label{int Gamma+}
		\begin{split}
			\int_{\R^3}\frac{1}{|v_1 - v|^e}\Gamma_{+}(f,f)(v_1)\,dv_1 &\leq C_b \iint_{\R^3\times \R^3\times \S^2} \frac{1}{|v_1 - v|^e} |v_1 - u_1 |^{\g} f(u'_1)f(v'_1)\,dv_1du_1d\sigma \\
			&= C_b\iint_{\R^3\times \R^3}  |v_1 - u_1|^{\g} f(u_1)f(v_1) \int_{\S^2} \frac{1}{|v'_1-v|^e}\,d\sigma\,dv_1du_1,
		\end{split}
	\end{equation}
	where we used $b(\sigma) \leq C_b$.	The $\sigma$ integral is computed by
    \begin{align*}
        \int_{\S^2} \frac{1}{|v'_1-v|^e}\,d\sigma&= \left(\frac{2}{|v_1-u_1|}\right)^{-e}\int_{\S^2} \left|\frac{v_1+u_1 - 2v}{|v_1-u_1|} + \sigma\right|^{-e}\,d\sigma\\
        &=\frac{2\pi}{(2-e)}\left(\frac{2}{|v_1-u_1|}\right)^{-e}\frac{1}{R}[(1+R)^{2-e} - |1-R|^{2-e}]
    \end{align*}
    for $e<2$, where $R = \frac{|v_1+u_1 - 2v|}{|v_1-u_1|}$. The final term is upper bounded by
    \begin{align*}
        \sup_{R>0}\frac{1}{R}[(1+R)^{2-e} - |1-R|^{2-e}]\leq \begin{cases}
            4\pi & 0\leq e<1,\\
            \frac{2^{3-e}\pi}{2-e} & 1\leq e<2.
        \end{cases}
    \end{align*}
	Applying this result to \eqref{int Gamma+}, we get
	\begin{equation*}
		\begin{split}
			&\int_{\R^3}\frac{1}{|v_1 - v|^e}\Gamma_{+}(f,f)(v_1)\,dv_1 \\
			&\leq C \int_{\R^3\times\R^3} \frac{1}{|v_1-u_1|^{e-\gamma} } f(u_1) f(v_1)\,dv_1du_1\\
			&\leq C \int_{\mathbb{R}^3} f(v_1)\,dv_1 \times \sup_{v_1}\int_{\mathbb{R}^3}  \frac{1}{|v_1-u_1|^{e-\gamma} } f(u_1)\,du_1
		\end{split}
	\end{equation*}
    for $0\leq e<2$.
\end{proof}

Let us recall $L(F) \geq \lambda_{F} \avv^{\g}$ from Proposition~\ref{prop: macro lb} of the previous section. From Lemma~\ref{lem_Carleman}, 
\begin{equation} \label{AB}
	\begin{split}
		&\frac{d}{dt}f(t, x+vt, v) + \lambda_{F}\langle v \rangle^{\color{black}{\g}} f(t, x+vt,v) \\
        &\leq \Gamma_{+}(f,f)(t, x+vt, v)\\
        &\leq C \underbrace{ \int_{\R^{3}_{\eta}} \frac{f(t, x+vt, \eta)}{|v-\eta|^{\color{black}2-\g}}\,d\eta }_{:=\mathcal{A}}
        \underbrace{\int_{w\in E_{v,\eta}} 
        f(t, x+vt, w) \mathbf{1}_{|w| \geq \frac{|v|}{\sqrt{2}}}\,dE_{w}}_{:=\mathcal{B}}.
	\end{split}
\end{equation}
We divide the integral into $\mathcal{A}$ and $\mathcal{B}$ as above and estimate each integral.

\begin{lemma}\label{lem unif est A} 
	$\mathcal{A}$ is uniformly bounded as
	\[
     \mathcal{A} \leq C_{F_0, \g, A_0,\varepsilon_L^{-1}}.
	\]
\end{lemma}
\begin{proof}
From the mild formulation and \eqref{lb_Lf}, we have
\begin{align*}
    f(t,x,v) &= f_0(x-vt,v)\exp\!\left[-\int_0^t L(F)\bigl(s,x-(t-s)v,v\bigr)\,ds\right] \\
    &\quad + \int_0^t \exp\!\left[-\int_s^t L(F)\bigl(\tau,x-(t-\tau)v,v\bigr)\,d\tau\right] \Gamma_+(f, f)(s,x-v(t-s),v)\,ds\\
    &\leq e^{-\varepsilon_{L}t} f_0(x-vt,v) + C\int_{0}^{t} e^{-\varepsilon_{L}(t-s)} \Gamma_+(f, f)(s,x-v(t-s),v)\,ds.
\end{align*}
Applying Lemma~\ref{lem_int Gamma} for $e = 2-\gamma$ to the estimate of $\int_{\R^{3}} \frac{\sup_x f(v_1)}{|v_1 - v|^{2-\gamma}}\,dv$, we get uniform estimate
	\begin{equation*}
		\begin{split}
			&\int_{\R^3} \frac{\sup_x f(t,x,v_1)}{|v_1 - v|^{\color{black}2-\g}}\,dv_1\\
            &\leq e^{-\varepsilon_{L}t} \int_{\R^3} \frac{\sup_x f_0(x - v_1 t, v_1)}{|v_1 - v|^{\color{black}2-\g}}\,dv_1 
			+ C\int_{0}^{t} e^{-\varepsilon_{L}(t-s)}
			\int_{\R^3} \frac{1}{|v_1 - v|^{\color{black}2-\g}} \Gamma_{+}(\sup_x f, \sup_x f) (s, x - v_1 (t-s), v_1)\,dv_1\,ds \\
			&\leq \Big( \|f_0\|_{L^\infty}\int_{|v_1-v|\leq 1} \frac{1}{|v_1-v|^{2-\g}} + \int_{|v_1 - v|\geq 1}\sup_x f_0(x,v)  \Big)  \\
			&\quad +
			\frac{C}{\varepsilon_L}  \sup_{s\leq t}\int_{\R^3} \sup_{x} f(s, x, v_1) dv_1 \times \sup_{s\leq t} \sup_{v_1}\int_{\R^3} {\color{black}\frac{1}{|v_1 - v_2|^{2-2\g}} } \sup_x f(s, x, v_2)\,dv_2  \\
			&\leq C_{F_0, \g, A_0} +
			C_{A_0, \varepsilon_L}\sup_{s\leq t} \sup_{v_1}\int_{\R^3} {\color{black}\frac{1}{|v_1 - v_2|^{2-2\g}} } \sup_x f(s, x, v_2) dv_2, \quad {\color{black} 0<\g \leq  1}.
		\end{split}
	\end{equation*}
    For $\gamma=1$, we are done. For the general case of $\g <1$, we repeat this process $m$-times so that $0\leq m\g -2<1$ for the first time. Then,
	\begin{equation} \label{induction est}
		\begin{split}			
			&\int_{\mathbb{R}^3} \frac{\sup_x f(t, x, v_1)}{|v_1 - v|^{\color{black}2-\g}}\,dv_1  \\	
			&\leq C_{F_0, \g, A_0} +  C_{A_0, \varepsilon_L^{-1}}\sup_{s\leq t} \sup_{v_1}\int_{\R^3} {\color{black}\frac{1}{|v_1 - v_2|^{2-2\g}} } \sup_x f(s, x, v_2)\,dv_2 \quad {\color{black} \g < 1} \\
			&\leq C_{F_0, \g, A_0} +  C_{A_0, \varepsilon_L^{-1}}
			\Big[
			C_{F_0, \g, A_0} +  C_{A_0, \varepsilon_L^{-1}}\sup_{s_1\leq t}\sup_{v_2}\int_{\R^3} {\color{black}\frac{1}{|v_2 - v_3|^{2-3\g}} } \sup_x f(s_1, x, v_3)\,dv_3
			\Big]  \\
			&\leq \cdots \quad \text{repeating $m-2$  more times}  \\
			&\leq C_{F_0, \g, A_0} +  C_{F_0, \g, A_0, \varepsilon_L^{-1}} \sup_{s_{m-1}\leq t} \sup_{v_m}\int_{\R^3} {\color{black}\frac{1}{|v_m - v_{m+1}|^{2-(m-1)\g}} }\sup_x f(s_{m-1}, x, v_{m+1})\,dv_{m+1}\\
			&\leq C_{F_0, \g, A_0} +  C_{F_0, \g, A_0, \varepsilon_L^{-1}} \sup_{s_m\leq t} \sup_{v_m}\int_{\R^3} {\color{black}\frac{1}{|v_m - v_{m+1}|^{2-(m-1)\g}} }\Gamma_+(\sup_x f, \sup_x f)(s_m, x, v_{m+1})\,dv_{m+1}\\
			&\leq C_{F_0, \g, A_0} + C_{F_0, \g, A_0, \varepsilon_L^{-1}} \sup_{s\leq t}\int_{\R^3\times \R^3}  |v_{m+1} - v_{m+2}|^{\color{black}m\g-2} \sup_x f(s, x , v_{m+1}) \sup_x f(s, x, v_{m+2})\,dv_{m+1} dv_{m+2}\\
			&\leq C_{F_0, \g, A_0} + C_{F_0, \g, A_0, \varepsilon_L^{-1}} \sup_{s\leq t}\Big[\int_{\R^3}  \langle v_{m+1} \rangle^{m\g-2} \sup_x f(s, x, v_{m+1})\,dv_{m+1} \Big]^2   \lesssim C_{F_0, \g, A_0, \varepsilon_L^{-1}}
		\end{split}
	\end{equation}
    as $m\gamma\geq 2$. In the second line from the bottom, we use the result from the proof of Lemma~\ref{lem_int Gamma}.
\end{proof}

For a given $w\in \R^{3}$ and a plane $E\subset \R^{3}$, define
\[
d_{E}(w) := \text{dist}(w, E).
\]
We approximate the measure $dE_w$ on the plane $E$ via $\phi^E_{j}(w)dw$ by
\[
\phi^{E}_{j}(w) = \sqrt{\frac{j}{\pi}} e^{-j d_{E}^2(w)}
\]
so that
\[
\lim_{j\rightarrow\infty} \int_{\R^{3}} \phi^{E}_{j}(w) \Gamma_{+}(f, f)(w)\,dw = \int_{E\subset\R^{3}} \Gamma_{+}(f,f)(w)\,dE_{w}.
\]
\begin{lemma}\label{lem unif est B} 
	For $\b >2$, we have
	{\color{black}
		\Be \label{est_int_E}
		\int_{w\in E_{v,\eta}} 
        f(t, x, w) \mathbf{1}_{|w| \geq \frac{|v|}{\sqrt{2}}}\,dE_{w} \leq \frac{C_{F_0,\b}}{\avv^{\b-2}} + \frac{C_{F_0, k, \g, A_0}}{\lambda_{F,\mathrm{pl}}\avv^{k+\g}}
		\Ee
	} 
	for almost everywhere $\eta$, $t$, and $x$. So, $\mathcal{B}$ is bounded by the right-hand side.
\end{lemma}
\begin{proof}
	Introduce cutoff $\psi$,
	\[
	\psi_{v}(w) := \begin{cases}
		0\quad |w| < \frac{|v|}{\sqrt{2}},\\
		1 \quad |w| \geq \frac{|v|}{\sqrt{2}}.
	\end{cases}
	\]
	On the support of $\psi_v(v_1)$, we have $\langle v_1\rangle\geq2^{-1/2}\langle v\rangle$. Set
	\[
	\lambda_{F,\mathrm{pl}}:=2^{-\gamma/2}\lambda_F>0.
	\]
	The coercivity estimate from Proposition~\ref{prop: macro lb} gives, for $0\leq s\leq t$,
	\begin{align*}
        \psi_v(v_1)\exp\!\left\{-\int_s^t L(F)(\tau,x-v_1(t-\tau),v_1)\,d\tau\right\}
        &\leq \psi_v(v_1)e^{-\lambda_F\langle v_1\rangle^\gamma(t-s)}\\
        &\leq \psi_v(v_1)e^{-\lambda_{F,\mathrm{pl}}\langle v\rangle^\gamma(t-s)}.
	\end{align*}
	In particular,
	\[
	\int_0^t e^{-\lambda_{F,\mathrm{pl}}\langle v\rangle^\gamma(t-s)}\,ds
	\leq\frac{1}{\lambda_{F,\mathrm{pl}}\langle v\rangle^\gamma}.
	\]
	Writing $E = E_{v,\eta}$, we apply $\int_{v_1} \phi^E_{j}(v_1) \psi_v(v_1)$ to the mild formulation of the solution to obtain
	\begin{equation} \label{integ_eq}
		\begin{split}
			&\int_{\R^3} \phi^{E}_{j}(v_1)\psi_{v}(v_1) f(t, x, v_1)\,dv_1  \\
			&\leq e^{-\lambda_{F,\mathrm{pl}} \avv^{\g} t}  \int_{\R^3}\phi^{E}_{j}(v_1)\psi_{{v}}(v_1) f_0(x-v_1 t, v_1)\,dv_1  \\
			&\quad + \int_{0}^{t}e^{-\lambda_{F,\mathrm{pl}} \avv^{\g} (t-s)} \int_{\R^3} \phi^{E}_{j}(v_1)\psi_{{v}}(v_1) \Gamma_{+}(f,f)(s, x -v_1 (t-s), v_1)\,dv_1 ds\\
			&\leq e^{-\lambda_{F,\mathrm{pl}} \avv^{\g} t} 
			\int_{\R^3}\phi^{E}_{j}(v_1) \psi_{{v}}(v_1) f_0(x-v_1 t, v_1)\,dv_1  \\
			&\quad + \int_{0}^{t} e^{-\lambda_{F,\mathrm{pl}} \avv^{\g} (t-s)} \int_{\R^3} \int_{\R^3} 
			\int_{\S^2} \phi^{E}_{j}(v'_1) \psi_{{v}}(v'_1) |v_1 - u_1|^{\color{black}\g} f(s, x - v'_1 (t-s), u_1)\\
            &\qquad \quad  \times f(s, x - v'_1 (t-s), v_1)\,d\sigma\,du_1 dv_1 ds
		\end{split}
	\end{equation}
	for fixed ${v}\in\R^{3}$, where we used change of variables $(v_1, u_1) \leftrightarrow (v'_1, u'_1)$. If both $|v_1|, |u_1| \leq \frac{|{v}|}{2}$, then $|v'_1|^2\leq |v_1|^2+|u_1|^2 \leq \frac{|{v}|^2}{2}$, so $\psi_{{v}}(v'_1)=0$. So at least one of $v_1$ and $u_1$ must be larger than $\frac{{v}}{2}$. Therefore, we have
    \begin{equation}\label{integ_eq_2}
        \begin{split}
            \eqref{integ_eq}&\leq e^{-\lambda_{F,\mathrm{pl}} \avv^{\g} t} 
            \int_{\R^3}\phi^{E}_{j}(v_1) \psi_{{v}}(v_1) \sup_x f_0(x, v_1)\,dv_1  \\
            &\quad + \int_{0}^{t} e^{-\lambda_{F,\mathrm{pl}} \avv^{\g} (t-s)} \int_{\R^3} \int_{\R^3} 
            |v_1 - u_1|^{\color{black}\g} \sup_x f(s, x, u_1)\sup_x f(s, x, v_1)\\
            & \qquad \quad \times \left(\mathbf{1}_{|v_1|\geq \frac{|v|}{2}} + \mathbf{1}_{|u_1|\geq \frac{|v|}{2}}\right)\underbrace{\int_{\S^2} \phi^{E}_{j}(v'_1)\,d\sigma}_{\eqref{max_length}}\,du_1 dv_1 ds.
        \end{split}
    \end{equation}
	
	We write $d = d_E\left(\frac{v_1+v_2}{2}\right)$ and $2R = |v_1-u_1|$. Then, by taking spherical coordinate centered at $\frac{v_1+v_2}{2}$ with $z$-axis perpendicular to the plane $E$, we have
	\begin{equation}\label{max_length}
	    \begin{split}
            \int_{\S^2} \phi^{E}_{j}(v'_1)\,d\sigma &= 2\pi \int_0^\pi \sqrt{\frac j\pi}e^{-j|d + R\cos\theta|^2} \sin\theta\,d\theta\\
            &=\frac{2\pi}{R} \int_{d-R}^{d+R} \sqrt{\frac j\pi}e^{-jt^2}\,dt\\
            &\le \frac{2\pi}{R} \le \frac{C}{|v_1 - u_1|}.
        \end{split}
	\end{equation}
    Now, the integrand is symmetric about $v_1$ and $u_1$, so we can choose $|v_1| \geq \frac{|{v}|}{2}$ by the symmetry. Therefore, we get
	\begin{equation}\label{integ_eq_3}
		\begin{split}
			\eqref{integ_eq_2} &\leq e^{-\lambda_{F,\mathrm{pl}} \avv^{\g} t} \int_{v_1\in\R^{3}}\phi^{E}_{j}(v_1)\psi_{{v}}(v_1) \sup_x f_0(x, v_1)\,dv_1 \\
			&\quad + C\int_{0}^{t} e^{-\lambda_{F,\mathrm{pl}} \avv^{\g} (t-s)}\underbrace{\int_{\R^{3}}\int_{\R^{3}}
			\frac{1}{|v_1-u_1|^{1-\g}}
			\sup_x f(s,x,u_1)\sup_x f(s,x,v_1)\mathbf{1}_{|v_1|\geq \frac{|v|}{2}}\,du_1 dv_1}_{(**)_{\g}}\,ds.
		\end{split}
	\end{equation}
    Using \eqref{apriori} or \eqref{apriori-poly}, the integral $(**)_{\g}$ can be bounded as follows.
	\begin{equation*}
		\begin{split}
			(**)_{\g} &\leq \int_{|v_1|\geq \frac{|{v}|}{2}}  \sup_{x} f(s, x, v_1)\,dv_1  \int_{\R^{3}} \frac{1}{|v_1-u_1|^{1-\g}} \sup_{x}f(s, x, u_1)\,du_1\\
			&\leq \frac{C_k}{\avv^k}\int_{|v_1|\geq \frac{|{v}|}{2}} \sup_x f(s, x, v_1)\langle v_1 \rangle^k\,dv_1 
			\times  \int_{\R^{3}} \frac{1}{|v_1-u_1|^{1-\g}} \sup_{x} f(s, x, u_1)\,du_1\\
			&\leq \frac{C_{F_0, k, \g, A_0}}{\avv^k},
		\end{split}
	\end{equation*}
    where we used \eqref{induction est} in the last step. Therefore,
    \begin{align*}
        &\int_{0}^{t} e^{-\lambda_{F,\mathrm{pl}} \avv^{\g} (t-s)} (**)_{\g}\,ds\leq \frac{C_{F_0, k, \g, A_0}}{\lambda_{F,\mathrm{pl}}\avv^{k+\g}}.
    \end{align*}
    For the initial term, we have
	\Be \label{initial_integ}
    	\int_{E} \psi_{{v}}(v_1) \sup_x f_0(x, v_1)\,dE(v_1) \leq \|wf_0\|_\infty \frac{2\pi}{\beta-2}\langle \frac{|v|}{\sqrt{2}}\rangle^{2-\beta} \leq \frac{C_{F_0,\b}}{\avv^{\b-2}} , \qquad \beta>2  
	\Ee
    from \eqref{initial} or \eqref{initial-poly}.
    
	Combining the estimates, we conclude that
	\begin{equation}
        \begin{split}
            \int_{\R^3} \phi^{E}_{j}(v_1) \psi_v(v_1) f(t, x, v_1)\,dv_1\leq \frac{C_{F_0,\b}}{\avv^{\b-2}} + \frac{C_{F_0, k, \g, A_0}}{\lambda_{F,\mathrm{pl}}\avv^{k+\g}}.
        \end{split}
	\end{equation}
    Since the upper bound is uniform about $j$, we take the limit $j\rightarrow\infty$ to obtain the results.
\end{proof}	

Applying Lemmas~\ref{lem unif est A} and \ref{lem unif est B} to \eqref{AB}, we perform the low-order estimate $\|\avv^\mathfrak{a} f\|_{L^\infty_{t,x,v}}$.
\begin{lemma} 
\label{lem_low order}
	For $\b > 2$, we have
	\[
		\|\avv^\mathfrak{a} f\|_{L^\infty_{t,x,v}} \leq \|\avv^\mathfrak{a} f_0\|_{L^\infty_{x,v}} + C_{F_0, \b, k, \g, A_0, \epsilon_L^{-1}, \lambda_{F,\mathrm{pl}}^{-1}} \quad \mathfrak{a}:= \min\{\b-2+\g, k+ 2\g\} \leq \b-1.
	\]
\end{lemma}
\begin{proof}
    By Lemmas~\ref{lem unif est A} and \ref{lem unif est B}, we have
    \begin{equation*}
        \begin{split}
            \Gamma_{+}(f,f)(s, x - v(t-s), v)&\leq C_{F_0,\gamma, A_0, \epsilon_L^{-1}}\left(\frac{C_{F_0,\b}}{\avv^{\b-2}} + \frac{C_{F_0, k, \g, A_0}}{\lambda_{F,\mathrm{pl}}\avv^{k+\g}}\right).
        \end{split}
    \end{equation*}
    Applying this bound to the mild formulation \eqref{AB}, we get
	\begin{equation*}
		\begin{split}
			f(t,x,v) &\leq e^{-\lambda_{F}\avv^{\g} t}f_0(x-vt,v) + \int_{0}^{t} e^{-\lambda_{F}\avv^{\g} (t-s)} \Gamma_{+}(f,f)(s, x - v(t-s), v)\,ds\\
			&\leq e^{-\lambda_{F}\avv^{\g} t}f_0(x-vt,v) \\
            &\quad + C_{F_0,\b, k, \g, A_0, \epsilon_L^{-1}, \lambda_{F,\mathrm{pl}}^{-1}}\int_{0}^{t} e^{-\lambda_{F}\avv^{\g} (t-s)} \left(\frac{1}{\avv^{\b-2}} + \frac{1}{\avv^{k+\g}}\right)\,ds.
		\end{split}
	\end{equation*}
	Using $\int_0^t e^{-\lambda_F\langle v\rangle^\gamma(t-s)}\,ds\leq\lambda_F^{-1}\langle v\rangle^{-\gamma}$, we get
		\begin{equation} \label{low order weight}
			\begin{split}
				\|\avv^\mathfrak{a} f\|_{L^\infty} &\leq \|\avv^\mathfrak{a} f_0\|_{L^\infty} + C_{F_0,\b, k, \g, A_0, \epsilon_L^{-1}, \lambda_{F,\mathrm{pl}}^{-1}}, \quad \mathfrak{a}:= \min\{\b-2+\g, k+ 2\g\} \leq \b-1.
			\end{split}
		\end{equation}
\end{proof}

\subsection{Decay estimates - II} Based on the previous result $\langle v\rangle^{\mathfrak{a}}f(t,x,v)\in L^\infty_{t,x,v}$, we refine the estimate of \eqref{AB} for each $\mathcal{A}$ and $\mathcal{B}$. As a result, we improve the decay of $f$ in the $L^\infty$ norm. We first bound $\mathcal{A}$ in the next lemma.
\begin{lemma}\label{lem_A decay} 
	Let $f(v)\in L^{1}_{k} \cap L^{\infty}_{\mathfrak{a}}$. Then for $0\leq e<2$, we have
		\[
		\int_{\R^{3}_{\eta}} \frac{f(\eta)}{|v-\eta|^e}\,d\eta \leq C_k\big( \|f\|_{L^\infty_{\mathfrak{a}}} + \|f\|_{L^{1}_{k}} \big) \frac{1}{\avv^{\mathfrak{b}}},\quad \mathfrak{b} := \min\Big\{e, \frac{e}{3}\mathfrak{a} + \frac{k(3-e)}{3} \Big\}.
		\]
        Therefore, for $f(t,x,v)\in L^\infty_{\mathfrak{a}}(t,x,v)$, we have
        \begin{align*}
            \mathcal{A}\leq C_k(\|f\|_{L^\infty_{\mathfrak{a}}(t,x,v)} + A_0)\frac{1}{\avv^{\mathfrak{b}}},\quad \mathfrak{b} := \min\Big\{2-\gamma, \frac{2-\gamma}{3}\mathfrak{a} + \frac{k(1+\g)}{3} \Big\}.
        \end{align*}
\end{lemma}
\begin{proof}
	If we prove the first inequality, the bound on $\mathcal{A}$ follows easily by setting $f = \sup_x f(x,v)$ and $e = 2-\g$. To prove the first inequality, We split
	\begin{equation*}
		\begin{split}
			O_1 &:= \{\eta : |\eta| \leq \frac{|v|}{2} \},  \\
			O_2 &:= \{ \eta : |\eta-v| \leq \avv^{-\theta} \},\quad \text{$\theta \geq -\frac{2}{3}$ to be chosen},  \\
			O_3 &:= \{ \eta : \R^{3} \backslash (O_1 \cup O_2)\}.
		\end{split}
	\end{equation*}
	(i) For $O_1 :=\{\eta : |\eta| \leq \frac{|v|}{2} \}$. Since $|v-\eta| \geq \frac{|v|}{2}$,
	\begin{equation*}
		\begin{split}
			\int_{O_1}  \frac{f(\eta)}{|v-\eta|^e}\,d\eta &\leq  \mathbf{1}_{|v|\geq 1} \frac{C}{|v|^e} \int_{O_1} f(\eta)\,d\eta + \mathbf{1}_{|v|\leq 1} \|f\|_{L^\infty} \int_{O_1}\frac{1}{|\eta-v|^e}\,d\eta   \\
			&\leq \mathbf{1}_{|v|\geq 1} \frac{C}{|v|^e} \|f\|_{L^{1}} + C \mathbf{1}_{|v|\leq 1} \|f\|_{L^\infty}  \\
			&\leq \frac{C}{\avv^e}( \|f\|_{L^{1}} + \|f\|_{L^\infty}).
		\end{split}
	\end{equation*}
	
	\noindent (ii) For $O_2 := \{ \eta : |\eta-v| \leq \avv^{-\theta} \}$, if $|v|\geq 10$ then we have
	$
	|\eta-v| \leq \avv^{-\theta} \leq \frac{|v|}{2}
	$
	which gives
	\[
	|\eta| \geq \frac{|v|}{2}\quad \text{so} \quad \frac{1}{\langle \eta \rangle} \leq \frac{C}{\langle v \rangle}\quad \text{for $|v|\geq 10$}.
	\]
	Then it holds that
	\begin{equation*}
		\begin{split}
			\int_{O_2}  \frac{f(\eta)}{|v-\eta|^e}\,d\eta &\leq  \mathbf{1}_{|v|\geq 10} \|f\|_{L^\infty_{\mathfrak{a}}} \int_{O_2} \frac{1}{\langle \eta \rangle^\mathfrak{a}} \frac{1}{|\eta - v|^e}\,d\eta +  \mathbf{1}_{|v|\leq 10} \|f\|_{L^\infty} \int_{O_2} \frac{1}{|\eta - v|^e}\,d\eta \\
			&\leq C\mathbf{1}_{|v|\geq 10} \|f\|_{L^\infty_{\mathfrak{a}}} \frac{1}{\avv^\mathfrak{a}} \int_{0}^{\avv^{-\theta}} \frac{1}{r^e} r^2\,dr + C \mathbf{1}_{|v|\leq 10} \|f\|_{L^\infty}  \\
			&\leq C\mathbf{1}_{|v|\geq 10} \|f\|_{L^\infty_{\mathfrak{a}}} \frac{1}{\avv^{\mathfrak{a} + \theta(3-e)}}  +  C \mathbf{1}_{|v|\leq 10} \|f\|_{L^\infty}  \\
			&\leq C \|f\|_{L^{\infty}_{\mathfrak{a}}} \frac{1}{\avv^{\mathfrak{a}+\theta(3-e)}}.
		\end{split}
	\end{equation*}
	
	\noindent (iii) For $O_3 $, we have $|\eta| \geq \frac{|v|}{2}$ and $|\eta-v|\geq \avv^{-\theta}$, so
	\begin{equation*}
		\begin{split}
			\int_{O_3}  \frac{f(\eta)}{|v-\eta|^e}\,d\eta &\leq \frac{1}{\avv^{-\theta e}} \int_{O_3} \langle \eta \rangle^k f(\eta) \frac{1}{\langle \eta \rangle^k}\,d\eta \\
			&\leq \frac{C_k}{ \avv^{k-\theta e}} \|f\|_{L^{1}_{k}}.
		\end{split}
	\end{equation*}
	For an optimal balance of $(ii)$ and $(iii)$, we choose 
	$\mathfrak{a}+\theta(3-e) = k-\theta e$,
	which yields 
	\[
	\theta = \frac{k}{3} - \frac{1}{3}\mathfrak{a}. 
	\] 
    Note that $\mathfrak{a} \leq k + 2\g$ implies $\theta \geq -\frac{2\g}{3} \geq -\frac{2}{3}$. With this choice of $\theta$, from above three estimates, we have decay rate : 
	\begin{equation} \label{decay rate}
		\begin{split}
			&\mathfrak{b} := \min\{e, \mathfrak{a} + \theta(3-e)\} = \min\{e, \mathfrak{a} + \Big(\frac{k}{3}- \frac{\mathfrak{a}}{3}\Big)(3-e)\}.
		\end{split}
	\end{equation}
\end{proof}

Next, we bound the part $\mathcal{B}$.
\begin{lemma}\label{lem_B decay} 
		\[
	\int_{E_{v,\eta}} f(t, x, w) \mathbf{1}_{|w| \geq \frac{|v|}{\sqrt{2}}}\,dE_{w}	\leq 
	\frac{C_{F_0, A_0}}{\avv^{\b-2}} +  \frac{C_{F_0, k,A_0}}{\langle {v} \rangle^{\mathfrak{b}_1 + k +\g}}\big( \|f\|_{L^\infty_{\mathfrak{a}}} + A_0 \big)
	\]
    for almost everywhere $\eta$, where $\mathfrak{b}_1$ is given by
    \[
	   \mathfrak{b}_1 := \min\Big\{1-\g, \frac{1-\g}{3}\mathfrak{a} + \frac{k(2+\g)}{3} \Big\}.
	\]
    So, $\mathcal{B}$ is bounded by the right-hand side.
\end{lemma}
\begin{proof}
	We start from the estimate \eqref{integ_eq_3} in the proof of Lemma~\ref{lem unif est B}. For $(**)_{\g}$, we apply Lemma~\ref{lem_A decay} for $e=1-\gamma$ to obtain
	\begin{equation}\label{est E_0^2}
		\begin{split}
			(**)_{\g} &\leq \int_{|v_1|\geq \frac{|{v}|}{2}} \int_{u_1\in\R^{3}} \frac{1}{|v_1-u_1|^{1-\g}} f(s, x - v'_1 (t-s), u_1) f(s, x - v'_1 (t-s), v_1)\,du_1 dv_1 \\
			&\leq \int_{|v_1|\geq |v|/2} \sup_x f(s,x, v_1) dv_1 \sup_{|v_1|\geq |v|/2}\int_{u_1} \sup_x f(s,x, u_1)\frac{1}{|u_1-v_1|^{1-\g}}\,du_1 \\
			&\leq \frac{C_k}{\avv^k} \int_{|v_1|\geq |v|/2} \sup_x f(s,x, v_1) \langle v_1 \rangle^k dv_1 \times  \big( \|f\|_{L^\infty_{\mathfrak{a}}} + A_0 \big) \frac{1}{\avv^{\mathfrak{b}_1}}.
		\end{split}
	\end{equation}
	Here,
	\[
		\mathfrak{b}_1 := \min\Big\{1-\g, \frac{1-\g}{3}\mathfrak{a} + \frac{k(2+\g)}{3} \Big\}.
	\]
	\\
	Applying \eqref{est E_0^2} to \eqref{integ_eq} for $E = E_{v,\eta}$, we get
    \begin{equation}\label{Beta est gamma}
        \begin{split}
            &\int_{\R^3} \phi^{E}_{j}(v_1)\psi_{v}(v_1)  f(t, x, v_1)\,dv_1  \\
            &\leq 
            e^{-\lambda_{F,pl} \langle {v}^{\g} \rangle t} \int_{\R^3} \phi^{E}_{j}(v_1)\psi_{v}(v_1) f_0(x - v_1 t, v_1)\,dv_1  +  
            \frac{C_k}{\avv^{k+\mathfrak{b}_1}}A_0 \big( \|f\|_{L^\infty_{\mathfrak{a}}} + A_0 \big) 
            \int_{0}^{t} e^{-\lambda_{F,pl} \avv^{\g} (t-s)}\,ds\\
            &\leq e^{-\lambda_{F,pl} \langle {v}^{\g} \rangle t} \int_{\R^3} \phi^{E}_{j}(v_1)\psi_{v}(v_1) f_0(x - v_1 t, v_1)\,dv_1  +  
            \frac{C_k}{\lambda_{F,pl} \avv^{k+\mathfrak{b}_1+\g}}A_0 \big( \|f\|_{L^\infty_{\mathfrak{a}}} + A_0 \big).
        \end{split}
    \end{equation}    
        Taking limit in $j$, we obtain
        \begin{equation*} 
            \int_{E} \psi_{{v}}(v_1) f(t, x, v_1) dE(v_1)\leq \frac{C_{F_0, A_0}}{\avv^{\b-2}} +  \frac{C_k A_0}{\lambda_{F,pl}\langle {v} \rangle^{\mathfrak{b}_1 + k + \g}}\big( \|f\|_{L^\infty_{\mathfrak{a}}} + A_0 \big),
		\end{equation*}
		where we used \eqref{initial_integ} for initial term treatment. \\
\end{proof}

Now, we are ready to prove Proposition~\ref{prop uniform bound}.
\begin{proof} [Proof of Proposition~\ref{prop uniform bound}]
	We go back to \eqref{AB}. We use Lemma~\ref{lem_A decay} and Lemma~\ref{lem_B decay} for $\mathcal{A}$ and $\mathcal{B}$ with uniform $L^\infty_{\mathfrak{a}}$ estimate of Lemma~\ref{lem_low order}. Then, we obtain
	\begin{equation*}
		\begin{split}
			&\frac{d}{dt}f(t, x+vt, v) + \lambda_{F}\langle v \rangle^{\g} f(t, x+vt,v) 
			\leq C_{F_0, k,A_0, \lambda_{F,pl}^{-1}}\frac{1}{\avv^\mathfrak{b}}
			\Big[
			\frac{1}{\avv^{\b-2}} +  \frac{1}{\langle {v} \rangle^{\mathfrak{b}_1 + k +\g}}\big( \|f\|_{L^\infty_{\mathfrak{a}}} + A_0 \big)
			\Big].
		\end{split}
	\end{equation*}
	Consequently, in the translated time variable,
	\[
	\operatorname*{ess\,sup}_{t\geq0}\|\langle v\rangle^{\beta_1}f(t)\|_\infty
	\leq C_{\infty,F}:=C_{\alpha,\kappa,\beta,k,\gamma}
	\big(F_0, A_0,\varepsilon_L,\lambda_F\big)
	\]
    with
    \begin{equation*}
		\begin{split}
			\b_1 &:= \min\{ \b, \mathfrak{b}+\g+\b-2, \mathfrak{b}+\mathfrak{b}_1 + k + 2\g \} \\
			&= \begin{cases}
				\min\{ \b, \mathfrak{b}+\g+\b-2\},\quad \b \leq \mathfrak{b}_1+k+\g + 2 \\
				\min\{ \b, \mathfrak{b}+\mathfrak{b}_1 + k +2\g \},\quad \b \geq \mathfrak{b}_1+k+\g + 2 ,
			\end{cases}
		\end{split}
	\end{equation*}
	where
	\begin{equation*} 
        \mathfrak{b} := \min\{2-\g, \frac{2-\g}{3}\mathfrak{a} + \frac{k(1+\g)}{3}\},\quad \mathfrak{b}_1 := \min\Big\{1-\g, \frac{1-\g}{3}\mathfrak{a} + \frac{k(2+\g)}{3} \Big\},
	\end{equation*}
	and $\mathfrak{a}:= \min\{\b-2+\g, k+ 2\g\} $. It is easy to check $\b_1 \geq 0$ once $\b \geq 2$.
\end{proof}

\begin{remark} \label{rmk k range-p}
	For some combinations of $\b$ and $k$, we can choose $\b_1=\b$. Suppose $\b \leq \mathfrak{b} + \g + \b - 2$ and $\b \leq \mathfrak{b}_1 + k + \g +2$. This is equivalent to
	\[
		0\leq \mathfrak{b} + \g - 2 = \min\{ 0, \frac{2-\g}{3}\mathfrak{a} + \frac{k(1+\g)}{3} + \g - 2  \} \ \ \text{and} \ \ 
		\b \leq \mathfrak{b}_1 + k + \g +2.
	\] 
	So, we must have
	\[
		\frac{2-\g}{3}\mathfrak{a} + \frac{k(1+\g)}{3} + \g - 2  \geq 0 \ \ \text{and} \ \ 
		\b \leq \min\Big\{1-\g, \frac{1-\g}{3}\mathfrak{a} + \frac{k(2+\g)}{3} \Big\} + k + \g +2.
	\]
	If $\mathfrak{a} = \b - 2 + \g$, i.e, $\b \leq k+ \g +2$, above two conditions yield
	\begin{equation*}
	    \begin{split}
	        5 - \g - \frac{k(1+\g)}{2-\g} \leq \b &\leq \min\{ k+\g+2, k+3, k + \frac{3k - \gamma^2 + 6\gamma + 4}{2+\gamma} \} = k+\g+2.
	    \end{split}
	\end{equation*}
    Therefore, the following is one of the conditions which we can choose $\b_1=\b$ :
    \[
        \max\{ 5 - \g - \frac{k(1+\g)}{2-\g}, 2 \} \leq \beta \leq k+\g+2.    
    \]
    For example, if $k = \b \geq \frac{10}{3}$, then we can choose $\b_1=\b$ for any $0 < \g \leq 1$.
   \\
\end{remark}

\section{Entropy decay}

To study decay of relative entropy $H(F \ \vert \ \mu)$, we split it into microscopic part and macroscopic part
\Be \notag
 	H(F \ \vert \ \mu) = H(F \ \vert \ M_F) + H(M_F \ \vert \ \mu), 
\Ee
where $M_F$ is the local Maxwellian of $F$.

\subsection{Microscopic entropy coercivity and time-integrated relaxation}

The study of $H(F\mid M_F)$ is closely related to Cercignani's conjecture for the spatially homogeneous Boltzmann equation. In particular, Villani~\cite{cercig} established almost exponential decay of $H(F\mid M_F)$ under suitable regularity assumptions and a uniform Gaussian lower bound. The key ingredient is an entropy--entropy dissipation inequality relating $H(F\mid M_F)$ to the entropy dissipation functional
\[
D(F)(t)=\frac14\iiint_{\R^3\times\R^3\times\mathbb S^2}
B(v-v_*,\sigma)(F'_*F'-F_*F)
\log\frac{F'_*F'}{F_*F}.
\]
In this subsection, we restrict our attention to the spatially homogeneous Boltzmann equation. By adapting the kernel-splitting argument of \cite{cercig}, one can obtain at least polynomial decay of $H(F\mid M_F)$ under a uniform Gaussian lower bound and a uniform bound in a suitably weighted $L^\infty$ space, without requiring the same degree of regularity. However, for the purpose of establishing exponential convergence in Theorem~\ref{thm:cond}, an optimal decay rate at this stage is unnecessary. We therefore employ a simpler argument to establish polynomial decay of $H(F\mid M_F)$.

\begin{proposition} \label{prop:micro decay} 
    Let $G$ satisfy 
    \Be \label{F range}
        e^{-a_g|v|^2 + c_g} \leq G(v) \leq C_G \avv^{-\b_2},\quad \forall v\in \R^3
    \Ee
    for some $a_g > 0$, $c_g\in \R$, and $C_G  > 0$. If $\b_2 > 11$, then the entropy dissipation functional of $G$ satisfies
    \Be \label{DG ode}
       D_{r^\g}(G) \geq C_{C_G, a_g, c_g, \b_2} K(G)^{\frac{3}{2}} H(G \ \vert \ M_{G})^{\frac{3}{2}}, 
    \Ee
    where 
    \[
        D_{r^{\g}}(G) := \int_{v, v_*, \sigma} |v-v_*|^\g (GG_* - G'G'_*)\log \frac{GG_*}{G'G'_*}
    \]
    and constant $K(G)$ depends on $K_B$, macroscopic $(\rho,T)$, and upper bound of $H(G)=\int_{v} G\log G$. 
\end{proposition}
\begin{proof}
First we recall Cercignani's conjecture for super quadratic kernel
\[
    B(v-v_*, \sigma) \geq K_{B}(1 + |v-v_*|^2)
\]
for some $K_B > 0$. Then entropy dissipation function $D(G)$ satisfies
\Be \label{super cercig}
    D_{1+r^2}(G)  \geq K(G)H(G \ \vert \ M_G) \quad \text{where}\quad r = |v-v_*|
\Ee
for a constant $K(G)$ depending only on $K_B$, macroscopic $(\rho, T)$ and upper bound of $H(G)$. We refer to \cite{cercig} for the proof.

For $p > 1$, applying H\"older inequality, 
\Be \label{D holder}
    D_{1+r^2}(G) \leq D_{r^\g}(G)^{\frac{1}{p}} D_{W_p(r)}(G)^{\frac{p-1}{p}},
\Ee
where $W_p(r) = (1+r^2)^{\frac{p}{p-1}} r^{-\frac{\g}{p-1}}$. From \eqref{F range}, we have
\[
    |\log G(v)| \leq C_{C_G, a_g, c_g, \b_2} \avv^{2}
\]
and thus
\begin{equation*}
    \begin{split}
        D_{W_p(r)}(G) &\leq C_{C_Ga_g, c_g, \b_2} \int_{v,v_*}W_p(|v-v_*|)G G_* (1 + \avv^2 + \avvv^2). 
    \end{split}
\end{equation*}
For small velocity,
\[
    W_p(r) \sim r^{-\frac{\g}{p-1}},
\]
so it is locally integrable for 
\Be \label{cond p}
    p >  1 +\frac{\g}{3}.
\Ee
For large velocity,
\[
    W_p(r) \sim r^{\frac{2p-\g}{p-1}}.
\]
To control the whole integration with the upper bound of \eqref{F range}, we need
\Be \label{cond b2}
    \b_2 > 3 + 2 + \frac{2p-\g}{p-1} = 7 + \frac{2-\g}{p-1} \quad \text{for}\quad 0\leq \g \leq 1.
\Ee
For simplicity, we choose $p = \frac{3}{2}$ and then for $\b_2 > 11$, uniform boundedness $D_{W_p(r)}(G) \leq C_{C_G, a_g, c_g, \b_2}$ holds.
Combining with \eqref{super cercig} and \eqref{D holder}, we obtain
\[
    D_{r^\g}(G) \geq C_{C_G, a_g, c_g, \b_2} K(G)^{\frac{3}{2}} H(G  \ \vert \ M_G)^{\frac{3}{2}}.
\]
As mentioned before, $K(G)$ depends on $K_B$, macroscopic $(\rho,T)$, and upper bound of $H(G)=\int_{v} G\log G$. 
\end{proof}

To apply Proposition~\ref{prop:micro decay} to the spatially inhomogeneous Boltzmann equation \eqref{BE} with initial datum \eqref{initial}, we need to verify that $F$ satisfies \eqref{F range} and that the constant $K(F)$ can be chosen independently of $x$.

We first consider the stretched-exponential setting. By Proposition~\ref{prop:glb}, we may assume, without loss of generality, that the initial datum $F_0$ admits a Gaussian lower bound. Combining Proposition~\ref{prop uniform bound}, \eqref{apriori}, and Proposition~\ref{prop: macro lb}, we obtain the macroscopic bounds required to apply the Gaussian lower-bound result of \cite{CPDE}. Consequently, we have the following proposition.

\begin{proposition} 
\label{prop:global glb}
There exist constants $a_g>0$ and $c_g\in\R$, depending on $F_0$ and $A_0$, such that
\[
F(t,x,v)\geq e^{-a_g|v|^2+c_g}
\]
for all $(t,x,v)\in[0,\infty)\times\mathbb T^3\times\R^3$.
\end{proposition}

Recall that, as explained at the end of Section~3, we may assume without loss of generality that the initial datum $F_0$ admits a Gaussian lower bound by shifting the initial time to a positive time at which such a lower bound has been established. This reduction applies to both the stretched-exponential and polynomially weighted settings. 

Moreover, Proposition~\ref{prop uniform bound} provides the uniform upper bound
\[
F(t,x,v)\leq C_{F_0,A_0}\avv^{-\b_2}.
\]
Hence, together with the Gaussian lower bound, $F$ satisfies \eqref{F range} on the corresponding time interval.

We next verify that $K(F)$ can be chosen independently of $x$. The macroscopic bounds obtained above yield
\[
0<\rho_-\leq\rho(t,x)\leq\rho_+,
\qquad T(t,x)\leq T_+.
\]
Furthermore, Proposition~\ref{prop uniform bound} and \eqref{F range} imply
\[
\sup_{t,x}\int_{\R^3}F(t,x,v)\log F(t,x,v)\,dv
\leq C_{F_0,A_0,a_g,c_g}.
\]
Applying Lemma~\ref{T lower}, we obtain a uniform positive lower bound on $T(t,x)$. Consequently, all the relevant macroscopic quantities are uniformly controlled, and $K(F)$ can be chosen independently of $x$.

We may therefore apply Proposition~\ref{prop:micro decay} to obtain
\[
\begin{aligned}
&C_{C_F,a_g,c_g,\b_2}K(F)^{3/2}
\int_0^T H(F\mid M_F)^{3/2}(t)\,dt\\
&\qquad\leq \int_0^T\mathcal D(F)(t)\,dt\\
&\qquad\leq H(F_0)-H(F(T))
\end{aligned}
\]
for every $T>0$ in the stretched-exponential setting. In the polynomially weighted setting, the same estimate holds with the lower limit of integration replaced by any fixed $\tau>0$ and the corresponding entropy difference.

Since the entropy is bounded from below, it follows that
\[
\int_\tau^\infty H(F\mid M_F)^{3/2}(t)\,dt<\infty.
\]
Consequently, for any sequence $t_n\to\infty$,
\[
\int_{t_n}^{t_n+1}H(F\mid M_F)^{3/2}(t)\,dt
\longrightarrow0
\qquad\text{as }n\to\infty.
\]
Finally, Hölder's inequality and the Csiszár--Kullback--Pinsker inequality
\[
\|F(t)-M_F(t)\|_{L^1_{x,v}}^2
\leq 2H(F\mid M_F)(t)
\]
give
\[
\begin{aligned}
&\int_{t_n}^{t_n+1}
\|F(t)-M_F(t)\|_{L^1_{x,v}}\,dt\\
&\qquad\leq
\sqrt{2}
\left(
\int_{t_n}^{t_n+1}
H(F\mid M_F)^{3/2}(t)\,dt
\right)^{1/3}
\longrightarrow0.
\end{aligned}
\]
Equivalently, writing $F_n(t,x,v)=F(t_n+t,x,v)$, we obtain
\begin{equation} \label{fMf conv}
\int_0^1\|F_n-M_{F_n}\|_{L^1_{x,v}}\,dt
\longrightarrow0
\qquad\text{as }n\to\infty.
\end{equation}
  
 \subsection{Decay of Macroscopic part}
 In this subsection, we prove $L^1$ decay and improve it to the decay of relative entropy.
 \begin{proposition} 
 \label{prop:rel-entro decay}
 	We have decay of relative entropy :
 	\Be \label{entropy decay}
 		H(F \ \vert \ \mu)(t) \longrightarrow 0\quad\text{as}\quad t\longrightarrow \infty. 
 	\Ee
 \end{proposition}
 \begin{proof}
 We consider countable set of test functions (for $\mathbf{1}_{R}:= \mathbf{1}_{|v|\leq R}(v)$)	
 \[
 \mathcal{D} := \{ \mathbf{1}_{R}(v), v_1 \mathbf{1}_{R}(v), v_2 \mathbf{1}_{R}(v), v_3 \mathbf{1}_{R}(v), |v|^2 \mathbf{1}_{R}(v) \ : \ R\in \mathbb{N}\}.
 \]
 
 Using $L^1$ velocity averaging, countability of $\mathcal{D}$, and diagonal argument, we can derive one common subsequence $\{n_k\}$ such that
 \begin{equation}
 	\begin{split}
 		\rho_{n_k}(t,x) \rightarrow \rho_*(t,x)\quad L^1((0,1)\times \mathbb{T}^3),\\
 		q_{n_k}(t,x) \rightarrow q_*(t,x)\quad L^1((0,1)\times \mathbb{T}^3),\\
 		E_{n_k}(t,x) \rightarrow E_*(t,x)\quad L^1((0,1)\times \mathbb{T}^3),
 	\end{split}
 \end{equation}
 where 
 \[
    \rho_n := \int_v F_n dv, \ \ q_n :=  \int_v F_n v dv, \ \ \text{and} \ \ E_n := \int_v F_n|v|^2.
 \]
 In conclusion, abusing sequence labelling into $n$ again, we obtain
 \begin{equation*}
 	\begin{split}
 		\rho_n  &\longrightarrow \rho_* \quad \text{in}\quad L^{1}([0,1]\times\mathbb{T}^3),\\
 		\rho_n U_n &\longrightarrow q_*:= \rho_* U_* \quad \text{in}\quad L^{1}([0,1]\times\mathbb{T}^3),\\
 		\rho_n |U_n|^2 + 3\rho_n T_n &\longrightarrow E_*:= \rho_* |U_*|^2 + 3\rho_* T_* \quad \text{in}\quad L^{1}([0,1]\times\mathbb{T}^3). \\
 	\end{split}
 \end{equation*}
 
 Recall the uniform bounds $\rho\geq\rho_->0$ and $\rho\leq\rho_+$, $|U|\leq U_+$, $T\leq T_+$ on $t\geq\tau$. These follow from the macroscopic coercivity and the local moment bounds in either of the two settings. Set $\sup_{t\geq 0}\|F(t)\|_{L^\infty_{x,v}} \leq M < \infty$. Choose $r>0$ sufficiently small that
 \[
 \int_{|v-U|<r}F\,dv\leq \sup_{t\geq 0} M |B_r|\leq\rho_-/2.
 \]
And then, 
 \[
 3\rho T = \int_{\mathbb{R}^3} |v-U|^2 F\,dv \geq r^2\int_{|v-U|\geq r} F\,dv  \geq r^2 \frac{\rho_{-}}{2},
 \]
 which gives 
 \[
 T(t,x) \geq T_{-} := \frac{r^2 \rho_{-}}{6\rho_{+}}.
 \]
 
From upper and lower bounds $\rho_{-}, \rho_{+}, U_{+}, T_{+}, T_{-}$, we obtain
\begin{equation*}
\begin{split}
    \rho_n  &\longrightarrow \rho_* \quad \text{in}\quad L^{1}([0,1]\times\mathbb{T}^3),\\
    U_n &\longrightarrow U_* := \frac{q_*}{\rho_*} \quad \text{in}\quad L^{1}([0,1]\times\mathbb{T}^3),\\
    T_n &\longrightarrow T_*:= \frac{E_{*}-\rho_* |U_*|^2}{3\rho_*} \quad \text{in}\quad L^{1}([0,1]\times\mathbb{T}^3),\\
\end{split}
\end{equation*}
which are well-defined. These convergences also hold almost everywhere in $(t,x)$ (up to subsequence if necessary).
Since $M_F(v)\le \frac{\rho_+}{(2\pi T_-)^{3/2}} \exp\!\left(\frac{U_+^2}{2T_+}\right) \exp\!\left(-\frac{|v|^2}{4T_+}\right)$, the  dominated convergence theorem yields 
\Be \label{Mf M* converge}
    M_{F_n} \longrightarrow M_* := \frac{\rho_*}{\sqrt{2\pi T_*}^3} e^{-\frac{|v-U_*|^2}{2T_*}}	\quad \text{in}\quad L^{1}([0,1]\times\mathbb{T}^3_{x}\times \R^{3}_{v}).
\Ee
Combining with \eqref{fMf conv}, we obtain
\Be \label{fM* conv}
    \| F_n - M_* \|_{L^{1}([0,1]\times\mathbb{T}^3_{x}\times \R^{3}_{v})} \longrightarrow 0.
\Ee
We now apply the time-shift identification proved in Proposition~\ref{Desvillettes}.
Together with \eqref{fM* conv} and the conservation laws for $F$, it implies $M_*=\mu$. 
Up to subsequence,
\[
    \| F_{n} - \mu\|_{L^{1}([0,1]\times\mathbb{T}^3_{x}\times \R^{3}_{v})} \longrightarrow 0.
\]
So, there exists a sequence $\{s_n\} \subset  [0,1]$ such that 
\[
    \| F(t_n + s_n) - \mu\|_{L^{1}(\mathbb{T}^3_{x}\times \R^{3}_{v})} \longrightarrow 0.
\]
From \eqref{uniform est}, Lemma \ref{lem_reverseCKP}, and Kullback-Pinsker inequality, with $\sup_{t\geq 0}\|F(t)\|_{L^\infty_{x,v}} \leq M < \infty$,
\begin{equation}\label{L1entropy}
    \begin{split}
        H(F\mid\mu)(s)
        &\leq C_{M}\|F(s)-\mu\|_{L^1_{x,v}(\langle v\rangle^2)}\\
        &\leq C_{M}
        \|F(s)-\mu\|_{L^1_{x,v}(\langle v\rangle^4)}^{1/2}
        \|F(s)-\mu\|_{L^1_{x,v}}^{1/2}\\
        &\leq C_{M}C_{F,\infty}\|F(s)-\mu\|_{L^1_{x,v}}^{1/2}.
    \end{split}
\end{equation}
At the selected times $s=t_n+s_n$, the right-hand side tends to zero. Finally, the entropy inequality
\[
    H(F(t)\mid\mu)+\int_s^t\mathcal D(F)(u)\,du
    \leq H(F(s)\mid\mu),\qquad t\geq s
\]
implies that $H(F(t)\mid\mu)\to0$ as $t\to\infty$.
\end{proof}

\bigskip

\section{Exponential convergence to global equilibrium}
In this section, we solve large amplitude problem with stretched exponential weight. 

\subsection{Perturbative setting}
Setting $F(t,x,v)=\mu(v) +f(t,x,v)$, the equation \eqref{BE} becomes the perturbed Boltzmann equation:
\begin{align} \label{FPBE}
	\partial_tf+v \cdot \nabla_x f +L(f)=Q(f,f).
\end{align}
On $\mathbb{T}^3$, the mass, momentum, energy conservations hold:
\begin{align} \label{conservationlaw}
	\int_{\mathbb{T}^3 \times \R^3} f(t,x,v)\begin{pmatrix}
		1  \\
		v  \\
		|v|^2
	\end{pmatrix} dvdx =0.
\end{align}

We will consider two perturbative settings. First of all, we set $F=\mu+\mu^{1/2}u$. (This setting will be used only in Section \ref{Lineardecayanalysis}.)  Then the Boltzmann equation \eqref{BE} can be written in terms of $u$ as (Here, Since the following perturbed equation is used only in the linear theory, $\Gamma$ term is unnecessary.)
\begin{align*}
	\partial_t u+v \cdot \nabla_x u +\mathcal{L}u=\Gamma(u,u).
\end{align*} 
Here, the linear operator $\mathcal{L}$ can be decomposed into $\mathcal{L}=\nu-\mathcal{K}$, where
\begin{align*}
	&\nu(v)u= \mu^{-1/2}Q^-(\mu^{1/2}u,\mu) \quad (\text{i.e.} \quad \nu(v)=\int_{\R^3} \int_{\S^2} B(v-v_*,\sigma)\mu(v_*) d\sigma dv_* ),\\
	& \mathcal{K}u = \mu^{-1/2}Q^+(\mu^{1/2}u,\mu)+ \mu^{-1/2}Q^+(\mu,\mu^{1/2}u)-\mu^{-1/2}Q^-(\mu,\mu^{1/2}u),
\end{align*}
and
\begin{align*}
	\Gamma(u,u)=\mu^{-1/2}Q(\mu^{1/2}u,\mu^{1/2}u).
\end{align*}
By Carleman representation \cite{Villani2002,Carleman1957}, we have
\begin{align*}
	\mathcal{K}u(v) = \int_{\R^3} \mathbf{k}(v,v_*)u(v_*)dv_*.
\end{align*}
It is well-known that from Guo \cite{Guo2010ARMA}, the kernel $\mathbf{k}$ is bounded by
\begin{align*}
	|\mathbf{k}(v,v') | \le C \{|v-v'|+|v-v'|^{-1}\}e^{-\frac{1}{8}|v-v'|^2-\frac{1}{8}\frac{||v|^2-|v'|^2|^2}{|v-v'|^2}}.
\end{align*}
Next, we set $F =\mu+f$. (This perturbation will be used in Section \ref{smallperturbationscheme}, \ref{largeamplitudescheme}.) Then the Boltzmann equation \eqref{BE} can be written in terms of $f$ as 
\begin{align} \label{perturbedBE} 
	\partial_t f+v \cdot \nabla_x f +Lf=Q(f,f).
\end{align} 
Here, the linear operator $L$ can be decomposed into $L=\nu-K$, where
\begin{align*}
	&\nu(v)f= Q^-(f,\mu) \quad (\text{i.e.} \quad \nu(v)=\int_{\R^3} \int_{\S^2} B(v-v_*,\sigma)\mu(v_*) d\sigma dv_* ),\\
	& Kf = Q^+(f,\mu)+ Q^+(\mu , f)-Q^-(\mu,f),
\end{align*}
We define the constant $\nu_0$ by
\begin{align} \label{nuzerodef}
\nu_0 := \inf_{v} \frac{\nu(v)}{\langle v \rangle^\gamma}.
\end{align}
By Carleman representation \cite{Villani2002,Carleman1957}, we have
\begin{align*}
	Kf(v) = \int_{\R^3} k(v,v_*)f(v_*)dv_*.
\end{align*}
We use the weight function $w(v)=\langle v \rangle^\beta e^{\alpha |v|^\kappa}$ and set $h=wf$. Then the perturbed equation \eqref{perturbedBE} for $h$ becomes
\begin{align*}
	\partial_t h+v \cdot \nabla_x h +\nu(v)h=K_w h+ wQ(f,f).
\end{align*}
Here, $K_w$ has a kernel $k_w$:
\begin{align*}
	k_w(v,v_*) = k(v,v_*)\frac{w(v)}{w(v_*)}.
\end{align*}
By rotation, we also have
\begin{align*}
	\int_{\R^3} k_w(v,v_*) dv_* &= 2\int_{\R^3}\int_{\S^2} b(\cos \theta)|v-v_*|^\gamma \frac{\langle  v \rangle^\beta e^{\alpha |v|^\kappa}}{\langle v'\rangle^\beta e^{\alpha |v'|^\kappa}} e^{-\frac{|v_*'|^2}{2}}  d\sigma d v_* \\
	& \quad - \int_{\R^3}\int_{\S^2} b(\cos \theta)|v-v_*|^\gamma \frac{\langle  v \rangle^\beta e^{\alpha |v|^\kappa}}{\langle v_*\rangle^\beta e^{\alpha |v_*|^\kappa}} e^{-\frac{|v|^2}{2}}  d\sigma d v_*.
\end{align*}

\bigskip

We state a Theorem for large amplitude problem in stretched exponential weight.

\begin{theorem} \label{Mainresultlargeamp}
	Let \(w(v)=(1+|v|^2)^{\beta/2}e^{\alpha|v|^{\kappa}}\) with $\beta\ge0$, $\alpha>0$, and $1<\kappa<2$. Assume that $F_0(x,v) = \mu(v) +  f_0(x,v) \ge 0$ satisfying 
	\begin{align*}
		\int_{\T^3 \times \R^3}f_0 
		\begin{pmatrix}
			1\\
			v\\
			|v|^2
		\end{pmatrix}
		dxdv =0.
	\end{align*}
	For any $M_0 \ge 1$, there exists a constant $\varepsilon_0>0$, depending only on $M_0$ and $\bar{M}$, such that if
	\begin{align*}
		\|wf_0\|_{L^\infty_{x,v}}\le M_0, \quad \sup_{t\geq 0}\|wf(t)\|_{L^\infty_{x,v}}\leq \bar{M}, \quad H(F_0 \ \vert \ \mu) \le  \varepsilon_0,
	\end{align*}
	then the solution $F(t,x,v) =\mu(v) + f(t,x,v) \ge 0 $ of the Boltzmann equation \eqref{BE} with initial data $F_0$ satisfies that 
	\begin{align*}
		\|wf(t)\|_{L^\infty_{x,v}} \le C_0(1+\bar{M}) M_0e^{-\lambda_\infty t}
	\end{align*}
	for all $t \ge 0$, where $C_0\ge 1$ is a generic constant. 
\end{theorem}
\begin{remark}
    In Proposition \ref{prop uniform bound}, we obtained a uniform bound $C_{\infty, F}$ for $\|wf(t)\|_{L^\infty_{x,v}}$. Therefore, we choose $M_0= \bar{M} = C_{\infty,F}$ in the proof of Theorem \ref{thm:cond} later. 
\end{remark}

\bigskip

To prove Theorem \ref{Mainresultlargeamp}, we first establish the following small-amplitude stability result. We will then show that, under the assumptions of the theorem, the solution enters the small amplitude regime after a finite time.

\begin{proposition} \label{Mainresultsmall} 
	Let \(w(v)=(1+|v|^2)^{\beta/2}e^{\alpha|v|^{\kappa}}\) with $\beta \ge 0 ,\alpha >0 ,0<\kappa<2$. Assume that $F_0(x,v) = \mu(v) +  f_0(x,v) \ge 0$ satisfying 
	\begin{align*}
		\int_{\T^3 \times \R^3}f_0 
		\begin{pmatrix}
			1\\
			v\\
			|v|^2
		\end{pmatrix}
		dxdv =0.
	\end{align*}
	Then there exists $\eta_0 >0$ such that for any $f_0$ in $L^\infty_{x,v}(w)$ with $\|wf_0\|_{L^\infty_{x,v}} \le \eta_0^{1+p}$ for some $0<p<1$, there is a unique solution $F(t,x,v) =\mu(v) +  f(t,x,v) \ge 0 $ of the Boltzmann equation \eqref{BE} with initial data $F_0$, which satisfies the mass, momentum, energy conservation and there exist $C_\infty,\lambda_\infty >0$ such that
	\begin{align*}
		\|wf(t)\|_{L_{x,v}^\infty}\le C_\infty e^{-\lambda_\infty t}(\|wf_0\|_{L_{x,v}^\infty}+\eta_0^{1+p})
	\end{align*}
	for all $t \ge 0$.\\
\end{proposition}

\bigskip

\subsection{Linear decay analysis} \label{Lineardecayanalysis}
In this section, we consider the (perturbed) linearized Boltzmann equation of $f \in L^2_{x,v}(\mu^{-1/2})$:
\begin{align} \label{LBE}
	\partial_tf+v \cdot \nabla_x f +Lf=0.
\end{align}
Setting $u(t,x,v) =\mu^{-1/2}(v)f(t,x,v) \in L^2_{x.v}$, the equation \eqref{LBE} can be written in terms of $u$:
\begin{align} \label{wLBE}
	\partial_t u+v \cdot \nabla_x u +\mathcal{L}u=0. 
\end{align}
Here, we can decompose $\mathcal{L} = \nu -\mathcal{K}$, and as in \cite{Guo2010ARMA}, the linear operator $\mathcal{K}$ can be expressed as
\begin{align*}
	\mathcal{K}u(v) = \int_{\R^3} \mathbf{k}(v,v_*)u(v_*)dv_*.
\end{align*}
We define the $L^2_{x,v}$ projection $P_{\mathcal{L}}f$ of $f$ corresponding to $\mathcal{L}$ as
\begin{align*}
	P_{\mathcal{L}}f = a(t,x)\mu^{1/2}(v)+b(t,x)\cdot v\mu^{1/2}(v)+c(t,x)\frac{|v|^2-3}{\sqrt{6}} \mu^{1/2}(v),
\end{align*}
where
\begin{align*}
	&a(t,x) = \int_{\R^3}f(t,x,v)\mu^{1/2}(v)dv,\\
	&b(t,x) = \int_{\R^3}vf(t,x,v)\mu^{1/2}(v)dv,\\
	&c(t,x) = \int_{\R^3}\frac{|v|^2-3}{\sqrt{6}}f(t,x,v)\mu^{1/2}(v)dv.
\end{align*}

\bigskip

We recall the exponential $L^2$ decay estimate for the linearized Boltzmann equation. This estimate will serve as the starting point for the weight $L^\infty$ decay analysis.

\begin{proposition} \label{L2decay} \cite{Guo2010ARMA,Briant2016}
	Let $u_0(x,v)$ be in $L^2_{x,v}$ satisfying 
	\begin{align*}
		\int_{\T^3 \times \R^3} \mu^{1/2}(v)u_0(x,v)\begin{pmatrix}
			1  \\
			v  \\
			|v|^2
		\end{pmatrix}dvdx =0.
	\end{align*}
	Suppose that $u(t,x,v) \in L^2_{x,v}$ be a solution to \eqref{wLBE} with initial data $u_0$, and satisfying the mass, momentum, energy conservation
	\begin{align*}
		\int_{\T^3 \times \R^3} \mu^{1/2}(v)u(t,x,v)\begin{pmatrix}
			1  \\
			v  \\
			|v|^2
		\end{pmatrix}dvdx =0.
	\end{align*}
	Then there exist $C_{\mathcal{L}}$, $\lambda_{\mathcal{L}}>0$, independent of $u_0$ and $u$, such that
	\begin{align*}
		\|u(t)\|_{L^2_{x,v}} \le C_{\mathcal{L}}e^{-\lambda_{\mathcal{L}} t} \|u_0\|_{L^2_{x,v}}
	\end{align*}
	for all $t \ge 0$.
\end{proposition}

\bigskip

We denote by $G(t,s)$ the solution operator to the linearized Boltzmann equation \eqref{LBE}. The following lemma establishes the  $L^\infty$ exponential decay in time for solutions to the linearized Boltzmann equation.
\begin{proposition}\label{linearLinftydecay} \cite{Guo2010ARMA}
	Let $f_0 \in L^\infty_{x,v}(\langle v\rangle^\beta \mu^{-1/2})$, $\beta>3$ which satisfies
    \begin{align*}
		\int_{\T^3 \times \R^3}f_0(x,v)\begin{pmatrix}
			1  \\
			v  \\
			|v|^2
		\end{pmatrix}dvdx =0.
	\end{align*} There exists a unique solution $f(t)=G(t,0)f_0 \in L^\infty_{x,v}(\langle v\rangle^\beta \mu^{-1/2})$ to the linearized Boltzmann equation \eqref{LBE} with initial data $f_0$. Moreover, there exist $\lambda_L, C_L>0$ such that
	\begin{align*}
		\|G(t,0)f_0\|_{L_{x,v}^\infty(\langle v\rangle^\beta \mu^{-1/2})} \le C_L e^{-\lambda_L t}\|f_0\|_{L^\infty_{x,v}(\langle v\rangle^\beta \mu^{-1/2})}
	\end{align*}
	for all $t\ge 0$.
	The constants $C_L$ and $\lambda_L$ depend on the collision kernel $B$, $\lambda_\mathcal{L}$, and $C_\mathcal{L}$.
\end{proposition}

\bigskip

\subsection{Nonlinear $L^\infty_{x,v}$ decay in small perturbation problem} \label{smallperturbationscheme}

The analysis developed in this section is inspired by the framework of Gualdani-Mischler-Mouhot \cite{Gualdani2017}, together with the related perturbative approaches in \cite{Briant2016,Briant2016DCDS,Briant2017}. More specifically, we follow the spirit of \cite{Gualdani2017} by separating the dissipative mechanism of the linearized dynamics from the regularizing effect of a suitably truncated collision operator. In the present setting, we implement this strategy through a coupled decomposition of the nonlinear perturbation. We consider the full nonlinear perturbed equation
\begin{align} \label{BE-f}
	\partial_tf+v \cdot \nabla_x f +Lf=Q(f,f) 
\end{align}
In what follows, we consider the weight function $w(v)=(1+|v|^2)^{\beta/2}e^{\alpha|v|^{\kappa}}$.\\
We define the projection in $L^2_{x,v}(\mu^{-1/2})$ onto the kernel of $G_L:=-v\cdot \nabla_x-L$ by
\begin{align*}
	P_{G_{L}}(f) &= \left(\int_{\T^3 \times \mathbb{R}^3}f(x,v)dxdv\right)\mu(v)+\left(\int_{\T^3 \times \mathbb{R}^3}vf(x,v)dxdv\right)\cdot v\mu(v)\\
	& \quad +\left(\int_{\T^3 \times \mathbb{R}^3}\frac{|v|^2-3}{\sqrt{6}}f(x,v)dxdv\right)\frac{|v|^2-3}{\sqrt{6}}\mu(v).
\end{align*}
Note that $P_{G_{L}}(f)=0$ means that $f$ satisfies the mass, momentum, and energy conservation.\\
\indent We follow the decomposition proposed in \cite{Gualdani2017}. For $0<\delta<1$, we define $\Theta_\delta=\Theta_\delta(v,v_*,\sigma)$ in $C^\infty$ that is one on the set
\begin{align*}
	\{|v|\le \delta^{-1} \ \text{and} \ 2\delta \le |v-v_*|\le \delta^{-1} \ \text{and} \ |\cos \theta|\le 1-2\delta \},
\end{align*}
and whose support is contained in 
\begin{align*}
	\{|v|\le 2\delta^{-1} \ \text{and} \ \delta \le |v-v_*|\le 2\delta^{-1} \ \text{and} \ |\cos \theta|\le 1-\delta \}.
\end{align*}
This truncation localizes the collision operator to bounded velocities and relative velocities, away from vanishing relative velocity and the angular endpoints. Then we can split the operator $K$ into two operators with singular part and nonsingular part :
\begin{align*}
	K = K_{ns}^{(\delta)}+ K_s^{(\delta)}
\end{align*}
with
\begin{align*}
	K_{ns}^{(\delta)} f(v) = \int_{\R^3 \times \S^2} \Theta_\delta(v,v_*,\sigma) b(\cos \theta)|v-v_*|^\gamma \left[\mu(v_*')f(v')+\mu(v')f(v_*')-\mu(v)f(v_*)\right]d\sigma dv_*
\end{align*}
and
\begin{align*}
	K_s^{(\delta)} f(v) =\int_{\R^3 \times \S^2} \left(1-\Theta_\delta(v,v_*,\sigma) \right) b(\cos \theta)|v-v_*|^\gamma \left[\mu(v_*')f(v')+\mu(v')f(v_*')-\mu(v)f(v_*)\right]d\sigma dv_*.
\end{align*}
Note that $K_{ns}^{(\delta)}$ has a kernel compactly supported and a regularizing effect. (See Lemma \ref{Ksest}.)
From this decomposition, the equation \eqref{BE-f} with initial data $f_0$ can be decomposed into two equations :
\begin{equation} \label{Decompositionsys}
	\begin{aligned}
		&\partial_t f_1+v\cdot \nabla_x f_1 +\nu(v) f_1= K_s^{(\delta)}f_1+Q(f_1+f_2,f_1+f_2) \quad \text{with} \ f_1(0,x,v)=f_0(x,v), \\
		&\partial_t f_2+v\cdot\nabla_x f_2 +Lf_2 = K_{ns}^{(\delta)}f_1 \quad \text{with} \ f_2(0,x,v)=0,
	\end{aligned}
\end{equation}
where $f =f_1+f_2$.

We establish several estimates that will be useful in the subsequent analysis in the small perturbation regime. Their proofs are motivated by the following works \cite{Gualdani2017,Briant2016DCDS,Briant2017,Briant2016}.
\bigskip

We next establish a weighted $L^\infty$ estimate for $K_{ns}^{(\delta)}$. The regularizing structure of $K_{ns}^{(\delta)}$, together with the compact velocity support of its kernel, yields the following sharper weighted $L^\infty$ estimate.

\begin{lemma} \label{Ksest}
	Let $w(v)=(1+|v|^2)^{\beta/2}e^{\alpha|v|^{\kappa}}$ with $\beta \ge 0 ,\alpha >0 ,0<\kappa<2$  and $\delta \in (0,1)$. Then there exists $C_{K_{ns}}>0$ such that for all f in $L^{\infty}_{x,v}(w)$, 
	\[ \|K_{ns}^{(\delta)}(f)\|_{L^{\infty}_{x,v}(\langle v \rangle^\beta \mu^{-1/2})} \le C_{K_{ns}} \|f\|_{L^{\infty}_{x,v}(w)}. \]
	The constant $C_{K_{ns}}$ only depends on $\delta$.
\end{lemma}
\begin{proof}
	Thanks to Carleman representation \cite{Carleman1957}, the operator $K_{ns}^{(\delta)}$ can be written as
	\begin{align*}
		K_{ns}^{(\delta)}(f)(x,v) = \int_{\mathbb{R}^3}k_{ns}^{(\delta)}(v,v_*)f(x,v_*)dv_*.
	\end{align*}
	Moreover, by definition of $K_{ns}^{(\delta)}$, the kernel $k_{ns}^{(\delta)}$ has compact support which implies the desired estimate. 
\end{proof}
\bigskip

The following lemmas state the weighted $L^\infty$ estimate for the nonlinear operator $Q$. 

\begin{lemma} \label{Q-estimate}
	Let \(w(v)=(1+|v|^2)^{\beta/2}e^{\alpha|v|^{\kappa}}\) with $\beta \ge 0 ,\alpha >0 ,0<\kappa<2$. Then for all f, g in $L^{\infty}_{x,v}(w)$,
	\[\|Q^-(f,g)\|_{L^{\infty}_{x,v}(\nu^{-1} w)} \le C_{Q^-} \|f\|_{L^{\infty}_{x,v}(w)}\|g\|_{L^{\infty}_{x,v}(w)}. \]
\end{lemma}
\begin{proof}
	We know that $Q^-(f,g)=\int_{\mathbb{R}^3\times\mathbb{S}^2}B(v_*,v,\sigma)f(v)g(v_*)d\sigma dv_*$.\\
	\begin{align*}
		\left| \nu^{-1}(v)w(v)Q^-(f,g)\right| &\le \|f\|_{L^{\infty}_{x,v}(w)}\|g\|_{L^{\infty}_{x,v}(w)}\nu^{-1}(v)\left(C_b\int_{\mathbb{R}^3}\frac{|v-v_*|^{\gamma}}{w(v_*)}dv_*\right)\\
		&\le C\left(C_b\int_{\mathbb{R}^3}\frac{(1+|v_*|)^{\gamma}}{w(v_*)}dv_*\right) \|f\|_{L^{\infty}_{x,v}(w)}\|g\|_{L^{\infty}_{x,v}(w)},
	\end{align*}
	where we have used the inequality $|v-v_*|^{\gamma} \le C\nu(v)(1+|v_*|)^{\gamma}$ and the integration in the last term is finite. 
\end{proof}
\bigskip

The following two lemmas are technical lemmas for proving the main lemma~\ref{pointwiseQ+est}. The next lemma states an elementary inequality.
\begin{lemma}  Let $0<\alpha\leq 1$ and $x\geq y> 0$. Then,
	\begin{align*}
		(x+y)^\alpha\leq x^\alpha + \alpha x^{\alpha-1} y = x^\alpha + \alpha \left(\frac{x}{y}\right)^{\alpha-1} y^\alpha.
	\end{align*}
\end{lemma}
\begin{proof}
	Let $t= x/y\geq 1$. Then, the inequality is the same as
	\begin{align*}
		(t+1)^\alpha\leq t^\alpha + \alpha t^{\alpha-1}.
	\end{align*}
	Since $2^\alpha\leq 1+\alpha$ for $0\leq\alpha\leq 1$, by the derivative test, the inequality holds for $t\geq 1$.
\end{proof}
From this lemma, we have
\begin{align*}
	(|v'|^2+|v'_*|^2)^{\frac{\kappa}{2}}&\leq \begin{cases}
		|v'|^\kappa + \frac{\kappa}{2}\left(\frac{|v'|^2}{|v'_*|^2}\right)^{\frac{\kappa}{2}-1}(|v'_*|^2)^{\frac{\kappa}{2}} & |v'_*|\leq |v'|\\
		|v'_*|^\kappa + \frac{\kappa}{2}\left(\frac{|v'_*|^2}{|v'|^2}\right)^{\frac{\kappa}{2}-1}(|v'|^2)^{\frac{\kappa}{2}} & |v'|\leq |v'_*|
	\end{cases}\\
	&\leq \begin{cases}
		|v'|^\kappa + |v'_*|^\kappa + \left(\frac{\kappa}{2}\left(\frac{|v'|^2}{|v'_*|^2}\right)^{\frac{\kappa}{2}-1}-1\right)|v'_*|^\kappa & |v'_*|\leq |v'|\\
		|v'_*|^\kappa + |v'|^\kappa + \left(\frac{\kappa}{2}\left(\frac{|v'_*|^2}{|v'|^2}\right)^{\frac{\kappa}{2}-1}-1\right)|v'|^\kappa & |v'|\leq |v'_*|
	\end{cases}\\
	&\leq \begin{cases}
		|v'|^\kappa + |v'_*|^\kappa + \left(\frac{\kappa}{2}-1\right)|v'_*|^\kappa & |v'_*|\leq |v'|\\
		|v'_*|^\kappa + |v'|^\kappa + \left(\frac{\kappa}{2}-1\right)|v'|^\kappa & |v'|\leq |v'_*|.
	\end{cases}
\end{align*}
By the previous computation, we also have the dyadic decomposition of $(|v'|^2+|v'_*|^2)^{\frac{\kappa}{2}} = (|v|^2+|v_*|^2)^{\frac{\kappa}{2}}$ to $|v|^\kappa$ and $|v_*|^\kappa$. Combining these two inequalities, we get for $0<\kappa\leq 2$,
\begin{align}\label{eq:expbound}
	\begin{split}
		&e^{-\alpha(|v'|^\kappa + |v'_*|^\kappa - |v|^\kappa)}\\
		&\leq e^{-\alpha((|v'|^2+|v'_*|^2)^{\frac{\kappa}{2}} - |v|^\kappa)}\times\begin{cases}
			e^{-\alpha(1-\frac{\kappa}{2})|v'_*|^\kappa} & |v'_*|\leq |v'|,\\
			e^{-\alpha(1-\frac{\kappa}{2})|v'|^\kappa} & |v'|\leq |v'_*|
		\end{cases}\\
		&\leq \left(e^{-\alpha(2^{\kappa/2}-1)|v'+v'_*-v|^\kappa}\mathbf{1}_{|v'+v'_*-v|\geq |v|} + \sum_{n=0}^\infty e^{-\frac{\alpha\kappa/4}{2^{(1-\kappa/2)(n+1)}}|v_*|^\kappa}\mathbf{1}_{2^{-n-1}|v|^2\leq |v_*|^2\leq 2^{-n}|v|^2}\right)\\
		&\qquad \times \begin{cases}
			e^{-\alpha(1-\frac{\kappa}{2})|v'_*|^\kappa} & |v'_*|\leq |v'|,\\
			e^{-\alpha(1-\frac{\kappa}{2})|v'|^\kappa} & |v'|\leq |v'_*|.
		\end{cases}
	\end{split}
\end{align}
Note that $v_* = v'+v'_*-v$. It refines the crude estimate \eqref{weight interac} in Section 5.

\bigskip

The next lemma, from \cite{Park2026ARMA}, provides a sharp bound for an integral that appears in the proof of Lemma~\ref{pointwiseQ+est}.
\begin{lemma}\label{lem:cone} \cite{Park2026ARMA}
	Let $x,y,a\in\mathbb{R}^3$ satisfies
	\begin{align*}
		x\perp y,\quad |x|\leq |a|/4,\quad R\leq |a|/4,\quad |a+x+y|\leq R.
	\end{align*}
	Then the set
	\begin{align*}
		\{x:\{y:y\perp x\}\cap \{y:|a+x+y|\leq R\}\neq \emptyset,\quad |x|\leq |a|/4\}
	\end{align*}
	is contained in a set $C_{a, R}$ satisfying the following: for any isotropic function $h(x) = h(|x|)$,
	\begin{align*}
		\int_{\{|x|\leq |a|/4\}\cap C_{a, R}} h(x)\,dx \leq C\frac{R}{|a|}\int_{0}^{|a|/4}|x|^{2}h(|x|)\,d|x|
	\end{align*}
	for some dimension-dependent constant $C$.
\end{lemma}

\bigskip

Now, we prove the pointwise estimate of $Q^+$ under the weight $w(v) = \langle v\rangle^\beta e^{\alpha |v|^\k}$. It will be used to control the nonlinear part of the $Q^+$ operator in constructing the solution of the Boltzmann equation.
\begin{lemma}\label{pointwiseQ+est}
	Suppose $0<\kappa< 2$ and $0\leq \gamma\leq 1$. Then, for any $\beta$ and $e_1,e_2\geq 0$, we have
	\begin{align*}
		\begin{split}
			&\left|\langle v \rangle^\beta e^{\alpha |v|^\kappa}Q^+(f,f)(v)\right|\\
			&\leq C\min\left\{(1+|v|)^{-\frac{1}{2}+\gamma-\frac{\kappa}{2}+e_1}, (1+|v|)^{\gamma-\kappa + e_2}\right\}\|(1+|v|^2)^{-e_1/2}\langle v \rangle^\beta e^{\alpha |v|^\kappa}f\|_{L^2_v}\|(1+|v|^2)^{-e_2/2}\langle v \rangle^\beta e^{\alpha |v|^\kappa}f\|_{L^\infty_v}.
		\end{split}
	\end{align*}
	Here, $C$ depends on $\alpha,\kappa,\gamma$, $\beta$, $e_1$, and $e_2$. 
\end{lemma}

Interestingly, one can use a mild polynomial weight by $\langle v\rangle^{-e_1}$ or $\langle v\rangle^{-e_2}$ for $wf$ to bound $wQ^+(f,f)$. Note that $\beta$ can be negative.
\begin{remark}
	In what follows, we will mainly use the above pointwise estimate in the following simplified form:
	\begin{align*}
		|wQ^+(f, f)(v)| \leq C_{\alpha,\kappa,\beta,\gamma}  (1+|v|)^{\gamma-\kappa} \|wf\|_{L^\infty_v}\left(\int_{\R^3}\frac{1}{(1+|\eta|)^6}|wf(\eta)|^2d\eta \right)^{1/2}.
	\end{align*}
	Moreover, from the proof of Lemma \ref{pointwiseQ+est}, we obtain the $L^\infty_{x,v} $ nonlinear estimate for $Q^+$:
	\begin{align} \label{Q+estimate}
		\|wQ^+(f,g)\|_{L^\infty_{x,v}(\langle v \rangle ^{-\gamma+\kappa})}  \le C_{\alpha,\kappa,\beta,\gamma} \|wf\|_{L^\infty_{x,v}}\|wg\|_{L^\infty_{x,v}}.
	\end{align}
\end{remark}
\begin{proof}
	Let $h(t,x,v)=\langle v\rangle^\beta e^{\alpha |v|^\kappa}f(t,x,v)$.
	Using \eqref{eq:expbound} and the Carleman representation, we write
	\begin{align*}
		\begin{split}
			&\left|\langle v \rangle^\beta e^{\alpha |v|^\kappa}Q^+(f,f)(v)\right| \\
			&\leq (1+|v|^2)^{\beta/2}\int_{\mathbb{R}^3}\frac{|h(v')|}{|v-v'|(1+|v'|^2)^{\beta/2}}\int_{v_*'\in v + E_{v'-v}}\frac{|h(v_*')|}{|v-v_*|^{1-\gamma}(1+|v'_*|^2)^{\beta/2}}\,dv_*' dv'\\
			&\qquad \times \left(e^{-\alpha(2^{\kappa/2}-1)|v'+v'_*-v|^\kappa}\mathbf{1}_{|v'+v'_*-v|\geq |v|} + \sum_{n=0}^\infty e^{-\frac{\alpha\kappa/4}{2^{(1-\kappa/2)(n+1)}}|v'+v'_*-v|^\kappa}\mathbf{1}_{2^{-n-1}|v|^2\leq |v'+v'_*-v|^2\leq 2^{-n}|v|^2}\right)\\
			&\qquad \times \begin{cases}
				e^{-\alpha(1-\frac{\kappa}{2})|v'_*|^\kappa} & |v'_*|\leq |v'|\\
				e^{-\alpha(1-\frac{\kappa}{2})|v'|^\kappa} & |v'|\leq |v'_*|
			\end{cases}\,dv'_* dv'.
		\end{split}
	\end{align*}
	
	We divide the summation of the exponential terms by
	\begin{align*}
		e^{-\alpha(2^{\kappa/2}-1)|v'+v'_*-v|^\kappa}\mathbf{1}_{|v'+v'_*-v|\geq |v|} + \sum_{n=0}^7 e^{-\frac{\alpha\kappa/4}{2^{(1-\kappa/2)(n+1)}}|v'+v'_*-v|^\kappa}\mathbf{1}_{2^{-n-1}|v|^2\leq |v'+v'_*-v|^2\leq 2^{-n}|v|^2}
	\end{align*}
	and the remainder. We bound the first part in \textbf{Case 1} and the other part in \textbf{Case 2}.\\
	
	\noindent \textbf{Case 1} For the first case, we note that the exponential decay about $|v_*| = |v'+v'_*-v|$ remains valid. Also, $|v_*|\geq 2^{-4}|v|$ implies the $|v_*|$ decay even can produce some exponential $|v|$ decay. To use these ideas, we crudely take the upper bound of the exponential functions by taking $n=7$:
	\begin{align*}
		&e^{-\alpha(2^{\kappa/2}-1)|v'+v'_*-v|^\kappa}\mathbf{1}_{|v'+v'_*-v|\geq |v|} + \sum_{n=0}^7 e^{-\frac{\alpha\kappa/4}{2^{(1-\kappa/2)(n+1)}}|v'+v'_*-v|^\kappa}\mathbf{1}_{2^{-n-1}|v|^2\leq |v'+v'_*-v|^2\leq 2^{-n}|v|^2}\\
		&\leq e^{-\frac{\alpha\kappa/4}{2^{8-4\kappa}}|v'+v'_*-v|^\kappa}\mathbf{1}_{2^{-4}|v|\leq |v'+v'_*-v|}\\
		&\leq e^{-\frac{1}{2}\frac{\alpha\kappa/4}{2^{8-4\kappa}}(2^{-4}|v|)^\kappa}e^{-\frac{1}{2}\frac{\alpha\kappa/4}{2^{8-4\kappa}}|v'+v'_*-v|^\kappa}\mathbf{1}_{2^{-4}|v|\leq |v'+v'_*-v|}.
	\end{align*}
	Therefore, we readily get
	\begin{align}\label{eq:1}
		\begin{split}
			&(1+|v|^2)^{\beta/2}\int_{\mathbb{R}^3}\frac{|h(v')|}{|v-v'|(1+|v'|^2)^{\beta/2}}\int_{v_*'\in v + E_{v'-v}}\frac{|h(v'_*)|}{|v'-v'_*|^{1-\gamma}(1+|v'_*|^2)^{\beta/2}}\\
			&\qquad \times \left(e^{-\alpha(2^{\kappa/2}-1)|v_*|^\kappa}\mathbf{1}_{|v_*|\geq |v|} + \sum_{n=0}^7 e^{-\frac{\alpha\kappa/4}{2^{(1-\kappa/2)(n+1)}}|v'+v'_*-v|^\kappa}\mathbf{1}_{2^{-n-1}|v|^2\leq |v'+v'_*-v|^2\leq 2^{-n}|v|^2}\right)\,dv' dv'_*\\
			&\leq (1+|v|^2)^{\beta/2}e^{-\frac{1}{2}\frac{\alpha\kappa/4}{2^8}|v|^\kappa}\int_{\mathbb{R}^3}\frac{|h(v')|}{|v-v'|(1+|v'|^2)^{\beta/2}}\int_{v_*'\in v + E_{v'-v}}\frac{|h(v'_*)|}{|v'-v'_*|^{1-\gamma}(1+|v'_*|^2)^{\beta/2}}\\
			&\qquad \times e^{-\frac{1}{2}\frac{\alpha\kappa/4}{2^{(8-4\kappa)}}|v'+v'_*-v|^\kappa}\mathbf{1}_{2^{-4}|v|\leq |v'+v'_*-v|}\,dv' dv'_*\\
			&\leq (1+|v|^2)^{\beta/2}e^{-\frac{1}{2}\frac{\alpha\kappa/4}{2^8}|v|^\kappa}\|(1+|v|^2)^{-e_2/2}h(v)\|_{L^\infty}\int_{\mathbb{R}^3}\frac{h(v')}{|v-v'|(1+|v'|^2)^{\beta/2}}\\
			&\qquad \times \int_{v_*'\in v + E_{v'-v}}\frac{1}{|v'-v'_*|^{1-\gamma}(1+|v'_*|^2)^{(\beta-e_2)/2}}e^{-\frac{1}{2}\frac{\alpha\kappa/4}{2^{(8-4\kappa)}}|v'+v'_*-v|^\kappa}\mathbf{1}_{2^{-4}|v|\leq |v'+v'_*-v|}\,dv' dv'_*.
		\end{split}
	\end{align}
	The exponential term $e^{-\frac{1}{2}\frac{\alpha\kappa/4}{2^8}|v|^\kappa}$ can absorb any polynomial growth of the $|v|$, so we will perform only crude estimates if the result of the integral is bounded by a polynomial of $|v|$. As a result, any results in \textbf{Case 1} will end with the decay $e^{-\frac{1}{4}\frac{\alpha\kappa/4}{2^8}|v|^\kappa}$.
	
	We divide the domain of $v'$ and $v'_*$ by (1) $|v'|\geq 2|v|$ or $|v'_*|\geq 2|v|$ and (2) $|v'|\leq 2|v|$ and $|v'_*|\leq 2|v|$.\\
	
	\noindent \textbf{Case 1-1} Suppose $|v'|\geq 2|v|$ or $|v'_*|\geq 2|v|$. By the symmetry between $v'$ and $v'_*$, it is enough to assume $|v'|\geq 2|v|$. If $|v'|\geq 2|v|$, using $(v'_*-v)\cdot (v'-v) = 0$,
	\begin{align*}
		|v'+v'_*-v|^2 &= (v'_*-v)^2 + 2v'\cdot (v'_*-v) + (v')^2 = (v'_*-v)^2 + 2v\cdot (v'_*-v) + (v')^2\\
		&=\left(v'_*-v + v\right)^2 -v^2 + (v')^2\geq (v'_*)^2 +\frac{3}{4}(v')^2.
	\end{align*}
	Therefore, using $(x+y)^{\kappa/2}\geq \frac{\kappa}{2}(x^{\kappa/2} + y^{\kappa/2})$ for $0<\kappa\leq 2$, we have
	\begin{align*}
		e^{-\frac{1}{2}\frac{\alpha\kappa/4}{2^{(8-4\kappa)}}|v'+v'_*-v|^\kappa}\leq e^{-\frac{1}{2}\frac{\alpha\kappa/4}{2^{(8-4\kappa)}}\frac{3}{4}\left((v')^2+(v'_*)^2\right)^{\frac{\kappa}{2}}}\leq e^{-\frac{1}{2}\frac{\alpha\kappa/4}{2^{(8-4\kappa)}}\frac{3\kappa}{8}\left(|v'|^{\kappa} + |v'_*|^{\kappa}\right)}.
	\end{align*}
	
	Using the inequality, we get
	\begin{align}\label{eq:2}
		\begin{split}
			&\int_{|v'|\geq 2|v|}\frac{|h(v')|}{|v-v'|(1+|v'|^2)^{\beta/2}}\int_{v_*'-v\in E_{v'-v}}\frac{1}{(|v'-v|^2 + |v'_* - v|^2)^{\frac{1-\gamma}{2}}(1+|v'_*|^2)^{(\beta-e_2)/2}}\\
			&\qquad \times e^{-\frac{1}{2}\frac{\alpha\kappa/4}{2^{(8-4\kappa)}}|v'+v'_*-v|^\kappa}\mathbf{1}_{2^{-4}|v|\leq |v'+v'_*-v|}\,dv'_* dv'\\
			&\leq \int_{|v'|\geq 2|v|}\frac{|h(v')|}{|v-v'|(1+|v'|^2)^{\beta/2}}e^{-\frac{1}{2}\frac{\alpha\kappa/4}{2^{(8-4\kappa)}}\frac{3\kappa}{8}|v'|^{\kappa}}\\
			&\qquad \times \int_{v_*'-v\in E_{v'-v}}\frac{1}{(|v'-v|^2 + |v'_* - v|^2)^{\frac{1-\gamma}{2}}(1+|v'_*|^2)^{(\beta-e_2)/2}} e^{-\frac{1}{2}\frac{\alpha\kappa/4}{2^{(8-4\kappa)}}\frac{3\kappa}{8}|v'_*|^{\kappa}}\,dv'_* dv'.
		\end{split}
	\end{align}
    Since $0\leq \gamma\leq 1$, it is enough to bound
    \begin{align*}
        \int_{v_*'-v\in E_{v'-v}}\frac{1}{|v'_* - v|^{1-\gamma}(1+|v'_*|^2)^{(\beta-e_2)/2}} e^{-\frac{1}{2}\frac{\alpha\kappa/4}{2^{(8-4\kappa)}}\frac{3\kappa}{8}|v'_*|^{\kappa}}\,dv'_*
    \end{align*}
    for $v'_*$ integral. This integral is uniformly bounded by a constant. Indeed,
    \begin{align*}
        &\int_{v_*'-v\in E_{v'-v}}\frac{1}{|v'_* - v|^{1-\gamma}(1+|v'_*|^2)^{(\beta-e_2)/2}} e^{-\frac{1}{2}\frac{\alpha\kappa/4}{2^{(8-4\kappa)}}\frac{3\kappa}{8}|v'_*|^{\kappa}}\,dv'_*\\
        &\leq \int_{v_*'-v\in E_{v'-v}, |v_*'-v|\leq 1}\frac{1}{|v'_* - v|^{1-\gamma}(1+|v'_*|^2)^{(\beta-e_2)/2}} e^{-\frac{1}{2}\frac{\alpha\kappa/4}{2^{(8-4\kappa)}}\frac{3\kappa}{8}|v'_*|^{\kappa}}\,dv'_*\\
        &\quad + \int_{v_*'-v\in E_{v'-v}, |v_*'-v|\geq 1}\frac{1}{|v'_* - v|^{1-\gamma}(1+|v'_*|^2)^{(\beta-e_2)/2}} e^{-\frac{1}{2}\frac{\alpha\kappa/4}{2^{(8-4\kappa)}}\frac{3\kappa}{8}|v'_*|^{\kappa}}\,dv'_*\\
        &\leq \int_{v_*'-v\in E_{v'-v}, |v_*'-v|\leq 1}\frac{1}{|v'_* - v|^{1-\gamma}} \,dv'_*+ \int_{v_*'-v\in E_{v'-v}}\frac{1}{(1+|v'_*|^2)^{(\beta-e_2)/2}} e^{-\frac{1}{2}\frac{\alpha\kappa/4}{2^{(8-4\kappa)}}\frac{3\kappa}{8}|v'_*|^{\kappa}}\,dv'_*\\
        &\leq C
    \end{align*}
    for some constant $C$. Therefore, we get
	\begin{align}\label{eq:case11}
        \begin{split}
            \eqref{eq:2}&\leq C\int_{|v'|\geq 2|v|}\frac{|h(v')|}{|v-v'|(1+|v'|^2)^{\beta/2}} e^{-\frac{1}{2}\frac{\alpha\kappa/4}{2^{(8-4\kappa)}}\frac{3\kappa}{8}|v'|^{\kappa}}dv'\\
    		&\leq C(\text{polynomial of $|v|$})\|(1+|v|^2)^{-e_1/2}h(v)\|_{L^2}.
        \end{split}
	\end{align}
	\\
	
	\noindent \textbf{Case 1-2} Now, suppose $|v'|\leq 2|v|$ and $|v'_*|\leq 2|v|$. Since $0\leq \gamma\leq 1$, we get
	\begin{align*}
		&\int_{v_*'-v\in E_{v'-v}}\frac{1}{(|v'-v|^2 + |v'_* - v|^2)^{\frac{1-\gamma}{2}}(1+|v'_*|^2)^{(\beta-e_2)/2}}e^{-\frac{1}{2}\frac{\alpha\kappa/4}{2^{(8-4\kappa)}}|v'+v'_*-v|^\kappa}\mathbf{1}_{|v'|\text{ and }|v'_*|\leq 2|v|}\,dv'_*\\
		&\leq C(1+|v|)^{\max\{-\beta+e_2, 0\}}\int_{v_*'-v\in E_{v'-v}, |v'_*-v|\leq 3|v|}\frac{1}{|v'_*-v|^{1-\gamma}}\,dx\\
		&\leq C(1+|v|)^{\max\{-\beta+e_2, 0\}}|v|^{1+\gamma}.
	\end{align*}
	Therefore,
	\begin{align}\label{eq:case12}
        \begin{split}
            &\int_{\mathbb{R}^3}\frac{|h(v')|}{|v-v'|(1+|v'|^2)^{\beta/2}}e^{-\frac{1}{2}\frac{\alpha\kappa/4}{2^{(8-4\kappa)}}\frac{3\kappa}{8}|v'|^{\kappa}}\\
    		&\qquad \times \int_{v_*'-v\in E_{v'-v}}\frac{1}{(|v'-v|^2 + |v'_* - v|^2)^{\frac{1-\gamma}{2}}(1+|v'_*|^2)^{\beta/2}} e^{-\frac{1}{2}\frac{\alpha\kappa/4}{2^{(8-4\kappa)}}\frac{3\kappa}{8}|v'_*|^{\kappa}}\mathbf{1}_{|v'|\text{ and }|v'_*|\leq 2|v|}\,dv'_* dv'\\
    		&\leq C(1+|v|)^{\max\{-\beta+e_2, 0\}}|v|^{1+\gamma}\int_{|v'|\leq 2|v|}\frac{|h(v')|}{|v-v'|(1+|v'|^2)^{\beta/2}}e^{-\frac{1}{2}\frac{\alpha\kappa/4}{2^{(8-4\kappa)}}\frac{3\kappa}{8}|v'|^{\kappa}}\,dv'\\
    		&\leq C(\text{polynomial of $|v|$})\|(1+|v|^2)^{-e_1/2}h(v)\|_{L^2}.
        \end{split}
	\end{align}
	
	Merging the results \eqref{eq:1}, \eqref{eq:case11}, and \eqref{eq:case12}, we get
	\begin{align}\label{eq:case1}
		\begin{split}
			&(1+|v|^2)^{\beta/2}e^{-\frac{1}{2}\frac{\alpha\kappa/4}{2^8}|v|^\kappa}\int_{\mathbb{R}^3}\frac{|h(v')|}{|v-v'|(1+|v'|^2)^{\beta/2}}\int_{v_*'\in v + E_{v'-v}}\frac{|h(v'_*)|}{|v'-v'_*|^{1-\gamma}(1+|v'_*|^2)^{\beta/2}}\\
			&\qquad \times \left(e^{-\alpha(2^{\kappa/2}-1)|v_*|^\kappa}\mathbf{1}_{|v_*|\geq |v|} + \sum_{n=0}^7 e^{-\frac{\alpha\kappa/4}{2^{(1-\kappa/2)(n+1)}}|v'+v'_*-v|^\kappa}\mathbf{1}_{2^{-n-1}|v|^2\leq |v'+v'_*-v|^2\leq 2^{-n}|v|^2}\right)\,dv' dv'_*\\
			&\leq Ce^{-\frac{1}{4}\frac{\alpha\kappa/4}{2^8}|v|^\kappa}\|(1+|v|^2)^{-e_2/2}h\|_{L^\infty}\|(1+|v|^2)^{-e_1/2}h\|_{L^2}.
		\end{split}
	\end{align}
	Since the exponential decay dominates the polynomial growth, we obtain an exponential decay of $|v|$.\\
	
	\noindent \textbf{Case 2} We move on to the second part. First, we have
	\begin{align*}
		e^{-\frac{\alpha\kappa/4}{2^{(1-\kappa/2)(n+1)}}|v'+v'_*-v|^\kappa}\mathbf{1}_{2^{-n-1}|v|^2\leq |v'+v'_*-v|^2\leq 2^{-n}|v|^2}&\leq e^{-\frac{\alpha\kappa/4}{2^{(1-\kappa/2)(n+1)}}(2^{-(n+1)/2}|v|)^\kappa}\mathbf{1}_{2^{-n-1}|v|^2\leq |v'+v'_*-v|^2\leq 2^{-n}|v|^2}\\
		&\leq e^{-\frac{\alpha\kappa/4}{2^{n+1}}|v|^\kappa}\mathbf{1}_{2^{-n-1}|v|^2\leq |v'+v'_*-v|^2\leq 2^{-n}|v|^2}.
    	\end{align*}
	As $|v_*|\leq \frac{1}{2^{n/2}}|v|$,
	\begin{align*}
		\left(1-\frac{1}{2^{n/2}}\right)|v|\leq |v-v_*|=|v'-v'_*|\leq \left(1+\frac{1}{2^{n/2}}\right)|v|.
	\end{align*}
	Also, it implies $|v'|\leq 2|v|$.
	
	Therefore, we write
	\begin{align*}
		&(1+|v|^2)^{\beta/2}\int_{\mathbb{R}^3}\frac{|h(v')|}{|v-v'|(1+|v'|^2)^{\beta/2}}\int_{v_*'\in v + E_{v'-v}}\frac{|h(v'_*)|}{|v'-v'_*|^{1-\gamma}}\frac{1}{(1+|v'_*|^2)^{\beta/2}}\\
		&\qquad \times \sum_{n=8}^\infty e^{-\frac{\alpha\kappa/4}{2^{(1-\kappa/2)(n+1)}}|v'+v'_*-v|^\kappa}\mathbf{1}_{2^{-n-1}|v|^2\leq |v'+v'_*-v|^2\leq 2^{-n}|v|^2}\times \begin{cases}
			e^{-\alpha(1-\frac{\kappa}{2})|v'_*|^\kappa} & |v'_*|\leq |v'|\\
			e^{-\alpha(1-\frac{\kappa}{2})|v'|^\kappa} & |v'|\leq |v'_*|
		\end{cases}\,dv'_* dv'.\\
		&\leq C(1+|v|^2)^{\beta/2} \sum_{n=8}^\infty\frac{e^{-\frac{\alpha\kappa/4}{2^{n+1}}|v|^\kappa}}{|v|^{1-\gamma}}\int_{|v'|\leq 2|v|}\frac{|h(v')|}{|v-v'|(1+|v'|^2)^{\beta/2}}\int_{v_*'\in v + E_{v'-v}}\frac{|h(v'_*)|}{(1+|v'_*|^2)^{\beta/2}}\\
		&\qquad \times  \mathbf{1}_{2^{-n-1}|v|^2\leq |v'+v'_*-v|^2\leq 2^{-n}|v|^2}\times \begin{cases}
			e^{-\alpha(1-\frac{\kappa}{2})|v'_*|^\kappa} & |v'_*|\leq |v'|\\
			e^{-\alpha(1-\frac{\kappa}{2})|v'|^\kappa} & |v'|\leq |v'_*|
		\end{cases}\,dv'_* dv'.
	\end{align*}
	
	Using $v'\leftrightarrow v'_*$ symmetry, we have freedom to choose $|v'|\leq |v'_*|$ or $|v'|\geq |v'_*|$. Each case has a slightly different $(1+|v|)$ decay, so we perform computations for each case. \textbf{Case 2-1} is for $|v'|\leq |v'_*|$ case and \textbf{Case 2-2} is for the other case. The critical case is the case $v' \simeq v$ and $v'_* \simeq 0$ in \textbf{Case 2-2}. Therefore, we divide the subcases reflecting this nature.\\
	
	\noindent \textbf{Case 2-1} When $|v'|\leq |v'_*|$, as $|v'_*|\geq \frac{1}{\sqrt{2}}|v|$ and $|v'|\leq \frac{4}{5}|v|$, we can rearrange the integral by
	
	\begin{align}\label{eq:case211}
		\begin{split}
			&\sum_{n=8}^\infty (1+|v|^2)^{\beta/2}\frac{e^{-\frac{\alpha\kappa/4}{2^{n+1}}|v|^\kappa}}{|v|^{1-\gamma}}\int_{|v'|\leq 2|v|}\frac{|h(v')|}{|v-v'|(1+|v'|^2)^{\beta/2}}e^{-\alpha(1-\frac{\kappa}{2})|v'|^\kappa}\int_{v_*'\in v + E_{v'-v}}\frac{|h(v'_*)|}{(1+|v'_*|^2)^{\beta/2}}\\
			&\qquad \times \mathbf{1}_{2^{-n-1}|v|^2\leq |v'+v'_*-v|^2\leq 2^{-n}|v|^2}\,dv'_* dv'\\
			&\leq C\sum_{n=8}^\infty \|(1+|v|^2)^{-e_2/2}h\|_{L^\infty}\frac{e^{-\frac{\alpha\kappa/4}{2^{n+1}}|v|^\kappa}}{|v|^{1-\gamma}(1+|v|^2)^{-e_2/2}}\\
			&\qquad \times \int_{|v'|\leq \frac{4}{5}|v|}\frac{|h(v')|}{|v-v'|(1+|v'|^2)^{(\beta-e_2)/2}}e^{-\alpha(1-\frac{\kappa}{2})|v'|^\kappa}\int_{v_*'\in v + E_{v'-v}}\mathbf{1}_{2^{-n-1}|v|^2\leq |v'+v'_*-v|^2\leq 2^{-n}|v|^2}\,dv'_* dv'.
		\end{split}
	\end{align}
	Seeing $v'_*-v$ as a plane variable and $\mathbf{1}_{2^{-n-1}|v|^2\leq |v'+v'_*-v|^2\leq 2^{-n}|v|^2}\mathbf{1}_{|v'|\leq |v'_*|}$ as a ball centered at $-v'$ with radius from $2^{-\frac{n+1}{2}}|v|$ to $2^{-\frac{n}{2}}|v|$, the integral is the area of the intersection between $E_{v'-v}$ and the ball. Therefore, it is bounded by
	\begin{align}\label{eq:v'*intbound}
		\int_{v_*'\in v + E_{v'-v}}\mathbf{1}_{2^{-n-1}|v|^2\leq |v'+v'_*-v|^2\leq 2^{-n}|v|^2}\mathbf{1}_{|v'|\leq |v'_*|} \,dv'_*\leq C (2^{-n/2}|v|)^{2}.
	\end{align}
	Now, we need to bound
	\begin{align*}
		\sum_{n=8}^\infty \frac{e^{-\frac{\alpha\kappa/4}{2^{n+1}}|v|^\kappa}}{|v|^{1-\gamma}}(2^{-n/2}|v|)^{2}\int_{|v'|\leq \frac{4}{5}|v|}\frac{|h(v')|}{|v-v'|(1+|v'|^2)^{\beta/2}}e^{-\alpha(1-\frac{\kappa}{2})|v'|^\kappa}dv'.
	\end{align*}
	By the H{\"o}lder's inequality,
	\begin{align*}
		&\int_{|v'|\leq \frac{4}{5}|v|}\frac{(1+|v'|^2)^{-e_1/2}|h(v')|}{|v-v'|(1+|v'|^2)^{(\beta-e_1)/2}}e^{-\alpha(1-\frac{\kappa}{2})|v'|^\kappa}dv'\\
		&\leq \|(1+|v'|^2)^{-e_1/2}h\|_{L^2}\left\|\frac{1}{|v-v'|(1+|v'|^2)^{(\beta-e_1)/2}}e^{-\alpha(1-\frac{\kappa}{2})|v'|^\kappa}\mathbf{1}_{|v'|\leq \frac{4}{5}|v|}\right\|_{L^2}\\
		&\leq \frac{C\min\{|v|^{\frac{1}{2}}, 1\}}{1+|v|}\|(1+|v|^2)^{-e_1/2}h\|_{L^2}.
	\end{align*}
	Hence,
	\begin{align}\label{eq:case212}
		\begin{split}
			&\sum_{n=8}^\infty C\frac{e^{-\frac{\alpha\kappa/4}{2^{n+1}}|v|^\kappa}}{|v|^{1-\gamma}}\int_{|v'|\leq 2|v|}\frac{|h(v')|}{|v-v'|(1+|v'|^2)^{\beta/2}}e^{-\alpha(1-\frac{\kappa}{2})|v'|^\kappa}\int_{v_*'\in v + E_{v'-v}}\mathbf{1}_{2^{-n-1}|v|^2\leq |v'+v'_*-v|^2\leq 2^{-n}|v|^2}\,dv'_* dv'\\
			&\leq C\sum_{n=8}^\infty \frac{e^{-\frac{\alpha\kappa/4}{2^{n+1}}|v|^\kappa}}{|v|^{1-\gamma}(1+|v|)^{1-e_2}}(2^{-n/2}|v|)^{2} \min\{|v|^{\frac{1}{2}}, 1\}\|(1+|v|^2)^{-e_1/2}h\|_{L^2}\\
			&\leq C\sum_{n=8}^\infty e^{-\frac{\alpha\kappa/4}{2^{n+1}}|v|^\kappa}2^{-n}\frac{|v|^{1+\gamma}\min\{|v|^{\frac{1}{2}}, 1\}}{(1+|v|)^{1-e_2}}\|(1+|v|^2)^{-e_1/2}h\|_{L^2}.
		\end{split}
	\end{align}
	

When $|v|\leq 2$, then the summation is well bounded, so we assume $|v|\geq 2$. Choose $M\in\mathbb{N}$ such that $2^{M-1}\leq |v|< 2^M$. When $n\leq \lfloor \kappa M\rfloor$, using $e^{-x}\leq \frac{1}{x^2}$,
	\begin{align*}
		\sum_{n=8}^\infty e^{-\frac{\alpha\kappa/4}{2^{n+1}}|v|^\kappa}2^{-n}\frac{|v|^{1+\gamma}}{(1+|v|)^{1-e_2}}&\leq C\sum_{n=8}^\infty \frac{2^{2n}}{|v|^{2\kappa}}2^{-n}\frac{|v|^{1+\gamma}}{(1+|v|)^{1-e_2}}\leq C\sum_{n=8}^{\lfloor \kappa M\rfloor} (1+|v|)^{\gamma-2\kappa + e_2} 2^n\\
		&\leq C (1+|v|)^{\gamma-\kappa + e_2}.
	\end{align*}
	When $n> \lfloor \kappa M\rfloor$, then we ignore $e^{-\frac{\alpha\kappa/4}{2^{n+1}}|v|^\kappa}$ by $1$, so
	\begin{align*}
		\sum_{\lfloor \kappa M\rfloor + 1}^\infty e^{-\frac{\alpha\kappa/4}{2^{n+1}}|v|^\kappa}2^{-n}\frac{|v|^{1+\gamma}}{(1+|v|)^{1-e_2}}&\leq (1+|v|)^{\gamma+e_2}\sum_{\lfloor \kappa M\rfloor + 1}^\infty 2^{-n} \leq 2(1+|v|)^{\gamma-\kappa + e_2}.
	\end{align*}
	
    Applying the summation to \eqref{eq:case212} and \eqref{eq:case211}, we get
	\begin{align}\label{eq:4}
		\begin{split}
			&\sum_{n=8}^\infty (1+|v|^2)^{\beta/2}\frac{e^{-\frac{\alpha\kappa/4}{2^{n+1}}|v|^\kappa}}{|v|^{1-\gamma}}\int_{|v'|\leq 2|v|}\frac{|h(v')|}{|v-v'|(1+|v'|^2)^{\beta/2}}\int_{v_*'\in v + E_{v'-v}}\frac{|h(v'_*)|}{(1+|v'_*|^2)^{\beta/2}}e^{-\alpha(1-\frac{\kappa}{2})|v'|^\kappa}\\
			&\qquad \times \mathbf{1}_{2^{-n-1}|v|^2\leq |v'+v'_*-v|^2\leq 2^{-n}|v|^2}\,dv'_* dv'\\
			&\leq C(1+|v|)^{\gamma-\kappa + e_2}\|(1+|v|^2)^{-e_1/2}h\|_{L^2}\|(1+|v|^2)^{-e_2/2}h\|_{L^\infty}
		\end{split}
	\end{align}
	when $|v'|\leq |v'_*|$.\\
	
	\noindent \textbf{Case 2-2} When $|v'_*|\leq |v'|$, as $|v'|\geq \frac{1}{\sqrt{2}}|v|$, we can rearrange the integral by
	\begin{align*}
		&\sum_{n=8}^\infty C\|(1+|v|^2)^{-e_2/2}h\|_{L^\infty}\frac{e^{-\frac{\alpha\kappa/4}{2^{n+1}}|v|^\kappa}}{|v|^{1-\gamma}}\int_{|v'|\leq 2|v|}\frac{|h(v')|}{|v-v'|}\\
		&\qquad \times \int_{v_*'\in v + E_{v'-v}}\frac{e^{-\alpha(1-\frac{\kappa}{2})|v'_*|^\kappa}}{(1+|v'_*|^2)^{(\beta-e_2)/2}}\mathbf{1}_{2^{-n-1}|v|^2\leq |v'+v'_*-v|^2\leq 2^{-n}|v|^2}\,dv'_* dv'.
	\end{align*}
	
	We divide the integral region of $v'$ by $|v'-v|\geq \frac{|v|}{4}$, $2^{-n/2+2}|v|\leq |v'-v|\leq \frac{|v|}{4}$, and the other part. Note that $|v'|\leq 2|v|$ since $|v_*|\leq 2^{-4}|v|$. As noted before, the critical case is when $v' \simeq v$ and $v'_* \simeq 0$.
	
	\noindent \textbf{Case 2-2-1} When $|v'-v|\geq \frac{|v|}{4}$, then $|v'_*|\geq \frac{|v|}{8}$. It means $e^{-\alpha(1-\frac{\kappa}{2})|v'_*|^\kappa}\leq e^{-\frac{\alpha}{8^\kappa}(1-\frac{\kappa}{2})|v|^\kappa}$. Following the previous computation \eqref{eq:v'*intbound},
	\begin{align*}
		&\sum_{n=8}^\infty \frac{e^{-\frac{\alpha\kappa/4}{2^{n+1}}|v|^\kappa}}{|v|^{1-\gamma}}\int_{|v'-v|\geq \frac{|v|}{4}, |v'|\leq 2|v|}\frac{|h(v')|}{|v-v'|}\int_{v_*'\in v + E_{v'-v}}\frac{e^{-\alpha(1-\frac{\kappa}{2})|v'_*|^\kappa}}{(1+|v'_*|^2)^{(\beta-e_2)/2}}\mathbf{1}_{2^{-n-1}|v|^2\leq |v'+v'_*-v|^2\leq 2^{-n}|v|^2}\,dv'_* dv'\\
		&\leq C\max\{(1+|v|^2)^{(e_2-\beta)/2}, 1\}e^{-\frac{\alpha}{8^\kappa}(1-\frac{\kappa}{2})|v|^\kappa}\sum_{n=8}^\infty \frac{e^{-\frac{\alpha\kappa/4}{2^{n+1}}|v|^\kappa}}{|v|^{1-\gamma}}\\
		&\qquad \times\int_{|v'-v|\geq \frac{|v|}{4}, |v'|\leq 2|v|}\frac{|h(v')|}{|v-v'|}\int_{v_*'\in v + E_{v'-v}}\mathbf{1}_{2^{-n-1}|v|^2\leq |v'+v'_*-v|^2\leq 2^{-n}|v|^2}\,dv'_* dv'\\
		&\leq C\max\{(1+|v|^2)^{(e_2-\beta)/2}, 1\}e^{-\frac{\alpha}{8^\kappa}(1-\frac{\kappa}{2})|v|^\kappa}\sum_{n=8}^\infty \frac{e^{-\frac{\alpha\kappa/4}{2^{n+1}}|v|^\kappa}}{|v|^{1-\gamma}}(2^{-n/2}|v|)^{2}\int_{|v'-v|\geq \frac{|v|}{4}, |v'|\leq 2|v|}\frac{|h(v')|}{|v-v'|}dv'\\
		&\leq C\max\{(1+|v|^2)^{(e_2-\beta)/2}, 1\}(1+|v|^2)^{e_1/2}e^{-\frac{\alpha}{8^\kappa}(1-\frac{\kappa}{2})|v|^\kappa}\\
		&\qquad \times \sum_{n=8}^\infty 2^{-n}|v|^{1+\gamma}|v|^{\frac{1}{2}}\|(1+|v|^2)^{-e_1/2}h\|_{L^2}.
	\end{align*}
	Because of the exponential decay of $|v|$, it is bounded uniformly about $|v|$. Therefore, we have
	\begin{align}\label{eq:case221}
		\begin{split}
			&\sum_{n=8}^\infty \frac{e^{-\frac{\alpha\kappa/4}{2^{n+1}}|v|^\kappa}}{|v|^{1-\gamma}}\int_{|v'-v|\geq \frac{|v|}{4}, |v'|\leq 2|v|}\frac{|h(v')|}{|v-v'|}\int_{v_*'\in v + E_{v'-v}}\frac{e^{-\alpha(1-\frac{\kappa}{2})|v'_*|^\kappa}}{(1+|v'_*|^2)^{(\beta-e_2)/2}}\mathbf{1}_{2^{-n-1}|v|^2\leq |v'+v'_*-v|^2\leq 2^{-n}|v|^2}\,dv'_* dv'\\
			&\leq Ce^{-\frac{1}{2}\frac{\alpha}{8^\kappa}(1-\frac{\kappa}{2})|v|^\kappa}\|(1+|v|^2)^{-e_1/2}h\|_{L^2}.
		\end{split}
	\end{align}
	\\
	
	\noindent \textbf{Case 2-2-2} When $2^{-n/2+2}|v|\leq |v'-v|\leq \frac{|v|}{4}$, then $|v'_*|\geq \frac{3}{4}|v'-v|\geq \frac{3}{4}2^{-n/2+2}|v|$. In this case, the size of the set of $v'$ such that $E_{v'-v}\cap |v'+'_*-v|\leq 2^{-n/2}|v|$ is proportional to $2^{-n/2}$ according to Lemma \ref{lem:cone} (for $x = v'-v$, $y = v'_*-v$, and $a = v$), so
	\begin{align}\label{eq:coneapply}
		\begin{split}
			&\int_{2^{-n/2+2}|v|\leq |v'-v|\leq \frac{|v|}{4}}\frac{|h(v')|}{|v-v'|} \int_{v_*'\in v + E_{v'-v}} \mathbf{1}_{2^{-n-1}|v|^2\leq |v'+v'_*-v|^2\leq 2^{-n}|v|^2}\,dv'_* dv'\\
			&=\int_{\{2^{-n/2+2}|v|\leq |v'-v|\leq \frac{|v|}{4}\}\cap C_{v, 2^{-n/2}|v|}}\frac{|h(v')|}{|v-v'|} \int_{v_*'\in v + E_{v'-v}} \mathbf{1}_{2^{-n-1}|v|^2\leq |v'+v'_*-v|^2\leq 2^{-n}|v|^2}\,dv'_* dv'
		\end{split}
	\end{align}
	Therefore, bounding $\frac{e^{-\frac{\alpha}{2}(1-\frac{\kappa}{2})|v'_*|^\kappa}}{(1+|v'_*|^2)^{(\beta-e_2)/2}}<\infty$ and using \eqref{eq:v'*intbound},
	\begin{align}\label{eq:case2221}
		\begin{split}
			&\sum_{n=8}^\infty \frac{e^{-\frac{\alpha\kappa/4}{2^{n+1}}|v|^\kappa}}{|v|^{1-\gamma}}\int_{2^{-n/2+2}|v|\leq |v'-v|\leq \frac{|v|}{4}}\frac{|h(v')|}{|v-v'|} \int_{v_*'\in v + E_{v'-v}} \frac{e^{-\alpha(1-\frac{\kappa}{2})|v'_*|^\kappa}}{(1+|v'_*|^2)^{(\beta-e_2)/2}}\mathbf{1}_{2^{-n-1}|v|^2\leq |v'+v'_*-v|^2\leq 2^{-n}|v|^2}\,dv'_* dv'\\
			&\leq C\sum_{n=8}^\infty \frac{e^{-\frac{\alpha\kappa/4}{2^{n+1}}|v|^\kappa}}{|v|^{1-\gamma}}\int_{\{2^{-n/2+2}|v|\leq |v'-v|\leq \frac{|v|}{4}\}\cap C_{v, 2^{-n/2}|v|}}\frac{|h(v')|}{|v-v'|}e^{-\frac{\alpha}{2}(1-\frac{\kappa}{2})(\frac{3}{4}|v'-v|)^\kappa}\\
			&\qquad \times \int_{v_*'\in v + E_{v'-v}}\mathbf{1}_{|v'+v'_*-v|^2\leq \frac{|v'-v|^2}{4^2}}\,dv_*' dv'\\
			&\leq C\sum_{n=8}^\infty \frac{e^{-\frac{\alpha\kappa/4}{2^{n+1}}|v|^\kappa}}{|v|^{1-\gamma}}(1+|v|^2)^{e_1/2}\\
			&\qquad \times \int_{\{2^{-n/2+2}|v|\leq |v'-v|\leq \frac{|v|}{4}\}\cap C_{v, 2^{-n/2}|v|}}(1+|v'|^2)^{-e_1/2}|h(v')|e^{-\frac{\alpha}{2}(1-\frac{\kappa}{2})(\frac{3}{4}|v'-v|)^\kappa} |v-v'|^{1} dv'.
		\end{split}
	\end{align}
	We can compute the $L^2$ norm of the integrand of $v'$ integral using Lemma \ref{lem:cone}.
	\begin{align*}
		&\int_{\{2^{-n/2+2}|v|\leq |v'-v|\leq \frac{|v|}{4}\}\cap C_{v,2^{-n/2}|v|}}|v-v'|^{2}e^{-\alpha(1-\frac{\kappa}{2})(\frac{3}{4}|v'-v|)^\kappa}\,dv'\\
		&\leq C 2^{-n/2}\int_{2^{-n/2+2}|v|}^{|v|/4} e^{-\alpha(1-\frac{\kappa}{2})\frac{3^\kappa}{4^\kappa} r^\kappa}r^{4}\,dr\\
		&\leq C 2^{-n/2}e^{-\alpha(1-\frac{\kappa}{2})\frac{3^\kappa}{4^\kappa}(2^{-n/2+2}|v|)^\kappa} \min\{|v|^{5}, 1\}(1+2^{-n/2+2}|v|)^{5-\kappa},
		\intertext{and}
		&\int_{\{2^{-n/2+2}|v|\leq |v'-v|\leq \frac{|v|}{4}\}\cap C_{v,2^{-n/2}|v|}}|v-v'|^{1}e^{-\frac{\alpha}{2}(1-\frac{\kappa}{2})(\frac{3}{4}|v'-v|)^\kappa}\,dv'\\
		&\leq C 2^{-n/2}\int_{2^{-n/2+2}|v|}^{|v|/4} e^{-\frac{\alpha}{2}(1-\frac{\kappa}{2})\frac{3^\kappa}{4^\kappa} r^\kappa}r^{3}\,dr\\
		&\leq C 2^{-n/2} e^{-\frac{\alpha}{2}(1-\frac{\kappa}{2})\frac{3^\kappa}{4^\kappa}(2^{-n/2+2}|v|)^\kappa} \min\{|v|^{4}, 1\}(1+2^{-n/2+2}|v|)^{4-\kappa}.
	\end{align*}
	Using H{\"o}lder's inequality, we get
	\begin{align}\label{eq:case2222}
		\begin{split}
			&\int_{\{2^{-n/2+2}|v|\leq |v'-v|\leq \frac{|v|}{4}\}\cap C_{v, 2^{-n/2}|v|}}(1+|v'|^2)^{-e_1/2}|h(v')|e^{-\frac{\alpha}{2}(1-\frac{\kappa}{2})(\frac{3}{4}|v'-v|)^\kappa} |v-v'|^{1} dv'\\
			&\leq C \left(2^{-n/2} e^{-\alpha(1-\frac{\kappa}{2})\frac{3^\kappa}{4^\kappa}(2^{-n/2+2}|v|)^\kappa} \min\{|v|^{5}, 1\}(1+2^{-n/2+2}|v|)^{5-\kappa}\right)^{1/2}\|(1+|v|^2)^{-e_1/2}h\|_{L^2}.
		\end{split}
	\end{align}
	Combining \eqref{eq:case2221} and \eqref{eq:case2222}, for $L^2$ case, we obtain
	\begin{align*}
		&\sum_{n=8}^\infty \frac{e^{-\frac{\alpha\kappa/4}{2^{n+1}}|v|^\kappa}}{|v|^{1-\gamma}}\int_{2^{-n/2+2}|v|\leq |v'-v|\leq \frac{|v|}{4}}\frac{|h(v')|}{|v-v'|} \int_{v_*'\in v + E_{v'-v}} \frac{e^{-\alpha(1-\frac{\kappa}{2})|v'_*|^\kappa}}{(1+|v'_*|^2)^{(\beta-e_2)/2}}\mathbf{1}_{2^{-n-1}|v|^2\leq |v'+v'_*-v|^2\leq 2^{-n}|v|^2}\,dv'_* dv'\\
		&\leq C\sum_{n=8}^\infty \frac{e^{-\frac{\alpha\kappa/4}{2^{n+1}}|v|^\kappa}}{|v|^{1-\gamma}}(1+|v|^2)^{e_1/2} \left(2^{-n/2} e^{-\alpha(1-\frac{\kappa}{2})\frac{3^\kappa}{4^\kappa}(2^{-n/2+2}|v|)^\kappa} \min\{|v|^{5}, 1\}(1+2^{-n/2+2}|v|)^{5-\kappa}\right)^{1/2}\\
		&\qquad \times \|(1+|v|^2)^{-e_1/2}h\|_{L^2}.
	\end{align*} $|v|^{\frac{3}{2}+\gamma}$
	If $|v|\leq 2$, then $|v|^{\frac{3}{2}+\gamma}$ is well-bounded, so the summation is well-defined. For $|v|\geq 2$, ignoring $e^{-\frac{\alpha\kappa/4}{2^{n+1}}|v|^\kappa}\leq 1$, we have
	\begin{align*}
		&\sum_{n=8}^\infty \frac{e^{-\frac{\alpha\kappa/4}{2^{n+1}}|v|^\kappa}}{|v|^{1-\gamma}}(1+|v|^2)^{e_1/2} \left(2^{-n/2} e^{-\alpha(1-\frac{\kappa}{2})\frac{3^\kappa}{4^\kappa}(2^{-n/2+2}|v|)^\kappa} (1+2^{-n/2+2}|v|)^{5-\kappa}\right)^{1/2}\\
		&=\sum_{n=8}^\infty e^{-\frac{\alpha}{2}(1-\frac{\kappa}{2})\frac{3^\kappa}{4^\kappa}(2^{-n/2+2}|v|)^\kappa}(1+|v|)^{e_1 - 1 + \gamma} (1+2^{-n/2+2}|v|)^{\frac{5}{2}-\frac{\kappa}{2}} 2^{-n/4}
	\end{align*}
	
	Suppose $2^{M-1}\leq |v|<2^M$ for some $M\in\mathbb{N}$. When $n\leq 2M$, using $x^{\frac{5-\frac{\kappa}{2}}{\kappa}} e^{-x}\leq C_\kappa$ for some constant $C_\kappa$ for $0<\kappa$,
	\begin{align*}
		&\sum_{n=8}^{2M} e^{-\frac{\alpha}{2}(1-\frac{\kappa}{2})\frac{3^\kappa}{4^\kappa}(2^{-n/2+2}|v|)^\kappa}(1+|v|)^{e_1 - 1 + \gamma} (1+2^{-n/2+2}|v|)^{\frac{5}{2}-\frac{\kappa}{2}} 2^{-n/4}\\
		&\leq C\sum_{n=8}^{2M} e^{-\frac{\alpha}{2}(1-\frac{\kappa}{2})\frac{3^\kappa}{4^\kappa}(2^{-n/2+2}|v|)^\kappa}(1+|v|)^{e_1+\gamma+\frac{3}{2}-\frac{\kappa}{2}} 2^{-(\frac{3}{2}-\frac{\kappa}{4})n}\\
		&\leq C\sum_{n=8}^{2M} \left(\frac{2^{n/2}}{|v|}\right)^{5-\frac{\kappa}{2}} (1+|v|)^{e_1+\gamma+\frac{3}{2}-\frac{\kappa}{2}} 2^{-(\frac{3}{2}-\frac{\kappa}{4})n} = C\sum_{n=8}^{2M} 2^n (1+|v|)^{-\frac{7}{2}+\gamma+e_1}\\
		&\leq C(1+|v|)^{-\frac{3}{2}+\gamma+e_1}.
	\end{align*}
	When $n>2M$, bounding $e^{-\alpha(1-\frac{\kappa}{2})\frac{3^\kappa}{4^\kappa}\left(2^{-n/2}|v|\right)^\kappa}\leq 1$,
	\begin{align*}
		&\sum_{n=2M+1}^\infty e^{-\frac{\alpha}{2}(1-\frac{\kappa}{2})\frac{3^\kappa}{4^\kappa}(2^{-n/2+2}|v|)^\kappa}(1+|v|)^{e_1 - 1 + \gamma} (1+2^{-n/2+2}|v|)^{\frac{5}{2}-\frac{\kappa}{2}} 2^{-n/4}\\
		&\leq (1+|v|)^{e_1 - 1 + \gamma} \sum_{n=2M+1}^\infty 2^{-n/4}\leq C(1+|v|)^{-\frac{3}{2}+\gamma+e_1}.
	\end{align*}
	
	It proves
	\begin{align}\label{eq:case2223}
		\begin{split}
			&\sum_{n=8}^\infty  \frac{e^{-\frac{\alpha\kappa/4}{2^{n+1}}|v|^\kappa}}{|v|^{1-\gamma}}\int_{2^{-n/2+2}|v|\leq |v'-v|\leq \frac{|v|}{4}}\frac{|h(v')|}{|v-v'|} \int_{v_*'\in v + E_{v'-v}} \frac{e^{-\alpha(1-\frac{\kappa}{2})|v'_*|^\kappa}}{(1+|v'_*|^2)^{(\beta-e_2)/2}}\mathbf{1}_{2^{-n-1}|v|^2\leq |v'+v'_*-v|^2\leq 2^{-n}|v|^2}\,dv'_* dv'\\
			&\leq C (1+|v|)^{-\frac{3}{2}+\gamma+e_1}\|(1+|v|^2)^{-e_1/2}h\|_{L^2}.
		\end{split}
	\end{align}
	\\
	
	\noindent \textbf{Case 2-2-3} Here, we consider the $L^2$ case. When $|v'-v|\leq 2^{-n/2+2}|v|$, there is no hope to extract some decay about $|v|$ or $|v'-v|$ from $e^{-\alpha(1-\frac{\kappa}{2})|v'_*|^\kappa}$ since the integral region on $v'_*$ can hit $v'_* = 0$. However, it can provide the $L^1$ bound of the $v'_*$ integral.
	\begin{align*}
		&\int_{v_*'\in v + E_{v'-v}}\frac{e^{-\alpha(1-\frac{\kappa}{2})|v'_*|^\kappa}}{(1+|v'_*|^2)^{(\beta-e_2)/2}} \mathbf{1}_{2^{-n-1}|v|^2\leq |v'+v'_*-v|^2\leq 2^{-n}|v|^2}\,dv'_*\\
		&\leq C\int_{v_*'\in v + E_{v'-v}}e^{-\frac{\alpha}{2}(1-\frac{\kappa}{2})|v'_*|^\kappa} \mathbf{1}_{|v'_*|^2\leq 2^{-n+5}|v|^2}\,dv'_*\\
		&\leq \begin{cases}
			C(2^{-n/2}|v|)^{2} & 2^{-n/2}|v|\leq 1\\
			C & 2^{-n/2}|v|\geq 1
		\end{cases}
	\end{align*}
	Furthermore, we can apply \eqref{eq:coneapply} for the same reason. So,
	\begin{align}\label{eq:case2231}
		\begin{split}
			&\sum_{n=8}^\infty \frac{e^{-\frac{\alpha\kappa/4}{2^{n+1}}|v|^\kappa}}{|v|^{1-\gamma}}\int_{|v'-v|\leq 2^{-n/2+2}|v|}\frac{|h(v')|}{|v-v'|} \int_{v_*'\in v + E_{v'-v}} \frac{e^{-\alpha(1-\frac{\kappa}{2})|v'_*|^\kappa}}{(1+|v'_*|^2)^{(\beta-e_2)/2}}\mathbf{1}_{2^{-n-1}|v|^2\leq |v'+v'_*-v|^2\leq 2^{-n}|v|^2}\,dv'_* dv'\\
			&\leq C\sum_{n=8}^\infty \frac{e^{-\frac{\alpha\kappa/4}{2^{n+1}}|v|^\kappa}}{|v|^{1-\gamma}}\int_{\{|v'-v|\leq 2^{-n/2+2}|v|\}\cap C_{v,2^{-n/2}|v|}}\frac{|h(v')|}{|v'-v|}\,dv'\times \begin{cases}
				C(2^{-n/2}|v|)^{2} & 2^{-n/2}|v|\leq 1,\\
				C & 2^{-n/2}|v|\geq 1.
			\end{cases}
		\end{split}
	\end{align}
    Using Lemma \ref{lem:cone}, we have
	\begin{align*}
		\int_{\{|v'-v|\leq 2^{-n/2+2}|v|\}\cap C_{v,2^{-n/2}|v|}}\frac{1}{|v'-v|^2}\,dv'&\leq C 2^{-n/2}\int_0^{2^{-n/2+2}|v|} \,dr \leq C2^{-n/2}(2^{-n/2+2}|v|)^{1}.
	\end{align*}
    It shows
	\begin{align*}
		\int_{\{|v'-v|\leq 2^{-n/2+2}|v|\}\cap C_{v,2^{-n/2}|v|}}\frac{|h(v')|}{|v'-v|}\,dv'\leq C\left(2^{-n/2}(2^{-n/2+2}|v|)^{1}\right)^{1/2}(1+|v|^2)^{e_1/2}\|(1+|v|^2)^{-e_1/2}h\|_{L^2}
	\end{align*}
	%
	
	Combining these inequalities with \eqref{eq:case2231}, we obtain
	\begin{align}\label{eq:case2232}
		\begin{split}
			&\sum_{n=8}^\infty \frac{e^{-\frac{\alpha\kappa/4}{2^{n+1}}|v|^\kappa}}{|v|^{1-\gamma}}\int_{|v'-v|\leq 2^{-n/2+2}|v|}\frac{|h(v')|}{|v-v'|} \int_{v_*'\in v + E_{v'-v}} \frac{e^{-\alpha(1-\frac{\kappa}{2})|v'_*|^\kappa}}{(1+|v'_*|^2)^{(\beta-e_2)/2}}\mathbf{1}_{2^{-n-1}|v|^2\leq |v'+v'_*-v|^2\leq 2^{-n}|v|^2}\,dv'_* dv'\\
			&\leq C\sum_{n=8}^\infty \frac{e^{-\frac{\alpha\kappa/4}{2^{n+1}}|v|^\kappa}}{|v|^{1-\gamma}}\left(2^{-n/2}(2^{-n/2+2}|v|)^{1}\right)^{1/2}(1+|v|^2)^{e_1/2}\|(1+|v|^2)^{-e_1/2}h\|_{L^2}\\
			&\qquad \times \begin{cases}
				C(2^{-n/2}|v|)^{2} & 2^{-n/2}|v|\leq 1\\
				C & 2^{-n/2}|v|\geq 1
			\end{cases}
		\end{split}
	\end{align}

	When $|v|\leq 2$, it is again uniformly bounded, so we assume $|v|\geq 2$. To fully use $e^{-\frac{\alpha\kappa/4}{2^{n+1}}|v|^\kappa}$ decay, we again divide $n\leq \lfloor \kappa M\rfloor$ and $n>\lfloor \kappa M\rfloor$ for $2^{M-1}\leq |v|\leq 2^M$. Following a similar computation for \eqref{eq:case2231}, we have
	\begin{align*}
		&\sum_{n=8}^\infty \frac{e^{-\frac{\alpha\kappa/4}{2^{n+1}}|v|^\kappa}}{|v|^{1-\gamma}}\left(2^{-n/2}(2^{-n/2+2}|v|)^{1}\right)^{1/2}\times \begin{cases}
			C(2^{-n/2}|v|)^{2} & 2^{-n/2}|v|\leq 1\\
			C & 2^{-n/2}|v|\geq 1,
		\end{cases}\\
		&\leq C(1+|v|)^{-\frac{1}{2} + \gamma -\frac{1}{2}\kappa}.
	\end{align*}

	As a result, we obtain
	\begin{align}\label{eq:case2234}
		\begin{split}
			&\sum_{n=8}^\infty \frac{e^{-\frac{\alpha\kappa/4}{2^{n+1}}|v|^\kappa}}{|v|^{1-\gamma}}\int_{|v'-v|\leq 2^{-n/2+2}|v|}\frac{|h(v')|}{|v-v'|} \int_{v_*'\in v + E_{v'-v}} \frac{e^{-\alpha(1-\frac{\kappa}{2})|v'_*|^\kappa}}{(1+|v'_*|^2)^{(\beta-e_2)/2}}\mathbf{1}_{2^{-n-1}|v|^2\leq |v'+v'_*-v|^2\leq 2^{-n}|v|^2}\,dv'_* dv'\\
			&\leq C(1+|v|)^{-\frac{1}{2} + \gamma -\frac{1}{2}\kappa + e_1}\|(1+|v|^2)^{-e_1/2}h\|_{L^2}.
		\end{split}
	\end{align}
	
	Combining \eqref{eq:case1}, \eqref{eq:4}, \eqref{eq:case221}, \eqref{eq:case2223}, and \eqref{eq:case2234}, we conclude the result of the proposition.
\end{proof}

\bigskip

\begin{corollary} \label{Qestimate}
	Let \(w(v)=(1+|v|^2)^{\beta/2}e^{\alpha|v|^{\kappa}}\) with $\beta \ge 0 ,\alpha >0 ,0<\kappa<2$. Then for all f, g in $L^{\infty}_{x,v}(w)$,
	\[\|Q(f,g)\|_{L^{\infty}_{x,v}(\nu^{-1} w)} \le C_{Q} \|f\|_{L^{\infty}_{x,v}(w)}\|g\|_{L^{\infty}_{x,v}(w)}. \]
\end{corollary}

\bigskip

The following lemma shows that $K_s^{(\delta)}$ is small relative to the collision frequency. This smallness allows its contribution to be absorbed in the weighted $L^\infty$ estimates based on the Duhamel formula.

\begin{lemma} \label{Knsest}
	Let \(w(v)=(1+|v|^2)^{\beta/2}e^{\alpha|v|^{\kappa}}\) with $\beta \ge 0 ,\alpha >0 ,0<\kappa<2$. Then for all f in $L^{\infty}_{x,v}(w)$,
	\[ \|K_s^{(\delta)}f\|_{L^{\infty}_{x,v}(\nu^{-1} w)} \le C_K(\delta) \|f\|_{L^{\infty}_{x,v}(w)}, \]
	where 
	\[ C_K(\delta) \rightarrow 0 \quad \text{as} \quad \delta\rightarrow 0.\]
\end{lemma}

\begin{proof}
	We split $K_s^{\delta}$ into three parts as follow:
	\begin{align} \label{Knsest11}
		\left|K_s^{(\delta)}f\right| \le & \int_{\mathbb{R}^3\times \mathbb{S}^2} B(v,v_*,\theta)\mathbf{1} _{|v|\ge R} \left[\mu(v_*')f(v')+\mu(v')f(v_*')\right]d\sigma dv_*\nonumber\\
		&+C_b\int_{\mathbb{R}^3\times \mathbb{S}^2}{1} _{|v|\le R}(1-\Theta_{\delta})|v-v_*|^{\gamma} \left[\mu(v_*')f(v')+\mu(v')f(v_*')\right]d\sigma dv_*\nonumber\\
		&+C_b\int_{\mathbb{R}^3\times \mathbb{S}^2}(1-\Theta_{\delta})|v-v_*|^{\gamma} \mu(v)f(v_*)d\sigma dv_*.
	\end{align}
	First of all, let us consider last two terms in \eqref{Knsest11}. We use $|v_*-v|^{\gamma} \le C \nu(v) \langle v_* \rangle^{\gamma}$ for the last one in \eqref{Knsest11}:
	\begin{align*}
		&w(v)\nu(v)^{-1}\int_{\mathbb{R}^3\times \mathbb{S}^2}(1-\Theta_{\delta})|v-v_*|^{\gamma} \mu(v)f(v_*)d\sigma dv_*\\
		&\le Cw(v)\mu(v) \mathbf{1}_{|v| \ge R}\left(\int_{\mathbb{R}^3\times \mathbb{S}^2}\frac{\langle v_* \rangle^{\gamma}}{w(v_*)}d\sigma dv_*\right) \|f\|_{L^{\infty}_{x,v}(w)}\\
		&\quad +Cw(R)\mathbf{1}_{|v| \le R}\left(\int_{\mathbb{R}^3\times \mathbb{S}^2}(1-\Theta_{\delta})\frac{\langle v_* \rangle^{\gamma}}{w(v_*)}d\sigma dv_*\right) \|f\|_{L^{\infty}_{x,v}(w)}.
	\end{align*}
	Note that
	\begin{align*}
		|v_*-v|^{\gamma}\frac{\mu(v')}{w(v_*')} \le C\frac{\nu(v)\langle v_*\rangle^{\gamma}}{w(v')w(v_*')} \le C\frac{\nu(v)\langle v_*\rangle^{\gamma}}{w\left(\frac{|v_*|}{\sqrt{2}}\right)}.
	\end{align*}
	Here, $|v'| \ge \frac{1}{\sqrt{2}}|v_*|$ or $|v_*'| \ge \frac{1}{\sqrt{2}}|v_*|$ since $|v_*'|^2+|v'|^2 \ge |v_*|^2$, and $\frac{1}{w(v')} \le \frac{1}{w\left(\frac{|v_*|}{\sqrt{2}}\right)}$ or $\frac{1}{w(v_*')} \le \frac{1}{w\left(\frac{|v_*|}{\sqrt{2}}\right)}$.\\
	As interchanging $v_*'$ and $v'$ to the second one in \eqref{Knsest11},
	\begin{align*}
		&w(v)\nu(v)^{-1}\int_{\mathbb{R}^3\times \mathbb{S}^2}\mathbf{1}_{|v| \le R}(1-\Theta_{\delta})|v-v_*|^{\gamma} \left[\mu(v_*')f(v')+\mu(v')f(v_*')\right]d\sigma dv_*\\
		&\le Cw(R)\left(\int_{\mathbb{R}^3\times \mathbb{S}^2}\mathbf{1}_{|v| \le R}(1-\Theta_{\delta})\frac{\langle v_* \rangle^{\gamma}}{w\left(\frac{|v_*|}{\sqrt{2}}\right)}d\sigma dv_*\right) \|f\|_{L^{\infty}_{x,v}(w)}.
	\end{align*}
	For $\frac{1}{\delta} \ge 2R$, it holds that
	\begin{align*}
		\mathbf{1}_{|v| \le R}(1-\Theta_{\delta}) &\le \mathbf{1}_{|v| \le R}\mathbf{1}_{|v-v_*| \ge \frac{1}{\delta}}+\mathbf{1}_{|v_*-v| \le 2\delta}+\mathbf{1}_{|\cos\theta| \ge 1-2\delta}\\
		&\le \mathbf{1}_{|v_*| \ge \frac{1}{2\delta}}+\mathbf{1}_{|v_*-v| \le 2\delta}+\mathbf{1}_{|\cos\theta| \ge 1-2\delta}.
	\end{align*}
	Then we have
	\begin{align*}
		w(v)\mu(v)\mathbf{1}_{|v|\ge R} \rightarrow 0 \quad \text{as} \quad R \rightarrow \infty,
	\end{align*}
	and
	\begin{align*}
		\left(\int_{\mathbb{R}^3\times \mathbb{S}^2}\mathbf{1}_{|v| \le R}(1-\Theta_{\delta})\frac{\langle v_* \rangle^{\gamma}}{w\left(\frac{|v_*|}{\sqrt{2}}\right)}d\sigma dv_*\right) \rightarrow 0 \quad \text{as} \quad \delta \rightarrow 0 \quad \text{for fixed }R
	\end{align*}
	since $\frac{\langle v_* \rangle^{\gamma}}{w\left(\frac{|v_*|}{\sqrt{2}}\right)}$ is integrable on $\mathbb{R}^3$. Thus we can write
	\begin{align*}
		w(v)\nu(v)^{-1}\left|K_s^{(\delta)}f\right| \le &  w(v)\nu(v)^{-1}\int_{\mathbb{R}^3\times \mathbb{S}^2} B(v,v_*,\theta)\mathbf{1} _{|v|\ge R} \left[\mu(v_*')f(v')+\mu(v')f(v_*')\right]d\sigma dv_*\\
		&+\left(\varepsilon_w(R)+w(R)\varepsilon_w(\delta) \right)\|f\|_{L^{\infty}_{x,v}(w)},
	\end{align*}
	where
	\begin{align*}
		&\varepsilon_w(R) \rightarrow 0 \quad \text{as} \quad R \rightarrow \infty,\quad  \text{and} \quad \varepsilon_w(\delta) \rightarrow 0 \quad  \text{as} \quad \delta  \rightarrow 0.
	\end{align*} 
	For the first one in \eqref{Knsest11}, we use the estimate \eqref{Q+estimate}. Then we have
	\begin{align*}
		\|Q^+(\mu,f) +Q^+(f,\mu) \|_{L^{\infty}_{x,v}(\nu^{-1}w)} \le \varepsilon_w(R)\|\mu\| _{L^{\infty}_{x,v}(w)}\|f\|_{L^{\infty}_{x,v}(w)} \le \varepsilon_w(R)\|f\|_{L^{\infty}_{x,v}(w)}.
	\end{align*}
    and 
    \begin{align*}
		\varepsilon_w(R) \le C\frac{\mathbf{1}_{|v| \ge R}}{\langle v \rangle^\kappa} \le \frac{C}{(1+R^2)^{\kappa/2 }} \rightarrow 0 \quad \text{as} \quad R\rightarrow \infty, 
	\end{align*}
	Thus we deduce that
	\begin{align*}
		w(v)\nu(v)^{-1}\left|K_s^{(\delta)}f\right| \le \left(\varepsilon_w(R)+w(R)\varepsilon_w(\delta) \right)\|f\|_{L^{\infty}_{x,v}(w)}.
	\end{align*}
	Since $R$ can be arbitrarily large, we obtain
	\begin{align*}
		\|K_s^{(\delta)}f\|_{L^{\infty}_{x,v}(\nu^{-1} w)} \le C_K(\delta) \|f\|_{L^{\infty}_{x,v}(w)},
	\end{align*}
	where $C_K(\delta) \rightarrow 0$ as $\delta \rightarrow 0.$
\end{proof}

Now, we construct a solution of the full nonlinear equation by using the above estimates.

Before proceeding, we study the nonlinear equation of $f_1$, treating $g$ as a prescribed function. The following proposition establishes global existence in the space $L^\infty_{x,v}(w)$, provided that $\delta$, the initial datum, and $g$ satisfy suitable smallness conditions. In particular, the solution $f_1$ inherits the exponential decay rate of $g$.

\begin{proposition} \label{T2}
	Let \(w(v)=(1+|v|^2)^{\beta/2}e^{\alpha|v|^{\kappa}}\) with $\beta \ge 0 ,\alpha >0 ,0<\kappa<2$. Let \(f_0\) be in \(L^{\infty}_{x,v}(w)\) and \(g(t,x,v)\) in \(L^{\infty}_{x,v}(w)\). Then there exists \(\delta_w>0\) such that for any \(\delta\) in \((0,\delta_w]\), there exist $C_s$, $\eta_1$, and $\lambda_1>0$ such that if                  
	\[\|f_0\|_{L^{\infty}_{x,v}(w)}\le \eta_1^{1+p} \quad \text{and} \quad \|g(t)\|_{L^{\infty}_{x,v}(w)}\le \eta_1 e^{-\lambda_1t},\]
	then there exists a solution $f_1$ to the equation
	\[\partial_t f_1+v\cdot \nabla_x f_1 +\nu(v) f_1= K_s^{(\delta)}f_1+Q(f_1+g,f_1+g)\] 
	with \(f_0\)  as initial data. Moreover, this solution satisfies
	\[\|f_1(t)\|_{L_{x,v}^\infty(w)}\le C_se^{-\lambda_1t}(\|f_0\|_{L_{x,v}^\infty(w)}+\eta_1^{1+p})\]
	for all $t \ge 0$. The constant $C_s$, $\eta_1$, and $\lambda_1$ are constructive and only depend on $w$ and $\delta$.
\end{proposition}

\begin{remark}
	In the proof of proposition \ref{T2}, note that we choose $\lambda_1 = \frac{\varepsilon}{2}\nu_0 $. This will be crucial in the proof of proposition \ref{T1}.
\end{remark}

\begin{proof}
	By Duhamel's principle, we obtain
    \begin{align*}
    	f_1(t,x,v) &= e^{-\nu(v)t}f_0(x-vt,v) + \int_0^t e^{-\nu(v)(t-s)}K_s^{(\delta)}f_1(s,x-v(t-s),v)ds\\
    	& \quad + \int_0^t e^{-\nu(v)(t-s)}Q(f_1+g,f_1+g)(s,x-v(t-s),v)ds.
    \end{align*}
    Set $h_0 = 0$ and we consider the iterative system as follows:\\
    \begin{align*}
    h_{m+1}(t,x,v) &= e^{-\nu(v)t}f_0(x-vt,v)+ \int_0^t e^{-\nu(v)(t-s)}K_s^{(\delta)}(h_{m+1})(s,x-v(t-s),v)ds\\ 
    & \quad +\int_0^te^{-\nu(v)(t-s)}Q(h_m+g,h_m+g)(s,x-v(t-s),v)ds 
    \end{align*}
    with initial data $h_{m+1}(0,x,v) = f_0(x,v)$. Then
    \begin{align*}
    	\left|h_{m+1}(t,x,v)-h_m(t,x,v)\right| 
    	&\le \int_0^t \left|e^{-\nu(v)(t-s)}K_s^{(\delta)}(h_{m+1}-h_m)(s)\right|ds\\
    	&\quad +\int_0^t  \Big| e^{-\nu(v)(t-s)}\left(Q(h_m+g,h_m+g)-Q(h_{m-1}+g,h_{m-1}+g)\right)(s)\Big|ds\\
    	&= \int_0^t  \left|e^{-\nu(v)(t-s)}K_s^{(\delta)}(h_{m+1}-h_m)(s)\right|ds\\
    	&\quad +\int_0^t  \Big|e^{-\nu(v)(t-s)}\left[Q(h_m-h_{m-1},h_m+g)+Q(h_{m-1}+g,h_m-h_{m-1})(s)\right]\Big|ds
    \end{align*}
    Fix $0<\varepsilon<1/2$. First of all, we consider the bound of the integration for $K_s^{(\delta)}$. Then for $\varepsilon$ in $(0,1)$,
    \begin{align} \label{I11}
    	 &w(v)\int_0^t \left|e^{-\nu(v)(t-s)}   K_s^{(\delta)}(h_{m+1}-h_m)(s)\right|ds\nonumber\\
    	 &\le \frac{1}{1-\varepsilon}e^{-\frac{\varepsilon\nu_0}{2}t}\sup_{0\le s \le t} \left[e^{\varepsilon \nu_0s}\left\|K_s^{(\delta)}(h_{m+1}-h_m)(s)\right\|_{L^{\infty}_{x,v}(\nu^{-1}w)}\right]\nonumber\\
    	 &\le \frac{C_K(\delta)}{1-\varepsilon}e^{-\frac{\varepsilon\nu_0}{2}t}\sup_{0\le s \le t} \left[e^{\frac{\varepsilon\nu_0}{2}s}\left\|(h_{m+1}-h_m)(s)\right\|_{L^{\infty}_{x,v}(w)}\right],
    \end{align}
    where we applied Lemma \ref{Knsest}. Next, we consider the bound for the integration for the nonlinear term $Q$. Using Corollary \ref{Qestimate}, we have
   \begin{align} \label{I14}
   		&w(v)\int_0^t  \Big| e^{-\nu(v)(t-s)} \left[Q(h_m-h_{m-1},h_m+g)+Q(h_{m-1}+g,h_m-h_{m-1})(s)\right]\Big|ds\nonumber\\
   		&\le \frac{C_w}{1-\varepsilon} \sup_{0\le s \le t} \left[ e^{-\frac{\frac{\varepsilon\nu_0}{2}}{2}(t-s)}\|\left[Q(h_m-h_{m-1},h_m+g)+Q(h_{m-1}+g,h_m-h_{m-1})(s)\right]\|_{L^\infty_{x,v}(w\nu^{-1})} \right]\nonumber \\
   		&\le \frac{C_w}{1-\varepsilon} \left(\|h_m\|_{L_{[0,t]}^{\infty}L^{\infty}_{x,v}(w)}+\|h_{m-1}\|_{L_{[0,t]}^{\infty}L^{\infty}_{x,v}(w)}+\|g\|_{L_t^{\infty}L^{\infty}_{x,v}(w)}\right)\nonumber\\
   		&\quad \times e^{-\frac{\varepsilon}{2}\nu_0t} \sup_{0\le s \le t} \left[e^{\frac{\varepsilon}{2}\nu_0s}\left\|(h_m-h_{m-1})(s)\right\|_{L^{\infty}_{x,v}(w)}\right].
   \end{align}
   Combining the bounds \eqref{I11} and \eqref{I14}, for any $T_0 >0$ and $t \in [0,T_0]$,
   \begin{align*}
   	\sup_{0\le s \le t} \left[e^{\frac{\varepsilon}{2}\nu_0s}\left\|(h_{m+1}-h_{m})(s)\right\|_{L^{\infty}_{x,v}(w)}\right] 
   	&\le \frac{C_K(\delta)}{1-\varepsilon}\sup_{0\le s \le t} \left[e^{\frac{\varepsilon}{2}\nu_0s}\left\|(h_{m+1}-h_m)(s)\right\|_{L^{\infty}_{x,v}(w)}\right]\\
   	&\quad + \frac{C_w}{1-\varepsilon}\left(\|h_m\|_{L_{[0,t]}^{\infty}L^{\infty}_{x,v}(w)}+\|h_{m-1}\|_{L_{[0,t]}^{\infty}L^{\infty}_{x,v}(w)}+\|g\|_{L_t^{\infty}L^{\infty}_{x,v}(w)}\right)\\
   	&\quad \quad \times \sup_{0\le s \le t} \left[e^{\frac{\varepsilon}{2}\nu_0s}\left\|(h_m-h_{m-1})(s)\right\|_{L^{\infty}_{x,v}(w)}\right].
   \end{align*}
   We can choose $\delta_w>0$ such that
   \begin{align*}
   	\frac{C_K(\delta)}{1-\varepsilon}<\frac{1}{2}
   \end{align*}
  for $0<\delta \le \delta_w$ since $C_K(\delta) \rightarrow 0$ as $\delta \rightarrow 0$. From our choice of constants, we can deduce
   \begin{align*}
   	\sup_{0\le s \le t} \left[e^{\frac{\varepsilon}{2}\nu_0s}\left\|(h_{m+1}-h_{m})(s)\right\|_{L^{\infty}_{x,v}(w)}\right] 
   	&\le C_w\left(\|h_m\|_{L_{[0,t]}^{\infty}L^{\infty}_{x,v}(w)}+\|h_{m-1}\|_{L_{[0,t]}^{\infty}L^{\infty}_{x,v}(w)}+\|g\|_{L_t^{\infty}L^{\infty}_{x,v}(w)}\right)\\
   	&\quad\times \sup_{0\le s \le t} \left[e^{\frac{\varepsilon}{2}\nu_0s}\left\|(h_m-h_{m-1})(s)\right\|_{L^{\infty}_{x,v}(w)}\right]
   \end{align*}
   for all $0\le t\le T_0$. Since $C_w$ is independent of $T_0$, we can extend the inequality to all time $t \ge 0$.
   We claim that $\left\|h_{l+1}\right\|_{L^\infty_{[0,t]}  L^{\infty}_{x,v}(w)}$ is uniformly bounded. In a similar argument to the above, we can easily deduce that
   \begin{align} \label{I21}
   	w(v)\int_0^t \left|e^{-\nu(v)(t-s)}  K_s^{(\delta)}h_{m+1}(s)\right|ds \le C e^{-\frac{\varepsilon\nu_0}{2}t}\sup_{0\le s \le t} \left[e^{\frac{\varepsilon\nu_0}{2}s}\left\|h_{m+1}(s)\right\|_{L^{\infty}_{x,v}(w)}\right].
   \end{align}
   We also obtain 
   \begin{align*}
   	& w(v)\int_0^t \Big|e^{-\nu(v)(t-s)} Q(h_m+g,h_m+g)(s)\Big|ds\\
    &\le C_\varepsilon e^{-\frac{\varepsilon}{2}\nu_0t}\sup_{0\le s \le t} \left[e^{\frac{\varepsilon}{2}\nu_0s}\left\|Q(h_m+g,h_m+g)(s)\right\|_{L^{\infty}_{x,v}(\nu^{-1}w)}\right]\\
   	&\le C_\varepsilon e^{-\frac{\varepsilon}{2}\nu_0t}\bigg\{\left(\|h_m\|_{L_{[0,t]}^{\infty}L^{\infty}_{x,v}(w)}+\|g\|_{L_{[0,t]}^{\infty}L^{\infty}_{x,v}(w)}\right)\sup_{0\le s \le t} \left[e^{\frac{\varepsilon}{2}\nu_0s}\left\|h_m(s)\right\|_{L^{\infty}_{x,v}(w)}\right]\\
   	&\quad + \|g\|_{L_{[0,t]}^{\infty}L^{\infty}_{x,v}(w)}\sup_{0\le s \le t} \left[e^{\frac{\varepsilon}{2}\nu_0s}\left\|g(s)\right\|_{L^{\infty}_{x,v}(w)}\right]\bigg\}.
   \end{align*}
   Take $\lambda_1 = \frac{\varepsilon}{2}\nu_0$ and assume $\|g(t)\|_{L^{\infty}_{x,v}(w)}\le \eta_1 e^{-\lambda_1t}$.
   Then we get
   \begin{align} \label{I22}
   	&w(v)\int_0^t \Big|e^{-\nu(v)(t-s)} Q(h_m+g,h_m+g)(s)\Big|ds\nonumber\\
   	&\le C_\varepsilon e^{-\frac{\varepsilon}{2}\nu_0t}\bigg\{\left(\|h_m\|_{L_{[0,t]}^{\infty}L^{\infty}_{x,v}(w)}+\eta_1\right)\sup_{0\le s \le t} \left[e^{\frac{\varepsilon}{2}\nu_0s}\left\|h_m(s)\right\|_{L^{\infty}_{x,v}(w)}\right] + \eta_1^2\bigg\}.
   \end{align}
   Write $A_m = \sup_{0\le s \le t} \left[e^{\frac{\varepsilon}{2}\nu_0s}\left\|h_m(s)\right\|_{L^{\infty}_{x,v}(w)}\right]$. By \eqref{I21} and \eqref{I22},  there exist $C_0^{(1)}$, $C_0^{(2)}>0$ such that
   \begin{align*}
   	A_{m+1} \le C_0^{(1)}\|f_0\|_{L^{\infty}_{x,v}(w)} + C_0^{(2)}\eta_1^2+C_0^{(2)}\left(\|h_m\|_{L_{[0,t]}^{\infty}L^{\infty}_{x,v}(w)}+\eta_1\right)A_m.
   \end{align*}
   By using an induction, we will show that 
   \begin{align*}
   	\sup_m A_m \le \left(1+C_0^{(1)}\right)\|f_0\|_{L^{\infty}_{x,v}(w)} + \left(1+C_0^{(1)}\right)\eta_1^{1+p}.
   \end{align*}
   We deduce
   \begin{align*}
   	A_{m+1} &\le C_0^{(1)}\|f_0\|_{L^{\infty}_{x,v}(w)} + C_0^{(2)}\eta_1^2\\
   	&\quad + C_0^{(2)}\left\{\left(2+C_0^{(1)}\right)\eta_1+\left(1+C_0^{(1)}\right)\eta_1^{1+p}\right\}\times\left\{\left(1+C_0^{(1)}\right)\|f_0\|_{L^{\infty}_{x,v}(w)}+\left(1+C_0^{(1)}\right)\eta_1^{1+p}\right\}\\
   	&\le C_0^{(1)}\|f_0\|_{L^{\infty}_{x,v}(w)} + C_0^{(2)}\eta_1^2\\
   	&\quad + \|f_0\|_{L^{\infty}_{x,v}(w)} + C_0^{(2)}\left\{\left(2+C_0^{(1)}\right)\eta_1+\left(1+C_0^{(1)}\right)\eta_1^{1+p}\right\}\left(1+C_0^{(1)}\right)\eta_1^{1+p}\\
   	&\le \left(1+C_0^{(1)}\right)\|f_0\|_{L^{\infty}_{x,v}(w)} + \left(1+C_0^{(1)}\right)\eta_1^{1+p}
   \end{align*}
   by choosing sufficiently small $\eta_1 \ll 1$ such that
   \begin{align*}
   	& C_0^{(2)}\left\{\left(2+C_0^{(1)}\right)\eta_1+\left(1+C_0^{(1)}\right)\eta_1^{1+p}\right\}\left(1+C_0^{(1)}\right) \le 1, \\
   	& C_0^{(2)}\left\{\left(2+C_0^{(1)}\right)\eta_1+\left(1+C_0^{(1)}\right)\eta_1^{1+p}\right\}\left(1+C_0^{(1)}\right)\eta_1^{1+p}+ C_0^{(2)}\eta_1^2 \le \left(1+C_0^{(1)}\right)\eta_1^{1+p}.
   \end{align*}
   It follows that
   \begin{align*}
   	\sup_m A_m \le \left(1+C_0^{(1)}\right)\|f_0\|_{L^{\infty}_{x,v}(w)} + \left(1+C_0^{(1)}\right)\eta_1^{1+p} \le 2\left(1+C_0^{(1)}\right)\eta_1^{1+p}\le 2\left(1+C_0^{(1)}\right)\eta_1.
   \end{align*}
   Thus, we obtain the uniform bound of $\left\|h_{m+1}\right\|_{L^\infty_{[0,t]}  L^{\infty}_{x,v}(w)}$ as follows:
   \begin{align*}
   	\left\|h_{m+1}\right\|_{L^\infty_{[0,t]}  L^{\infty}_{x,v}(w)} \le 2\left(1+C_0^{(1)}\right)\eta_1.
   \end{align*}
   Hence we conclude that
   \begin{align*}
   		\sup_{0\le s \le t} \left[e^{\frac{\varepsilon}{2}\nu_0s}\left\|(h_{m+1}-h_{m})(s)\right\|_{L^{\infty}_{x,v}(w)}\right] &\le 6\left(1+C_0^{(1)}\right)\eta_1C \sup_{0\le s \le t} \left[e^{\frac{\varepsilon}{2}\nu_0s}\left\|(h_m-h_{m-1})(s)\right\|_{L^{\infty}_{x,v}(w)}\right].
   \end{align*}
   Choosing sufficiently small $\eta_1$ such that $6\left(1+C_0^{(1)}\right)\eta_1 C < 1$, $\left(h_m\right)_{m\in \mathbb{N}}$ is Cauchy in $L_{t}^\infty L_{x,v}^\infty(w)$, so that  $\left(h_m\right)_{m\in \mathbb{N}}$ converges to $f_1$ in $L_{t}^\infty L_{x,v}^\infty(w)$. Hence $f_1$ is a solution to the equation
   \begin{align*}
   		\partial_t f_1+v\cdot \nabla_x f_1 +\nu(v) f_1= K_s^{(\delta)}f_1+Q(f_1+g,f_1+g),
   \end{align*}
   and we get the following exponential decay
   \begin{align*}
   	\|f_1(t)\|_{L_{x,v}^\infty(w)}\le C_se^{-\lambda_1t}(\|f_0\|_{L_{x,v}^\infty(w)}+\eta_1^{1+p}).
   \end{align*}
   for all $t \ge 0$.
\end{proof}

\bigskip
We next consider the linear equation for $f_2$ driven by the term $K_{ns}^{(\delta)}g$. The following proposition establishes the global existence and uniqueness in the stronger weight space. Moreover, under the conservation assumption, the solution $f_2$ inherits the expoential decay rate of $g$. 

\begin{proposition} \label{T3}
	Let $w(v)=(1+|v|^2)^{\beta/2}e^{\alpha|v|^{\kappa}}$ with $\beta \ge 0 ,\alpha >0 ,0<\kappa<2$. Let $g=g(t,x,v)$ be in \(L^\infty_tL^{\infty}_{x,v}(w)\). Then there exists a unique function $f_2$ in \(L^\infty_{loc,t}L^{\infty}_{x,v}(\langle v \rangle^{\beta_3} \mu^{-1/2})\) with $\beta_3 > \max\{3,\beta\}$ such that
	\[\partial_t f_2+v\cdot\nabla_x f_2 +Lf_2 = K_{ns}^{(\delta)}g \quad and \quad f_2(0,x,v)=0.\]
    If $P_{G_L}(f_2+g)=0$ and there exist $\lambda_g$, $\eta_g>0$ such that $\lambda_g <\lambda_L$, where $\lambda_L$ is defined in proposition \ref{linearLinftydecay}, and for all $t \ge 0$ 
	\[\|g(t)\|_{L_{x,v}^\infty(w)}\le \eta_ge^{-\lambda_g t},\]
	then there exists $C_{ns}>0$ such that
	\[\|f_2(t)\|_{L_{x,v}^\infty(\langle v \rangle^{\beta_3} \mu^{-1/2})}\le C_{ns}\eta_ge^{-\lambda_g t}\]
	for all $t \ge 0$. This guarantees $f_2 \in L^\infty_{t}L^{\infty}_{x,v}(\langle v \rangle^{\beta_3} \mu^{-1/2})$. The constant $C_{ns}$ only depends on $\lambda_L, \lambda_g$.
\end{proposition}

\begin{proof}
	From Lemma \ref{Ksest} and $g \in L^\infty_tL^{\infty}_{x,v}(w)$, we obtain $K_{ns}^{(\delta)}(g) \in L^\infty_{t}L^{\infty}_{x,v}(\langle v \rangle^\beta \mu^{-1/2}).$
	By proposition \ref{linearLinftydecay}, there exists a unique function $f_2$ in \(L^\infty_{\text{loc},t}L^{\infty}_{x,v}(\langle v \rangle^{\beta_3} \mu^{-1/2})\) with $\beta_3>\max\{3,\beta\}$ such that
	\[\partial_t f_2+v\cdot\nabla_x f_2 +Lf_2 = K_{ns}^{(\delta)}(g) \quad \text{with} \ f_2(0,x,v)=0,\]
	and it is given by
	\begin{align*}
		f_2(t,x,v) = \int_0^tG(t,s)K_{ns}^{(\delta)}(g)(s)ds.
	\end{align*}
	Assume $P_{G_L}(f_2+g)=0$ and there exist $\lambda_g$, $\eta_g>0$ such that $\|g(t)\|_{L_{x,v}^\infty(w)}\le \eta_ge^{-\lambda_g t}$ for all $t\ge 0$. Then for all $t\ge 0$, we have
	\begin{align} \label{est32}
		\|P_{G_L}(f_2)(t)\|_{L_{x,v}^\infty(\langle v \rangle^{\beta_3} \mu^{-1/2})} &= \|P_{G_L}(g)(t)\|_{L_{x,v}^\infty(\langle v \rangle^{\beta_3} \mu^{-1/2})}\nonumber\\ 
		&\le C \|g(t)\|_{L_{x,v}^\infty(w)}\\
		&\le C\eta_ge^{-\lambda_g t}. \nonumber
	\end{align}
	We claim $\partial_t\left[P_{G_L}^\perp\left(f_2\right)\right] = G_L\left[P_{G_L}^\perp\left(f_2\right)\right]+P_{G_L}^\perp\left[K_{ns}^{(\delta)}(g)\right]$. Taking $P_{G_L}$ to both sides of $\partial_t f_2 = G_L(f_2) + K_{ns}^{(\delta)}(g)$, we obtain
	\begin{align} \label{est31}
		\partial_t\left[P_{G_L} \left(f_2\right)\right] = P_{G_L} \left[ \partial_t f_2 \right] = P_{G_L} \left[ G_L(f_2) + K_{ns}^{(\delta)}(g) \right] = P_{G_L} \left[G_L(f_2)\right] + P_{G_L} \left[K_{ns}^{(\delta)}(g)\right].
	\end{align}
	Since $(L(f_2),\mu)_{L^2_{v}(\mu^{-1/2})}=0$, we have
	\begin{align*}
		P_{G_L} \left[G_L(f_2)\right] = 0.
	\end{align*}
	From $G_L\left[P_{G_L}\left(f_2\right)\right]=0$,  \eqref{est31} becomes
	\begin{align*}
		\partial_t\left[P_{G_L} \left(f_2\right)\right] = P_{G_L} \left[K_{ns}^{(\delta)}(g)\right]= G_L\left[P_{G_L}\left(f_2\right)\right] +P_{G_L} \left[K_{ns}^{(\delta)}(g)\right],
	\end{align*}
	and it follows that 
	\begin{align*}
\partial_t\left[P_{G_L}^\perp\left(f_2\right)\right] = G_L\left[P_{G_L}^\perp\left(f_2\right)\right]+P_{G_L}^\perp\left[K_{ns}^{(\delta)}(g)\right].
	\end{align*}
	By Duhamel's principle, we get
	\begin{align*}
		P_{G_L}^\perp\left(f_2\right) = \int_0^t G(t,s)P_{G_L}^\perp\left(K_{ns}^{(\delta)}(g)(s)\right)ds.
	\end{align*}
	Thanks to proposition \ref{linearLinftydecay} and Lemma \ref{Ksest}, we deduce that
	\begin{align}\label{est33}
		\|P_{G_L}^{\perp}(f_2)(t)\|_{L_{x,v}^\infty(\langle v \rangle^{\beta_3} \mu^{-1/2})} 
		&\le C_L \int_0^t e^{-\lambda_L(t-s)}\|K_{ns}^{(\delta)}(g)(s)\|_{L_{x,v}^\infty(\langle v \rangle^{\beta_3} \mu^{-1/2})}ds\nonumber\\ 
		&\le C \int_0^t e^{-\lambda_L(t-s)}\|g(s)\|_{L_{x,v}^\infty(w)}ds\\
		&\le C \int_0^t e^{-\lambda_L(t-s)}\eta_g e^{-\lambda_g s}ds \nonumber\\ & = C\eta_g\frac{1}{\lambda_L-\lambda_g}\left(e^{-\lambda_g t}-e^{-\lambda_L t}\right)\nonumber\\
		&\le C\eta_g e^{-\lambda_g t}.\nonumber
	\end{align}
	Hence we conclude that
	\begin{align*}
		\|f_2(t)\|_{L_{x,v}^\infty(\langle v \rangle^{\beta_3} \mu^{-1/2})}\le C\eta_ge^{-\lambda_g t} 
	\end{align*} 
	for all $t \ge 0$.
    Therefore, this guarantees $f_2 \in L^\infty_{t}L^{\infty}_{x,v}(\langle v \rangle^{\beta_3} \mu^{-1/2})$.
\end{proof}

\bigskip

The following propositions (Propositon \ref{T1}, \ref{UT}, \ref{positivityprop}) establish the small-amplitude result (Proposition \ref{Mainresultsmall}). They assert the global existence, stability, uniqueness, positivity of solutions to the Boltzmann equation \eqref{BE} for sufficiently small initial data $f_0$ in the weighted space $L^\infty_{x,v}(w)$.

\begin{proposition} 
\label{T1}
	Let \(w(v)=(1+|v|^2)^{\beta/2}e^{\alpha|v|^{\kappa}}\) with $\beta \ge 0 ,\alpha >0 ,0<\kappa<2$. Then there exists $\eta_0 >0$ such that for any $f_0$ in $L^\infty_{x,v}(w)$ with $P_{G_L}(f_0)=0$ and $\|f_0\|_{L^\infty_{x,v}(w)} \le \eta_0^{1+p}$ for some $0<p<1$, there exists at least one solution f of \eqref{BE-f} with initial data $f_0$. Moreover, f satisfies the mass, momentum, energy conservation and there exist $C_\infty,\lambda_\infty >0$ such that
	\begin{align*}
		\|f(t)\|_{L_{x,v}^\infty(w)}\le C_\infty e^{-\lambda_\infty t}(\|f_0\|_{L_{x,v}^\infty(w)}+\eta_0^{1+p})
	\end{align*}
	for all $t \ge 0$.
\end{proposition}

\begin{proof}
	We will use the sequence argument $f_2^{(n)} \rightarrow f_1^{(n+1)} \rightarrow f_2^{(n+1)} \rightarrow \cdots$ and we will deduce the following estimates: 
	\begin{align*}
		&\sup_n\|f_1^{(n)}(t)\|_{L_{x,v}^\infty(w)} \le C_se^{-\lambda t}(\|f_0\|_{L_{x,v}^\infty(w)}+\eta_0^{1+p}),\\
		&\sup_n\|f_2^{(n)}(t)\|_{L_{x,v}^\infty(\langle v \rangle^{\beta_3} \mu^{-1/2})} \le C_sC_{ns}e^{-\lambda t}(\|f_0\|_{L_{x,v}^\infty(w)}+\eta_0^{1+p}).
	\end{align*} 
	Set $f_1^{(0)} = f_2^{(0)} = 0$ and we use the iterative system as follows:
	\begin{equation} \label{T1d4}
		\begin{aligned}
			&\partial_t f_1^{(n+1)}+v\cdot\nabla_x f_1^{(n+1)} +\nu(v)f_1^{(n+1)}=K_s^{(\delta)}f_1^{(n+1)}+Q(f_1^{(n+1)}+f_2^{(n)},f_1^{(n+1)}+f_2^{(n)})\quad  \text{with}\quad f_1^{(n+1)}(0,x,v)=f_0,\\
			&\partial_t f_2^{(n+1)}+v\cdot\nabla_x f_2^{(n+1)} +Lf_2^{(n+1)} = K_{ns}^{(\delta)}f_1^{(n+1)}\quad \text{with}\quad f_2^{(n+1)}(0,x,v)=0.
		\end{aligned}
	\end{equation}
	Assume that $\|f_0\|_{L^\infty_{x,v}(w)} \le \eta_0^{1+p}$ for some $0<p<1$.\\
	First, we consider the step $f_2^{(n)} \rightarrow f_1^{(n+1)}$ and assume $\|f_2^{(n)}(t)\|_{L_{x,v}^\infty(\langle v \rangle^{\beta_3} \mu^{-1/2})} \le C_sC_{ns}e^{-\lambda t}(\|f_0\|_{L_{x,v}^\infty(w)}+\eta_0^{1+p})$ for all $t\ge 0$.\\
	In proposition \ref{T2}, putting $g = f_2^{(n)}$ and $\tilde{C}= \sup_{x,v}\left(\langle v \rangle^{-{\beta_3}} \mu^{1/2} w \right)$, we get
	\begin{align*}
		\|g(t)\|_{L_{x,v}^\infty(w)} = \|f_2^{(n)}(t)\|_{L_{x,v}^\infty(w)} &\le \tilde{C}\|f_2^{(n)}(t)\|_{L_{x,v}^\infty(\langle v \rangle^{\beta_3} \mu^{-1/2})} \le \tilde{C}C_sC_{ns}e^{-\lambda t}(\|f_0\|_{L_{x,v}^\infty(w)}+\eta_0^{1+p})\\
		&\le \tilde{C}C_sC_{ns}e^{-\lambda t}(2\eta_0^{1+p}) \le \eta_0 e^{-\lambda t},
	\end{align*}
	where, in the last inequality, we choose $\eta_0$ sufficiently small such that $\tilde{C}C_sC_{ns}(2\eta_0^{1+p}) \le \eta_0$. Therefore, $g=f_2^{(n)}$ satisfies the smallness assumption of Proposition \ref{T2}. Applying Proposition \ref{T2}, we obtain
	\begin{align*}
		\|f_1^{(n+1)}(t)\|_{L_{x,v}^\infty(w)} \le C_se^{-\lambda t}(\|f_0\|_{L_{x,v}^\infty(w)}+\eta_0^{1+p}).
	\end{align*}
	for all $t \ge 0.$
	Next, we consider the step $f_1^{(n+1)} \rightarrow f_2^{(n+1)}$ and assume $\|f_1^{(n)}(t)\|_{L_{x,v}^\infty(w)} \le C_se^{-\lambda t}(\|f_0\|_{L_{x,v}^\infty(w)}+\eta_0^{1+p})$ for all $t \ge 0$.\\
	From $\partial_t \left(f_1^{(n+1)}+f_2^{(n+1)}\right)=G_L\left(f_1^{(n+1)}+f_2^{(n+1)}\right)+Q\left(f_1^{(n+1)}+f_2^{(n)},f_1^{(n+1)}+f_2^{(n)}\right)$, we can deduce $\partial_t P_{G_L}\left(f_1^{(n+1)}+f_2^{(n+1)}\right)=0$, which implies
	\begin{align} \label{sequenceconservation}
		P_{G_L}\left(f_1^{(n+1)}+f_2^{(n+1)}\right)=P_{G_L}\left(f_0\right)=0.
	\end{align}
	In proposition \ref{T3}, taking $g=f_1^{(n+1)}$, $\eta_g = C_s\left(\|f_0\|_{L_{x,v}^\infty(w)}+\eta_0^{1+p}\right)$, and $\lambda_g = \lambda_1=\frac{\varepsilon}{2} \nu_0$, we get
	\begin{align*}
		\|f_2^{(n+1)}(t)\|_{L_{x,v}^\infty(\langle v \rangle^{\beta_3} \mu^{-1/2})} \le C_{ns}\eta_g e^{-\lambda_g t} = C_sC_{ns}e^{-\lambda t}\left(\|f_0\|_{L_{x,v}^\infty(w)}+\eta_0^{1+p}\right),
	\end{align*}
	where we choose enough small $\varepsilon >0$ such that $\lambda_g  = \frac{\varepsilon}{2} \nu_0 < \lambda_L$.\\
	Now, we consider the sequence $f_2^{(n+1)}-f^{(n)}_2$. This sequence satisfies 
	\begin{align*}
		\partial_t (f_2^{(n+1)}-f^{(n)}_2)+v\cdot\nabla_x (f_2^{(n+1)}-f^{(n)}_2) +L(f_2^{(n+1)}-f^{(n)}_2) = K_{ns}^{(\delta)}(f_1^{(n+1)}-f_1^{(n)}).
	\end{align*}
	From \eqref{sequenceconservation}, \eqref{est32}, and \eqref{est33}, we have
	\begin{align*}
		\|P_{G_L}(f_2^{(n+1)}-f^{(n)}_2)(t)\|_{L_{x,v}^\infty(\langle v \rangle^{\beta_3} \mu^{-1/2})} &\le C_{\beta_3}\sup_{t\ge 0}\|(f_1^{(n+1)}-f_1^{(n)})(t)\|_{L_{x,v}^\infty(w)},\\
		\|P_{G_L}^{\perp}(f_2^{(n+1)}-f^{(n)}_2)(t)\|_{L_{x,v}^\infty(\langle v \rangle^{\beta_3} \mu^{-1/2})} 
		&\le C_\beta \int_0^t e^{-\lambda_L(t-s)}\|(f_1^{(n+1)}-f_1^{(n)})(s)\|_{L_{x,v}^\infty(w)}ds\\
		& \le C_{\beta,\lambda_L} \sup_{t\ge 0}\|(f_1^{(n+1)}-f_1^{(n)})(t)\|_{L_{x,v}^\infty(w)},
	\end{align*}
	which implies that 
	\begin{align} \label{f2Cauchy}
		\sup_{t\ge 0 }\|(f_2^{(n+1)}-f^{(n)}_2)(t)\|_{L_{x,v}^\infty(w)}&\le C_{\alpha,\kappa}\sup_{t\ge 0 }\|(f_2^{(n+1)}-f^{(n)}_2)(t)\|_{L_{x,v}^\infty(\langle v \rangle^{\beta_3} \mu^{-1/2})}\nonumber\\&\le C_{\alpha,\beta,\kappa,\lambda_L} \sup_{t\ge 0}\|(f_1^{(n+1)}-f_1^{(n)})(t)\|_{L_{x,v}^\infty(w)}.
	\end{align}
	Hence, from the sequence argument, we get
	\begin{equation} \label{T1d1}
		\begin{aligned}
			&\sup_n\|f_1^{(n)}(t)\|_{L_{x,v}^\infty(w)} \le C_se^{-\lambda t}(\|f_0\|_{L_{x,v}^\infty(w)}+\eta_0^{1+p}),\\
			&\sup_n\|f_2^{(n)}(t)\|_{L_{x,v}^\infty(\langle v \rangle^{\beta_3} \mu^{-1/2})} \le C_sC_{ns}e^{-\lambda t}(\|f_0\|_{L_{x,v}^\infty(w)}+\eta_0^{1+p}).
		\end{aligned}
	\end{equation}
	Combining the scheme in proposition \ref{T2} for $f_1^{(n)}$ and \eqref{f2Cauchy}, \eqref{T1d1}, $\left(f_1^{(n)}\right)_{n\in \mathbb{N}}$ is Cauchy, and converges strongly to $f_1$ in $L^\infty_tL^\infty_{x,v}(w)$. It follows from \eqref{f2Cauchy} that $\left(f_2^{(n)}\right)_{n\in \mathbb{N}}$ is Cauchy, and converges strongly to $f_2$ in $L^\infty_tL^\infty_{x,v}(w)$. By \eqref{T1d1}, the following holds
	\begin{align*} 
		\|f_1(t)\|_{L_{x,v}^\infty(w)} \le C_se^{-\lambda t}(\|f_0\|_{L_{x,v}^\infty(w)}+\eta_0^{1+p}),\\ 
		\|f_2(t)\|_{L_{x,v}^\infty(w)} \le C_{\alpha,\kappa} C_sC_{ns}e^{-\lambda t}(\|f_0\|_{L_{x,v}^\infty(w)}+\eta_0^{1+p})
	\end{align*}
	for all $t\ge 0$.
	Therefore, taking the limit to the iterative system \eqref{T1d4}, $f=f_1+f_2$ is a solution to \eqref{BE-f} and we conclude that
	\begin{align*}
		P_{G_L}(f) = 0 \quad \text{and} \quad 
		\quad\|f(t)\|_{L_{x,v}^\infty(w)} \le Ce^{-\lambda t}(\|f_0\|_{L_{x,v}^\infty(w)}+\eta_0^{1+p})
	\end{align*}
	for all $t \ge 0$.
\end{proof}

\bigskip

\begin{proposition} [Uniqueness] \label{UT}
	Let \(w(v)=(1+|v|^2)^{\beta/2}e^{\alpha|v|^{\kappa}}\) with $\beta \ge 0 ,\alpha >0 ,0<\kappa<2$. There exists \(\eta_0 >0 \) such that for any \(f_0 \) in \(L^{\infty}_{x,v}(w)\) with \(\|f_0\|_{L^{\infty}_{x,v}(w)}\le \eta_0^{1+p}\) for \( 0<p<1 \) and \( P_{G_L}(f_0)=0 \), there exists at most one solution \(f(t,x,v)\) in \(L^\infty_tL^{\infty}_{x,v}(w)\) to \eqref{BE-f} with \(f_0\)  as initial data. 
\end{proposition}

\begin{proof}
	Let $f$ be a solution constructed in Proposition \ref{T1} of the equation \eqref{BE-f} with initial data $f_0$. Suppose that there exists a solution $\tilde{f}$ in $L^\infty_{t}L^\infty_{x,v}(w)$ to the equation \eqref{BE-f} with same initial data $f_0$. Let $\tilde{w}(v) =\langle v \rangle^\beta e^{\frac{\alpha}{2}|v|^\kappa}$. We define
	\begin{align*}
		M_\infty:= 1+ \sup_{0\le t< \infty} \|\tilde{f}(t)\|_{L^\infty_{x,v}(w)}.
	\end{align*}
	Assume that $\tilde{f}(t_0)=f(t_0)$ for some $t_0 \ge 0$.
	For $t\ge 0$, Duhamel principle gives
	\begin{align*}
		\tilde{w}\tilde{f}(t_0+t,x,v) =e^{-\nu(v)t} \tilde{w}\tilde{f}(t_0,x-vt,v)+\int_0^t e^{-\nu(v)(t-s)} \left[ \tilde{w}K\tilde{f}+\tilde{w}Q(\tilde{f},\tilde{f})\right](t_0+s,x-v(t-s),v)ds.
	\end{align*}
	By Lemma \ref{Ksest} and Lemma \ref{Knsest}, we obtain
	\begin{align*}
		\left|\int_0^t e^{-\nu(v)(t-s)}\tilde{w}K\tilde{f}(t_0+s,x-v(t-s),v)ds  \right| \le C_{\tilde{w}}t M_\infty,
	\end{align*}
	and from Lemma \ref{Qestimate}, we also have
	\begin{align*}
		\left|\int_0^t e^{-\nu(v)(t-s)}\tilde{w}Q(\tilde{f},\tilde{f})(t_0+s,x-v(t-s),v)ds  \right| \le C_{\tilde{w}} t M_\infty^2.
	\end{align*} 
	Thus, we derive
	\begin{align} \label{uni1}
		\|\tilde{w}\tilde{f}(t_0+t)\|_{L^\infty_{x,v}} \le \|\tilde{w}\tilde{f}(t_0)\|_{L^\infty_{x,v}} +C_w t M_\infty^2.
	\end{align}
	for $t\ge 0$. By Theorem \ref{T1},  $f$ has the following bound:
	\begin{align*}
		\sup_{t\ge 0}\|f(t)\|_{L^\infty_{x,v}(w)} \le 2C_\infty \eta_0^{1+p}.
	\end{align*}
	This gives
	\begin{align*}
		\|\tilde{f}(t_0)\|_{L^\infty_{x,v}(\tilde{w})} = \|f(t_0)\|_{L^\infty_{x,v}(\tilde{w})} \le 2C_\infty \eta_0^{1+p}.
	\end{align*}
	We choose sufficiently small $t_*>0$ so that $C_w t_* M_\infty^2 \le 2C_\infty \eta_0^{1+p}$, and it follows from \eqref{uni1} that
	\begin{align*}
		\sup_{t_0 \le t \le t_0+t_*} \|\tilde{f}(t)\|_{L^\infty_{x,v}(\tilde{w})} \le 4C_\infty \eta_0^{1+p}. 
	\end{align*}
		Subtracting the equations satisfied by $f$ and $\tilde{f}$, we get
	\begin{align*}
		\partial_t (f-\tilde{f})+v \cdot \nabla_x (f-\tilde{f}) +L(f-\tilde{f})=Q(f-\tilde{f},f)+Q(\tilde{f},f-\tilde{f})
	\end{align*}
	and decompose the linear operator $L$ to $K_{ns}^{(\delta)}$ and $K_s^{(\delta)}$
	\begin{align*}
		\partial_t (f-\tilde{f})+v \cdot \nabla_x (f-\tilde{f}) +\nu(v)(f-\tilde{f})=K_{ns}^{(\delta)}(f-\tilde{f})+K_s^{(\delta)}(f-\tilde{f})+Q(f-\tilde{f},f)+Q(\tilde{f},f-\tilde{f}).
	\end{align*}
	By Duhamel's principle,
	\begin{align*}
		(f-\tilde{f})(t_0+t,x,v) &= \int^t_0 e^{-\nu(v)(t-s)}\left[K_s^{(\delta)}(f-\tilde{f})(t_0+s)+Q(f-\tilde{f},f)(t_0+s)+Q(\tilde{f},f-\tilde{f})(t_0+s)\right]ds\\
		&\quad + \int^t_0 e^{-\nu(v)(t-s)}K_{ns}^{(\delta)}(f-\tilde{f})(t_0+s)ds
	\end{align*}
	 for $0\le t \le t_*$. By the proof of proposition \ref{T2} and choosing enough small $\delta >0$, we obtain
	\begin{align*}
		\left|\int^t_0 e^{-\nu(v)(t-s)} \tilde{w}(v)  \left[K_s^{(\delta)}(f-\tilde{f})(t_0+s)\right]ds\right| &\le\left(1-C_s \right) \|f-\tilde{f}\|_{L^{\infty}_{[t_0,t_0+t_*]}L^{\infty}_{x,v}(\tilde{w})}.
	\end{align*}
	Similarly,
	\begin{align*}
		&\left|\int^t_0 e^{-\nu(v)(t-s)}\tilde{w}(v) \left[Q(f-\tilde{f},f)(t_0+s)+Q(\tilde{f},f-\tilde{f})(t_0+s)\right]ds\right| \\
		& \quad \le C_Q\left(\|f\|_{L^{\infty}_{[t_0,t_0+t_*]}L^{\infty}_{x,v}(\tilde{w})}+\|\tilde{f}\|_{L^{\infty}_{[t_0,t_0+t_*]}L^{\infty}_{x,v}(\tilde{w})}\right) \|f-\tilde{f}\|_{L^{\infty}_{[t_0,t_0+t_*]}L^{\infty}_{x,v}(\tilde{w})}.
	\end{align*}
	By using Lemma \ref{Ksest},
	\begin{align*}
		\left|\int^t_0 e^{-\nu(v)(t-s)} \tilde{w}(v) \left[K_{ns}^{(\delta)}(f-\tilde{f})(t_0+s)\right]ds\right| 
		&\le \int^t_0 \|K_{ns}^{(\delta)}(f-\tilde{f})(t_0+s)\|_{L^{\infty}_{x,v}(\tilde{w})}ds\\
		& \le C\int^t_0 \|K_{ns}^{(\delta)}(f-\tilde{f})(t_0+s)\|_{L^{\infty}_{x,v}(\langle v \rangle^\beta \mu^{-1/2})}ds\\
		& \le C_{ns} t_* \|f-\tilde{f}\|_{L^{\infty}_{[t_0,t_0+t_*]}L^{\infty}_{x,v}(\tilde{w})}.
	\end{align*}
	Thus, we obtain
	\begin{align*}
		\|f-\tilde{f}\|_{L^{\infty}_{[t_0,t_0+t_*]}L^{\infty}_{x,v}(\tilde{w})} \le \left[(1-C_s)+C_{ns}t_*+6C_\infty C_Q\eta_0^{1+p} \right]\|f-\tilde{f}\|_{L^{\infty}_{[t_0,t_0+t_*]}L^{\infty}_{x,v}(\tilde{w})}.
	\end{align*}
	We choose $\eta_0>0$ such that $(1-\frac{3}{4}C_s)+6C_\infty C_Q\eta_0^{1+p} < 1$, and choose $t_*>0$ such that $C_{ns}t_* \le C_s/4$.
	This implies that $f\equiv \tilde{f}$ on $[t_0,t_0+t_*]$. Starting at $t_0=0$, we can iterate the process on $[0,t_*], [t_*,2t_*], \cdots.$ since $t_*$ is independent of $t_0$. Hence, we can conclude that $f\equiv \tilde{f}$ on $\mathbb{R}^+$.

\end{proof}

\bigskip

\begin{proposition} \label{positivityprop} 
	Let $w(v)=(1+|v|^2)^{\beta/2}e^{\alpha|v|^{\kappa}}$ with $\beta \ge 0 ,\alpha >0 ,0<\kappa<2$. Let $f_0$ in $L^\infty_{x,v}(w)$ with $P_{G_L}(f_0)=0$ and $\|f_0\|_{L^\infty_{x,v}(w)} \le \eta_0^{1+p}$ for some $0<p<1$, where $\eta_0$ is chosen so that proposition \ref{T1} and proposition \ref{UT} hold	and assume f is a unique solution of \eqref{BE-f} with initial data $f_0$. If $F_0 = \mu + f_0 \ge 0$, then $F=\mu+f \ge 0$.
\end{proposition}
\begin{proof}
	We decompose the collision operator $Q$ into two parts:
	\begin{align*}
		Q(F,F) = Q^{+}(F,F) - Q^{-}(F,F),
	\end{align*}
	where
	\begin{align*}
		&Q^{+}(F,F)(v) =\int_{\mathbb{R}^3 \times \mathbb{S}^2} B(v-v_*,\sigma)F(v')F(v_*')d\sigma dv_*,\\
		&Q^{-}(F,F)(v) =\int_{\mathbb{R}^3 \times \mathbb{S}^2} B(v-v_*,\sigma)F(v)F(v_*)d\sigma dv_* = F(v) R(F)(v).
	\end{align*}
	Now, we construct an iterative system for $F^{(n)}$ as follows:
	\begin{align} \label{PIS}
		\partial_t F^{(n+1)}+ v\cdot \nabla_x F^{(n+1)} + R(F^{(n)})F^{(n+1)} = Q^{+}(F^{(n)},F^{(n)}),
	\end{align}
	with $F^{(0)} = \mu$, $F^{(n+1)}|_{t=0} = \mu + f_0$. 	We define $f^{(n)} = F^{(n)}-\mu$. Then we can derive from \eqref{PIS} that
	\begin{align*}
		\partial_tf^{(n+1)} + v\cdot \nabla_xf^{(n+1)} + \nu(v) f^{(n+1)} = K(f^{(n)}) + Q^{+}(f^{(n)},f^{(n)})-R(f^{(n)})f^{(n+1)}.
	\end{align*}
	By a similar argument to proposition \ref{UT}, we obtain $f^{(n)}$ is a Cauchy sequence in $L^\infty_{x,v}(w)$, in local time $[0,T_0]$, where $T_0$ depends on $\|f_0\|_{L^\infty_{x,v}(w)}$, and so converges to $f$ in $L^\infty_{x,v}(w)$.\\
	We use an induction to show $F^{(n)} \ge 0$ for all $n\in \mathbb{N}$ over $[0,T_0]$. Assume that $F^{(n)} \ge 0$. For $(n+1)$th iteration, the Duhamel formulation gives
	\begin{align*}
		F^{(n+1)}(t,x,v) = e^{-\int_0^t R(F^{(n)})(\tau)d\tau}F_0(x-vt,v) + \int_0^t e^{-\int_s^t R(F^{(n)})(\tau)d\tau}Q^{+}(F^{(n)},F^{(n)})(s)ds \ge 0
	\end{align*}
	since $F_0 \ge 0$ and $Q^{+}(F^{(n)},F^{(n)}) \ge 0$. By induction, $F^{(n)} \ge 0$ for all $n\in \mathbb{N}$ over $[0,T_0]$, and thus $F \ge 0$ over $[0,T_0]$. Hence, $F \ge 0$ for all $t\ge 0$ by repeating $[0, T_0]$, $[T_0, 2T_0]$, ..., $[kT_0, (k + 1)T_0]$.
\end{proof}

%
%
%

\subsection{Useful Lemmas for Large amplitude problem} \label{largeamplitudescheme}
The proof of the following lemma is motivated by \cite{Cao2022JFA}. The following lemma provides a velocity decay estimate for a weighted collision integral involving a stretched-exponential weight. The prescribed range of $\varepsilon$ ensures decay as $|v|\to\infty$.

\begin{lemma} \label{K2estimate} 	Let $0\le \gamma \le 1$, $\beta\ge 0$, $\alpha>0$, and $1<\kappa<2$. For any $0\le \varepsilon < -\gamma+\frac{\kappa(\gamma+3)}{4}$, there is a constant $C_{\beta,\alpha,\kappa,\gamma,\varepsilon}$, depending on $\gamma, \kappa, \beta, \alpha$, and $\varepsilon$, so that
	\begin{align*}
		\int_{\R^3}\int_{\S^2} |v-v_*|^\gamma \frac{\langle  v \rangle^\beta e^{\alpha |v|^\kappa}}{\langle v'\rangle^\beta e^{\alpha |v'|^\kappa}} e^{-\frac{|v_*'|^2}{2}} |v-v_*|^\varepsilon d\sigma d v_* \le C_{\beta,\alpha,\kappa,\gamma,\varepsilon} \langle v \rangle^{\gamma-\frac{\kappa(\gamma+3)}{4}+\varepsilon}.
	\end{align*}
\end{lemma}
\begin{proof}
	We set
	\begin{align*}
		I:= \int_{\R^3}\int_{\S^2} |v-v_*|^\gamma \frac{\langle  v \rangle^\beta e^{\alpha |v|^\kappa}}{\langle v'\rangle^\beta e^{\alpha |v'|^\kappa}} e^{-\frac{|v_*'|^2}{2}} |v-v_*|^\varepsilon d\sigma d v_*.
	\end{align*}
	Since $\gamma-1\le 0$, we have
	\begin{align*}
		I&\le 4\int_{\R^3} \frac{1}{|v-v'|} \frac{\langle  v \rangle^\beta e^{\alpha |v|^\kappa}}{\langle v'\rangle^\beta e^{\alpha |v'|^\kappa}} \int_{\{w\, :\, w\cdot (v'-v)=0\}} \left(|v-v'|^2+|w|^2 \right)^\frac{\gamma-1+\varepsilon}{2} e^{-\frac{|v+w|^2}{2}}  dw d v'\\
		& \le 4 \int_{\R^3} \frac{1}{|v-v'|^\frac{3-\gamma}{2}}\frac{\langle  v \rangle^\beta e^{\alpha |v|^\kappa}}{\langle v'\rangle^\beta e^{\alpha |v'|^\kappa}} \int_{\{w\, :\, w\cdot (v'-v)=0\}} |w|^{\frac{\gamma-1}{2}+\varepsilon} e^{-\frac{|v+w|^2}{2}} dwdv'\\
		& \quad +4 \int_{\R^3} \frac{1}{|v-v'|^{\frac{3-\gamma}{2}-\varepsilon}}\frac{\langle  v \rangle^\beta e^{\alpha |v|^\kappa}}{\langle v'\rangle^\beta e^{\alpha |v'|^\kappa}} \int_{\{w\, :\, w\cdot (v'-v)=0\}} |w|^{\frac{\gamma-1}{2}} e^{-\frac{|v+w|^2}{2}} dwdv'\\
		&=:I_1+I_2,
	\end{align*}
	where we have used $(|v-v'|^2+|w|^2)^\frac{\varepsilon}{2} \le |v-v'|^\varepsilon+|w|^\varepsilon$ in the last inequality.\\
	\indent First of all, let us estimate the term $I_1$.
	Applying the following decomposition:
	\begin{align*}
		&v = v_\perp + v_\shortparallel , \quad v_\perp = \frac{v \cdot (v'-v)}{|v-v'|^2}(v'-v), \quad |v|^2= |v_\perp|^2+ |v_\shortparallel |^2, \\
		&v_\perp \perp v_\shortparallel ,\quad |v+w|^2 = |w+v_\shortparallel|^2 + |v_\perp|^2,
	\end{align*}
	we deduce that
	\begin{align*}
		I_1 \le 4 \int_{\R^3} \frac{1}{|v-v'|^\frac{3-\gamma}{2}}\frac{\langle  v \rangle^\beta e^{\alpha |v|^\kappa}}{\langle v'\rangle^\beta e^{\alpha |v'|^\kappa}}e^{-\frac{|v_\perp|^2}{2}} \int_{\{w\, :\, w\cdot (v'-v)=0\}} |w|^{\frac{\gamma-1}{2}+\varepsilon} e^{-\frac{|v_\shortparallel+w|^2}{2}} dwdv'.
	\end{align*}
	Since $\frac{\gamma-1}{2}+\varepsilon>-2$, it holds that
	\begin{align*}
		\int_{\{w\, :\, w\cdot (v'-v)=0\}} |w|^{\frac{\gamma-1}{2}+\varepsilon} e^{-\frac{|v_\shortparallel+w|^2}{2}} dw = \int_{\{w\, :\, w\cdot (v'-v)=0\}} |w-v_\shortparallel|^{\frac{\gamma-1}{2}+\varepsilon} e^{-\frac{|w|^2}{2}} dw \le C_{\gamma,\varepsilon} \langle v_\shortparallel \rangle^{\frac{\gamma-1}{2}+\varepsilon},
	\end{align*}
	which implies that
	\begin{align} \label{K2termI1}
		I_1 \le C_{\gamma,\varepsilon} \int_{\R^3} \frac{1}{|v-v'|^\frac{3-\gamma}{2}}\frac{\langle  v \rangle^\beta e^{\alpha |v|^\kappa}}{\langle v'\rangle^\beta e^{\alpha |v'|^\kappa}}e^{-\frac{|v_\perp|^2}{2}}\langle v_\shortparallel \rangle^{\frac{\gamma-1}{2}+\varepsilon}dv'.
	\end{align}
	If $|v|\le 1$, \eqref{K2termI1} becomes
	\begin{align*}
		I_1 \le C_{\alpha,\kappa,\beta,\gamma,\varepsilon} \int_{\R^3} \frac{|v-v'|^{\frac{\gamma-3}{2}}}{\langle v'\rangle^\beta e^{\alpha |v'|^\kappa}}dv' \le C_{\alpha,\kappa,\beta,\gamma,\varepsilon} \le C_{\alpha,\kappa,\beta,\gamma,\varepsilon} \langle v \rangle^{-1}.
	\end{align*}
	We consider the case $|v|\ge 1$.\\
	\newline
	$\mathbf{Case\ 1:\ } |v_\perp|\ge \frac{|v|}{2}$.\\
	From \eqref{K2termI1}, we can easily compute
	\begin{align*}
		I_1 & \le C_{\gamma,\varepsilon} \int_{\R^3} |v-v'|^\frac{\gamma-3}{2}\frac{\langle  v \rangle^\beta e^{\alpha |v|^\kappa}}{\langle v'\rangle^\beta e^{\alpha |v'|^\kappa}}e^{-\frac{|v|^2}{16}}dv'\\
		& \le C_{\alpha,\kappa,\beta,\gamma,\varepsilon} e^{-\frac{|v|^2}{32}} \int_{\R^3} \frac{|v-v'|^\frac{\gamma-3}{2}}{\langle v'\rangle^\beta e^{\alpha |v'|^\kappa}}dv\\
		& \le C_{\alpha,\kappa,\beta,\gamma,\varepsilon} e^{-\frac{|v|^2}{32}}.
	\end{align*}
	\newline
	$\mathbf{Case\ 2:\ } |v_\perp|< \frac{|v|}{2}$.\\
	In this case, we have $|v_\shortparallel| \ge \frac{|v|}{2}$, and it follows that $\langle v_\shortparallel \rangle^{\frac{\gamma-1}{2}} \le C \langle v \rangle^{\frac{\gamma-1}{2}}$. Then \eqref{K2termI1} is bounded by
	\begin{align*}
		I_1 \le C_{\gamma,\varepsilon} \langle v \rangle^\frac{\gamma-1}{2}  \int_{\R^3} |v-v'|^\frac{\gamma-3}{2}\frac{\langle  v \rangle^\beta e^{\alpha |v|^\kappa}}{\langle v'\rangle^\beta e^{\alpha |v'|^\kappa}}e^{-\frac{|v_\perp|^2}{2}} \langle v_\shortparallel \rangle^\varepsilon  dv'.
	\end{align*}
	We make a change of variables 
	\begin{align*}
		&r = |v-v'|, \quad v \cdot(v'-v) = r|v|\cos \theta, \quad dv'= r^2\sin\theta dr d\theta d\phi,\\
		&|v_\perp|^2 = |v|^2 \cos^2\theta, \quad |v_\shortparallel|^2 = |v|^2 \sin^2\theta
	\end{align*}
	to obtain
	\begin{align*}
		I_1 &\le C_{\gamma,\varepsilon} \langle v \rangle^\frac{\gamma-1}{2} \int_0^\infty \int_0^\pi  r^{\frac{\gamma+1}{2}} \frac{(1+|v|^2)^{\beta/2}}{\left( 1+(|v|^2+r^2+2r|v|\cos \theta)  \right)^{\beta/2}}\\
		& \qquad \times e^{\alpha|v|^\kappa - \alpha (|v|^2+r^2+2r|v|\cos \theta)^{\kappa/2}} e^{-\frac{|v|^2\cos^2\theta}{2}} \sin\theta (1+|v|^2\sin^2\theta)^{\varepsilon/2}d\theta dr\\
		& \le C_{\gamma,\varepsilon} \langle v \rangle^{\frac{\gamma-1}{2}+\varepsilon} \int_0^\infty \int_0^\pi  r^{\frac{\gamma+1}{2}} \frac{(1+|v|^2)^{\beta/2}}{\left( 1+(|v|^2+r^2+2r|v|\cos \theta)  \right)^{\beta/2}}\\
		& \qquad \times e^{\alpha|v|^\kappa - \alpha (|v|^2+r^2+2r|v|\cos \theta)^{\kappa/2}} e^{-\frac{|v|^2\cos^2\theta}{2}} \sin\theta d\theta dr,
	\end{align*}
	where we have used $|v'|^2 = |v-v'|^2+|v|^2 +2v\cdot(v'-v)$.\\
	To control the term $\frac{(1+|v|^2)^{\beta/2}}{\left( 1+(|v|^2+r^2+2r|v|\cos \theta)  \right)^{\beta/2}}$, we decompose the range of $\theta$ into $\{\frac{\pi}{4} \le \theta \le \frac{3\pi}{4} \}$ and $\{ 0\le \theta < \frac{\pi}{4}\ \text{or}\ \frac{3\pi}{4}<\theta \le \pi\}$. For the case $\{\frac{\pi}{4} \le \theta \le \frac{3\pi}{4} \}$, we have $\cos\theta \ge - \frac{1}{\sqrt{2}}$, and implies that
	\begin{align*}
		|v|^2+r^2+2r|v|\cos \theta &\ge |v|^2+r^2-\sqrt{2}r|v|\\
		& =\left(r-\frac{|v|}{\sqrt{2}}\right)^2 + \frac{|v|^2}{2}.
	\end{align*}
	From the above inequality, we can derive
	\begin{align} \label{ineq3}
		\frac{(1+|v|^2)^{\beta/2}}{\left( 1+(|v|^2+r^2+2r|v|\cos \theta)  \right)^{\beta/2}} \le \frac{(1+|v|^2)^{\beta/2}}{\left( 1+\frac{|v|^2}{2} \right)^{\beta/2}} \le C_\beta.
	\end{align}
	For the case $\{ 0\le \theta < \frac{\pi}{4}\ \text{or}\ \frac{3\pi}{4}<\theta \le \pi\}$, we have $\cos^2\theta \ge 1/2$. Then we obtain
	\begin{align} \label{ineq4}
		\frac{(1+|v|^2)^{\beta/2}}{\left( 1+(|v|^2+r^2+2r|v|\cos \theta)  \right)^{\beta/2}} e^{-\frac{|v|^2\cos^2\theta}{2}} &\le \frac{(1+|v|^2)^{\beta/2}}{\left( 1+(|v|^2+r^2+2r|v|\cos \theta)  \right)^{\beta/2}} e^{-\frac{|v|^2}{8}} e^{-\frac{|v|^2\cos^2\theta}{4}}\nonumber\\
		& \le C_{\beta} e^{-\frac{|v|^2\cos^2\theta}{4}}.
	\end{align}
	From \eqref{ineq3} and \eqref{ineq4}, we deduce
	\begin{align} \label{ineq5}
		I_1 &\le C_{\beta,\gamma,\varepsilon} \langle v \rangle^{\frac{\gamma-1}{2}+\varepsilon} \int_0^\infty \int_0^\pi  r^{\frac{\gamma+1}{2}} e^{\alpha|v|^\kappa - \alpha (|v|^2+r^2+2r|v|\cos \theta)^{\kappa/2}} e^{-\frac{|v|^2\cos^2\theta}{4}} \sin\theta d\theta dr.
	\end{align}
	Since $1<\kappa<2$, there exists a constant $C_{\alpha,\kappa}>0$ such that
	\begin{align*}
		\frac{|x|^2}{8}+C_{\alpha,\kappa} \ge 4\alpha|x|^\kappa \quad \text{for all}\ x\in \R,
	\end{align*}
	and it follows that
	\begin{align} \label{ineq1}
		\frac{|v|^2 \cos^2 \theta}{8}+C_{\alpha,\kappa} \ge 4\alpha(|v|^2\cos^2 \theta)^{\kappa/2}.
	\end{align}
	We also have
	\begin{align}  \label{ineq2} 
		\alpha (|v|^2+r^2+2r|v|\cos \theta)^{\kappa/2} + 4\alpha (|v|^2 \cos^2\theta)^{\kappa/2}& \ge \alpha (|v|^2+r^2 +2r|v|\cos \theta + 4|v|^2 \cos^2\theta)^{\kappa/2} \nonumber\\
		& \ge \alpha \left( |v|^2 +\frac{r^2}{2} \right)^{\kappa/2}.
	\end{align}
	By \eqref{ineq1} and \eqref{ineq2}, the inequality \eqref{ineq5} becomes
	\begin{align} \label{ineq6}
		I_1 \le C_{\alpha,\kappa,\beta,\gamma,\varepsilon} \langle v \rangle^{\frac{\gamma-1}{2}+\varepsilon} \int_0^\infty \int_0^\pi  r^{\frac{\gamma+1}{2}} e^{\alpha|v|^\kappa - \alpha \left(|v|^2+\frac{r^2}{2}\right)^{\kappa/2}} e^{-\frac{|v|^2\cos^2\theta}{8}} \sin\theta d\theta dr.
	\end{align}
	Making a change of variables
	\begin{align*}
		r=|v|x, \quad y =|v|\cos \theta, \quad dy =|v|\sin \theta d\theta,
	\end{align*}
	the inequality \eqref{ineq6} implies
	\begin{align} \label{ineq100}
		I_1 &\le C_{\alpha,\kappa,\beta,\gamma,\varepsilon} \langle v \rangle^{\gamma+\varepsilon} \int_0^\infty  x^{\frac{\gamma+1}{2}} e^{\alpha|v|^\kappa - \alpha |v|^\kappa \left(1+\frac{x^2}{2}\right)^{\kappa/2}}dx \int_{-|v|}^{|v|} e^{-\frac{y^2}{8}} dy\nonumber\\
		& \le C_{\alpha,\kappa,\beta,\gamma,\varepsilon} \langle v \rangle^{\gamma+\varepsilon} \int_0^\infty  x^{\frac{\gamma+1}{2}} e^{\alpha|v|^\kappa - \alpha |v|^\kappa \left(1+\frac{x^2}{2}\right)^{\kappa/2}}dx.
	\end{align}
	We split it into two cases $\{x\ge 2\}$ and $\{0\le x< 2\}$. For $\{x\ge 2\}$, it holds that $\left(1+\frac{x^2}{2} \right)^{\kappa/2}-1 \ge C_\kappa |x|^\kappa $ for some $C_\kappa>0$. Then we have
	\begin{align} \label{ineq7}
		\int_2^\infty  x^{\frac{\gamma+1}{2}} e^{\alpha|v|^\kappa - \alpha |v|^\kappa \left(1+\frac{x^2}{2}\right)^{\kappa/2}}dx &\le \int_2^\infty  x^{\frac{\gamma+1}{2}} e^{-C_\kappa \alpha|v|^\kappa |x|^\kappa}dx\nonumber\\
		& \le e^{-\frac{\alpha}{2}C_\kappa|v|^\kappa} \int_0^\infty x^{\frac{\gamma+1}{2}} e^{-\frac{\alpha}{2}C_\kappa|x|^\kappa}dx\nonumber\\
		&\le C_{\alpha,\kappa,\gamma}e^{-\frac{\alpha}{2}C_\kappa|v|^\kappa}.
	\end{align}
	For $\{0\le x< 2\}$, it holds that $\left(1+\frac{x^2}{2} \right)^{\kappa/2}-1 \ge C_\kappa |x|^2$ for some $C_\kappa>0$. We make a change of variable $z=|v|^{\kappa/2}x$ to derive
	\begin{align} \label{ineq8}
		\int_0^2  x^{\frac{\gamma+1}{2}} e^{\alpha|v|^\kappa - \alpha |v|^\kappa \left(1+\frac{x^2}{2}\right)^{\kappa/2}}dx &\le \int_0^\infty  x^{\frac{\gamma+1}{2}} e^{-C_\kappa \alpha|v|^\kappa x^2}dx\nonumber\\
		& \le |v|^{-\frac{\kappa(\gamma+3)}{4}} \int_0^\infty  z^{\frac{\gamma+1}{2}} e^{-C_\kappa \alpha z^2}dz\nonumber\\
		& \le C_{\alpha,\kappa,\gamma} |v|^{-\frac{\kappa(\gamma+3)}{4}}.
	\end{align}
	Combining two cases \eqref{ineq7} and \eqref{ineq8}, it follows from \eqref{ineq100} that 
	\begin{align} \label{ineq9}
		I_1 \le C_{\alpha,\kappa,\beta,\gamma,\varepsilon}  \langle v  \rangle^{\gamma-\frac{\kappa(\gamma+3)}{4}+\varepsilon}.
	\end{align}
	By a similar argument to estimate $I_1$, we deduce
	\begin{align} \label{ineq10}
		I_2 \le C_{\alpha,\kappa,\beta,\gamma,\varepsilon}  \langle v  \rangle^{\gamma-\frac{\kappa(\gamma+3)}{4}+\left(1-\frac{\kappa}{2}\right)\varepsilon}.
	\end{align}
	Gathering \eqref{ineq9} and \eqref{ineq10}, we conclude that
	\begin{align*}
		I \le C_{\alpha,\kappa,\beta,\gamma,\varepsilon}  \langle v  \rangle^{\gamma-\frac{\kappa(\gamma+3)}{4}+\varepsilon}.
	\end{align*}	
\end{proof}

\bigskip

The following lemma provides a weighted estimate for the kernel $k$, yielding the velocity decay for a suitable range of $\varepsilon$.

\begin{lemma} \label{KESTIMATE}
	Let $0\le \gamma \le 1$, $\beta\ge0$, $\alpha>0$, and $1<\kappa<2$. For any $0\le \varepsilon < -\gamma+\frac{\kappa(\gamma+3)}{4}$, there is a constant $C_{\beta,\alpha,\kappa,\gamma,\varepsilon}$, depending on $\gamma, \kappa, \beta, \alpha$, and $\varepsilon$, so that
	\begin{align*}
		\int_{\R^3}|k(v,v_*)| \frac{\langle  v \rangle^\beta e^{\alpha |v|^\kappa}}{\langle v_*\rangle^\beta e^{\alpha |v_*|^\kappa}} |v-v_*|^\varepsilon  d v_* \le C_{\beta,\alpha,\kappa,\gamma,\varepsilon} \langle v \rangle^{\gamma-\frac{\kappa(\gamma+3)}{4}+\varepsilon}.
	\end{align*}
\end{lemma}

Note that if $1<\kappa<2$ and $0\le \gamma \le 1$, it holds that $\gamma-\frac{\kappa(\gamma+3)}{4}<0$. For simplicity, we take $\varepsilon = \frac{1}{2}\left( -\gamma+\frac{\kappa(\gamma+3)}{4}\right)$.

\begin{proof}
We have
\begin{align*}
	\int_{\R^3}|k(v,v_*)| \frac{\langle  v \rangle^\beta e^{\alpha |v|^\kappa}}{\langle v_*\rangle^\beta e^{\alpha |v_*|^\kappa}} |v-v_*|^\varepsilon  d v_* 
	&\le C\int_{\R^3}\int_{\S^2} |v-v_*|^\gamma \frac{\langle  v \rangle^\beta e^{\alpha |v|^\kappa}}{\langle v'\rangle^\beta e^{\alpha |v'|^\kappa}} e^{-\frac{|v_*'|^2}{2}} |v-v_*|^\varepsilon d\sigma d v_* \\
	& \quad + C\int_{\R^3}\int_{\S^2} |v-v_*|^\gamma \frac{\langle  v \rangle^\beta e^{\alpha |v|^\kappa}}{\langle v_*\rangle^\beta e^{\alpha |v_*|^\kappa}} e^{-\frac{|v|^2}{2}} |v-v_*|^\varepsilon d\sigma d v_*\\
	& =:I_1+I_2.
\end{align*}
	 By Lemma \ref{K2estimate}, we obtain
	 \begin{align*}
	 	I_1 \le C_{\beta,\alpha,\kappa,\gamma,\varepsilon} \langle v \rangle^{\gamma-\frac{\kappa(\gamma+3)}{4}+\varepsilon}.
	 \end{align*}
	 Let us estimate $I_2$.
	 \begin{align*}
	 	I_2 &\le C \langle v \rangle^\beta e^{\alpha |v|^\kappa} e^{-\frac{|v|^2}{2}} \int_{\R^3} |v-v_*|^{\gamma+\varepsilon}  \frac{1}{\langle v_*\rangle^\beta e^{\alpha |v_*|^\kappa}}dv_*\\
	 	& \le C_{\beta,\alpha,\kappa,\gamma,\varepsilon}  \langle v \rangle^{\beta+\gamma+\varepsilon} e^{\alpha |v|^\kappa} e^{-\frac{|v|^2}{2}}\\
	 	& \le C_{\beta,\alpha,\kappa,\gamma,\varepsilon}  e^{-\frac{|v|^2}{4}}.
	 \end{align*}
	 Hence we complete the proof of this lemma.
\end{proof}

\bigskip

The following lemma provides an $L^2$ estimate for the weighted kernel $k$.
\begin{lemma} \label{KL2estimate}
	Let $0\le \gamma \le 1$, $\beta\ge0$, $\alpha>0$, and $1<\kappa<2$. There is a constant $C_{\beta,\alpha,\kappa,\gamma}$, depending on $\gamma, \kappa, \beta, \alpha$, so that
	\begin{align*}
		\int_{\R^3} |k(v,v_*)|^2 \left(\frac{\langle  v \rangle^\beta e^{\alpha |v|^\kappa}}{\langle v_*\rangle^\beta e^{\alpha |v_*|^\kappa}}\right)^2   d v_* \le C_{\alpha,\kappa,\beta,\gamma}  \langle v  \rangle^{2\gamma-2-\frac{\kappa}{2}}.
	\end{align*}
\end{lemma}
\begin{proof}
	As in Lemma \ref{K2estimate} and \ref{KESTIMATE}, the Carleman representation yields
	\begin{align*} 
		&\int_{\R^3} |k(v,v_*)|^2 \left(\frac{\langle  v \rangle^\beta e^{\alpha |v|^\kappa}}{\langle v_*\rangle^\beta e^{\alpha |v_*|^\kappa}}\right)^2 dv_*\nonumber\\
		& \le C \int_{\R^3}\left| \frac{1}{|v-v'|} \frac{\langle  v \rangle^\beta e^{\alpha |v|^\kappa}}{\langle v'\rangle^\beta e^{\alpha |v'|^\kappa}}\int_{\{w : w\cdot (v'-v)=0\}}(|v-v'|^2 + |w|^2)^{\frac{\gamma-1}{2}}e^{-\frac{|v+w|^2}{2}}dw \right|^2 dv'\nonumber\\
		&\quad + C\int_{\R^3}|v-v_*|^{2\gamma} \frac{\langle  v \rangle^{2\beta} e^{2\alpha |v|^\kappa}}{\langle v_*\rangle^{2\beta} e^{2\alpha |v_*|^\kappa}} e^{-|v|^2}dv_*\\
		&=:I_1+I_2.
	\end{align*}
	Since $\gamma \le 1$, we have
	\begin{align*}
		\int_{\{w : w\cdot (v'-v)=0\}}(|v-v'|^2 + |w|^2)^{\frac{\gamma-1}{2}}e^{-\frac{|v+w|^2}{2}}dw  \le \int_{\R^2} |w|^{\gamma-1}e^{-\frac{|v+w|^2}{2}}dw  
	\end{align*}
	For $I_1$, applying the following decomposition:
	\begin{align*}
		&v = v_\perp + v_\shortparallel , \quad v_\perp = \frac{v \cdot (v'-v)}{|v-v'|^2}(v'-v), \quad |v|^2= |v_\perp|^2+ |v_\shortparallel |^2, \\
		&v_\perp \perp v_\shortparallel, \quad |v+w|^2 = |w+v_\shortparallel|^2 + |v_\perp|^2,
	\end{align*}
	we deduce that
	\begin{align*}
		I_1 \le C\int_{\R^3} \left|\frac{1}{|v-v'|} \frac{\langle  v \rangle^{\beta} e^{\alpha |v|^\kappa}}{\langle v'\rangle^{\beta} e^{\alpha |v'|^\kappa}}e^{-\frac{|v_\perp|^2}{2}} \int_{\{w\, :\, w\cdot (v'-v)=0\}} |w|^{\gamma-1} e^{-\frac{|v_\shortparallel+w|^2}{2}} dw\right|^2dv'.
	\end{align*}
	Since $\gamma-1>-2$, it holds that
	\begin{align*}
		\int_{\{w\, :\, w\cdot (v'-v)=0\}} |w|^{\gamma-1} e^{-|v_\shortparallel+w|^2} dw = \int_{\{w\, :\, w\cdot (v'-v)=0\}} |w-v_\shortparallel|^{\gamma-1} e^{-|w|^2} dw \le C_{\gamma} \langle v_\shortparallel \rangle^{\gamma-1},
	\end{align*}
	which implies that
	\begin{align} \label{K2termI11}
		I_1 \le C \int_{\R^3} \frac{1}{|v-v'|^2} \frac{\langle  v \rangle^{2\beta} e^{2\alpha |v|^\kappa}}{\langle v'\rangle^{2\beta} e^{2\alpha |v'|^\kappa}}e^{-|v_\perp|^2}\langle v_\shortparallel \rangle^{2\gamma-2}dv'.
	\end{align}
	If $|v|\le 1$, \eqref{K2termI11} becomes
	\begin{align*}
		I_1 \le C_{\alpha,\kappa,\beta,\gamma} \int_{\R^3} \frac{|v-v'|^{-2}}{\langle v'\rangle^{2\beta} e^{2\alpha |v'|^\kappa}}dv' \le C_{\alpha,\kappa,\beta,\gamma} .
	\end{align*}
	We consider the case $|v|\ge 1$.\\
	\newline
	$\mathbf{Case\ 1:\ } |v_\perp|\ge \frac{|v|}{2}$.\\
	From \eqref{K2termI11}, we can easily compute
	\begin{align*}
		I_1 & \le C_{\gamma} \int_{\R^3} |v-v'|^{-2}\frac{\langle  v \rangle^{2\beta} e^{2\alpha |v|^\kappa}}{\langle v'\rangle^{2\beta} e^{2\alpha |v'|^\kappa}}e^{-\frac{|v|^2}{4}}dv'\\
		& \le C_{\alpha,\kappa,\beta,\gamma} e^{-\frac{|v|^2}{8}} \int_{\R^3} \frac{|v-v'|^{-2}}{\langle v'\rangle^{2\beta} e^{2\alpha |v'|^\kappa}}dv\\
		& \le C_{\alpha,\kappa,\beta,\gamma} e^{-\frac{|v|^2}{8}}.
	\end{align*}
	\newline
	$\mathbf{Case\ 2:\ } |v_\perp|< \frac{|v|}{2}$.\\
	In this case, we have $|v_\shortparallel| \ge \frac{|v|}{2}$, and it follows that $\langle v_\shortparallel \rangle^{2\gamma-2} \le C \langle v \rangle^{2\gamma-2}$. Then \eqref{K2termI11} is bounded by
	\begin{align*}
		I_1 \le C \langle v \rangle^{2\gamma-2}  \int_{\R^3} |v-v'|^{-2}\frac{\langle  v \rangle^{2\beta} e^{2\alpha |v|^\kappa}}{\langle v'\rangle^{2\beta} e^{2\alpha |v'|^\kappa}}e^{-|v_\perp|^2}   dv'.
	\end{align*}
	We make a change of variables 
	\begin{align*}
		&r = |v-v'|, \quad v \cdot(v'-v) = r|v|\cos \theta, \quad dv'= r^2\sin\theta dr d\theta d\phi,\\
		&|v_\perp|^2 = |v|^2 \cos^2\theta, \quad |v_\shortparallel|^2 = |v|^2 \sin^2\theta
	\end{align*}
	to obtain
	\begin{align*}
		I_1 &\le C_{\gamma} \langle v \rangle^{2\gamma-2} \int_0^\infty \int_0^\pi   \frac{(1+|v|^2)^{\beta}}{\left( 1+(|v|^2+r^2+2r|v|\cos \theta)  \right)^{\beta}}\\
		& \qquad \times e^{2\alpha|v|^\kappa - 2\alpha (|v|^2+r^2+2r|v|\cos \theta)^{\kappa/2}} e^{-|v|^2\cos^2\theta} \sin\theta d\theta dr,
	\end{align*}
	where we have used $|v'|^2 = |v-v'|^2+|v|^2 +2v\cdot(v'-v)$.\\
	To control the term $\frac{(1+|v|^2)^{\beta}}{\left( 1+(|v|^2+r^2+2r|v|\cos \theta)  \right)^{\beta}}$, we decompose the range of $\theta$ into $\{\frac{\pi}{4} \le \theta \le \frac{3\pi}{4} \}$ and $\{ 0\le \theta < \frac{\pi}{4}\ \text{or}\ \frac{3\pi}{4}<\theta \le \pi\}$. For the case $\{\frac{\pi}{4} \le \theta \le \frac{3\pi}{4} \}$, we have $\cos\theta \ge - \frac{1}{\sqrt{2}}$, and implies that
	\begin{align*}
		|v|^2+r^2+2r|v|\cos \theta &\ge |v|^2+r^2-\sqrt{2}r|v|\\
		& =\left(r-\frac{|v|}{\sqrt{2}}\right)^2 + \frac{|v|^2}{2}.
	\end{align*}
	From the above inequality, we can derive
	\begin{align} \label{ineq31}
		\frac{(1+|v|^2)^{\beta}}{\left( 1+(|v|^2+r^2+2r|v|\cos \theta)  \right)^{\beta}} \le \frac{(1+|v|^2)^{\beta}}{\left( 1+\frac{|v|^2}{2} \right)^{\beta}} \le C_\beta.
	\end{align}
	For the case $\{ 0\le \theta < \frac{\pi}{4}\ \text{or}\ \frac{3\pi}{4}<\theta \le \pi\}$, we have $\cos^2\theta \ge 1/2$. Then we obtain
	\begin{align} \label{ineq41}
		\frac{(1+|v|^2)^{\beta}}{\left( 1+(|v|^2+r^2+2r|v|\cos \theta)  \right)^{\beta}} e^{-|v|^2\cos^2\theta} &\le \frac{(1+|v|^2)^{\beta}}{\left( 1+(|v|^2+r^2+2r|v|\cos \theta)  \right)^{\beta}} e^{-\frac{|v|^2}{4}} e^{-\frac{|v|^2\cos^2\theta}{2}}\nonumber\\
		& \le C_{\beta} e^{-\frac{|v|^2\cos^2\theta}{2}}.
	\end{align}
	From \eqref{ineq31} and \eqref{ineq41}, we deduce
	\begin{align} \label{ineq51}
		I_1 &\le C_{\beta,\gamma} \langle v \rangle^{2\gamma-2} \int_0^\infty \int_0^\pi  e^{2\alpha|v|^\kappa - 2\alpha (|v|^2+r^2+2r|v|\cos \theta)^{\kappa/2}} e^{-\frac{|v|^2\cos^2\theta}{2}} \sin\theta d\theta dr.
	\end{align}
	Since $1<\kappa<2$, there exists a constant $C_{\alpha,\kappa}>0$ such that
	\begin{align*}
		\frac{|x|^2}{4}+C_{\alpha,\kappa} \ge 8\alpha|x|^\kappa \quad \text{for all}\ x\in \R,
	\end{align*}
	and it follows that
	\begin{align} \label{ineq11}
		\frac{|v|^2 \cos^2 \theta}{4}+C_{\alpha,\kappa} \ge 8\alpha(|v|^2\cos^2 \theta)^{\kappa/2}.
	\end{align}
	We also have
	\begin{align}  \label{ineq21} 
		2\alpha (|v|^2+r^2+2r|v|\cos \theta)^{\kappa/2} + 8\alpha (|v|^2 \cos^2\theta)^{\kappa/2}& \ge 2\alpha (|v|^2+r^2 +2r|v|\cos \theta + 4|v|^2 \cos^2\theta)^{\kappa/2} \nonumber\\
		& \ge 2\alpha \left( |v|^2 +\frac{r^2}{2} \right)^{\kappa/2}.
	\end{align}
	By \eqref{ineq11} and \eqref{ineq21}, the inequality \eqref{ineq51} becomes
	\begin{align} \label{ineq61}
		I_1 \le C_{\alpha,\kappa,\beta,\gamma} \langle v \rangle^{2\gamma-2} \int_0^\infty \int_0^\pi   e^{2\alpha|v|^\kappa - 2\alpha \left(|v|^2+\frac{r^2}{2}\right)^{\kappa/2}} e^{-\frac{|v|^2\cos^2\theta}{4}} \sin\theta d\theta dr.
	\end{align}
	Making a change of variables
	\begin{align*}
		r=|v|x, \quad y =|v|\cos \theta, \quad dy =|v|\sin \theta d\theta,
	\end{align*}
	the inequality \eqref{ineq61} implies
	\begin{align} \label{ineq91}
		I_1 &\le C_{\alpha,\kappa,\beta,\gamma} \langle v \rangle^{2\gamma-2} \int_0^\infty  e^{2\alpha|v|^\kappa - 2\alpha |v|^\kappa \left(1+\frac{x^2}{2}\right)^{\kappa/2}}dx \int_{-|v|}^{|v|} e^{-\frac{y^2}{4}} dy\nonumber\\
		& \le C_{\alpha,\kappa,\beta,\gamma} \langle v \rangle^{2\gamma-2} \int_0^\infty   e^{2\alpha|v|^\kappa - 2\alpha |v|^\kappa \left(1+\frac{x^2}{2}\right)^{\kappa/2}}dx.
	\end{align}
	We split it into two cases $\{x\ge 2\}$ and $\{0\le x< 2\}$. For $\{x\ge 2\}$, it holds that $\left(1+\frac{x^2}{2} \right)^{\kappa/2}-1 \ge C_\kappa |x|^\kappa $ for some $C_\kappa>0$. Then we have
	\begin{align} \label{ineq71}
		\int_2^\infty  e^{2\alpha|v|^\kappa - 2\alpha |v|^\kappa \left(1+\frac{x^2}{2}\right)^{\kappa/2}}dx &\le \int_2^\infty   e^{-C_\kappa \alpha|v|^\kappa |x|^\kappa}dx\nonumber\\
		& \le e^{-\frac{\alpha}{2}C_\kappa|v|^\kappa} \int_0^\infty  e^{-\frac{\alpha}{2}C_\kappa|x|^\kappa}dx\nonumber\\
		&\le C_{\alpha,\kappa,\gamma}e^{-\frac{\alpha}{2}C_\kappa|v|^\kappa}.
	\end{align}
	For $\{0\le x< 2\}$, it holds that $\left(1+\frac{x^2}{2} \right)^{\kappa/2}-1 \ge C_\kappa |x|^2$ for some $C_\kappa>0$. We make a change of variable $z=|v|^{\kappa/2}x$ to derive
	\begin{align} \label{ineq81}
		\int_0^2   e^{2\alpha|v|^\kappa - 2\alpha |v|^\kappa \left(1+\frac{x^2}{2}\right)^{\kappa/2}}dx &\le \int_0^\infty  e^{-C_\kappa \alpha|v|^\kappa x^2}dx\nonumber\\
		& \le |v|^{-\frac{\kappa}{2}} \int_0^\infty   e^{-C_\kappa \alpha z^2}dz\nonumber\\
		& \le C_{\alpha,\kappa,\gamma} |v|^{-\frac{\kappa}{2}} .
	\end{align}
	Combining two cases \eqref{ineq71} and \eqref{ineq81}, it follows from \eqref{ineq91} that 
	
	\begin{align*}
		I_1 \le C_{\alpha,\kappa,\beta,\gamma}  \langle v  \rangle^{2\gamma-2-\frac{\kappa}{2}} .
	\end{align*}
	Let us estimate $I_2$.
	\begin{align*}
		I_2 &\le C\langle  v \rangle^{2\beta+\gamma} e^{2\alpha |v|^\kappa}e^{-|v|^2}\int_{\R^3}\frac{\langle v_*\rangle^{2\gamma}}{\langle v_*\rangle^{2\beta} e^{2\alpha |v_*|^\kappa}} dv_*\le C_{\alpha,\beta,\kappa,\gamma}e^{-\frac{|v|^2}{2}}.
	\end{align*}
	This completes the proof of this lemma.
\end{proof}

\bigskip

\begin{lemma} \label{pointwiseQ-est}
	Suppose $\beta \ge 0$, $\alpha>0$ and $0<\kappa<2$. Then, for any $\beta$, we have
	\begin{align*}
		|wQ^-(f,f)(v)| \le C_{\alpha,\kappa,\gamma} \langle v\rangle^\gamma \|wf\|_{L^\infty_v}\left(\int_{\R^3} e^{-\frac{1}{2}\alpha|v_*|^\kappa} |wf(v_*)|^2dv_*\right)^{1/2}.
	\end{align*}
\end{lemma}
\begin{proof}
	We use the Cauchy-Schwarz inequality to obtain
	\begin{align*}
		|wQ^-(f,f)(v)| &\le \langle v\rangle^\beta e^{\alpha|v|^\kappa} \left|\int_{\R^3}\int_{\S^2}b(\cos \theta) |v-v_*|^\gamma  f(v)f(v_*) d\sigma dv_* \right|\\
		& \le C_\gamma \langle v\rangle^\gamma \|wf\|_{L^\infty_v}\left(\int_{\R^3}\langle v_*\rangle^\gamma \frac{1}{w(v_*)} |wf(v_*)|dv_*\right)\\
		& \le C_{\alpha,\kappa,\gamma} \langle v\rangle^\gamma \|wf\|_{L^\infty_v}\left(\int_{\R^3} e^{-\frac{1}{2}\alpha|v_*|^\kappa} |wf(v_*)|^2dv_*\right)^{1/2}.
	\end{align*}
\end{proof}

\bigskip

\subsection{Nonlinear $L^\infty_{x,v}$ estimate} 
Setting $h=w_{\beta,\alpha,\kappa}f=\langle v\rangle^\beta e^{\alpha|v|^\kappa}f$, the equation \eqref{FPBE} can be written as
\begin{align}\label{WFPBE2}
	\partial_t h + v \cdot \nabla_x h + \nu(v) h=K_wh + wQ\left(f ,f \right).
\end{align}
To show the global-in-time existence to the equation \eqref{FPBE}, we need to deal with the nonlinear $L^\infty_{x,v}$ estimate for the equation \eqref{WFPBE2}. This subsection is devoted to deriving the nonlinear $L^\infty_{x,v}$ estimate.

Recall the definition for relative entropy $H(F \ \vert \ \mu)$ from \eqref{def rel entropy} :
\begin{align*}
	H(F \ \vert \ \mu) = \int_{\T^3 \times \mathbb{R}^3} \left(\frac{F}{\mu}\log\frac{F}{\mu}-\frac{F}{\mu}+1\right) \mu dxdv.
\end{align*}
\bigskip
Since $F$ satisfies entropy inequality \cite{DL2} and physical conservations,
\Be \label{relativedecrease}
H(F_t \ \vert \ \mu) \le H(F_0 \ \vert \ \mu) \quad \text{for all } t\ge 0.
\Ee

With help of \eqref{relativedecrease}, the following lemma states initial relative entropy provides the bound for the $L^1$ and $L^2$ norms of $F-\mu$ on the different domains.
\begin{lemma} \label{L1L2cont} \cite{Guo2010Bounded}
	Assume $F$ satisfies the Boltzmann equation \eqref{BE}. We have
	\begin{align*}
		\int_{\T^3\times \mathbb{R}^3} \frac{1}{4\mu}|F-\mu|^2 \mathbf{1}_{\{|F-\mu|\le \mu\}} dxdv + \int_{\T^3 \times \mathbb{R}^3} \frac{1}{4}|F-\mu| \mathbf{1}_{\{|F-\mu| > \mu\}} dxdv \le H(F_0 \ \vert \ \mu)
	\end{align*}
	for all $t \ge 0$.
\end{lemma}

\begin{lemma} \label{LALinftyestimate} 
	Let $\beta\ge0$, $\alpha>0$, and $1<\kappa<2$. Under the a priori assumption, there exists a generic constant $C_3 \ge 1$, depending on $\alpha,\beta,\kappa$, and $\gamma$, such that
	\begin{align} \label{LALinftyestimate1}
		\|h(t)\|_{L^\infty_{x,v}} &\le C_3\left(1+ \bar{M} \right)\|h_0\|_{L^\infty_{x,v}}e^{-\frac{\nu_0}{2} t} \nonumber\\
		&\quad +C_3\left(N^{\max\{-1,\frac{1}{2}\left(\gamma-\frac{\kappa(\gamma+3)}{4}\right)\}}+ \varepsilon C_{\alpha,\beta,\kappa,\gamma,N}\right) \left[\bar{M}+\bar{M}^2+\bar{M}^3\right]\nonumber\\
		& \quad +  C_{\alpha,\beta,\kappa,\gamma,N,\varepsilon}\left[1 + \bar{M}^{1/2}+\bar{M}+\bar{M}^{3/2}+\bar{M}^2+\bar{M}^{5/2}\right]H(F_0 \ \vert \ \mu)^{1/2}
	\end{align}	
	for all $t \ge 0$, where $\varepsilon>0$ can be arbitrarily small and $N>0$ can be arbitrary large.
\end{lemma}
\begin{proof}
	Let $(t,x,v) \in [0,\infty) \times \T^3 \times \R^3$. We apply the Duhamel principle to the weighted Boltzmann equation \eqref{WFPBE2}:
	\begin{equation} \label{LA10}
		\begin{aligned}
			h(t,x,v) &= e^{-\nu(v)t} h_0(x-vt,v)\\
			&\quad + \int_0^t e^{-\nu(v)(t-s)}K_wh(s,x-v(t-s),v) ds\\
			& \quad + \int_0^t e^{-\nu(v)(t-s)}wQ(f,f)(s,x-v(t-s),v) ds\\
			& =: I_1+I_2+I_3.
		\end{aligned}
	\end{equation}
	First of all, we can easily compute that
	\begin{align} \label{LA62}
		|I_1| \le e^{-\nu_0t}\|h_0\|_{L^\infty_{x,v}}.
	\end{align}
	Now, let us estimate the term $I_2$ in \eqref{LA10}:\
	\begin{align}\label{LA11}
		|I_2|&\le \int_0^{t} e^{-\nu(v)(t-s)}\int_{\R^3}|k_w(v,v')||h(s,x-v(t-s),v')|dv'ds\nonumber\\
		& \le \mathbf{1}_{|v|\ge N}\int_0^{t} e^{-\nu(v)(t-s)}\int_{\R^3}|k_w(v,v')||h(s,x-v(t-s),v')|dv'ds\nonumber\\
		& \quad +  \mathbf{1}_{|v|\le  N}\int_0^{t} e^{-\nu(v)(t-s)}\int_{\R^3}|k_w(v,v')||h(s,x-v(t-s),v')|dv'ds\nonumber\\
		& \le C_{\alpha,\beta,\kappa,\gamma}N^{\gamma-\frac{\kappa(\gamma+3)}{4}} \bar{M}\nonumber\\
		& \quad + \mathbf{1}_{|v|\le  N}\int_0^{t} e^{-\nu(v)(t-s)}\int_{\R^3}|k_w(v,v')||h(s,x-v(t-s),v')|dv'ds,
	\end{align}
	where $\gamma-\frac{\kappa(\gamma+3)}{4}<0$, and we apply the Duhamel principle to $|h(s,x-v(t-s),v')|$ in the last term of the above:
	\begin{align*}
		&\mathbf{1}_{|v|\le  N}\int_0^{t} e^{-\nu(v)(t-s)}\int_{\R^3}|k_w(v,v')|e^{-\nu(v')s}|h_0(x-v(t-s)-v's,v')|dv'ds\\
		& +\mathbf{1}_{|v|\le  N}\int_0^{t} e^{-\nu(v)(t-s)}\int_{\R^3}|k_w(v,v')|\int_0^s e^{-\nu(v')(s-s')}|K_wh(s',x-v(t-s)-v'(s-s'),v')|ds'dv'ds \\
		& + \mathbf{1}_{|v|\le  N}\int_0^{t} e^{-\nu(v)(t-s)}\int_{\R^3}|k_w(v,v')|\int_0^s e^{-\nu(v')(s-s')}|wQ(f,f)(s',x-v(t-s)-v'(s-s'),v')|ds'dv'ds\\
		& =: I_{21}+I_{22}+I_{23}.
	\end{align*}
	For $I_{21}$, we use Lemma \ref{KESTIMATE} to obtain
	\begin{align} \label{LA12}
		|I_{21}| &\le \mathbf{1}_{|v|\le  N}\int_0^{t} e^{-\nu(v)(t-s)}\int_{\R^3}|k_w(v,v')|e^{-\nu(v')s}|h_0(x-v(t-s)-v's,v')|dv'ds\nonumber\\
		& \le C_{\alpha,\beta,\kappa,\gamma} \langle v \rangle^{\gamma-\frac{\kappa(\gamma+3)}{4}} t e^{-\nu_0 t}\|h_0\|_{L^\infty_{x,v}}\nonumber\\
		& \le C_{\alpha,\beta,\kappa,\gamma}e^{-\frac{\nu_0}{2} t}\|h_0\|_{L^\infty_{x,v}},
	\end{align}
	where $\gamma-\frac{\kappa(\gamma+3)}{4}<0$.\\
	For $I_{22}$, the term can be controlled by
	\begin{align*}
		&\mathbf{1}_{|v|\le  N}\int_0^{t}e^{-\nu(v)(t-s)}\int_{\R^3}|k_w(v,v')|\\
		& \quad  \times\int_0^s e^{-\nu(v')(s-s')}\int_{\R^3}|k_w(v',v'')||h(s',x-v(t-s)-v'(s-s'),v'')|dv''ds'dv'ds,
	\end{align*}
	and we divide this term into two cases.\\
	\newline
	$\mathbf{Case\ 1:}$ $|v|\le N, |v'|\ge 2N$ or $|v'|\le 2N, |v''|\ge 3N$.\\
	Either of the following is valid:
	\begin{equation} \label{class0}
		\begin{aligned}
			&|k_w(v,v')| \lesssim N^{\frac{1}{2}\left(\gamma-\frac{\kappa(\gamma+3)}{4}\right)}|k_w(v,v')||v-v'|^{-\frac{1}{2}\left(\gamma-\frac{\kappa(\gamma+3)}{4}\right)}\\
			\text{or} \quad &|k_w(v',v'')| \lesssim N^{\frac{1}{2}\left(\gamma-\frac{\kappa(\gamma+3)}{4}\right)}|k_w(v',v'')||v'-v''|^{-\frac{1}{2}\left(\gamma-\frac{\kappa(\gamma+3)}{4}\right)},
		\end{aligned}
	\end{equation}
	From Lemma \ref{KESTIMATE}, either of the following is valid:
	\begin{equation} \label{class1}
		\begin{aligned}
			&\int_{|v'| \ge 2N}|k_w(v,v')||v-v'|^{-\frac{1}{2}\left(\gamma-\frac{\kappa(\gamma+3)}{4}\right)}dv' \le C_{\alpha,\beta,\kappa,\gamma}\langle v \rangle^{\frac{1}{2}\left(\gamma-\frac{\kappa(\gamma+3)}{4}\right)}\\
			\text{or} \quad & \int_{|v''|\ge 3N}|k_w(v',v'')||v'-v''|^{-\frac{1}{2}\left(\gamma-\frac{\kappa(\gamma+3)}{4}\right)}dv'' \le C_{\alpha,\beta,\kappa,\gamma}\langle v' \rangle^{\frac{1}{2}\left(\gamma-\frac{\kappa(\gamma+3)}{4}\right)}.
		\end{aligned}
	\end{equation}
	From \eqref{class0} and \eqref{class1}, we obtain
	\begin{align}\label{LA37}
		I_{22}^1 &\le \mathbf{1}_{|v|\le N}\int_{0}^t\int_{0}^s e^{-\nu_0(t-s')}   \int_{|v'| \ge 2N}\int_{\R^3} |k_w(v,v')||k_w(v',v'')|\nonumber \\
		& \qquad \times |h(s',x-v(t-s)-v'(s-s'),v'')|dv''dv'ds'ds\nonumber\\
		& \quad + \int_{0}^t\int_{0}^s e^{-\nu_0(t-s')}  \int_{|v'|\le 2N}\int_{|v''| \ge 3N} |k_w(v,v')||k_w(v',v'')| \nonumber\\
		& \qquad \quad  \times |h(s',x-v(t-s)-v'(s-s'),v'')|dv''dv'ds'ds\nonumber\\
		& \le C_{\alpha,\beta,\kappa,\gamma}N^{\frac{1}{2}\left(\gamma-\frac{\kappa(\gamma+3)}{4}\right)} \bar{M}.
	\end{align}
	\newline
	$\mathbf{Case\ 2:}$ $|v|\le N, |v'|\le 2N, |v''|\le 3N$.\\
	Using the Cauchy-Schwarz inequality and Lemma \ref{KL2estimate}, we deduce
	\begin{align*}
		I_{22}^2 
		& \le \int_0^t \int_{s-\varepsilon}^{s} e^{-\nu_0(t-s')}\mathbf{1}_{|v|\le N} \int_{|v'|\le 2N}\int_{|v''|\le 3N}|k_w(v,v')||k_w(v',v'')|\\
		& \qquad \times|h(s',x-v(t-s)-v'(s-s'),v'')|dv''dv'ds'ds\\
		& \quad + \int_{0}^{t}\int_{0}^{s-\varepsilon} e^{-\nu_0(t-s')}\mathbf{1}_{|v|\le N}\left(\int_{|v'|\le 2N}\int_{|v''|\le 3N}|k_w(v,v')|^2|k_w(v',v'')|^2dv''dv'\right)^{1/2} \\
		& \qquad \times \left(\int_{|v'|\le 2N}\int_{|v''|\le 3N}|h(s',x-v(t-s)-v'(s-s'),v'')|^2dv''dv'  \right)^{1/2}ds'ds\\
		& \le \varepsilon C_{\alpha,\beta,\kappa,\gamma} \bar{M}       \\
		& \quad + C_{\alpha,\beta,\kappa,\gamma,N} \int_{0}^{t}\int_{0}^{s-\varepsilon} e^{-\nu_0(t-s')}\mathbf{1}_{|v|\le N}\\
		& \qquad \times\left(\int_{|v'|\le 2N}\int_{|v''|\le 3N}|h(s',x-v(t-s)-v'(s-s'),v'')|^2dv''dv'  \right)^{1/2} ds'ds.
	\end{align*}
	Making a change of variables $v' \mapsto y=x-v(t-s)-v'(s-s')$ with $\left|\frac{dy}{dv'}\right| = (s-s')^3$, the last term is bounded by
	\begin{align*}
		C_{\alpha,\beta,\kappa,\gamma,N,\varepsilon}\int_{0}^{t}\int_{0}^{s} e^{-\nu_0(t-s')}\left(\int_{\T^3}\int_{|v''|\le 3N}|h(s',y,v'')|^2dv''dy  \right)^{1/2}ds'ds.
	\end{align*}
	Here, we use Lemma \ref{L1L2cont} and the Cauchy-Schwarz inequality to compute
	\begin{equation} \label{RTIRE}
		\begin{aligned}
			&\int_{\T^3}\int_{|v'|\le 2N}|h(s,y,v')|^2dv'dy\\
			& \le C_N \int_{\T^3}\int_{|v'|\le 2N}\frac{1}{4\mu(v')}|f(s,y,v')|^2 \mathbf{1}_{\{|f|\le \mu\}} dv'dy +C_N\int_{\T^3}\int_{|v'|\le 2N}|h(s,y,v')|^2 \mathbf{1}_{\{|f|> \mu\}}dv'dy\\
			&\le C_N H(F_0 \ \vert \ \mu) + C_N\|h(s)\|_{L^\infty_{x,v}}\left(\int_{\T^3}\int_{|v'|\le 2N}|h(s,y,v')| \mathbf{1}_{\{|f|> \mu\}}dv'dy\right) \\
			&\le C_N H(F_0 \ \vert \ \mu) + C_{N}\|h(s)\|_{L^\infty_{x,v}}\left(\int_{\T^3}\int_{|v'|\le 2N}|f(s,y,v')| \mathbf{1}_{\{|f|> \mu\}}dv'dy\right)\\
			&\le C_N H(F_0 \ \vert \ \mu) + C_{N}\|h(s)\|_{L^\infty_{x,v}}H(F_0 \ \vert \ \mu).
		\end{aligned}
	\end{equation}
	It follows from \eqref{RTIRE} that 
	\begin{align}\label{LA38}
		I_{22}^2  &\le \varepsilon C_{\alpha,\beta,\kappa,\gamma} \bar{M}+C_{\alpha,\beta,\kappa,\gamma,\varepsilon,N}\int_{0}^{t}\int_{0}^{s} e^{-\nu_0(t-s')} \left[H(F_0 \ \vert \ \mu)^{1/2} + \|h(s)\|_{L^\infty_{x,v}}^{1/2}H(F_0 \ \vert \ \mu)^{1/2}\right]ds'ds\nonumber\\
		&\le\varepsilon C_{\alpha,\beta,\kappa,\gamma} \bar{M}+  C_{\alpha,\beta,\kappa,\gamma,\varepsilon,N}\left[H(F_0 \ \vert \ \mu)^{1/2} + \bar{M}^{1/2}H(F_0 \ \vert \ \mu)^{1/2}\right].
	\end{align}
	\newline
	Gathering \eqref{LA37} and  \eqref{LA38}, we derive
	\begin{align} \label{LA39}
		I_{22} \le C_{\alpha,\beta,\kappa,\gamma}N^{\frac{1}{2}\left(\gamma-\frac{\kappa(\gamma+3)}{4}\right)} \bar{M} +\varepsilon C_{\alpha,\beta,\kappa,\gamma} \bar{M}+  C_{\alpha,\beta,\kappa,\gamma,\varepsilon,N}\left[H(F_0 \ \vert \ \mu)^{1/2} + \bar{M}^{1/2}H(F_0 \ \vert \ \mu)^{1/2}\right].
	\end{align}
	\newline
	For $I_{23}$, by Lemma \ref{pointwiseQ-est} and Lemma \ref{pointwiseQ+est}, it is bounded by
	\begin{align*}
		&\mathbf{1}_{|v|\le  N}C_{\alpha,\beta,\kappa,\gamma}\int_0^{t} e^{-\nu(v)(t-s)}\int_{\R^3}|k_w(v,v')|\int_0^s e^{-\nu(v')(s-s')}\langle v'\rangle^\gamma \|h(s')\|_{L^\infty_{x,v}}\\
		& \quad \times  \left(\int_{\R^3}\frac{1}{(1+|\eta|)^{6}}|h(s',x-v(t-s)-v'(s-s'),\eta)|^2 d\eta \right)^{1/2}ds'dv'ds\\
	\end{align*}
	and we divide this term into three cases.\\
	\newline
	$\mathbf{Case\ 1:}$ $|v|\le N, |v'|\ge 2N$.\\
	Using \eqref{class0} and \eqref{class1}, we deduce that
	\begin{align} \label{LA43}
		I_{23}^1 & \le C_{\alpha,\beta,\kappa,\gamma} \bar{M}^2 \int_{0}^{t} e^{-\nu_0(t-s)}\mathbf{1}_{|v|\le N} \int_{|v'|\ge 2N}|k_w(v,v')|\int_0^s e^{-\nu(v')(s-s')}\langle v'\rangle^\gamma dv'ds'ds\nonumber\\
		& \le  C_{\alpha,\beta,\kappa,\gamma}N^{\frac{1}{2}\left(\gamma-\frac{\kappa(\gamma+3)}{4}\right)}\bar{M}^2.
	\end{align}	
	\newline
	$\mathbf{Case\ 2:}$ $|v|\le N, |v'|\le 2N, |\eta| \ge N$.\\
	We apply Lemma \ref{KESTIMATE} to obtain 
	\begin{align}\label{LA44}
		I_{23}^2 &\le \frac{C_{\alpha,\beta,\kappa,\gamma}}{N} \bar{M}^2 \int_{0}^{t} e^{-\nu_0(t-s)}\mathbf{1}_{|v|\le N} \int_{|v'|\le 2N}|k_w(v,v')|\int_0^s e^{-\nu(v')(s-s')}\langle v'\rangle^\gamma dv'ds'ds\nonumber\\
		& \le \frac{C_{\alpha,\beta,\kappa,\gamma}}{N}\bar{M}^2.
	\end{align}
	\newline
	$\mathbf{Case\ 3:}$ $|v|\le N, |v'|\le 2N, |\eta| \le N$.\\
	Using the Cauchy-Schwarz inequality and Lemma \ref{KESTIMATE}, we get
	\begin{align*}
		I_{23}^3 
		&\le C_{\alpha,\beta,\kappa,\gamma}\int_0^{t} e^{-\nu(v)(t-s)}\int_{|v'|\le 2N}|k_w(v,v')|\int_{s-\varepsilon}^s e^{-\nu(v')(s-s')}\langle v'\rangle^\gamma \|h(s')\|_{L^\infty_{x,v}}\\
		& \qquad \times  \left(\int_{|\eta|\le N}\frac{1}{(1+|\eta|)^{6}}|h(s',x-v(t-s)-v'(s-s'),\eta)|^2 d\eta \right)^{1/2}ds'dv'ds\\
		& \quad +C_{\alpha,\beta,\kappa,\gamma}\int_0^{t} e^{-\nu(v)(t-s)}\int_{|v'|\le 2N}|k_w(v,v')|\int_0^{s-\varepsilon} e^{-\nu(v')(s-s')}\langle v'\rangle^\gamma \|h(s')\|_{L^\infty_{x,v}}\\
		& \qquad \times  \left(\int_{|\eta|\le N}\frac{1}{(1+|\eta|)^{6}}|h(s',x-v(t-s)-v'(s-s'),\eta)|^2 d\eta \right)^{1/2}ds'dv'ds\\
		& \le \varepsilon C_{\alpha,\beta,\kappa,\gamma,N} \bar{M}^2\\
		& \quad +C_{\alpha,\beta,\kappa,\gamma,N} \bar{M} \int_0^t \int_0^{s-\varepsilon} e^{-\nu_0(t-s')}    \int_{|v'|\le 2N}|k_w(v,v')|   \\
		& \qquad \times \left(\int_{|\eta|\le N}|h(s',x-v(t-s)-v'(s-s'),\eta)|^2 d\eta \right)^{1/2}dv'ds'ds\\
		& \le \varepsilon C_{\alpha,\beta,\kappa,\gamma,N} \bar{M}^2\\
		&\quad + C_{\alpha,\beta,\kappa,\gamma,N}\bar{M} \int_{0}^{t}\int_{0}^{s-\varepsilon} e^{-\nu_0(t-s')}\\
		& \qquad \times\left(\int_{|v'|\le 2N}\int_{|\eta|\le N}|h(s',x-v(t-s)-v'(s-s'),\eta)|^2 d\eta dv' \right)^{1/2} ds'ds.
	\end{align*}
	Making a change of variables $v' \mapsto y=x-v(t-s)-v'(s-s')$ with $\left|\frac{dy}{dv'}\right| = (s-s')^3$, the last term is bounded by
	\begin{align*}
		C_{\alpha,\beta,\kappa,\gamma,N,\varepsilon} \bar{M} \int_{0}^{t}\int_{0}^{s} e^{-\nu_0(t-s')}\left(\int_{\T^3}\int_{|\eta|\le N}|h(s',y,\eta)|^2d\eta dy  \right)^{1/2}ds'ds.
	\end{align*}
	It follows from \eqref{RTIRE} that 
	\begin{align}\label{LA45}
		I_{23}^3 \le  \varepsilon C_{\alpha,\beta,\kappa,\gamma,N} \bar{M}^2+C_{\alpha,\beta,\kappa,\gamma,N,\varepsilon} \left[\bar{M}H(F_0 \ \vert \ \mu)^{1/2}+\bar{M}^{3/2}H(F_0 \ \vert \ \mu)^{1/2}\right] .
	\end{align}
	Gathering \eqref{LA43}, \eqref{LA44}, and \eqref{LA45}, we derive
	\begin{align} \label{LA46}
		I_{23} &\le C_{\alpha,\beta,\kappa,\gamma}N^{\max\{-1,\frac{1}{2}\left(\gamma-\frac{\kappa(\gamma+3)}{4}\right)\}}\bar{M}^2+ \varepsilon C_{\alpha,\beta,\kappa,\gamma,N} \bar{M}^2\nonumber\\
		&\quad +C_{\alpha,\beta,\kappa,\gamma,N,\varepsilon} \left[\bar{M}H(F_0 \ \vert \ \mu)^{1/2}+\bar{M}^{3/2}H(F_0 \ \vert \ \mu)^{1/2}\right].
	\end{align}
	\newline
	Inserting \eqref{LA11}, \eqref{LA12}, \eqref{LA39}, and \eqref{LA46} into $I_{2}$, we derive
	\begin{align} \label{LA52}
		|I_2| &\le C_{\alpha,\beta,\kappa,\gamma}e^{-\frac{\nu_0}{2} t}\|h_0\|_{L^\infty_{x,v}}\nonumber\\
		&\quad +\left(C_{\alpha,\beta,\kappa,\gamma}N^{\max\{-1,\frac{1}{2}\left(\gamma-\frac{\kappa(\gamma+3)}{4}\right)\}}+ \varepsilon C_{\alpha,\beta,\kappa,\gamma,N}\right) \left[\bar{M}+\bar{M}^2\right]\nonumber\\
		& \quad +  C_{\alpha,\beta,\kappa,\gamma,N,\varepsilon}\left[H(F_0 \ \vert \ \mu)^{1/2} + \bar{M}^{1/2}H(F_0 \ \vert \ \mu)^{1/2}\right]\nonumber\\
		&\quad +C_{\alpha,\beta,\kappa,\gamma,N,\varepsilon} \left[\bar{M}H(F_0 \ \vert \ \mu)^{1/2}+\bar{M}^{3/2}H(F_0 \ \vert \ \mu)^{1/2}\right].
	\end{align}
	Next, let us estimate $I_{3}$. From Lemma \ref{pointwiseQ+est}, it is bounded by
	\begin{align} \label{LA53}
		|I_{3}| 
		&\le C_{\alpha,\beta,\kappa,\gamma}\int_0^t e^{-\nu(v)(t-s)} \langle v\rangle^\gamma \|h(s)\|_{L^\infty_{x,v}}\nonumber\\
		&\qquad \times    \left(\left(\int_{|\eta|\ge N}+\int_{|\eta|\le N}\right)\frac{1}{(1+|\eta|)^{6}}|h(s,x-v(t-s),\eta)|^2 d\eta \right)^{1/2} ds \nonumber\\
		& \le \frac{C_{\alpha,\beta,\kappa,\gamma}}{N}\bar{M}^2\nonumber\\
		& \quad + C_{\alpha,\beta,\kappa,\gamma} \bar{M}\int_{0}^{t} e^{-\nu(v)(t-s)} \langle  v\rangle^\gamma \left(\int_{|\eta|\le N}\frac{1}{(1+|\eta|)^{6}}|h(s,x-v(t-s),\eta)|^2 d\eta \right)^{1/2} ds.
	\end{align}
	We apply the Duhamel principle to bound the term $|h(s,x-v(t-s),\eta)|$ in the last term of the above
	\begin{align*}
		&\left(\int_{|\eta|\le N}\frac{1}{(1+|\eta|)^{6}}|h(s,x-v(t-s),\eta)|^2 d\eta \right)^{1/2} \\
		& \le \left(\int_{|\eta|\le N}\frac{1}{(1+|\eta|)^{6}}|e^{-\nu(\eta)s}h_0(x-v(t-s)-\eta s,\eta)|^2 d\eta \right)^{1/2}\\
		& \quad +\left(\int_{|\eta|\le N}\frac{1}{(1+|\eta|)^{6}}\left|\int_{0}^{s}e^{-\nu(\eta)(s-s')}K_wh(s',x-v(t-s)-\eta(s-s'),\eta)ds'\right|^2 d\eta \right)^{1/2}\\
		& \quad +\left(\int_{|\eta|\le N}\frac{1}{(1+|\eta|)^{6}}\left|\int_{0}^{s}e^{-\nu(\eta)(s-s')}wQ(f,f)(s',x-v(t-s)-\eta(s-s'),\eta)ds'\right|^2 d\eta \right)^{1/2}\\
		&=: \tilde{L}_1+\tilde{L}_2+\tilde{L}_3.
	\end{align*}
	We define
	\begin{align*}
		L_{i} = C_{\alpha,\beta,\kappa,\gamma} \bar{M}\int_{0}^{t} e^{-\nu(v)(t-s)} \langle v\rangle^\gamma \tilde{L}_ids.
	\end{align*}
	For $L_1$, we easily compute 
	\begin{align} \label{LA54}
		L_1 &\le C_{\alpha,\beta,\kappa,\gamma} \bar{M} \|h_0\|_{L^\infty_{x,v}} e^{-\frac{\nu_0}{2}t}\int_{0}^{t} e^{-\frac{\nu(v)}{2}(t-s)} \langle v\rangle^\gamma ds\nonumber\\
		& \le C_{\alpha,\beta,\kappa,\gamma} \bar{M} \|h_0\|_{L^\infty_{x,v}} e^{-\frac{\nu_0}{2}t}.
	\end{align}
	For $L_2$, we first need to estimate $\tilde{L}_2$. We use Lemma \ref{KESTIMATE}, \eqref{class0}, and \eqref{class1} to obtain
	\begin{align*}
		\tilde{L}_2 &\le C_{\alpha,\beta,\kappa,\gamma}\left(\int_{|\eta|\le N}\frac{1}{(1+|\eta|)^{6}}\int_{0}^{s}e^{-2\nu_0(s-s')}\left|K_wh(s',x-v(t-s)-\eta(s-s'),\eta)\right|^2ds' d\eta \right)^{\frac{1}{2}}\\
		& \le C_{\alpha,\beta,\kappa,\gamma} N^{\frac{1}{2}\left(\gamma-\frac{\kappa(\gamma+3)}{4}\right)} \bar{M}\\
		& \quad + C_{\alpha,\beta,\kappa,\gamma,N}\left(\left(\int_{s-\varepsilon^2}^{s}+\int_0^{s-\varepsilon^2}\right)e^{-\nu_0(s-s')}\int_{|\eta|\le N}\int_{|v'|\le 2N}|h(s',x-v(t-s)-\eta(s-s'),v')|^2dv' d\eta ds'  \right)^{1/2}\\
		& \le C_{\alpha,\beta,\kappa,\gamma} N^{\frac{1}{2}\left(\gamma-\frac{\kappa(\gamma+3)}{4}\right)} \bar{M}+\varepsilon C_{\alpha,\beta,\kappa,\gamma,N}\bar{M}+ C_{\alpha,\beta,\kappa,\gamma,N,\varepsilon} H(F_0 \ \vert \ \mu)^{1/2}+C_{\alpha,\beta,\kappa,\gamma,N,\varepsilon}\bar{M}^{1/2}H(F_0 \ \vert \ \mu)^{1/2},
	\end{align*}
	where we have used \eqref{RTIRE} in the last term.
	It follow that
	\begin{align} \label{LA57}
		L_2 &\le C_{\alpha,\beta,\kappa,\gamma} N^{\frac{1}{2}\left(\gamma-\frac{\kappa(\gamma+3)}{4}\right)}\bar{M}^2+\varepsilon C_{\alpha,\beta,\kappa,\gamma,N}\bar{M}^2+C_{\alpha,\beta,\kappa,\gamma,N,\varepsilon} \bar{M}H(F_0 \ \vert \ \mu)^{1/2}+ C_{\alpha,\beta,\kappa,\gamma,N,\varepsilon}\bar{M}^{3/2}H(F_0 \ \vert \ \mu)^{1/2}.
	\end{align}
	For $L_3$, we first need to estimate $\tilde{L}_3$. We use Lemma \ref{pointwiseQ-est} and Lemma \ref{pointwiseQ+est} to deduce
	\begin{align*}
		\tilde{L}_3
		& \le C_{\alpha,\beta,\kappa,\gamma}\bar{M}\Bigg(\int_{|\eta|\le N} \frac{1}{(1+|\eta|)^{6}} \Bigg| \int_{0}^{s}e^{-\nu(\eta)(s-s')}\langle \eta\rangle^\gamma\\
		& \qquad \times \Bigg( \left(\int_{|\eta'|\ge 2N}+\int_{|\eta'|\le 2N}\right) \frac{1}{(1+|\eta'|)^{6}}|h(s',x-v(t-s)-\eta(s-s'),\eta')|^2d\eta'\Bigg)^{1/2}ds' \Bigg|^2 d\eta \Bigg)^{1/2}  \\
		& \le \frac{C_{\alpha,\beta,\kappa,\gamma}}{N} \bar{M}^2 \\
		& \quad +C_{\alpha,\beta,\kappa,\gamma,N}\bar{M}\Bigg(\int_{|\eta|\le N} \frac{1}{(1+|\eta|)^{6}} \Bigg| \left(\int_{s-\varepsilon}^{s}+\int_0^{s-\varepsilon}\right)e^{-\nu_0(s-s')}\\
		& \qquad \quad \times \Bigg( \int_{|\eta'|\le 2N} \frac{1}{(1+|\eta'|)^{6}}|h(s',x-v(t-s)-\eta(s-s'),\eta')|^2d\eta'\Bigg)^{1/2}ds' \Bigg|^2 d\eta \Bigg)^{1/2} \\
		& \le \frac{C_{\alpha,\beta,\kappa,\gamma}}{N}\bar{M}^2+\varepsilon C_{\alpha,\beta,\kappa,\gamma,N}\bar{M}^2 \\
		& \quad +  C_{\alpha,\beta,\kappa,\gamma,N}\bar{M}\Bigg(\int_0^{s-\varepsilon}e^{-\nu_0(s-s')}\int_{|\eta|\le N}\int_{|\eta'|\le 2N} |h(s',x-v(t-s)-\eta(s-s'),\eta')|^2d\eta'd\eta ds'\Bigg)^{1/2}\\
		& \le \frac{C_{\alpha,\beta,\kappa,\gamma}}{N}\bar{M}^2+\varepsilon C_{\alpha,\beta,\kappa,\gamma,N}\bar{M}^2 +C_{\alpha,\beta,\kappa,\gamma,N,\varepsilon}\left[\bar{M}H(F_0 \ \vert \ \mu)^{1/2} + \bar{M}^{3/2}H(F_0 \ \vert \ \mu)^{1/2}\right],
	\end{align*}
	where we have used \eqref{RTIRE} and it follows that
	\begin{align} \label{LA58}
		L_3 \le \frac{C_{\alpha,\beta,\kappa,\gamma}}{N}\bar{M}^3+\varepsilon C_{\alpha,\beta,\kappa,\gamma,N}\bar{M}^3 +C_{\alpha,\beta,\kappa,\gamma,N,\varepsilon}\left[\bar{M}^2H(F_0 \ \vert \ \mu)^{1/2} + \bar{M}^{5/2}H(F_0 \ \vert \ \mu)^{1/2}\right].
	\end{align}
	Inserting \eqref{LA54}, \eqref{LA57}, and \eqref{LA58} into \eqref{LA53}, we derive
	\begin{align} \label{LA59}
		I_{3} 
		& \le C_{\alpha,\beta,\kappa,\gamma} \bar{M} \|h_0\|_{L^\infty_{x,v}} e^{-\frac{\nu_0}{2}t}+ \left(C_{\alpha,\beta,\kappa,\gamma}N^{\max\{-1,\frac{1}{2}\left(\gamma-\frac{\kappa(\gamma+3)}{4}\right)\}}+\varepsilon C_{\alpha,\beta,\kappa,\gamma,N} \right) \left[\bar{M}^2 +\bar{M}^3 \right] \nonumber\\
		& \quad +C_{\alpha,\beta,\kappa,\gamma,N,\varepsilon}\left[ \bar{M}H(F_0 \ \vert \ \mu)^{1/2}+ \bar{M}^{3/2}H(F_0 \ \vert \ \mu)^{1/2}+\bar{M}^2H(F_0 \ \vert \ \mu)^{1/2} + \bar{M}^{5/2}H(F_0 \ \vert \ \mu)^{1/2}\right].
	\end{align}
	Plugging \eqref{LA62}, \eqref{LA52}, and \eqref{LA59} into \eqref{LA10}, we can derive
	\begin{align*}
		\|h(t)\|_{L^\infty_{x,v}} &\le \left(1+C_{\alpha,\beta,\kappa,\gamma}+C_{\alpha,\beta,\kappa,\gamma} \bar{M} \right)\|h_0\|_{L^\infty_{x,v}}e^{-\frac{\nu_0}{2} t} \nonumber\\
		&\quad +\left(C_{\alpha,\beta,\kappa,\gamma}N^{\max\{-1,\frac{1}{2}\left(\gamma-\frac{\kappa(\gamma+3)}{4}\right)\}}+ \varepsilon C_{\alpha,\beta,\kappa,\gamma,N}\right) \left[\bar{M}+\bar{M}^2+\bar{M}^3\right]\nonumber\\
		& \quad +  C_{\alpha,\beta,\kappa,\gamma,N,\varepsilon}\left[1 + \bar{M}^{1/2}+\bar{M}+\bar{M}^{3/2}+\bar{M}^2+\bar{M}^{5/2}\right]H(F_0 \ \vert \ \mu)^{1/2}.
	\end{align*}
\end{proof}

\bigskip

\subsection{Proof of Theorem \ref{Mainresultlargeamp}}
\begin{proof}[Proof of Theorem \ref{Mainresultlargeamp}]
Recall that we choose $0<\eta_0\leq1$ from Proposition~\ref{Mainresultsmall} with its exponent $0<p<1$.
	Take
	\begin{align} \label{LA2}
		T_0 := \frac{2}{\nu_0}\left[\log\left(4C_3 \left(1+\bar{M}\right)M_0\right)+|\log\eta_0^{1+p}|\right].
	\end{align}
	By \eqref{LALinftyestimate1}, we get
	\begin{align} \label{LA3}
		\|h(t)\|_{L^\infty_{x,v}} &\le C_3\left(1+ \bar{M} \right)\|h_0\|_{L^\infty_{x,v}}e^{-\frac{\nu_0}{2} t} +\tilde{E}
	\end{align}
	for all $t \ge 0$, where
	\begin{align*}
		\tilde{E}:&= C_3\left(N^{\max\{-1,\frac{1}{2}\left(\gamma-\frac{\kappa(\gamma+3)}{4}\right)\}}+ \varepsilon C_{\alpha,\beta,\kappa,\gamma,N}\right) \left[\bar{M}+\bar{M}^2+\bar{M}^3\right]\nonumber\\
		& \quad +  C_{\alpha,\beta,\kappa,\gamma,N,\varepsilon}\left[1 + \bar{M}^{1/2}+\bar{M}+\bar{M}^{3/2}+\bar{M}^2+\bar{M}^{5/2}\right]H(F_0 \ \vert \ \mu)^{1/2}.
	\end{align*}
	We first choose $N>0$ sufficiently large, then choose $\varepsilon >0$ small enough, and assume $H(F_0 \ \vert \ \mu) \le \varepsilon_0(\eta_0, \bar{M})$ with small enough $\varepsilon_0(\eta_0, \bar{M})$ so that
	\begin{align*}
		\tilde{E} \le \frac{\eta_0^{1+p}}{4},
	\end{align*}
	and it follows that
	\begin{align} \label{LA6}
		\|h(t)\|_{L^\infty_{x,v}} \le C_3\left(1+ \bar{M} \right)\|h_0\|_{L^\infty_{x,v}}e^{-\frac{\nu_0}{2} t} + \frac{\eta_0^{1+p}}{4}
	\end{align}
	for all $t \ge 0$.	\\	
	Assume $H(F_0 \ \vert \ \mu) \le \varepsilon_0$. From \eqref{LA6}, we get
	\begin{align*}
		\|h(T_0)\|_{L^\infty_{x,v}} \le C_3\left(1+ \bar{M} \right)\|h_0\|_{L^\infty_{x,v}}e^{-\frac{\nu_0}{2} T_0} + \frac{\eta_0^{1+p}}{4}  \le \frac{\eta_0^{1+p}}{2}.
	\end{align*}
	It remains to show the exponential decay of the Boltzmann solution $f(t)$ in the weighted $L^\infty_{x,v}$ space. Assume $H(F_0 \ \vert \ \mu) \le \varepsilon_0$. By proposition \ref{Mainresultsmall}, for all $t \ge T_0$,
	\begin{align} \label{LA9}
		\|h(t)\|_{L^\infty_{x,v}} \le C_\infty (\|h(T_0)\|_{L^\infty_{x,v}}+\eta_0^{1+p}) e^{-\lambda_\infty(t-T_0)} \le 2 C_\infty \eta_0^{1+p} e^{-\lambda_\infty(t-T_0)}.
	\end{align}
	Since $\lambda_\infty < \frac{\nu_0}{2}$, it follows from \eqref{LA6} and \eqref{LA9} that 
	\begin{align*}
		\|h(t)\|_{L^\infty_{x,v}} \le \max\left\{2,8C_\infty\right\}C_3(1+\bar{M}) M_0e^{-\lambda_\infty t}
	\end{align*}
	for all $t \ge 0$.
	Therefore, we complete the proof of this theorem.
\end{proof}

\subsection{Large amplitude problem with polynomial weight}
We state a Theorem for large amplitude problem in polynomial weight.
\begin{theorem} \label{Mainresultlargeamp-p}
	Let \(w(v)=(1+|v|^2)^{\beta/2}\) with $\beta>\beta_\infty:=\max\left\{1+\gamma+\frac{4C_b}{c_b},\left(\frac{4C_1}{\nu_0}\right)^{\frac{4}{\gamma+3}}\right\}$, where
	\begin{align} \label{def C1-p}
		C_1=8\sqrt{2}C_b \left(2+\frac{4}{\gamma+3}+\frac{1}{e}\right) \left(\frac{32\pi}{3+\gamma} e^{-3/4} +8\pi \right),
	\end{align}
    and $\nu_0$ is defined in \eqref{nuzerodef}.\\
	Assume that $F_0(x,v) = \mu(v) +  f_0(x,v) \ge 0$ satisfying 
	\begin{align*}
		\int_{\T^3 \times \R^3}f_0 
		\begin{pmatrix}
			1\\
			v\\
			|v|^2
		\end{pmatrix}
		dxdv =0.
	\end{align*}
	For any $M_0 \ge 1$, there exists a constant $\varepsilon_0>0$, depending  only on $\eta_0$, $M_0$, and $\bar{M}$, such that if
	\begin{align*}
		\|wf_0\|_{L^\infty_{x,v}}\le M_0, \quad \sup_{t\geq 0}\|wf(t)\|_{L^\infty_{x,v}}\leq \bar{M}, \quad H(F_0 \ \vert \ \mu) \le  \varepsilon_0,
	\end{align*}
	then the solution $F(t,x,v) =\mu(v) +  f(t,x,v) \ge 0 $ of the Boltzmann equation \eqref{BE} with initial data $F_0$ satisfies 
	\begin{align*}
		\|wf(t)\|_{L^\infty_{x,v}} \le C_0(1+\bar{M}) M_0e^{-\lambda_\infty t}
	\end{align*}
	for all $t \ge 0$, where $C_0\ge 1$ is a generic constant.
\end{theorem}
The overall approach for the proof of Theorem \ref{Mainresultlargeamp-p} is similar to the proof of Theorem \ref{Mainresultlargeamp}. However, unlike in the case of stretched exponential weights, significant differences arise in the estimates of the linearized Boltzmann collision operator and its $v$-growth. We defer the detailed proof to the Appendix.  

\bigskip

\section{Proof of Main Theorems}

The uniqueness assertion in Theorem~\ref{thm:cond} was established at the end of Section~2. The same argument also applies to Corollary~\ref{coro:cond}, Theorem~\ref{thm:cond-p}, and Corollary~\ref{coro:cond-p}. Therefore, in this section, we focus exclusively on proving exponential convergence to equilibrium.

\subsection{Proof of Theorem \ref{thm:cond} and Corollary \ref{coro:cond}}

\begin{proof} [Proof of Theorem \ref{thm:cond}]
	Since $F$ satisfies \eqref{apriori}, $F$ is the unique global dissipative mild solution by Lemma \ref{lem_unique}. Moreover, $F\in \mathcal{D}_{DL}(F_0)$ since $\emptyset \neq \mathcal{D}_{DL}\subset \mathcal{U}_{\infty}(F_0)$. So, $F$ satisfies entropy inequality \eqref{entropy ineq}. Obviously, $F$ satisfies \eqref{macro_conserv} either by \eqref{apriori}. And, from Proposition \ref{prop:glb}, we assume $F_0$ has a uniform Gaussian lower bound without loss of generality by choosing $T_{glb} < T_{loc}$ as our new initial time, where $T_{loc}$ is local existence time from Proposition \ref{prop_loc} with $\|\avv^{\b}e^{\a|v|^\k}F(T_{loc})\|_{L^\infty_{x,v}} \leq 2\|\avv^{\b}e^{\a|v|^\k}F_0\|_{L^\infty_{x,v}}$.  
	
 Once we have \eqref{apriori} and initial Gaussian lower bound, we obtain a coercivity estimate for $L(F)$ as in Proposition \ref{prop: macro lb}. Then applying to Proposition \ref{prop uniform bound}, we get a uniform weighted $L^\infty$ bound. This also implies that $F$ satisfies entropy equality either. Let $\beta_1\ge0$ be the exponent specified in
Proposition \ref{prop uniform bound} and set 
	\begin{align*}
		M_0 = \max\left\{1, \sup_{t\ge 0}\|\avv^{\b_1}e^{\a|v|^\k}(F(t)-\mu)\|_{L^\infty_{x,v}}\right\}<\infty.
	\end{align*}
	Then we choose $\varepsilon_0=\varepsilon_0(M_0)$ as in Theorem \ref{Mainresultlargeamp} where we choose $M_0=\bar{M}=C_{\infty,F}$ in Theorem \ref{Mainresultlargeamp}. By entropy decay in Proposition \ref{prop:rel-entro decay}, there exists a finite time $T_{ent}(\varepsilon_0) < \infty$ such that
	\[
		H(F \ \vert \ \mu)(T_{ent}) \leq \varepsilon_0
	\]
	where $F = \mu + f$. We apply Theorem \ref{Mainresultlargeamp} with the weight exponent $\beta_1$, taking $T_{ent}$ as the new initial time. Therefore, we obtain decay estimate of Theorem \ref{thm:cond}. The decay rate $\lambda_\infty$ is constructive and independent of $F$.


\end{proof}

\begin{proof} [Proof of Corollary \ref{coro:cond}]
If we assume \eqref{cond Coro}, the initial $F$ satisfies \eqref{macro_conserv} and entropy equality obviously. Lemma \ref{lem_unique} also hold and we can apply Proposition \ref{prop:rel-entro decay} directly. Similar as in the proof of Theorem \ref{thm:cond}, there exists a finite time $T_{ent} < \infty$ from which large amplitude problem regime begins as in the proof of Theorem \ref{thm:cond}. Then we obtain \eqref{decay F coro} by Theorem \ref{Mainresultlargeamp}. In particular, we do not lose velocity decay rate $\avv^{-\b_1}e^{-\a|v|^\k}$ of a priori assumption \eqref{cond Coro}.
\end{proof}

\subsection{Proof of Theorem \ref{thm:cond-p} and Corollary \ref{coro:cond-p}}

\begin{proof} [Proof of Theorem \ref{thm:cond-p}]
	Since $F$ satisfies \eqref{apriori-poly} for $k\geq 2+\gamma$, $F$ is the unique global dissipative mild solution by Lemma \ref{lem_unique}. Moreover, $F\in \mathcal{D}_{DL}(F_0)$ since $\emptyset \neq \mathcal{D}_{DL}\subset \mathcal{U}_{\infty}(F_0)$. So, $F$ satisfies entropy inequality \eqref{entropy ineq}. Obviously, $F$ satisfies \eqref{macro_conserv} either by \eqref{apriori-poly}. By initial assumption \ref{ini profile-poly}, we replace \eqref{lower bound in time} into the following argument : we have
    \begin{equation} \notag
		u^{(a)}(t,x)  \geq \int_{|v_*|\leq R} |a-v_*|^{\g} F dv_* \geq e^{-A_0\langle R \rangle^{\g} t}\int_{|v_*|\leq R} |a-v_*|^{\g} \varphi(v) dv_*
	\end{equation}
    and the last integral on RHS is continuous in $a$ and positive for each $a$. From compactness, 
    \[
        \inf_{|a|\leq R} u^{(a)}(t,x) = \inf_{|a|\leq R} e^{-A_0\langle R \rangle^{\g} t}\int_{|v_*|\leq R} |a-v_*|^{\g} \varphi(v) dv_* > 0.    
    \]
    Therefore, we obtain a coercivity estimate for $L(F)$ as in Proposition \ref{prop: macro lb} with help of \eqref{apriori-poly}. Then applying to Proposition \ref{prop uniform bound} (which holds for $F$ as mentioned in Remark \ref{rmk unif-p}), we get a uniform bound in $L^\infty(\avv^{\b_1})$. And, from the definition of $\b_1$ of Proposition \ref{prop uniform bound}, $\b_1$ is sufficiently large so that $\b_1 > \max\{11, \beta_\infty\}$ if $\b$ and $k$ are large enough. So there exists $\b_*$ and $k_*$ such that if $\b \geq \b_*$ and $k \geq k_*$ then, $\b_1 > \max\{11, \beta_\infty\}$ holds so that Proposition \ref{prop:rel-entro decay} also hold.	
	We set 
	\begin{align*}
		M_0 = \max\left\{1, \sup_{t\ge 0}\|\avv^{\b_1}(F(t)-\mu)\|_{L^\infty_{x,v}}\right\}<\infty, \quad \beta_1 \ge \max\{11,\beta_\infty\}.	
    \end{align*}
	Then we choose $\varepsilon_0$  as in $\varepsilon_0(M_0)$ of Theorem \ref{Mainresultlargeamp-p}. By entropy decay in Proposition \ref{prop:rel-entro decay}, there exists a finite time $T_{ent}(\varepsilon_0) < \infty$ such that
	\[
		H(F \ \vert \ \mu)(T_{ent}) \leq \varepsilon_0
	\]
	where $F = \mu + f$.  We apply Theorem \ref{Mainresultlargeamp-p} with the weight exponent $\beta_1$, taking $T_{ent}$ as the new initial time. Therefore, we obtain decay estimate of Theorem \ref{thm:cond-p}.

\end{proof}

\begin{proof} [Proof of Corollary \ref{coro:cond-p}]
Let us choose $\b_* > \max\{11, \beta_\infty\}$. If we assume \eqref{cond Coro-p}, the initial $F$ satisfies \eqref{macro_conserv} and entropy equality obviously. Lemma \ref{lem_unique} also hold and we can apply Proposition \ref{prop:rel-entro decay} directly since $\b \geq 11$. In this case, we do not need to assume \eqref{ini profile-poly}, since a Gaussian lower bound appears near initial time. Using \eqref{cond Coro-p} again, we can assume $F_0$ has a Gaussian lower bound, without loss of generality. Now, similar as in the proof of Theorem \ref{thm:cond-p}, there exists a finite time $T_{ent} < \infty$ from which large amplitude problem regime begins. Then we obtain \eqref{decay F coro-p} by Theorem \ref{Mainresultlargeamp-p} since $\b_1 \geq \beta_\infty$. In particular, we do not lose decay rate $\avv^{-\b_1}$ of \eqref{cond Coro-p}.


\end{proof}

\subsection{Proof of Theorem \ref{thm:cond-low 2}}

\begin{proof} [Proof of Theorem \ref{thm:cond-low 2}]
	Since we assume 
    \Be \notag
    \sup_{0\le t\le T}
    \|\langle v\rangle^{\beta_1}F(t)\|_{L^\infty_{x,v}}<\infty
    \qquad\text{for every } \ T<\infty,
    \Ee
    for sufficiently large $\b_1$, we obtain conservation of mass, momentum, and energy, and entropy inequality. Moreover, $\|F\|_\infty$ is uniformly bounded for very short time near initial time. Therefore, we can apply argument of Section 3 and there exists $T_{glb} > 0$ where $F(T_{glb})$ admits a uniform Gaussian lower bound as in Proposition \ref{prop:glb}. Therefore, we can treat $T_{glb}$ as our new initial time with a Gaussian lower bound $\varphi(v) > 0$. Then, we can proceed with the argument of Section 4. We have a weak convergence
	\begin{equation}\label{eq:r1-average-limit2}
		\int F_n\psi\,dv\longrightarrow\int m\psi\,dv
		\quad\hbox{strongly in }L^1(\Omega_T)
	\end{equation}
	as in \eqref{eq:r1-average-limit}, but $m \neq \mu$ in general, because \eqref{eq:r1-energy-tightness} is not assumed. Following Lemma \ref{lem:r1-Maxwellian-rigidity} and Proposition \ref{Desvillettes}, $m(v)$ is a stationary global Maxwellian. Moreover, \eqref{M2} gives mass and momentum velocity tightness, so we obtain $m(v) = e^{a + c|v|^2}$ for some $a\in \R$ and $c < 0$. Then we can follow the same argument of Section 4, replacing $\mu$ with $m(v)$ since \eqref{eq:r1-energy-tightness} is only used to identify $\mu=m$.
    
    	Now, we consider Proposition \ref{prop uniform bound}. If we choose $k=2$, we obtain
	\[
		\mathfrak{a}=2+2\g,\ \ \theta = -\frac{2\g}{3}, \ \ \mathfrak{b}_1=1-\g,\ \ \mathfrak{b}=2-\g,
	\]
	by direct computation. This gives $\mathfrak{b}_1 + k + \g + 2 = 5$ and therefore, we have uniform-in-time bound
	\[
		\sup_{t\geq 0}\|\avv^5 F\|_{L^\infty_{x,v}} \leq K_5 < \infty,
	\]  
	from Proposition \ref{prop uniform bound}. To improve $L^\infty_5$ to $L^\infty_{\b_1}$, we consider Carleman representation
	\[
		Q_{+}(F)(t,x,v) \leq C_b\int_{u\in\R^3}\int_{w\in\R^3}\sqrt{|u-v|^2 + |w-v|^2}^{\g-1} \delta((u-v)\cdot(w-v)) F(u)F(w) dudw.
	\]
	Using symmetric property of $u$ and $w$, for $t\leq T < \infty$,
	\begin{equation} \label{Q gain Carleman}
	\begin{split}
		\avv^{\b_1} Q_{+}(F)(t,x,v) &\leq C_{B}\avv^{\b_1} \int_{u\in\R^3} \frac{F(t,x,u)}{|u-v|^{2-\g}}\int_{w\in E_{v,u}}F(t,x,w) 1_{|w|\geq |u|} dE_{w} du \\
		&\leq C_B \sup_{t\leq T}\|\avv^{\b_1}F(t)\|_{L^\infty_{x,v}} \avv^{\b_1} \int_{u\in\R^3} \frac{F(t,x,u)}{|u-v|^{2-\g}}\int_{w\in E_{v,u}} \langle w\rangle^{-\b_1} 1_{|w|\geq |u|} dE_{w} du,
	\end{split}
	\end{equation}
	where $C_{B}$ is some uniform constant, depending on kernel $B$, in particular $C_b$. For the last planar integral on $E_{v,u}$, we can compute directly. Set $d_{E_{v,u}}=\text{dist}(0,E_{v,u})$. When $|u|\geq d_{E_{v,u}}$,  
	\[
	\begin{split}
		\int_{w\in E_{v,u}} \langle w\rangle^{-\b_1} 1_{|w|\geq |u|} dE_{w} 1_{|u|\geq d_{E_{v,u}}} 
		&= \int_{\sqrt{|u|^2-d^2}}^{\infty} \frac{2\pi r}{(1 + d^2 + r^2)^{\frac{\b_1}{2}}} dr \\
		&= \frac{2\pi}{\b_1-2} (1+|u|^2)^{1-\b_1/2},
	\end{split}
	\]
	where we used abbreviation $d= d_{E_{v,u}}$. Similarly, if $|u|\leq d_{E_{v,u}}$, we obtain 
	\[
		\int_{w\in E_{v,u}} \langle w\rangle^{-\b_1} 1_{|w|\geq |u|} dE_{w} 1_{|u|\leq d_{E_{v,u}}} \\
		= \int_{0}^{\infty} \frac{2\pi r}{(1 + d^2 + r^2)^{\frac{\b_1}{2}}} dr \\
		= \frac{2\pi}{\b_1-2} (1+d_{E_{v,u}}^2)^{1-\b_1/2}.
	\]  
	Therefore, we obtain
	\begin{equation} \label{plane int}
		\int_{w\in E_{v,u}} \langle w\rangle^{-\b_1} 1_{|w|\geq |u|} dE_{w} \leq \frac{2\pi}{\b_1-2} \Big(1+\max\{d_{E_{v,u}}^2,|u|^2\}\Big)^{1-\b_1/2}.
	\end{equation}
	
	First, when $|v|\leq 4\sqrt{\b_1}$, we can easily obtain 
	\[
		|\avv^{\b_1}F(t,x,v)| \leq K_5 \avv^{\b_1-5} \leq C_{\b_1} K_5 < \infty
	\]
	for some constant $C_{\b_1} < \infty$. Hence we consider $|v|\geq 4\sqrt{\b_1}$ case only. We perform estimates of $\avv^{\b_1} Q_{+}(F)$ for the following three regions. \\

	\textit{Region 1.} Let us consider the region $|u| \leq \frac{1}{4\sqrt{\b_1}}|v|$ on RHS of \eqref{Q gain Carleman}. In this region,
	\[
		d_{E_{v,u}} \geq |v|\sqrt{1 - \frac{1}{16\b_1}} 
		\geq \frac{|v|}{4\sqrt{\b_1}} \geq |u|, \ \ \text{for} \ \ \b_1 \geq \frac{1}{8}.
	\]
	Using \eqref{plane int}, 
	\begin{equation} \label{region 1 est}
		\begin{split}
			\text{RHS of \eqref{Q gain Carleman}} &\leq \frac{C_B}{\b_1-2} \sup_{t\leq T}\|\avv^{\b_1}F(t)\|_{L^\infty_{x,v}} \avv^{2} \int_{|u|\leq \frac{1}{4\sqrt{\b_1}}|v|} \frac{F(t,x,u)}{|u-v|^{2-\g}}du \\
			&\leq \frac{C_B}{\b_1-2} \sup_{t\leq T}\|\avv^{\b_1}F(t)\|_{L^\infty_{x,v}} \int_{|u|\leq \frac{1}{4\sqrt{\b_1}}|v|} |v-u|^\g F(t,x,u) du \\
			&\leq \frac{C_B}{\b_1-2} \sup_{t\leq T}\|\avv^{\b_1}F(t)\|_{L^\infty_{x,v}} L(F)(t,x,v), 
		\end{split}
	\end{equation}
	since $|v-u| \approx |v| \approx \avv$ uniformly in this region. In particular, $C_B$ is independent of $\b_1$. \\
	
	\textit{Region 2.} Let us consider the region $\frac{1}{4\sqrt{\b_1}}|v| \leq |u|\leq 2|v|$ on RHS of \eqref{Q gain Carleman}. In this region, we have $d^2(E_{v,u}) + |u|^2 \geq |v|^2$. (If $|u|\leq |v|$, $d^2+|u|^2 \geq |v|^2$. Otherwise, if $|u|\geq |v|$, $d^2+|u|^2 \geq |v|^2$.) Applying to \eqref{plane int},  
	\begin{equation} \notag
		\int_{w\in E_{v,u}} \langle w\rangle^{-\b_1} 1_{|w|\geq |u|} dE_{w} \leq C_{\b_1} \avv^{2-\b_1}.
	\end{equation}
	Therefore,
	\begin{equation} \label{region 2 est}
		\begin{split}
			\text{RHS of \eqref{Q gain Carleman}} &\leq C_{\b_1} \sup_{t\leq T}\|\avv^{\b_1}F(t)\|_{L^\infty_{x,v}} K_5 \avv^{2} \int_{\frac{1}{4\sqrt{\b_1}}|v| \leq |u|\leq 2|v|} \frac{1}{\langle u \rangle^5} \frac{1}{|u-v|^{2-\g}} du \\
			&\leq  C_{\b_1} \sup_{t\leq T}\|\avv^{\b_1}F(t)\|_{L^\infty_{x,v}} K_5 \avv^{\g-2}.
		\end{split}
	\end{equation}
	
	\textit{Region 3.} Let us consider the region $|u|\geq 2|v|$ on RHS of \eqref{Q gain Carleman}. In this region, $d^2(E_{v,u}) + |u|^2 \geq |v|^2$ and
	\begin{equation} \notag
		\int_{w\in E_{v,u}} \langle w\rangle^{-\b_1} 1_{|w|\geq |u|} dE_{w} \leq C_{\b_1} \avv^{2-\b_1}
	\end{equation}
	still hold. Therefore,
	\begin{equation} \label{region 3 est}
		\begin{split}
			\text{RHS of \eqref{Q gain Carleman}} &\leq C_{\b_1} \sup_{t\leq T}\|\avv^{\b_1}F(t)\|_{L^\infty_{x,v}} K_5 \avv^{2} \int_{|u|\geq 2|v|} \frac{1}{\langle u \rangle^5} \frac{1}{|u-v|^{2-\g}} du \\
			&\leq  C_{\b_1} \sup_{t\leq T}\|\avv^{\b_1}F(t)\|_{L^\infty_{x,v}} K_5 \avv^{\g-2}.
		\end{split}
	\end{equation}
	
	From \eqref{region 1 est}, \eqref{region 2 est}, and \eqref{region 3 est}, we obtain
	\begin{equation} \label{final est Q gain}
		\avv^{\b_1} Q_{+}(F)(t,x,v) \leq \sup_{t\leq T}\|\avv^{\b_1}F(t)\|_{L^\infty_{x,v}} \Big[ \frac{C_B}{\b_1-2} L(F)(t,x,v) + C_{B,\b_1} K_5 \avv^{\g-2} \Big].
	\end{equation}
	In particular, $C_B$ does not depend on $\b_1$. Choose $\b_1 > \max\{4, 2+4C_B\}$ and for a chosen $\b_1$, we define
	\[
		R_0 := \max\big\{ 4\sqrt{\b_1}, \big( \frac{4C_{B,\b_1}K_5}{\lambda_F}\big)^{\half} \big\}.
	\]
	Then for $|v|\geq R_0$, we have 
	\[
		\frac{C_B}{\b_1-2} \leq \frac{1}{4} \ \ \text{and} \ \ C_{B,\b_1}K_5\avv^{\g-2} \leq \frac{1}{4}L(F)(t,x,v).
	\]
	Applying to \eqref{final est Q gain}, for $|v|\geq R_0$ and $t\leq T < \infty$,
	\begin{equation} 
		\begin{split}
			&\avv^{\b_1} F(t,x,v) \\
			&\leq e^{-\int_{0}^{t} L(F)(s, x-v(t-s), v) ds} \avv^{\b_1} F_0(x-vt, v) \\
			&\quad + \int_{0}^{t} e^{-\int_{s}^{t} L(F)(\tau, x-v(t-\tau), v) d\tau} \avv^{\b_1} Q_{+}(F,F)(s, x-v(t-s), v) ds \\
			&\leq \avv^{\b_1} F_0(x-vt, v) + \int_{0}^{t} e^{-\int_{s}^{t} L(F)(\tau, x-v(t-\tau), v) d\tau} \avv^{\b_1} Q_{+}(F,F)(s, x-v(t-s), v) ds \\
			&\leq \|\avv^{\b_1} F_0\|_{L^\infty_{x,v}} + \sup_{0\leq t\leq T}\|\avv^{\b_1} F(t)\|_{L^\infty_{x,v}} \int_{0}^{t} e^{-\int_{s}^{t} L(F)(\tau, x-v(t-\tau), v) d\tau} \frac{1}{2}L(F)(s, x-v(t-s), v) ds .\\
		\end{split}
	\end{equation}

    For $|v|\leq R_0$, we easily get $|\avv^{\b_1} F(t,x,v)| \leq \langle R_0 \rangle^{\b_1-5} K_5$, and hence
\[
\sup_{0\leq t\leq T}\|\avv^{\b_1} F(t)\|_{L^\infty_{x,v}} \leq
\max\left\{
K_5\langle R_0\rangle^{\beta_1-5},
\|F_0\|_{L^\infty_{\beta_1}} + \frac{1}{2} \sup_{0\leq t\leq T}\|\avv^{\b_1} F(t)\|_{L^\infty_{x,v}}
\right\}.
\]
Since $\sup_{0\leq t\leq T}\|\avv^{\b_1} F(t)\|_{L^\infty_{x,v}} < \infty$, we obtain uniform-in-time boundedness
\[
\sup_{0\leq t\leq T}\|\avv^{\b_1} F(t)\|_{L^\infty_{x,v}} \leq
\max\left\{
K_5\langle R_0\rangle^{\beta_1-5},
2\|F_0\|_{L^\infty_{\beta_1}}
\right\}.
\]
Finally, set
\[
\b_*:=\max\{11,\b_\infty,2+4C_B\},
\]
where $\b_\infty$ is defined in Theorem~\ref{Mainresultlargeamp-p} and $C_B$ is the constant appearing in \eqref{Q gain Carleman}. If $\b_1>\b_*$, the preceding argument based on the Carleman representation yields a uniform-in-time bound for $F$ in $L^\infty_{\b_1}$, namely,
\[
F\in L^\infty([0,\infty);L^\infty_{\b_1}).
\]
This uniform bound allows us to identify the weak limit $m(v)$ with $\mu(v)$. Consequently, both Propositions~\ref{prop:rel-entro decay} and~\ref{prop:micro decay} apply. In particular, the relative entropy becomes sufficiently small after some finite time, allowing us to invoke Theorem~\ref{Mainresultlargeamp-p} to conclude the proof.
\end{proof}

\bigskip

\noindent{\bf Data availability:} No data was used for the research described in the article.
\newline

\noindent{\bf Conflict of interest:} The authors declare that they have no conflict of interest.\newline

\noindent{\bf Disclaimer on the use of AI tools} During the preparation of this manuscript, the authors utilized generative AI (Chat GPT) exclusively for English language editing and cross-checking the validity of mathematical derivations. All core mathematical ideas, formulations, and results originated entirely from the authors. The authors take full responsibility for the contents and ultimate integrity of this work. \newline

\noindent{\bf Acknowledgement} The research of L.-B. He was supported by NSF of China under Grant No.11771236 and New Cornerstone Investigator Program 100001127.
J.-C. Jiang is supported by  National Sci-Tech Grant 
MOST 115-2115-M-007-014-MY3 and NSTC 115-2124-M-002 -009.  J. Kim, D. Lee, and S. Park are supported by the National Research Foundation of Korea(NRF) grant funded by the Korea government(MSIT)(No.RS-2023-00212304 and No.RS-2023-00219980). S. Park is also supported by the Basic Science Research Institute Fund(No.RS-2021-NR060139). \newline

\bigskip

\appendix

\section{Local theory}

\begin{proposition} \label{prop_loc}
	Let $w(v)=\langle v\rangle^\beta e^{\alpha|v|^\kappa}$ with $\beta\ge0$, $\alpha>0$, and $\gamma\le \kappa<2$.  Assume that $F_0(x,v) \ge 0$ and $ \|\langle v\rangle^\beta e^{\alpha|v|^\kappa}F_0\|_{L^\infty_{x,v}}<\infty$. Then there exists a positive time $\hat{t}_0:= \left(
	\hat{C}\|\langle v\rangle^\beta e^{\alpha|v|^\kappa}F_0\|_{L^\infty_{x,v}}\right)^{-1}$ such that there is a unique solution $F(t,x,v) \ge 0 $ on $[0,\hat{t}_0]$ to the Boltzmann equation \eqref{BE} with initial data $F(0,x,v)=F_0(x,v)$ satisfying
	\begin{align*}
		\sup_{0\le t \le \hat{t}_0}\|\langle v\rangle^\beta e^{\alpha|v|^\kappa}F(t)\|_{L^\infty_{x,v}} \le 2\|\langle v\rangle^\beta e^{\alpha|v|^\kappa}F_0\|_{L^\infty_{x,v}}.
	\end{align*}
\end{proposition}

\begin{proof}
	To consider the local existence of solutions for the Boltzmann equation \eqref{BE}, we start from the iteration system that for $m=0,1,2,\cdots$,
	\begin{align} \label{iterativesys1}
		\partial_t F^{(m+1)} + v \cdot \nabla_x F^{(m+1)} + F^{(m+1)} \int_{\R^3} \int_{\S^2} B(v-v_*, \omega) F^{(m)}(t,x,v_*)d\omega dv_*= Q^+(F^{(m)},F^{(m)})
	\end{align}
	with initial data $F^{(m+1)}(t,x,v) |_{t=0} =F_0(x,v) \ge 0$, where we set $F^{(0)}(t,x,v) = 0$. We denote
	\begin{align*}
		\nu(F^{(m)}):=\int_{\R^3} \int_{\S^2} B(v-v_*, \omega) F^{(m)}(t,x,v_*)d\omega dv_*.
	\end{align*}
	Then the iterative system \eqref{iterativesys1} for $F^{(m+1)}$ with weight $\langle v\rangle^\beta e^{\alpha|v|^\kappa}$ can be written as, for $m=0,1,2,\cdots$,
	\begin{align} \label{iterativesys2}
		&\partial_t (\langle v\rangle^\beta e^{\alpha|v|^\kappa}F^{(m+1)}) + v \cdot \nabla_x (\langle v\rangle^\beta e^{\alpha|v|^\kappa}F^{(m+1)}) + \nu(F^{(m)})\langle v\rangle^\beta e^{\alpha|v|^\kappa}F^{(m+1)} \\\nonumber
		&=  \langle v\rangle^\beta e^{\alpha|v|^\kappa}Q^+\left(F^{(m)},F^{(m)}\right)	
	\end{align}
	with initial data $F^{(m+1)}(t,x,v) |_{t=0} =F_0(x,v)$ and $F^{(0)}(t,x,v) \equiv 0$.\\
	\newline
	(1) $\mathit{(Uniform\ boundedness\ and\ Positivity)}$\\
	We use the induction argument on $m=0,1,2,\cdots$. Now, we show the uniform bound and positivity for solutions to the system \eqref{iterativesys1}:
	\begin{align} \label{uniformbound}
		\|\langle v\rangle^\beta e^{\alpha|v|^\kappa}F^{(m)}(t)\|_{L^\infty_{x,v}} \le 2 \|\langle v\rangle^\beta e^{\alpha|v|^\kappa}F_0\|_{L^\infty_{x,v}}, \quad F^{(m)}(t,x,v) \ge 0 
	\end{align} 
	for all $m$ and $0 \le t \le \hat{t}_1$, where $\hat{t}_1$ is determined later. \\
	We denote the function $I^{(m)}(t,s)$ by
	\begin{align*}
		I^{(m)}(t,s) = \exp\left\{-\int_s^t \nu(F^{(m)})(\tau,x-v(t-\tau),v)d\tau\right\}. 
	\end{align*}
	for $m=0,1,2,\cdots$.
	First of all, for $m=0$, we can easily get
	\begin{align} \label{local1}
		\|\langle v\rangle^\beta e^{\alpha|v|^\kappa}F^{(1)}(t)\|_{L^\infty_{x,v}} \le 2 \|\langle v\rangle^\beta e^{\alpha|v|^\kappa}F_0\|_{L^\infty_{x,v}}
	\end{align}
	for all $t \ge 0$. We also have
	\begin{align*} 
		F^{(1)}(t,x,v) &= I^{(0)}(t,0) F_0(x-vt,v) =F_0(x-vt,v) \ge 0.
	\end{align*}
	Suppose that \eqref{uniformbound} holds for $m=0,1,2, \cdots, n-1$. We consider \eqref{uniformbound} for $m =n$. Using the Duhamel principle, the system \eqref{iterativesys1} implies that
	\begin{align*} 
		\langle v\rangle^\beta e^{\alpha|v|^\kappa}F^{(n+1)}(t,x,v) &= I^{(n)}(t,0)\langle v\rangle^\beta e^{\alpha|v|^\kappa}F_0 + \int_0^t I^{(n)}(t,s)\langle v\rangle^\beta e^{\alpha|v|^\kappa}Q^+ \left(F^{(n)},F^{(n)}\right)(s)ds.
	\end{align*}
	Note that $I^{(n)}(t,s) \le 1$ since $F^{(n)}\ge 0$ by the induction assumption.
	From \eqref{uniformbound} and Lemma \ref{pointwiseQ+est}, we obtain
	\begin{align*}
		\|\langle v\rangle^\beta e^{\alpha|v|^\kappa}F^{(n+1)}(t)\|_{L^\infty_{x,v}} &\le \|\langle v\rangle^\beta e^{\alpha|v|^\kappa}F_0\|_{L^\infty_{x,v}}+\int_0^t   \left\|\langle v\rangle^\beta e^{\alpha|v|^\kappa}Q^+ \left(F^{(n)},F^{(n)}\right)(s)\right\|_{L^\infty_{x,v}} ds\\
		& \le \|\langle v\rangle^\beta e^{\alpha|v|^\kappa}F_0\|_{L^\infty_{x,v}} + C_{\alpha,\kappa,\beta,\gamma}\int_0^t \|\langle v\rangle^\beta e^{\alpha|v|^\kappa}F^{(n)}(s)\|_{L^\infty_{x,v}}^2ds\\
		& \le \|\langle v\rangle^\beta e^{\alpha|v|^\kappa}F_0\|_{L^\infty_{x,v}}\left(1+\hat{C}_1t\|\langle v\rangle^\beta e^{\alpha|v|^\kappa}F_0\|_{L^\infty_{x,v}}  \right),
	\end{align*}
	where $\hat{C}_1$ is a constant depending on $\alpha,\kappa,\beta,\gamma$.
	Taking $\hat{t}_1:= \left(\hat{C}_1\|\langle v\rangle^\beta e^{\alpha|v|^\kappa}F_0\|_{L^\infty_{x,v}}\right)^{-1}$, it holds that
	\begin{align*}
		\|\langle v\rangle^\beta e^{\alpha|v|^\kappa}F^{(n+1)}(t)\|_{L^\infty_{x,v}} \le 2\|\langle v\rangle^\beta e^{\alpha|v|^\kappa}F_0\|_{L^\infty_{x,v}}
	\end{align*}
	for $0\le t \le \hat{t}_1$. We conclude the uniform bound for solutions to the system \eqref{iterativesys2}.\\
	We now show the positivity for $F^{(n+1)}$. We have
	\begin{align*}
		F^{(n+1)}(t,x,v) &=  I^{(n)}(t,0)F_0(x-vt,v)\\
		&\quad +\int_0^t  I^{(n)}(t,s)Q^+(F^{(n)},F^{(n)})(s,x-v(t-s),v) ds\\
		& \ge 0,
	\end{align*}
	Here, since $F^{(n)}\ge 0$, $\nu(F^{(n)})\ge 0$ and $Q^+(F^{(n)},F^{(n)}) \ge 0$.
	Thus, $F^{(m+1)}(t,x,v) \ge 0$ over $[0,\hat{t}_1]$. By the induction, \eqref{uniformbound} holds for all $m=0,1,2, \cdots$.\\
	\newline
	(2) $(Existence)$\\
	We denote $\hat{w}(v) = \frac{\langle v\rangle^\beta e^{\alpha|v|^\kappa}}{\sqrt{1+|v|^2}}$ and
	\begin{align*}
		\hat{F}^{(m)}:= \frac{\langle v\rangle^\beta e^{\alpha|v|^\kappa}}{\sqrt{1+|v|^2}}F^{(m)}(t,x,v).
	\end{align*}
	\newline
	By Duhamel principle, for each $m$, we have
	\begin{align*}
		\hat{F}^{(m+1)}(t,x,v) &= I^{(m)}(t,0) \hat{F}_0(x-vt,v) + \int_0^t I^{(m)}(t,s)\hat{w}Q^+ \left(\frac{\hat{F}^{(m)}}{\hat{w}},\frac{\hat{F}^{(m)}}{\hat{w}}\right)(s,x-v(t-s),v)ds.	
	\end{align*}
	By subtracting $\hat{F}^{(m+2)} - \hat{F}^{(m+1)}$, we obtain
	\begin{equation} \label{local8}
		\begin{aligned}
			&\left|\left(\hat{F}^{(m+2)} - \hat{F}^{(m+1)} \right)(t,x,v)\right|\\
			& \le |I^{(m+1)}(t,0)- I^{(m)}(t,0) ||\hat{F}_0(x-vt,v)|\\
			& \quad +\int_0^t |I^{(m+1)}(t,s)- I^{(m)}(t,s) | \left| \hat{w}Q^+ \left(\frac{\hat{F}^{(m+1)}}{\hat{w}},\frac{\hat{F}^{(m+1)}}{\hat{w}}\right)(s,x-v(t-s),v)\right|ds\\
			& \quad + \int_0^t I^{(m)}(t,s)\left|\left( \hat{w}Q^+ \left(\frac{\hat{F}^{(m+1)}}{\hat{w}},\frac{\hat{F}^{(m+1)}}{\hat{w}}\right)-\hat{w}Q^+ \left(\frac{\hat{F}^{(m)}}{\hat{w}},\frac{\hat{F}^{(m)}}{\hat{w}}\right) \right)(s,x-v(t-s),v)\right|ds\\
			& =:J_1^{(m)} +J_2^{(m)}+J_3^{(m)}.
		\end{aligned}
	\end{equation}
	First of all, let us estimate $J_1^{(m)}$. By a direct calculation, we have
	\begin{equation} \label{local2}
		\begin{aligned}
			&|I^{(m+1)}(t,s)- I^{(m)}(t,s) |\\ &\le \left|\int_s^t \left[\nu(F^{(m+1)})-\nu(F^{(m)})\right](\tau,x-v(t-\tau),v) d\tau \right| \\
			& \le C \nu(v) \int_s^t  \int_{\R^3} \nu(v_*)\langle v_*\rangle^{-\beta+1} e^{-\alpha|v_*|^\kappa} \left|\left(\hat{F}^{(m+1)}-\hat{F}^{(m)}\right)(\tau,x-v(t-\tau),v_*)
			\right|dv_*d\tau \\
			&\le C\nu(v) \int_s^t\left\|\left(\hat{F}^{(m+1)}-\hat{F}^{(m)}\right)(\tau)\right\|_{L^\infty_{x,v}}d\tau.
		\end{aligned}
	\end{equation}
	Then from \eqref{local2},
	\begin{align} \label{local3}
		\|J_1^{(m)}\|_{L^\infty_{x,v}} &\le C \|\langle v\rangle^\beta e^{\alpha|v|^\kappa}F_0\|_{L^\infty_{x,v}}\int_0^t\left\|\left(\hat{F}^{(m+1)}-\hat{F}^{(m)}\right)(\tau)\right\|_{L^\infty_{x,v}}d\tau\nonumber\\
		& \le C t\|\langle v\rangle^\beta e^{\alpha|v|^\kappa}F_0\|_{L^\infty_{x,v}}\sup_{0\le \tau \le t}\left\|\left(\hat{F}^{(m+1)}-\hat{F}^{(m)}\right)(\tau)\right\|_{L^\infty_{x,v}}.
	\end{align}
	Next, let us estimate $J_2^{(m)}$. By Lemma \ref{pointwiseQ+est} and \eqref{local2}, we obtain 
	\begin{align} \label{local4}
		\|J_2^{(m)}\|_{L^\infty_{x,v}} &\le C \int_0^t \left\| \langle v\rangle^\beta e^{\alpha|v|^\kappa}Q^+ \left(F^{(m+1)},F^{(m+1)}\right)(s)\right\|_{L^\infty_{x,v}}\nonumber\\
		&\qquad \times \int_s^t\left\|\left(\hat{F}^{(m+1)}-\hat{F}^{(m)}\right)(\tau)\right\|_{L^\infty_{x,v}}d\tau ds \nonumber\\
		& \le C \int_0^t   \|\langle v\rangle^\beta e^{\alpha|v|^\kappa}F^{(m+1)}(s)\|_{L^\infty_{x,v}}^2  \int_s^t\left\|\left(\hat{F}^{(m+1)}-\hat{F}^{(m)}\right)(\tau)\right\|_{L^\infty_{x,v}}d\tau ds\nonumber\\
		& \le C t^2 \sup_{0\le s \le t}\|\langle v\rangle^\beta e^{\alpha|v|^\kappa}F^{(m+1)}(s)\|_{L^\infty_{x,v}}^2  \sup_{0\le \tau \le t}\left\|\left(\hat{F}^{(m+1)}-\hat{F}^{(m)}\right)(\tau)\right\|_{L^\infty_{x,v}}.
	\end{align}
	Combining \eqref{local3} and \eqref{local4}, and using \eqref{uniformbound}, we get for $t \le \hat{t}_1$,
	\begin{align} \label{local5}
		&\|J_1^{(m)}+J_2^{(m)}\|_{L^\infty_{x,v}}\nonumber \\ &\le Ct\left(\|\langle v\rangle^\beta e^{\alpha|v|^\kappa}F_0\|_{L^\infty_{x,v}}+t \sup_{0\le s \le t} \|\langle v\rangle^\beta e^{\alpha|v|^\kappa}F^{(m+1)}(s)\|_{L^\infty_{x,v}}^2  \right)\nonumber\\
		& \qquad \times\sup_{0\le \tau \le t}\left\|\left(\hat{F}^{(m+1)}-\hat{F}^{(m)}\right)(\tau)\right\|_{L^\infty_{x,v}}\nonumber\\
		& \le C t\|\langle v\rangle^\beta e^{\alpha|v|^\kappa}F_0\|_{L^\infty_{x,v}}\sup_{0\le \tau \le t}\left\|\left(\hat{F}^{(m+1)}-\hat{F}^{(m)}\right)(\tau)\right\|_{L^\infty_{x,v}}.
	\end{align}
	Note that
	\begin{align*}
		&\hat{w}Q^+ \left(\frac{\hat{F}^{(m+1)}}{\hat{w}},\frac{\hat{F}^{(m+1)}}{\hat{w}}\right)-\hat{w}Q^+ \left(\frac{\hat{F}^{(m)}}{\hat{w}},\frac{\hat{F}^{(m)}}{\hat{w}}\right) \\
		&= \hat{w}Q^+ \left(\frac{\hat{F}^{(m+1)}}{\hat{w}}-\frac{\hat{F}^{(m)}}{\hat{w}},\frac{\hat{F}^{(m+1)}}{\hat{w}}\right) +\hat{w}Q^+ \left(\frac{\hat{F}^{(m)}}{\hat{w}},\frac{\hat{F}^{(m+1)}}{\hat{w}}-\frac{\hat{F}^{(m)}}{\hat{w}} \right).
	\end{align*}
	For $J_3^{(m)}$, we use Lemma \ref{pointwiseQ+est} and the positivity for $F^{(m)}$ to obtain
	\begin{align} \label{local6}
		\|J_3^{(m)}\|_{L^\infty_{x,v}} \le C\int_0^t \left(\|\hat{F}^{(m)}(s)\|_{L^\infty_{x,v}} + \|\hat{F}^{(m+1)}(s)\|_{L^\infty_{x,v}}\right)\|(\hat{F}^{(m+1)}-\hat{F}^{(m)})(s)\|_{L^\infty_{x,v}}ds.
	\end{align}
	For $t\le \hat{t}_1$, we apply the uniform boundedness \eqref{uniformbound} to \eqref{local6} :
	\begin{align} \label{local7}
		\|J_3^{(m)}\|_{L^\infty_{x,v}} \le Ct \|\langle v\rangle^\beta e^{\alpha|v|^\kappa}F_0\|_{L^\infty_{x,v}}\sup_{0\le s\le t} \|(\hat{F}^{(m+1)}-\hat{F}^{(m)})(s)\|_{L^\infty_{x,v}}.
	\end{align}
	Plugging \eqref{local5} and \eqref{local7} into \eqref{local8}, we derive
	\begin{align} \label{local9}
		\sup_{0\le s  \le t}\|(\hat{F}^{(m+2)}-\hat{F}^{(m+1)})(s)\|_{L^\infty_{x,v}} \le \hat{C}_2 t \|\langle v\rangle^\beta e^{\alpha|v|^\kappa}F_0\|_{L^\infty_{x,v}}\sup_{0\le s\le t} \|(\hat{F}^{(m+1)}-\hat{F}^{(m)})(s)\|_{L^\infty_{x,v}}
	\end{align}
	for all $0\le t \le \hat{t}_1$.\\
	\newline
	We set $\hat{C}:=\max\{\hat{C}_1,2\hat{C}_2\}$ and choose
	\begin{align*}
		\hat{t}_0 = \frac{1}{\hat{C}\|\langle v\rangle^\beta e^{\alpha|v|^\kappa}F_0\|_{L^\infty_{x,v}}}.
	\end{align*}
	Then it holds that $\hat{t}_0 \le \hat{t}_1$, and it follows from \eqref{local9} that $0\le t \le \hat{t}_0$,
	\begin{align*}
		\sup_{0\le s  \le t}\|(\hat{F}^{(m+2)}-\hat{F}^{(m+1)})(s)\|_{L^\infty_{x,v}} \le \frac{1}{2} \sup_{0\le s\le t} \|(\hat{F}^{(m+1)}-\hat{F}^{(m)})(s)\|_{L^\infty_{x,v}}.
	\end{align*}
	Hence $\hat{F}^{(m)}$ is a Cauchy sequence, and $F^{(m)}$ is also a Cauchy sequence. Therefore, there exists a limit $F(t,x,v)$, which is a solution to the equation \eqref{BE}.\\
	\newline	
	(3) $\mathit{(Uniqueness)}$\\
	Assume that there exists another solution $\tilde{F}(t,x,v)$ to the Boltzmann equation \eqref{BE} with the same initial data as for $F$, and the bound $\sup_{0\le t \le \hat{t}_0}\|\langle v\rangle^\beta e^{\alpha|v|^\kappa}\tilde{F}(t)\|_{L^\infty_{x,v}}<\infty$.
	By the similar argument to the existence, we obtain
	\begin{align*}
		&\|(1+|v|^2)^{-1/2}(\langle v\rangle^\beta e^{\alpha|v|^\kappa}(F-\tilde{F}))(t)\|_{L^\infty_{x,v}}\\
		& \le C\left(\sup_{0\le t \le \hat{t}_0}\|\langle v\rangle^\beta e^{\alpha|v|^\kappa}\tilde{F}(t)\|_{L^\infty_{x,v}}+\sup_{0\le t \le \hat{t}_0}\|\langle v\rangle^\beta e^{\alpha|v|^\kappa}F(t)\|_{L^\infty_{x,v}} \right)\\
		& \qquad \times\int_0^t \|(1+|v|^2)^{-1/2}(\langle v\rangle^\beta e^{\alpha|v|^\kappa}(F-\tilde{F}))(s)\|_{L^\infty_{x,v}}ds.
	\end{align*}
	By the Gr\"onwall inequality, the uniqueness is guaranteed.
\end{proof}

\section{Lemmas}
Some useful and basic lemmas.
\begin{lemma} 
\label{lem_abslog}
	Assume that following mass, energy, and entropy bounds, i.e.,
	\[
	\iint_{\mathbb{T}^{3}\times\R^{3}} f(1+|v|^2+\log f) dvdx \leq C_0
	\]
	for all $t\geq 0$. Then
	\[
	\iint_{\mathbb{T}^{3}\times\R^{3}} f|\log f| dvdx \leq C_0^{\prime}.
	\]
\end{lemma}
\begin{proof}
	Decomposing $\{f \geq 1\}$ and $\{f < 1\}$, we get
	\begin{equation*}
		\begin{split}
			\iint_{\mathbb{T}^{3}\times\R^{3}} f|\log f| dvdx &= \iint_{f \geq 1} f\log f dvdx + \iint_{f < 1} f\log\frac{1}{f} dvdx \\
			&= \iint f\log f dvdx + 2 \iint_{f < 1} f\log\frac{1}{f} dvdx. \\
		\end{split}
	\end{equation*}
	To control the 2nd term on RHS, we use $x\log\frac{1}{x} \leq \sqrt{x}$ to obtain
	\begin{equation*}
		\begin{split}
			\iint_{f < 1} f\log\frac{1}{f} dvdx &\leq  \iint_{e^{-|v|^2} \leq f < 1} f\log\frac{1}{f} dvdx + \iint_{f < e^{-|v|^2}} f\log\frac{1}{f} dvdx  \\
			&\leq \iint f|v|^2 dvdx + \iint_{f < e^{-|v|^2}} \sqrt{f} dvdx < \infty.
		\end{split}
	\end{equation*}
\end{proof}

\begin{lemma} 
\label{lem_min entropy}
	Relative entropy is always non-negative. Moreover, local Maxwellian $M_{\rho, U, T}$ minimize entropy among all distributions which have macroscopic quantities $(\rho, U, T)$ of $F$.
\end{lemma}
\begin{proof}
	Since $s\log s - s + 1 \geq 0$ for $s\geq 0$, relative entropy \eqref{def rel entropy} should be non-negative.
	To prove entropy minimizing property of $M_{\rho, U, T}$, let us assume that a distribution $F$ has $(\rho, U, T)$ for its macroscopic quantities. Then
	\begin{equation*}
		\begin{split}
			0 &\leq H(F \ \vert \ M_{\rho, U, T}) = \iint F\log F - \iint F\log M_{\rho, U, T} \\
			&= \iint F\log F - \iint F \Big[ \log\rho - \frac{3}{2}\log(2\pi) -\frac{3}{2}\log T - \frac{|v-U|^2}{2T} \Big] \\
			&= \iint F\log F - \iint M_{\rho, U, T} \Big[ \log\rho - \frac{3}{2}\log(2\pi) -\frac{3}{2}\log T - \frac{|v-U|^2}{2T} \Big] \\
			&= H(F) - H(M_{\rho, U, T}).
		\end{split}
	\end{equation*}
\end{proof}

\begin{lemma} 
\label{lem_reverseCKP}
Assume that $F\geq0$ has the same total mass as $\mu$ and satisfies $\|F\|_{L^\infty_{x,v}}\leq M<\infty$. Then
\[
H(F\mid\mu)\leq C_M\|F-\mu\|_{L^1_2},
\]
where $C_M>0$ depends only on $M$ and $\mu$. 
\end{lemma}

\begin{proof}
Let $\Phi(z)=z\log z-z+1$, with the convention $0\log0=0$. Since $F$ and $\mu$ have the same total mass,
\[
H(F\mid\mu)=\iint\mu\Phi(F/\mu)\,dx\,dv.
\]
For $0\leq z\leq1$, we have $\Phi(z)\leq1-z$, while for $z\geq1$,
\[
\Phi(z)\leq(z-1)\log z.
\]
Consequently, on the set $\{F\leq\mu\}$,
\[
\mu\Phi(F/\mu)\leq\mu-F=|F-\mu|.
\]
On the other hand, on $\{F\geq\mu\}$, the uniform upper bound $F\leq M$ gives
\[
\begin{aligned}
\mu\Phi(F/\mu)
&\leq(F-\mu)\log(F/\mu)\\
&\leq(F-\mu)\log(M/\mu)\\
&\leq C_M\langle v\rangle^2|F-\mu|.
\end{aligned}
\]
Combining the two estimates and integrating over $(x,v)$ yields the desired inequality.
\end{proof}

\begin{lemma} \label{T lower}
	Assume
		\[
			\rho_{-} \leq \rho(t,x) \leq \rho_{+},\quad T(t,x) \leq T_{+},\quad \sup_{t,x}\int_{v} |F\log F| \leq H_1.
		\]
		Then, $T(t,x)$ has a uniform lower bound. 
\end{lemma}
\begin{proof}
	Note that
	\[
	3\rho T = \int_{v} F |v-U|^2 dv \geq \int_{|v-U| > r} F |v-U|^2 dv \geq r^2 \int_{|v-U| > r}F.
	\]
	For any $r>0$ we get
	\[
	\int_{|v-U| < r} F dv \leq \int_{|v-U| < r, F<j} F + \frac{1}{\log j}\int_{|v-U| < r, F>j} |F \log F| \leq j|B_r| + \frac{H_1}{\log j}.
	\]
	We choose $j=j(\rho_{-}, H_1)$ so that $\frac{H_1}{\log j}=\frac{\rho_{-}}{4}$ and $r=r(\rho_{-}, H_1)$ so that $j|B_r|=\frac{\rho_{-}}{4}$. Then,
	\[
	\int_{|v-U| < r} F \leq \frac{\rho_{-}}{2},\quad \int_{|v-U| > r} F  = \int_{v} F - \int_{|v-U| < r} F \geq \frac{\rho_{-}}{2}.
	\]
	Hence $3\rho T \geq r^2\frac{\rho_{-}}{2}$ and $T \geq \frac{r^2 \rho_{-}}{6\rho} \geq \frac{\rho_{-}}{6\rho_{+}}r^2$. 
\end{proof}

\section{Exponential convergence to global equilibrium with polynomial weight}
This section is devoted to proving Theorem \ref{Mainresultlargeamp-p}.
We state a Theorem for large amplitude problem in polynomial weight.

\begin{theorem} 
	Let \(w(v)=(1+|v|^2)^{\beta/2}\) with $\beta>\beta_\infty:=\max\left\{1+\gamma+\frac{4C_b}{c_b},\left(\frac{4C_1}{\nu_0}\right)^{\frac{4}{\gamma+3}}\right\}$, where
	\begin{align*}
		C_1=8\sqrt{2}C_b \left(2+\frac{4}{\gamma+3}+\frac{1}{e}\right) \left(\frac{32\pi}{3+\gamma} e^{-3/4} +8\pi \right),
	\end{align*}
    and $\nu_0$ is defined in \eqref{nuzerodef}.\\
Assume that $F_0(x,v) = \mu(v) +  f_0(x,v) \ge 0$ satisfying 
	\begin{align*}
		\int_{\T^3 \times \R^3}f_0 
		\begin{pmatrix}
			1\\
			v\\
			|v|^2
		\end{pmatrix}
	dxdv =0.
	\end{align*}
	For any $M_0 \ge 1$, there exists a constant $\varepsilon_0>0$, depending  only on $M_0$ and $\bar{M}$, such that if
	\begin{align*}
		\|wf_0\|_{L^\infty_{x,v}}\le M_0,\quad \sup_{t\ge 0}\|wf(t)\|_{L^\infty_{x,v}}\le \bar{M}, \quad H(F_0 \mid \mu ) \le  \varepsilon_0,
	\end{align*}
	then the solution $F(t,x,v) =\mu(v) +  f(t,x,v) \ge 0 $ of the Boltzmann equation \eqref{BE} with initial data $F_0$ satisfies 
    \begin{align*}
		 	\|wf(t)\|_{L^\infty_{x,v}} \le C_0(1+\bar{M}) M_0e^{-\lambda_\infty t}
	\end{align*}
	for all $t \ge 0$, where $C_0\ge 1$ is a generic constant.
\end{theorem}

\bigskip

\begin{proposition} \label{poMainresultsmall} 
	Let $w(v)=\left \langle v \right \rangle ^\beta=(1+|v|^2)^{\beta/2}$ with $\beta>1+\gamma+\frac{4 C_b}{c_b}$. Assume that $F_0(x,v) = \mu(v) +  f_0(x,v) \ge 0$ satisfying 
	\begin{align*}
		\int_{\T^3 \times \R^3}f_0 
		\begin{pmatrix}
			1\\
			v\\
			|v|^2
		\end{pmatrix}
	dxdv =0.
	\end{align*}
	 Then there exists $\eta_0 >0$ such that for any $f_0$ in $L^\infty_{x,v}(w)$ with $\|wf_0\|_{L^\infty_{x,v}} \le \eta_0^{1+p}$ for some $0<p<1$, there is a unique solution $F(t,x,v) =\mu(v) +  f(t,x,v) \ge 0 $ of the Boltzmann equation \eqref{BE} with initial data $F_0$, which satisfies the mass, momentum, energy conservation and there exist $C_\infty,\lambda_\infty >0$ such that
    \begin{align*}
    	\|wf(t)\|_{L_{x,v}^\infty}\le C_\infty e^{-\lambda_\infty t}(\|wf_0\|_{L_{x,v}^\infty}+\eta_0^{1+p})
    \end{align*}
    for all $t \ge 0$.\\
\end{proposition}
	
\subsection{Nonlinear $L^\infty_{x,v}$ decay in small perturbation problem}
\label{posmallperturbationscheme}
	 We establish several estimates that will be useful in the subsequent analysis in the small perturbation regime. Their proofs are motivated by the following works \cite{Gualdani2017,Briant2016,Briant2017,Briant2016DCDS,Cao2022JFA}.

\begin{lemma} \label{poKnsest} \cite{Briant2017}
Let $w(v)=\langle v\rangle^\beta$ with $\beta \ge 0 $ and $\delta \in (0,1)$. Then there exists $C_{K_{ns}}>0$ such that for all f in $L^{\infty}_{x,v}(w)$, 
\[ \|K_{ns}^{(\delta)}(f)\|_{L^{\infty}_{x,v}(\langle v \rangle^\beta \mu^{-1/2})} \le C_{K_{ns}} \|f\|_{L^{\infty}_{x,v}(w)}. \]
The constant $C_{K_{ns}}$ only depends on $\delta$.
\end{lemma}
\begin{proof}
	Thanks to Carleman representation, the operator $K_{ns}^{(\delta)}$ can be written as
	\begin{align*}
		K_{ns}^{(\delta)}(f)(x,v) = \int_{\mathbb{R}^3}k_{ns}^{(\delta)}(v,v_*)f(x,v_*)dv_*.
	\end{align*}
	Moreover, by definition of $K_{ns}^{(\delta)}$, the kernel $k_{ns}^{(\delta)}$ are of compact support which implies the desired estimate. 
\end{proof}

\bigskip

\begin{lemma} \label{poQ-estimate} 
	Let \(w(v)=\langle v\rangle^\beta\) with $\beta > 3+\gamma$. Then for all f, g in $L^{\infty}_{x,v}(w)$,
\[\|Q^-(f,g)\|_{L^{\infty}_{x,v}(\nu^{-1} w)} \le C_{Q^-} \|f\|_{L^{\infty}_{x,v}(w)}\|g\|_{L^{\infty}_{x,v}(w)}. \]
\end{lemma}
\begin{proof}
	We have
	\begin{align*}
		\left| \nu^{-1}(v)w(v)Q^-(f,g)\right| &\le \|f\|_{L^{\infty}_{x,v}(w)}\|g\|_{L^{\infty}_{x,v}(w)}\nu^{-1}(v)\left(4\pi C_b\int_{\mathbb{R}^3}\frac{|v-v_*|^{\gamma}}{w(v_*)}dv_*\right)\\
		&\le C\left(4\pi C_b\int_{\mathbb{R}^3}\frac{(1+|v_*|)^{\gamma}}{w(v_*)}dv_*\right) \|f\|_{L^{\infty}_{x,v}(w)}\|g\|_{L^{\infty}_{x,v}(w)},
	\end{align*}
where we have used the inequality $|v-v_*|^{\gamma} \le C\nu(v)(1+|v_*|)^{\gamma}$ and the integration in the last term is finite. 
\end{proof}
\bigskip
\begin{lemma} \cite{Cao2022JFA} \label{pogaintermestimate}
	For any $-3< \gamma \le 2$, for any constant $\alpha >\max\{3,3+\gamma\}$, we have
	\begin{align*}
		\int_{\R^3} \int_{\S^2} |v-v_*|^\gamma \frac{\langle v\rangle^\alpha}{\langle v'\rangle^\alpha \langle v_*'\rangle^\alpha} d\sigma dv_* \le C_\alpha \langle v\rangle^\gamma
	\end{align*}	
	for all $v \in \R^3$. Moreover, for $0\le \gamma \le 1$, we have a stronger estimate
	\begin{align*}
		\int_{\R^3} \int_{\S^2} |v-v_*|^\gamma \frac{\langle v\rangle^\alpha}{\langle v'\rangle^\alpha \langle v_*'\rangle^\alpha} d\sigma dv_* \le  \frac{C}{\alpha} \langle v\rangle^\gamma+ C_\alpha \frac{1}{\langle v\rangle^{2-\gamma}}.
	\end{align*}
\end{lemma}

\bigskip

\begin{lemma} \label{poQ+estimate}
Let \(w(v)= \langle v \rangle ^\beta\) with $\beta>3+\gamma$. Then for all f, g in $L^{\infty}_{x,v}(w)$,
\[ \|Q^+(f,g)\|_{L^{\infty}_{x,v}(\nu^{-1} w)} \le C_\beta \|f\|_{L^{\infty}_{x,v}(w)}\|g\|_{L^{\infty}_{x,v}(w)}. \]
\end{lemma}
\begin{proof}
	By Lemma \ref{pogaintermestimate}, we easily obtain
	\begin{align*}
		|\nu^{-1}(v)w(v)Q^+(f,g)| &\le \nu^{-1}(v)\|f\|_{L^{\infty}_{x,v}(w)}\|g\|_{L^{\infty}_{x,v}(w)}\int_{\R^3}\int_{\S^2}b(\cos\theta)|v-v_*|^\gamma \frac{\langle v\rangle^\beta}{\langle v'\rangle^\beta\langle v_*'\rangle^\beta} d\sigma dv_*\\
		& \le C_\beta \|f\|_{L^{\infty}_{x,v}(w)}\|g\|_{L^{\infty}_{x,v}(w)}.
	\end{align*}
	This completes the proof of the lemma.
\end{proof}

\bigskip

\begin{corollary} \label{poQestimate} 
	Let \(w(v)=\langle v\rangle^{\beta}\) with $\beta >3+\gamma$. Then for all f, g in $L^{\infty}_{x,v}(w)$,
\[\|Q(f,g)\|_{L^{\infty}_{x,v}(\nu^{-1} w)} \le C_{Q} \|f\|_{L^{\infty}_{x,v}(w)}\|g\|_{L^{\infty}_{x,v}(w)}. \]
\end{corollary}

\bigskip

\begin{lemma} \label{poKsest} 
Let \(w(v)= \langle v \rangle ^\beta\) with $\beta>3+\gamma$. Then for all f in $L^{\infty}_{x,v}(w)$,
\[ \|K_s^{(\delta)}(f)\|_{L^{\infty}_{x,v}(\nu^{-1} w)} \le C_K(\delta) \|f\|_{L^{\infty}_{x,v}(w)}, \]
where
\[C_K(\delta) \rightarrow \frac{4 C_b}{c_b} \frac{1}{\beta-1-\gamma}\quad as \quad \delta \rightarrow 0.\]
\end{lemma}

\begin{proof}
We split $K_s^{\delta}$ into three pieces as follow:
\begin{align} \label{poKnest11}
	\left|K_s^{(\delta)}(f)\right| \le & \int_{\mathbb{R}^3\times \mathbb{S}^2} B(v,v_*,\theta)\mathbf{1} _{|v|\ge R} \left[\mu(v_*')|f(v')|+\mu(v')|f(v_*')|\right]d\sigma dv_*\nonumber\\
	&+C_b\int_{\mathbb{R}^3\times \mathbb{S}^2}{1} _{|v|\le R}(1-\Theta_{\delta})|v-v_*|^{\gamma} \left[\mu(v_*')|f(v')|+\mu(v')|f(v_*')|\right]d\sigma dv_*\nonumber\\
	&+C_b\int_{\mathbb{R}^3\times \mathbb{S}^2}(1-\Theta_{\delta})|v-v_*|^{\gamma} \mu(v)|f(v_*)|d\sigma dv_*.
\end{align}
First of all, let us consider last two terms in \eqref{poKnest11}.
We use $|v_*-v|^{\gamma} \le C \nu(v) \langle v_* \rangle^{\gamma}$ for the last one in \eqref{poKnest11}:
\begin{align*}
	w&(v)\nu(v)^{-1}\int_{\mathbb{R}^3\times \mathbb{S}^2}(1-\Theta_{\delta})|v-v_*|^{\gamma} \mu(v)|f(v_*)|d\sigma dv_*\\
	\le &Cw(v)\mu(v) \mathbf{1}_{|v| \ge R}\left(\int_{\mathbb{R}^3\times \mathbb{S}^2}\frac{\langle v_* \rangle^{\gamma}}{w(v_*)}d\sigma dv_*\right) \|f\|_{L^{\infty}_{x,v}(w)}\\
	&+Cw(R)\left(\int_{\mathbb{R}^3\times \mathbb{S}^2}\mathbf{1}_{|v| \le R}(1-\Theta_{\delta})\frac{\langle v_* \rangle^{\gamma}}{w(v_*)}d\sigma dv_*\right) \|f\|_{L^{\infty}_{x,v}(w)}.
\end{align*}
Note that
\begin{align*}
	|v_*-v|^{\gamma}\frac{\mu(v')}{w(v_*')} \le C\frac{\nu(v)\langle v_*\rangle^{\gamma}}{w(v')w(v_*')} \le C\frac{\nu(v)\langle v_*\rangle^{\gamma}}{w\left(\frac{|v_*|}{\sqrt{2}}\right)}.
\end{align*}
Here, $|v'| \ge \frac{1}{\sqrt{2}}|v_*|$ or $|v_*'| \ge \frac{1}{\sqrt{2}}|v_*|$ since $|v_*'|^2+|v'|^2 \ge |v_*|^2$, and $\frac{1}{w(v')} \le \frac{1}{w\left(\frac{|v_*|}{\sqrt{2}}\right)}$ or $\frac{1}{w(v_*')} \le \frac{1}{w\left(\frac{|v_*|}{\sqrt{2}}\right)}$.\\
As interchanging $v_*'$ and $v'$ to the second one in \eqref{poKnest11},
\begin{align*}
	w&(v)\nu(v)^{-1}\int_{\mathbb{R}^3\times \mathbb{S}^2}\mathbf{1}_{|v| \le R}(1-\Theta_{\delta})|v-v_*|^{\gamma} \left[\mu(v_*')|f(v')|+\mu(v')|f(v_*')|\right]d\sigma dv_*\\
	\le &Cw(R)\left(\int_{\mathbb{R}^3\times \mathbb{S}^2}\mathbf{1}_{|v| \le R}(1-\Theta_{\delta})\frac{\langle v_* \rangle^{\gamma}}{w\left(\frac{|v_*|}{\sqrt{2}}\right)}d\sigma dv_*\right) \|f\|_{L^{\infty}_{x,v}(w)}.
\end{align*}
For $\frac{1}{\delta} \ge 2R$, it holds that
\begin{align*}
	\mathbf{1}_{|v| \le R}(1-\Theta_{\delta}) &\le \mathbf{1}_{|v| \le R}\mathbf{1}_{|v-v_*| \ge \frac{1}{\delta}}+\mathbf{1}_{|v_*-v| \le 2\delta}+\mathbf{1}_{|\cos\theta| \ge 1-2\delta}\\
	&\le \mathbf{1}_{|v_*| \ge \frac{1}{2\delta}}+\mathbf{1}_{|v_*-v| \le 2\delta}+\mathbf{1}_{|\cos\theta| \ge 1-2\delta}.
\end{align*}
Then we have
\begin{align*}
	w(v)\mu(v)\mathbf{1}_{|v|\ge R} \rightarrow 0 \quad \text{as} \quad R \rightarrow \infty.
\end{align*}
and
\begin{align*}
	\left(\int_{\mathbb{R}^3\times \mathbb{S}^2}\mathbf{1}_{|v| \le R}(1-\Theta_{\delta})\frac{\langle v_* \rangle^{\gamma}}{w\left(\frac{|v_*|}{\sqrt{2}}\right)}d\sigma dv_*\right) \rightarrow 0 \quad \text{as} \quad \delta \rightarrow 0 \quad \text{for fixed }R
\end{align*}
since $\frac{\langle v_* \rangle^{\gamma}}{w\left(\frac{|v_*|}{\sqrt{2}}\right)}$ is integrable on $\mathbb{R}^3$. Thus we can write
\begin{align*}
	w(v)\nu(v)^{-1}\left|K_s^{(\delta)}(f)\right| \le &  w(v)\nu(v)^{-1}\int_{\mathbb{R}^3\times \mathbb{S}^2} B(v,v_*,\theta)\mathbf{1} _{|v|\ge R} \left[\mu(v_*')|f(v')|+\mu(v')|f(v_*')|\right]d\sigma dv_*\\
	&+\left(\varepsilon_w(R)+w(R)\varepsilon_w(\delta) \right)\|f\|_{L^{\infty}_{x,v}(w)}.
\end{align*}
By the Carleman representation \cite{Carleman1957} and rotation, we obtain
\begin{align*}
	I&:=\int_{\mathbb{R}^3\times \mathbb{S}^2} B(v,v_*,\theta)\mathbf{1} _{|v|\ge R} \left[\mu(v_*')w^{-1}(v')+\mu(v')w^{-1}(v_*')\right]d\sigma dv_*\\
	&	 \le 8C_b\int_{\R^3} \mathbf{1} _{|v|\ge R}\frac{\mu(v')}{|v-v'|} \int_{\{w:w\cdot(v'-v)=0\}} (|v-v'|^2+|w|^2)^{\frac{\gamma-1}{2}}\langle v+w\rangle^{-\beta}dwdv',
\end{align*}
where $R>1$. Fix $0<\varepsilon<1$. We split the integral domain into two parts $|v'|\ge \varepsilon |v|$ and $|v'|\le \varepsilon|v|$.\\
On $|v'|\le \varepsilon |v|$, we have $|v-v'|\ge (1-\varepsilon)|v|$ and $|v-v'|^2+|w|^2 \ge (1-\varepsilon)^2 (|v|^2+|w|^2)$.
 Note that on $\{w:w\cdot(v'-v)=0\}$,
 \begin{align*}
 	\langle v+w\rangle^2 &= 1+|v|^2+|w|^2+2v\cdot w\\
 	& = 1+|v|^2+|w|^2+2v'\cdot w\\
 	& \ge 1+|v|^2+|w|^2 - 2\varepsilon|v||w|\\
 	& \ge (1-\varepsilon)(|v|^2+|w|^2).
 \end{align*} 
 Then we can estimate
\begin{align} \label{poKsIfirstpart}
	I \mathbf{1}_{|v'|\le \varepsilon |v|} &\le \mathbf{1} _{|v|\ge R}\frac{8C_b(1-\varepsilon)^{-2+\gamma-\frac{\beta}{2}}}{|v|} \int_{|v'|\le \varepsilon |v|} \mu(v')  \int_{\R^2} (|v|^2+|w|^2)^{\frac{-\beta-1+\gamma}{2}}dwdv'\nonumber\\
	& \le \mathbf{1} _{|v|\ge R}\frac{16\pi C_b (1-\varepsilon)^{-2+\gamma-\frac{\beta}{2}} }{\beta-1-\gamma} |v|^{\gamma-\beta}.
\end{align}
On $|v'|\ge \varepsilon |v|$, we have $(|v-v'|^2+|w|^2)^\frac{\gamma-1}{2} \le |v-v'|^{\gamma-1}$. This gives
\begin{align*}
	I \mathbf{1}_{|v'|\ge \varepsilon |v|} &\le C\mathbf{1}_{|v|\ge R}\int_{|v'|\ge \varepsilon |v|}\mu(v') |v-v'|^{\gamma-2}\int_{\{w:w\cdot(v'-v)=0\}} \langle v+w\rangle^{-\beta}dwdv'\\
	& \le C_\beta \mathbf{1}_{|v|\ge R}\int_{|v'|\ge \varepsilon |v|}\mu(v') |v-v'|^{\gamma-2} dv'\\
	& \le C_\beta \mathbf{1}_{|v|\ge R} e^{-\frac{\varepsilon^2|v|^2}{4}} \int_{\R^3}e^{-\frac{|v'|^2}{4}} |v-v'|^{\gamma-2} dv'.
\end{align*}
Here, we can compute
\begin{align*}
	\int_{\R^3}e^{-\frac{|v'|^2}{4}} |v-v'|^{\gamma-2} dv' &\le \int_{|v-v'|\le 1} |v-v'|^{\gamma-2}dv' + \int_{|v-v'|\ge 1} e^{-\frac{|v'|^2}{4}}dv'\\
	& \le 4\pi\int_0^1 r^{\gamma} + \int_{\R^3}e^{-\frac{|v'|^2}{4}}dv'\\
	& = \frac{4\pi}{\gamma+1}+(4\pi)^{3/2}.
\end{align*}
It follows that
\begin{align} \label{poKsIsecondpart}
	I \mathbf{1}_{|v'|\ge \varepsilon |v|} \le C_\beta \mathbf{1}_{|v|\ge R} e^{-\frac{\varepsilon^2|v|^2}{4}} .
\end{align}
Combining \eqref{poKsIfirstpart} and \eqref{poKsIsecondpart} and using the fact that $\frac{1}{\nu(v)} \le \frac{1}{ 4\pi c_b(|v|^{\gamma}-A)}$ for some constant $A$, we can derive
\begin{align*}
	w(v)\nu(v)^{-1}I &\le w(v)\nu(v)^{-1}\mathbf{1} _{|v|\ge R}\left(\frac{16\pi C_b (1-\varepsilon)^{-2+\gamma-\frac{\beta}{2}} }{\beta-1-\gamma} |v|^{\gamma-\beta}+C_\beta \mathbf{1}_{|v|\ge R} e^{-\frac{\varepsilon^2|v|^2}{4}}\right)\\
	& \le \mathbf{1} _{|v|\ge R}\left(\frac{16\pi C_b}{4\pi c_b(\beta-1-\gamma)} (1-\varepsilon)^{-2+\gamma-\frac{\beta}{2}} \frac{w(v)|v|^\gamma }{(|v|^\gamma-A)|v|^\beta} + C_\beta w(v)\nu(v)^{-1} e^{-\frac{\varepsilon^2|v|^2}{4}}\right).
\end{align*}
Taking $\varepsilon \rightarrow 0$, we have
\begin{align*}
	\limsup_{R\rightarrow \infty} \sup_{|v|\ge R}\left[w(v)\nu(v)^{-1}I \right]= \frac{4 C_b}{ c_b(\beta-1-\gamma)} 
\end{align*}
Therefore, we obtain
\begin{align*}
    \|K_s^{(\delta)}(f)\|_{L^{\infty}_{x,v}(\nu^{-1} w)} \le C_K(\delta) \|f\|_{L^{\infty}_{x,v}(w)},
\end{align*}
where $C_K(\delta) \rightarrow \frac{4 C_b}{ c_b} \frac{1}{\beta-1-\gamma}$ as $\delta \rightarrow 0.$\\
\end{proof}

Now, we construct a solution of the full nonlinear equation.
\begin{proposition} \label{poT2}
Let \(w(v)= \langle v \rangle ^\beta\) with $\beta>1+\gamma+\frac{4C_b}{c_b}$. Let \(f_0\) be in \(L^{\infty}_{x,v}(w)\) and \(g(t,x,v)\) in \(L^{\infty}_{x,v}(w)\). Then there exists \(\delta_w>0\) such that for any \(\delta\) in \((0,\delta_w]\), there exist $C_{s}$, $\eta_1$, and $\lambda_1>0$ such that if                  
\[\|f_0\|_{L^{\infty}_{x,v}(w)}\le \eta_1^{1+p} \quad \text{and} \quad \|g(t)\|_{L^{\infty}_{x,v}(w)}\le \eta_1 e^{-\lambda_1t},\]
then there exists a solution $f_1$ to the equation
\[\partial_t f_1+v\cdot \nabla_x f_1 +\nu(v)f_1=K_s^{(\delta)}f_1+Q(f_1+g,f_1+g)\] 
with \(f_0\)  as initial data.
Moreover, this solution satisfies
\[\|f_1(t)\|_{L_{x,v}^\infty(w)}\le C_{s}e^{-\lambda_1t}(\|f_0\|_{L_{x,v}^\infty(w)}+\eta_1^{1+p})\]
for all $t \ge 0$. The constant $C_{s}$, $\eta_1$, and $\lambda_1$ are constructive and only depend on w, $\delta$.
\end{proposition}

\begin{proof}
    By Duhamel's principle,
    \begin{align*}
    	f_1(t,x,v) &= e^{-\nu(v)t}f_0(x-vt,v) + \int_0^t e^{-\nu(v)(t-s)}K_s^{(\delta)}f_1(s,x-v(t-s),v)ds\\
    	& \quad + \int_0^t e^{-\nu(v)(t-s)}Q(f_1+g,f_1+g)(s,x-v(t-s),v)ds.
    \end{align*}
    Set $h_0 = 0$ and consider an iterative system as follows:
    \begin{align*}
    h_{l+1}(t,x,v) &= e^{-\nu(v)t}f_0(x-vt,v) + \int_0^t e^{-\nu(v)(t-s)}K_s^{(\delta)}h_{l+1}(s,x-v(t-s),v)ds \\
    & \quad +\int_0^t e^{-\nu(v)(t-s)} Q(h_l+g,h_l+g)(s,x-v(t-s),v)ds.
    \end{align*}
    with initial data $h_{l+1}(0,x,v) = f_0(x,v)$. Then
    \begin{align*}
    	&\left|h_{l+1}(t,x,v)-h_l(t,x,v)\right|\\ &\le \int_0^t e^{-\nu(v)(t-s)} \left|K_s^{(\delta)}(h_{l+1}-h_l)(s,x-v(t-s),v)\right|ds\\
    	&\quad +\int_0^t e^{-\nu(v)(t-s)} \Big|\left[Q(h_l+g,h_l+g)-Q(h_{l-1}+g,h_{l-1}+g)\right](s,x-v(t-s),v)\Big|ds\\
    	&= \int_0^t e^{-\nu(v)(t-s)} \left|K_s^{(\delta)}(h_{l+1}-h_l)(s,x-v(t-s),v)\right|ds\\
    	&\quad +\int_0^t e^{-\nu(v)(t-s)} \Big|\left[Q(h_l-h_{l-1},h_l+g)+Q(h_{l-1}+g,h_l-h_{l-1})\right](s,x-v(t-s),v)\Big|ds
    \end{align*}
   Fix $0<\varepsilon <1/2$. First of all, we consider the bound of the integration for $K_s^{(\delta)}$.
      By Lemma \ref{poKsest}, for $\varepsilon$ in $(0,1)$,
    \begin{align}
    	 &w(v)\int_0^t e^{-\nu(v)(t-s)} \left|K_s^{(\delta)}(h_{l+1}-h_l)(s,x-v(t-s),v)\right|ds\nonumber \\
    	 &\le \int_0^t \nu(v)e^{-\nu(v)(t-s)} \left\|K_s^{(\delta)}(h_{l+1}-h_l)(s)\right\|_{L^{\infty}_{x,v}(\nu^{-1}w)}ds\nonumber\\
    	 &\le C_K(\delta)\int_0^t \nu(v)e^{-\nu(v)(t-s)} \left\|(h_{l+1}-h_l)(s)\right\|_{L^{\infty}_{x,v}(w)}ds\nonumber\\
    	 &\le C_K(\delta)e^{-\frac{\varepsilon\nu_0}{2}t}\sup_{0\le s \le t} \left[e^{\frac{\varepsilon\nu_0}{2}s}\left\|(h_{l+1}-h_l)(s)\right\|_{L^{\infty}_{x,v}(w)}\right]\int_0^t \nu(v)e^{-(1-\varepsilon)\nu(v)(t-s)}ds\nonumber\\ \label{TI1}
    	 &\le \frac{C_K(\delta)}{1-\varepsilon}e^{-\frac{\varepsilon\nu_0}{2}t}\sup_{0\le s \le t} \left[e^{\frac{\varepsilon\nu_0}{2}s}\left\|(h_{l+1}-h_l)(s)\right\|_{L^{\infty}_{x,v}(w)}\right].
    \end{align} 
    Next, we consider the bound for the integration for the nonlinear term. By Corollary \ref{poQestimate}, for $\varepsilon$ in $(0,1)$,
   \begin{align}
   		&w(v)\int_0^t e^{-\nu(v)(t-s)} \Big|\left[Q(h_l-h_{l-1},h_l+g)+Q(h_{l-1}+g,h_l-h_{l-1})\right](s,x-v(t-s),v)\Big|ds \nonumber\\
   		&\le \int_0^t \nu(v)e^{-\nu(v)(t-s)} \left\|\left[Q(h_l-h_{l-1},h_l+g)+Q(h_{l-1}+g,h_l-h_{l-1})\right](s)\right\|_{L^{\infty}_{x,v}(\nu^{-1}w)}ds\nonumber\\
   		&\le e^{-\frac{\varepsilon\nu_0}{2}t}\sup_{0\le s \le t} \left[e^{\frac{\varepsilon\nu_0}{2}s}\left\|\left[Q(h_l-h_{l-1},h_l+g)+Q(h_{l-1}+g,h_l-h_{l-1})\right](s)\right\|_{L^{\infty}_{x,v}(\nu^{-1}w)}\right]\int_0^t \nu(v)e^{-(1-\varepsilon)\nu(v)(t-s)}ds\nonumber\\
   		&\le \frac{1}{1-\varepsilon} e^{-\frac{\varepsilon\nu_0}{2}t} \sup_{0\le s \le t} \left[e^{\frac{\varepsilon\nu_0}{2}s}\left\|\left[Q(h_l-h_{l-1},h_l+g)+Q(h_{l-1}+g,h_l-h_{l-1})\right](s)\right\|_{L^{\infty}_{x,v}(\nu^{-1}w)}\right]\nonumber\\
   		&\le \frac{C}{1-\varepsilon}\left(\|h_l\|_{L_{[0,t]}^{\infty}L^{\infty}_{x,v}(w)}+\|h_{l-1}\|_{L_{[0,t]}^{\infty}L^{\infty}_{x,v}(w)}+\|g\|_{L_t^{\infty}L^{\infty}_{x,v}(w)}\right)\nonumber\\ \label{TI2}
   		&\quad \times e^{-\frac{\varepsilon\nu_0}{2}t} \sup_{0\le s \le t} \left[e^{\frac{\varepsilon\nu_0}{2}s}\left\|(h_l-h_{l-1})(s)\right\|_{L^{\infty}_{x,v}(w)}\right].
   		\end{align}
   Combining \eqref{TI1} and \eqref{TI2}, for $T_0 >0$ and $t\in [0,T_0]$,
   \begin{align*}
   	\sup_{0\le s \le t} \left[e^{\frac{\varepsilon\nu_0}{2}s}\left\|(h_{l+1}-h_{l})(s)\right\|_{L^{\infty}_{x,v}(w)}\right] &\le \frac{C_K(\delta)}{1-\varepsilon}\sup_{0\le s \le t} \left[e^{\frac{\varepsilon\nu_0}{2}s}\left\|(h_{l+1}-h_l)(s)\right\|_{L^{\infty}_{x,v}(w)}\right]\\
   	&\quad + \frac{C}{1-\varepsilon}\left(\|h_l\|_{L_{[0,t]}^{\infty}L^{\infty}_{x,v}(w)}+\|h_{l-1}\|_{L_{[0,t]}^{\infty}L^{\infty}_{x,v}(w)}+\|g\|_{L_t^{\infty}L^{\infty}_{x,v}(w)}\right)\\
   	&\quad \quad \times \sup_{0\le s \le t} \left[e^{\frac{\varepsilon\nu_0}{2}s}\left\|(h_l-h_{l-1})(s)\right\|_{L^{\infty}_{x,v}(w)}\right].
   \end{align*}
   We first choose $\delta_w>0$ such that $0<C_K(\delta_w)<1$ because $C_K(\delta) \rightarrow \frac{4 C_b}{c_b} \frac{1}{\beta-1-\gamma} < 1$ as $\delta \rightarrow 0$, and we take $\varepsilon_0>0$ such that $\frac{C_K(\delta_w)}{1-\varepsilon_0} < 1$. We can deduce 
   \begin{align*}
   	\sup_{0\le s \le t} \left[e^{\frac{\varepsilon\nu_0}{2}s}\left\|(h_{l+1}-h_{l})(s)\right\|_{L^{\infty}_{x,v}(w)}\right] &\le C_w\left(\|h_l\|_{L_{[0,t]}^{\infty}L^{\infty}_{x,v}(w)}+\|h_{l-1}\|_{L_{[0,t]}^{\infty}L^{\infty}_{x,v}(w)}+\|g\|_{L_t^{\infty}L^{\infty}_{x,v}(w)}\right)\\
   	&\quad\times \sup_{0\le s \le t} \left[e^{\frac{\varepsilon\nu_0}{2}s}\left\|(h_l-h_{l-1})(s)\right\|_{L^{\infty}_{x,v}(w)}\right].
   \end{align*}
   for all $t \in [0,T_0]$. Since $C_w$ is independent of $T_0$, we can extend the inequality to all time $t\ge 0$.\\
   Now, we claim that $\left\|h_{l+1}\right\|_{L^\infty_{[0,t]}  L^{\infty}_{x,v}(w)}$ is uniformly bounded.\\
   In a similar argument to the above, we can easily deduce that 
   \begin{align} \label{TI3}
   	w(v)\int_0^t e^{-\nu(v)(t-s)} \left|K_s^{(\delta)}h_{l+1}(s,x-v(t-s),v)\right|ds \le \frac{C_K(\delta)}{1-\varepsilon}e^{-\frac{\varepsilon\nu_0}{2}t}\sup_{0\le s \le t} \left[e^{\frac{\varepsilon\nu_0}{2}s}\left\|(h_{l+1})(s)\right\|_{L^{\infty}_{x,v}(w)}\right].
   \end{align}
   We also obtain 
   \begin{align}
   	& w(v)\int_0^t e^{-\nu(v)(t-s)}\left|Q(h_l+g,h_l+g)(s,x-v(t-s),v)\right|ds \nonumber\\
   	&\le \int_0^t \nu(v)e^{-\nu(v)(t-s)} \left\|Q(h_l+g,h_l+g)(s)\right\|_{L^{\infty}_{x,v}(\nu^{-1}w)}ds\nonumber\\
   	&\le Ce^{-\frac{\varepsilon\nu_0}{2}t}\left(\int_0^t \nu(v)e^{-(1-\varepsilon)\nu(v)(t-s)}ds\right)\nonumber\\
   	&\quad \times \bigg\{\left(\|h_l\|_{L_{[0,t]}^{\infty}L^{\infty}_{x,v}(w)}+\|g\|_{L_{t}^{\infty}L^{\infty}_{x,v}(w)}\right)\sup_{0\le s \le t} \left[e^{\frac{\varepsilon\nu_0}{2}s}\left\|(h_l)(s)\right\|_{L^{\infty}_{x,v}(w)}\right]\nonumber\\
   	&\quad \quad \quad + \|g\|_{L_{t}^{\infty}L^{\infty}_{x,v}(w)}\sup_{0\le s \le t} \left[e^{\frac{\varepsilon\nu_0}{2}s}\left\|g(s)\right\|_{L^{\infty}_{x,v}(w)}\right]\bigg\}\nonumber\\
   	&\le \frac{C}{1-\varepsilon}e^{-\frac{\varepsilon\nu_0}{2}t}\bigg\{\left(\|h_l\|_{L_{[0,t]}^{\infty}L^{\infty}_{x,v}(w)}+\|g\|_{L_{t}^{\infty}L^{\infty}_{x,v}(w)}\right)\sup_{0\le s \le t} \left[e^{\frac{\varepsilon\nu_0}{2}s}\left\|(h_l)(s)\right\|_{L^{\infty}_{x,v}(w)}\right]\nonumber\\ \label{TI4}
   	&\quad + \|g\|_{L_{t}^{\infty}L^{\infty}_{x,v}(w)}\sup_{0\le s \le t} \left[e^{\frac{\varepsilon\nu_0}{2}s}\left\|g(s)\right\|_{L^{\infty}_{x,v}(w)}\right]\bigg\}.
   \end{align}
   Take $\lambda_1 = \frac{\varepsilon}{2} \nu_0$ and assume $\|g(t)\|_{L^{\infty}_{x,v}(w)}\le \eta_1 e^{-\lambda_1t}$. Then we get
   \begin{align*}
   	&w(v)\int_0^t e^{-\nu(v)(t-s)}\left|Q(h_l+g,h_l+g)(s,x-v(t-s),v)\right|ds\\
   	&\le \frac{C}{1-\varepsilon}e^{-\frac{\varepsilon\nu_0}{2}t}\bigg\{\left(\|h_l\|_{L_{[0,t]}^{\infty}L^{\infty}_{x,v}(w)}+\eta_1\right)\sup_{0\le s \le t} \left[e^{\frac{\varepsilon\nu_0}{2}s}\left\|(h_l)(s)\right\|_{L^{\infty}_{x,v}(w)}\right] + \eta_1^2\bigg\}.
   \end{align*}
   Write $A_l = \sup_{0\le s \le t} \left[e^{\frac{\varepsilon}{2}\nu_0s}\left\|(h_l)(s)\right\|_{L^{\infty}_{x,v}(w)}\right]$. By \eqref{TI3} and \eqref{TI4},  there exist $C_0^{(1)}$, $C_0^{(2)}>0$ such that
   \begin{align*}
   	A_{l+1} \le C_0^{(1)}\|f_0\|_{L^{\infty}_{x,v}(w)} + C_0^{(2)}\eta_1^2+C_0^{(2)}\left(\|h_l\|_{L_{[0,t]}^{\infty}L^{\infty}_{x,v}(w)}+\eta_1\right)A_l.
   \end{align*}
   By using induction, we show that 
   \begin{align*}
   	\sup_l A_l \le \left(1+C_0^{(1)}\right)\|f_0\|_{L^{\infty}_{x,v}(w)} + \left(1+C_0^{(1)}\right)\eta_1^{1+p}.
   \end{align*}
  We deduce
   \begin{align*}
   	A_{l+1} &\le C_0^{(1)}\|f_0\|_{L^{\infty}_{x,v}(w)} + C_0^{(2)}\eta_1^2\\
   	&\quad + C_0^{(2)}\left\{\left(2+C_0^{(1)}\right)\eta_1+\left(1+C_0^{(1)}\right)\eta_1^{1+p}\right\}\times\left\{\left(1+C_0^{(1)}\right)\|f_0\|_{L^{\infty}_{x,v}(w)}+\left(1+C_0^{(1)}\right)\eta_1^{1+p}\right\}\\
   	&\le C_0^{(1)}\|f_0\|_{L^{\infty}_{x,v}(w)} + C_0^{(2)}\eta_1^2\\
   	&\quad + \|f_0\|_{L^{\infty}_{x,v}(w)} + C_0^{(2)}\left\{\left(2+C_0^{(1)}\right)\eta_1+\left(1+C_0^{(1)}\right)\eta_1^{1+p}\right\}\left(1+C_0^{(1)}\right)\eta_1^{1+p}\\
   	&\le \left(1+C_0^{(1)}\right)\|f_0\|_{L^{\infty}_{x,v}(w)} + \left(1+C_0^{(1)}\right)\eta_1^{1+p}
   \end{align*}
   by choosing sufficiently small $\eta_1 \ll 1$ such that
   \begin{align*}
   	& C_0^{(2)}\left\{\left(2+C_0^{(1)}\right)\eta_1+\left(1+C_0^{(1)}\right)\eta_1^{1+p}\right\}\left(1+C_0^{(1)}\right) \le 1, \\
   	& C_0^{(2)}\left\{\left(2+C_0^{(1)}\right)\eta_1+\left(1+C_0^{(1)}\right)\eta_1^{1+p}\right\}\left(1+C_0^{(1)}\right)\eta_1^{1+p}+ C_0^{(2)}\eta_1^2 \le \left(1+C_0^{(1)}\right)\eta_1^{1+p}.
   \end{align*}
   It follows that
   \begin{align*}
   	\sup_l A_l \le \left(1+C_0^{(1)}\right)\|f_0\|_{L^{\infty}_{x,v}(w)} + \left(1+C_0^{(1)}\right)\eta_1^{1+p} \le 2\left(1+C_0^{(1)}\right)\eta_1^{1+p}\le 2\left(1+C_0^{(1)}\right)\eta_1.
   \end{align*}
   Thus, we obtain the uniform bound of $\left\|h_{l+1}\right\|_{L^\infty_{[0,t]}  L^{\infty}_{x,v}(w)}$ as follows:
   \begin{align*}
   	\left\|h_{l+1}\right\|_{L^\infty_{[0,t]}  L^{\infty}_{x,v}(w)} \le 2\left(1+C_0^{(1)}\right)\eta_1.
   \end{align*}
  Hence we conclude that
   \begin{align*}
   		\sup_{0\le s \le t} \left[e^{\frac{\varepsilon}{2}\nu_0s}\left\|(h_{l+1}-h_{l})(s)\right\|_{L^{\infty}_{x,v}(w)}\right] &\le 6\left(1+C_0^{(1)}\right)\eta_1C \sup_{0\le s \le t} \left[e^{\frac{\varepsilon}{2}\nu_0s}\left\|(h_l-h_{l-1})(s)\right\|_{L^{\infty}_{x,v}(w)}\right].
   \end{align*}
   Choosing sufficiently small $\eta_1$ such that $6\left(1+C_0^{(1)}\right)\eta_1 C < 1$, $\left(h_l\right)_{l\in \mathbb{N}}$ is Cauchy in $L_{t}^\infty L_{x,v}^\infty(w)$, so that  $\left(h_l\right)_{l\in \mathbb{N}}$ converges to $f_1$ in $L_{t}^\infty L_{x,v}^\infty(w)$. Hence, $f_1$ is a solution to the equation
   \begin{align*}
   	\partial_t f_1+v\cdot \nabla_x f_1 +\nu(v)f_1= K_s^{(\delta)}f_1+Q(f_1+g,f_1+g)
   \end{align*} 
   and we get the following exponential decay
   \begin{align*}
   	\|f_1(t)\|_{L_{x,v}^\infty(w)}\le C_{s}e^{-\lambda_1t}(\|f_0\|_{L_{x,v}^\infty(w)}+\eta_1^{1+p}).
   \end{align*}
   for all $t \ge 0$.
\end{proof}

\bigskip

\begin{proposition} \label{poT3}
Let $w(v)=(1+|v|)^{\beta}$ with $\beta > 5$. Let $g=g(t,x,v)$ be in \(L^\infty_tL^{\infty}_{x,v}(w)\). Then there exists a unique function $f_2$ in \(L^\infty_{loc,t}L^{\infty}_{x,v}(\langle v \rangle^{\beta_3} \mu^{-1/2})\) with $\beta_3> \max\{3,\beta\} $ such that
	\[\partial_t f_2+v\cdot\nabla_x f_2 +Lf_2 = K_{ns}^{(\delta)}g \quad and \quad f_2(0,x,v)=0.\]
    If $P_{G_L}(f_2+g)=0$ and there exist $\lambda_g$, $\eta_g>0$ such that $\lambda_g <\lambda_L$, where $\lambda_L$ is defined in proposition \ref{linearLinftydecay}, and for all $t \ge 0$ 
	\[\|g(t)\|_{L_{x,v}^\infty(w)}\le \eta_ge^{-\lambda_g t},\]
	then there exists $C_{ns}>0$ such that
	\[\|f_2(t)\|_{L_{x,v}^\infty(\langle v \rangle^{\beta_3} \mu^{-1/2})}\le C_{ns}\eta_ge^{-\lambda_g t}\]
	for all $t \ge 0$. This guarantees $f_2 \in L^\infty_{t}L^{\infty}_{x,v}(\langle v \rangle^\beta \mu^{-1/2})$. The constant $C_{ns}$ only depends on $\lambda_L, \lambda_g$.
\end{proposition}

\begin{proof}
	From Lemma \ref{poKnsest} and $g \in L^\infty_tL^{\infty}_{x,v}(w)$, we obtain $K_{ns}^{(\delta)}(g) \in L^\infty_{t}L^{\infty}_{x,v}(\langle v \rangle^{\beta} \mu^{-1/2}).$
	By proposition \ref{linearLinftydecay}, there exists a unique function $f_2$ in \(L^\infty_{\text{loc},t}L^{\infty}_{x,v}(\langle v \rangle^{\beta_3} \mu^{-1/2})\) such that
	\[\partial_t f_2+v\cdot\nabla_x f_2 +Lf_2 = K_{ns}^{(\delta)}(g) \quad \text{with} \ f_2(0,x,v)=0,\]
	and it is given by
	\begin{align*}
		f_2(t,x,v) = \int_0^tG(t,s)K_{ns}^{(\delta)}(g)(s)ds.
	\end{align*}
	Assume $P_{G_L}(f_2+g)=0$ and there exist $\lambda_g$, $\eta_g>0$ such that $\|g(t)\|_{L_{x,v}^\infty(w)}\le \eta_ge^{-\lambda_g t}$ for all $t\ge 0$. Then for all $t\ge 0$, we have
	\begin{align} \label{poest32}
		\|P_{G_L}(f_2)(t)\|_{L_{x,v}^\infty(\langle v \rangle^{\beta_3} \mu^{-1/2})} &= \|P_{G_L}(g)(t)\|_{L_{x,v}^\infty(\langle v \rangle^{\beta_3} \mu^{-1/2})}\nonumber\\ 
		&\le C \|g(t)\|_{L_{x,v}^\infty(w)}\\
		&\le C\eta_ge^{-\lambda_g t}. \nonumber
	\end{align}
	We claim $\partial_t\left[P_{G_L}^\perp\left(f_2\right)\right] = G_L\left[P_{G_L}^\perp\left(f_2\right)\right]+P_{G_L}^\perp\left[K_{ns}^{(\delta)}(g)\right]$. Taking $P_{G_L}$ to both sides of $\partial_t f_2 = G_L(f_2) + K_{ns}^{(\delta)}(g)$, we obtain
	\begin{align} \label{poest31}
		\partial_t\left[P_{G_L} \left(f_2\right)\right] = P_{G_L} \left[ \partial_t f_2 \right] = P_{G_L} \left[ G_L(f_2) + K_{ns}^{(\delta)}(g) \right] = P_{G_L} \left[G_L(f_2)\right] + P_{G_L} \left[K_{ns}^{(\delta)}(g)\right].
	\end{align}
	Since $(L(f_2),\mu)_{L^2_{v}(\mu^{-1/2})}=0$, we have
	\begin{align*}
		P_{G_L} \left[G_L(f_2)\right] = 0.
	\end{align*}
	From $G_L\left[P_{G_L}\left(f_2\right)\right]=0$,  \eqref{poest31} becomes
	\begin{align*}
		\partial_t\left[P_{G_L} \left(f_2\right)\right] = P_{G_L} \left[K_{ns}^{(\delta)}(g)\right]= G_L\left[P_{G_L}\left(f_2\right)\right] +P_{G_L} \left[K_{ns}^{(\delta)}(g)\right],
	\end{align*}
	and it follows that 
	\begin{align*}
\partial_t\left[P_{G_L}^\perp\left(f_2\right)\right] = G_L\left[P_{G_L}^\perp\left(f_2\right)\right]+P_{G_L}^\perp\left[K_{ns}^{(\delta)}(g)\right].
	\end{align*}
	By Duhamel's principle, we get
	\begin{align*}
		P_{G_L}^\perp\left(f_2\right) = \int_0^t G(t,s)P_{G_L}^\perp\left(K_{ns}^{(\delta)}(g)(s)\right)ds.
	\end{align*}
	Thanks to proposition \ref{linearLinftydecay} and Lemma \ref{poKsest}, we deduce that
	\begin{align}\label{poest33}
		\|P_{G_L}^{\perp}(f_2)(t)\|_{L_{x,v}^\infty(\langle v \rangle^{\beta_3} \mu^{-1/2})} 
		&\le C_L \int_0^t e^{-\lambda_L(t-s)}\|K_{ns}^{(\delta)}(g)(s)\|_{L_{x,v}^\infty(\langle v \rangle^{\beta_3} \mu^{-1/2})}ds\nonumber\\ 
		&\le C \int_0^t e^{-\lambda_L(t-s)}\|g(s)\|_{L_{x,v}^\infty(w)}ds\\
		&\le C \int_0^t e^{-\lambda_L(t-s)}\eta_g e^{-\lambda_g s}ds \nonumber\\ & = C\eta_g\frac{1}{\lambda_L-\lambda_g}\left(e^{-\lambda_g t}-e^{-\lambda_L t}\right)\nonumber\\
		&\le C\eta_g e^{-\lambda_g t}.\nonumber
	\end{align}
	Hence we conclude that
	\begin{align*}
		\|f_2(t)\|_{L_{x,v}^\infty(\langle v \rangle^{\beta_3} \mu^{-1/2})}\le C\eta_ge^{-\lambda_g t} 
	\end{align*} 
	for all $t \ge 0$.
    Therefore, this guarantees $f_2 \in L^\infty_{t}L^{\infty}_{x,v}(\langle v \rangle^{\beta_3} \mu^{-1/2})$.
\end{proof}

\bigskip

\begin{proposition} \label{poT1}
Let \(w(v)=\langle v\rangle^\beta\) with $\beta>1+\gamma+\frac{4C_b}{c_b}$. Then there exists $\eta_0 >0$ such that for any $f_0$ in $L^\infty_{x,v}(w)$ with $P_{G_L}(f_0)=0$ and $\|f_0\|_{L^\infty_{x,v}(w)} \le \eta_0^{1+p}$ for some $0<p<1$, there exists at least one solution f of \eqref{BE-f} with initial data $f_0$. Moreover, f satisfies the mass, momentum, energy conservation and there exist $C_\infty,\lambda_\infty >0$ such that
    \begin{align*}
    	\|f(t)\|_{L_{x,v}^\infty(w)}\le C_\infty e^{-\lambda_\infty t}(\|f_0\|_{L_{x,v}^\infty(w)}+\eta_0^{1+p})
    \end{align*}
    for all $t \ge 0$.
\end{proposition}

\begin{proof}
	We will use the sequence argument $f_2^{(n)} \rightarrow f_1^{(n+1)} \rightarrow f_2^{(n+1)} \rightarrow \cdots$ and we will deduce the following estimates: 
	\begin{align*}
		&\sup_n\|f_1^{(n)}(t)\|_{L_{x,v}^\infty(w)} \le C_se^{-\lambda t}(\|f_0\|_{L_{x,v}^\infty(w)}+\eta_0^{1+p}),\\
		&\sup_n\|f_2^{(n)}(t)\|_{L_{x,v}^\infty(\langle v \rangle^{\beta_3} \mu^{-1/2})} \le C_sC_{ns}e^{-\lambda t}(\|f_0\|_{L_{x,v}^\infty(w)}+\eta_0^{1+p}).
	\end{align*} 
	Set $f_1^{(0)} = f_2^{(0)} = 0$ and we use the iterative system as follows:
	\begin{equation} \label{poT1d4}
	\begin{aligned}
		&\partial_t f_1^{(n+1)}+v\cdot\nabla_x f_1^{(n+1)} +\nu(v)f_1^{(n+1)}=K_s^{(\delta)}f_1^{(n+1)}+Q(f_1^{(n+1)}+f_2^{(n)},f_1^{(n+1)}+f_2^{(n)})\quad  \text{with}\quad f_1^{(n+1)}(0,x,v)=f_0,\\
		&\partial_t f_2^{(n+1)}+v\cdot\nabla_x f_2^{(n+1)} +Lf_2^{(n+1)} = K_{ns}^{(\delta)}f_1^{(n+1)}\quad \text{with}\quad f_2^{(n+1)}(0,x,v)=0.
	\end{aligned}
	\end{equation}
	Assume that $\|f_0\|_{L^\infty_{x,v}(w)} \le \eta_0^{1+p}$ for some $0<p<1$.\\
	First, we consider the step $f_2^{(n)} \rightarrow f_1^{(n+1)}$ and assume $\|f_2^{(n)}(t)\|_{L_{x,v}^\infty(\langle v \rangle^{\beta_3} \mu^{-1/2})} \le C_sC_{ns}e^{-\lambda t}(\|f_0\|_{L_{x,v}^\infty(w)}+\eta_0^{1+p})$ for all $t\ge 0$.\\
	In proposition \ref{poT2}, putting $g = f_2^{(n)}$ and $\tilde{C}= \sup_{x,v}\left(\langle v \rangle^{-{\beta_3}} \mu^{1/2} w \right)$, we get
	\begin{align*}
		\|g(t)\|_{L_{x,v}^\infty(w)} = \|f_2^{(n)}(t)\|_{L_{x,v}^\infty(w)} &\le \tilde{C}\|f_2^{(n)}(t)\|_{L_{x,v}^\infty(\langle v \rangle^{\beta_3} \mu^{-1/2})} \le \tilde{C}C_sC_{ns}e^{-\lambda t}(\|f_0\|_{L_{x,v}^\infty(w)}+\eta_0^{1+p})\\
		&\le \tilde{C}C_sC_{ns}e^{-\lambda t}(2\eta_0^{1+p}) \le \eta_0 e^{-\lambda t},
	\end{align*}
	where in the last term we choose enough small $\eta_0$ such that $\tilde{C}C_sC_{ns}(2\eta_0^{1+p}) \le \eta_0$.\\
	By proposition \ref{poT2}, we obtain
	\begin{align*}
		\|f_1^{(n+1)}(t)\|_{L_{x,v}^\infty(w)} \le C_se^{-\lambda t}(\|f_0\|_{L_{x,v}^\infty(w)}+\eta_0^{1+p}).
	\end{align*}
	for all $t \ge 0.$
	Next, we consider the step $f_1^{(n+1)} \rightarrow f_2^{(n+1)}$ and assume $\|f_1^{(n)}(t)\|_{L_{x,v}^\infty(w)} \le C_se^{-\lambda t}(\|f_0\|_{L_{x,v}^\infty(w)}+\eta_0^{1+p})$ for all $t \ge 0$.\\
	From $\partial_t \left(f_1^{(n+1)}+f_2^{(n+1)}\right)=G_L\left(f_1^{(n+1)}+f_2^{(n+1)}\right)+Q\left(f_1^{(n+1)}+f_2^{(n)},f_1^{(n+1)}+f_2^{(n)}\right)$, we can deduce $\partial_t P_{G_L}\left(f_1^{(n+1)}+f_2^{(n+1)}\right)=0$, which implies
	\begin{align} \label{posequenceconservation}
		P_{G_L}\left(f_1^{(n+1)}+f_2^{(n+1)}\right)=P_{G_L}\left(f_0\right)=0.
	\end{align}
	In proposition \ref{poT3}, taking $g=f_1^{(n+1)}$, $\eta_g = C_s\left(\|f_0\|_{L_{x,v}^\infty(w)}+\eta_0^{1+p}\right)$, and $\lambda_g = \lambda_1=\frac{\varepsilon}{2} \nu_0$, we get
	\begin{align*}
		\|f_2^{(n+1)}(t)\|_{L_{x,v}^\infty(\langle v \rangle^{\beta_3} \mu^{-1/2})} \le C_{ns}\eta_g e^{-\lambda_g t} = C_sC_{ns}e^{-\lambda t}\left(\|f_0\|_{L_{x,v}^\infty(w)}+\eta_0^{1+p}\right),
	\end{align*}
	where we choose enough small $\varepsilon >0$ such that $\lambda_g  = \frac{\varepsilon}{2} \nu_0 < \lambda_L$.\\
	Now, we consider the sequence $f_2^{(n+1)}-f^{(n)}_2$. This sequence satisfies 
	\begin{align*}
		\partial_t (f_2^{(n+1)}-f^{(n)}_2)+v\cdot\nabla_x (f_2^{(n+1)}-f^{(n)}_2) +L(f_2^{(n+1)}-f^{(n)}_2) = K_{ns}^{(\delta)}(f_1^{(n+1)}-f_1^{(n)}).
	\end{align*}
	From \eqref{posequenceconservation}, \eqref{poest32}, and \eqref{poest33}, we have
	\begin{align*}
		\|P_{G_L}(f_2^{(n+1)}-f^{(n)}_2)(t)\|_{L_{x,v}^\infty(\langle v \rangle^{\beta_3} \mu^{-1/2})} &\le C_\beta\sup_{t\ge 0}\|(f_1^{(n+1)}-f_1^{(n)})(t)\|_{L_{x,v}^\infty(w)},\\
		\|P_{G_L}^{\perp}(f_2^{(n+1)}-f^{(n)}_2)(t)\|_{L_{x,v}^\infty(\langle v \rangle^{\beta_3} \mu^{-1/2})} 
		 &\le C_\beta \int_0^t e^{-\lambda_L(t-s)}\|(f_1^{(n+1)}-f_1^{(n)})(s)\|_{L_{x,v}^\infty(w)}ds\\
		 & \le C_{\beta,\lambda_L} \sup_{t\ge 0}\|(f_1^{(n+1)}-f_1^{(n)})(t)\|_{L_{x,v}^\infty(w)},
	\end{align*}
	which implies that 
	\begin{align} \label{pof2Cauchy}
		\sup_{t\ge 0 }\|(f_2^{(n+1)}-f^{(n)}_2)(t)\|_{L_{x,v}^\infty(w)}&\le \sup_{t\ge 0 }\|(f_2^{(n+1)}-f^{(n)}_2)(t)\|_{L_{x,v}^\infty(\langle v \rangle^{\beta_3} \mu^{-1/2})}\nonumber\\&\le C_{\beta,\lambda_L} \sup_{t\ge 0}\|(f_1^{(n+1)}-f_1^{(n)})(t)\|_{L_{x,v}^\infty(w)}.
	\end{align}
	 Hence, from the sequence argument, we get
	\begin{equation} \label{poT1d1}
	\begin{aligned}
		&\sup_n\|f_1^{(n)}(t)\|_{L_{x,v}^\infty(w)} \le C_se^{-\lambda t}(\|f_0\|_{L_{x,v}^\infty(w)}+\eta_0^{1+p}),\\
		&\sup_n\|f_2^{(n)}(t)\|_{L_{x,v}^\infty(\langle v \rangle^{\beta_3} \mu^{-1/2})} \le C_sC_{ns}e^{-\lambda t}(\|f_0\|_{L_{x,v}^\infty(w)}+\eta_0^{1+p}).
	\end{aligned}
	\end{equation}
	Combining the scheme in proposition \ref{poT2} for $f_1^{(n)}$ and \eqref{pof2Cauchy}, \eqref{poT1d1}, $\left(f_1^{(n)}\right)_{n\in \mathbb{N}}$ is Cauchy, and converges strongly to $f_1$ in $L^\infty_tL^\infty_{x,v}(w)$. It follows from \eqref{pof2Cauchy} that $\left(f_2^{(n)}\right)_{n\in \mathbb{N}}$ is Cauchy, and converges strongly to $f_2$ in $L^\infty_tL^\infty_{x,v}(w)$. By \eqref{poT1d1}, the following holds
	\begin{align*} 
		\|f_1(t)\|_{L_{x,v}^\infty(w)} \le C_se^{-\lambda t}(\|f_0\|_{L_{x,v}^\infty(w)}+\eta_0^{1+p}),\\ 
		\|f_2(t)\|_{L_{x,v}^\infty(w)} \le  C_sC_{ns}e^{-\lambda t}(\|f_0\|_{L_{x,v}^\infty(w)}+\eta_0^{1+p})
	\end{align*}
	for all $t\ge 0$.
	Therefore, taking the limit to the iterative system \eqref{poT1d4}, $f=f_1+f_2$ is a solution to \eqref{BE-f} and we conclude that
	\begin{align*}
		P_{G_L}(f) = 0 \quad \text{and} \quad 
		\quad\|f(t)\|_{L_{x,v}^\infty(w)} \le Ce^{-\lambda t}(\|f_0\|_{L_{x,v}^\infty(w)}+\eta_0^{1+p})
	\end{align*}
	for all $t \ge 0$.
\end{proof}

\bigskip

\begin{proposition} \label{poUT}
Let \(w(v)=\langle v\rangle^\beta\) with $\beta>1+\gamma+\frac{4C_b}{c_b}$. There exists \(\eta_0 >0 \) such that for any \(f_0 \) in \(L^{\infty}_{x,v}(w)\) with \(\|f_0\|_{L^{\infty}_{x,v}(w)}\le \eta_0^{1+p}\) for \( 0<p<1 \) and \( P_{G_L}(f_0)=0 \), there exists at most one solution \(f(t,x,v)\) in \(L^\infty_tL^{\infty}_{x,v}(w)\) to \eqref{BE-f} with \(f_0\)  as initial data. 
\end{proposition}

\begin{proof}
The proof of this theorem is similar to that of Theorem \ref{UT} and is therefore omitted.
\end{proof}

\bigskip

\begin{proposition}
	Let $w(v)=\langle v\rangle^\beta$ with $\beta>1+\gamma+\frac{4C_b}{c_b}$. Let $f_0$ in $L^\infty_{x,v}(w)$ with $P_{G_L}(f_0)=0$ and $\|f_0\|_{L^\infty_{x,v}(w)} \le \eta_0^{1+p}$ for some $0<p<1$, where $\eta_0$ is chosen so that Theorem \ref{poT1} and Theorem \ref{poUT} hold	and assume f is a unique solution of \eqref{BE-f} with initial data $f_0$. If $F_0 = \mu + f_0 \ge 0$, then $F=\mu+f \ge 0$.
\end{proposition}
\begin{proof}
	The proof of this theorem is similar to that of the stretched exponent case and is omitted.
\end{proof}	

\bigskip	

\subsection{Useful estimate for Large amplitude problem}

	The following lemma can be found in \cite{Cao2022JFA}. Here, we derive an explicit value for the constant $C_1$ appearing in its statement. The constant $C_1$ plays a key role in determining the admissible range of $\beta$ in Theorem \ref{Mainresultlargeamp-p}.

\begin{lemma} \label{poKESTIMATE} \cite[Lemma 6.3]{Cao2022JFA}
	Let $0\le \gamma \le 1$, $\beta>3+\gamma$. We have 
	\begin{align*}
		\int_{\R^3}|k_w(v,v_*)|
		   d v_* \le \frac{C_1}{\beta^{\frac{\gamma+3}{4}}}\langle v\rangle^\gamma+ \frac{C_\beta}{\langle v\rangle}
	\end{align*}
	for some constant $C_\beta>0$ and
	\begin{align*}
		C_1=8\sqrt{2}C_b \left(2+\frac{4}{\gamma+3}+\frac{1}{e}\right) \left(\frac{32\pi}{3+\gamma} e^{-3/4} +8\pi \right).
	\end{align*}
	 Moreover, for $\delta>0$ sufficiently small, we have
	\begin{align*}
		\int_{\{|v-v_*|>\frac{\langle v\rangle}{\delta}\} \cup \{|v-v_*|<\delta \langle v\rangle\}} |k_w(v,v_*)|dv_* \le C_{\beta,\delta}\langle v\rangle^\gamma,\quad \lim_{\delta \rightarrow 0 } C_{\beta,\delta} = 0.
	\end{align*}
\end{lemma}
\begin{proof}
	Note that 
\begin{align*}
	\int_{\R^3}k_w(v,v_*)   d v_*  &\le \frac{1}{(2\pi)^{3/2}} C_b \int_{\R^3} \int_{\S^2} |v-v_*|^\gamma \frac{\langle v\rangle^\beta}{\langle v_*\rangle^\beta}e^{-\frac{|v|^2}{2}} d\sigma dv_*\\
	&\quad  +  \frac{2}{(2\pi)^{3/2}} C_b \int_{\R^3} \int_{\S^2} |v-v_*|^\gamma \frac{\langle v\rangle^\beta}{\langle v'\rangle^\beta}e^{-\frac{|v_*'|^2}{2}} d\sigma dv_*\\
	& = I_1+I_2.
\end{align*}
We can easily obtain
\begin{align*}
	I_1 \le C_{\beta,\gamma} e^{-\frac{|v|^2}{4}}.
\end{align*}
Let us estimate $I_2$. By Lemma \ref{poCarlemansimple}, we have
\begin{align*}
	I_2 &= 8C_b \frac{1}{(2\pi)^{3/2}}\int_{\R^3} \frac{1}{|v-v'|} \frac{\langle v\rangle^\beta}{\langle v'\rangle^\beta} \int_{\{w:w \cdot (v'-v)=0\}} \frac{(\sqrt{|v-v'|^2+|w|^2})^\gamma}{\sqrt{|v-v'|^2+|w|^2}} e^{-\frac{|v+w|^2}{2}}dwdv'\\
	& \le 8C_b \frac{1}{(2\pi)^{3/2}}\int_{\R^3} \frac{1}{|v-v'|^{\frac{3-\gamma}{2}}} \frac{\langle v\rangle^\beta}{\langle v'\rangle^\beta} \int_{\{w:w \cdot (v'-v)=0\}} |w|^{\frac{\gamma-1}{2}} e^{-\frac{|v+w|^2}{2}}dwdv',
\end{align*}
where we have used $(|v-v'|^2+|w|^2)^{\frac{\gamma-1}{2}}\le |v-v'|^{\frac{\gamma-1}{2}}|w|^{\frac{\gamma-1}{2}}$. We split $v$ into two parts
\begin{align*}
	v= v_\perp+v_\shortparallel,\quad v_\perp =\frac{v\cdot (v'-v)}{|v-v'|^2}(v'-v),\quad v_{\perp} \perp v_\shortparallel. 
\end{align*}
Note that $|v|^2 = |v_\perp|^2+ |v_\shortparallel|^2$ and $|v+w|^2 = |w+v_\shortparallel|^2+|v_\perp|^2$.
Then we obtain
\begin{align*}
	I_2 &\le 8C_b \frac{1}{(2\pi)^{3/2}}\int_{\R^3} \frac{1}{|v-v'|^{\frac{3-\gamma}{2}}} e^{-\frac{|v_\perp|^2}{2}}\frac{\langle v\rangle^\beta}{\langle v'\rangle^\beta} \int_{\{w:w \cdot (v'-v)=0\}} |w|^{\frac{\gamma-1}{2}} e^{-\frac{|v_\shortparallel+w|^2}{2}}dwdv'.
\end{align*}
Here,
\begin{align} \label{pok11}
	\int_{\{w:w \cdot (v'-v)=0\}} |w|^{\frac{\gamma-1}{2}} e^{-\frac{|v_\shortparallel+w|^2}{2}}dw 
	&\le \int_{\R^2} |v_\shortparallel-w|^{\frac{\gamma-1}{2}}e^{-\frac{|w|^2}{2}}dw\nonumber\\
	&= \int_{|v_\shortparallel-w|\le \frac{\langle v_\shortparallel\rangle}{2}} |v_\shortparallel-w|^{\frac{\gamma-1}{2}}e^{-\frac{|w|^2}{2}}dw\nonumber\\
	& \quad +\int_{|v_\shortparallel-w|>\frac{\langle v_\shortparallel\rangle}{2}} |v_\shortparallel-w|^{\frac{\gamma-1}{2}}e^{-\frac{|w|^2}{2}}dw.
\end{align}
Since $|w|^2 \ge \frac{1}{2}|v_\shortparallel|^2 -|v_\shortparallel-w|^2 \ge\frac{1}{2}|v_\shortparallel|^2-\frac{\langle v_\shortparallel\rangle^2}{4}$ on $|v_\shortparallel-w|\le \frac{\langle v_\shortparallel\rangle}{2}$, the first term in \eqref{pok11} is bounded by
\begin{align*}
	\int_{|v_\shortparallel-w|\le \frac{\langle v_\shortparallel\rangle}{2}} |v_\shortparallel-w|^{\frac{\gamma-1}{2}}e^{-\frac{|w|^2}{2}}dw 
	&\le \int_{|v_\shortparallel-w|\le \frac{\langle v_\shortparallel\rangle}{2}} |v_\shortparallel-w|^{\frac{\gamma-1}{2}}e^{-\frac{|v_\shortparallel|^2}{4}}e^{\frac{\langle v_\shortparallel\rangle^2}{8}}dw\\
	&=2\pi  e^{-\frac{|v_\shortparallel|^2}{4}}e^{\frac{\langle v_\shortparallel\rangle^2}{8}} \int_0^{\frac{\langle v_\shortparallel \rangle}{2}} r^{\frac{\gamma+1}{2}}dr\\
	& \le \frac{4\pi}{3+\gamma}e^{-\frac{|v_\shortparallel|^2}{4}}e^{\frac{\langle v_\shortparallel\rangle^2}{8}}  \left(\frac{\langle v_\shortparallel\rangle}{2}\right)^{\frac{\gamma+3}{2}}\\
	& =  \frac{4\pi}{3+\gamma}\left(\frac{1}{2}\right)^{\frac{\gamma+3}{2}}e^{-\frac{|v_\shortparallel|^2}{4}}e^{\frac{\langle v_\shortparallel\rangle^2}{8}} \langle v_\shortparallel\rangle^2 \langle v_\shortparallel\rangle^{\frac{\gamma-1}{2}}\\
	& \le\frac{32\pi}{3+\gamma}\left(\frac{1}{2}\right)^{\frac{\gamma+3}{2}} e^{-3/4}\langle v_\shortparallel\rangle^{\frac{\gamma-1}{2}},
\end{align*}
where we have used $e^{-\frac{|v_\shortparallel|^2}{2}}e^{\frac{\langle v_\shortparallel\rangle^2}{4}} \langle v_\shortparallel\rangle^2 \le 8e^{-3/4}$ in the last term.\\
The last term in \eqref{pok11} is bounded by
\begin{align*}
	\int_{|v_\shortparallel-w|>\frac{\langle v_\shortparallel\rangle}{2}} |v_\shortparallel-w|^{\frac{\gamma-1}{2}}e^{-\frac{|w|^2}{2}}dw 
	&< \left(\frac{1}{2}\right)^{\frac{\gamma-1}{2}} \langle v_\shortparallel\rangle^{\frac{\gamma-1}{2}}\int_{\R^2} e^{-\frac{|w|^2}{2}}dw\\
	& = 2\pi \left(\frac{1}{2}\right)^{\frac{\gamma-1}{2}} \langle v_\shortparallel\rangle^{\frac{\gamma-1}{2}}.
\end{align*}
Thus
\begin{align*}
	I_2 \le 8C_b \frac{1}{(2\pi)^{3/2}}\left(\frac{32\pi}{3+\gamma}\left(\frac{1}{2}\right)^{\frac{\gamma+3}{2}} e^{-3/4} +2\pi \left(\frac{1}{2}\right)^{\frac{\gamma-1}{2}}\right) \int_{\R^3} \frac{1}{|v-v'|^{\frac{3-\gamma}{2}}} e^{-\frac{|v_\perp|^2}{2}}\frac{\langle v\rangle^\beta}{\langle v'\rangle^\beta}\langle v_\shortparallel\rangle^{\frac{\gamma-1}{2}} dv'.
\end{align*}
We now estimate as follows: 
\begin{align*}
	 J:=\int_{\R^3} \frac{1}{|v-v'|^{\frac{3-\gamma}{2}}} e^{-\frac{|v_\perp|^2}{2}}\frac{\langle v\rangle^\beta}{\langle v'\rangle^\beta}\langle v_\shortparallel\rangle^{\frac{\gamma-1}{2}} dv'.
\end{align*}
We split $v$-domain into two regions $|v|\le 1$ and $|v|\ge 1$.\\
\newline
$\textbf{Case 1: } |v|\le 1$ \\
We can compute
\begin{align*}
	J \mathbf{1}_{|v|\le 1} \le \int_{\R^3} \frac{1}{|v-v'|^{\frac{3-\gamma}{2}}} \frac{1}{\langle v'\rangle^\beta} dv' \le C_\beta \langle v\rangle^{\frac{\gamma-3}{2}}\le C_\beta \langle v\rangle^{-5},
\end{align*}
where $\beta>3$.\\
\newline
$\textbf{Case 2: } |v|\ge 1$ \\
We split this into two regions $|v_\perp|\ge \frac{|v|}{\beta}$ and $|v_\perp|\le \frac{|v|}{\beta}$.\\
\newline
$\textbf{Case 2-1: } |v|\ge 1$ and $|v_\perp|\ge \frac{|v|}{\beta}$\\
We can easily compute
\begin{align*}
	J \mathbf{1}_{|v|\ge 1 \cap |v_\perp|\ge \frac{|v|}{\beta} } \le \int_{\R^3} \frac{1}{|v-v'|^{\frac{3-\gamma}{2}}} e^{-\frac{|v|^2}{2\beta^2}}\frac{\langle v\rangle^\beta}{\langle v'\rangle^\beta} dv' \le C_\beta e^{-\frac{|v|^2}{4\beta^2}} \langle v\rangle^{\frac{\gamma-3}{2}} \le C_\beta \langle v\rangle^{-5}.
\end{align*}
\newline
$\textbf{Case 2-2: } |v|\ge 1$ and $|v_\perp|\le \frac{|v|}{\beta}$\\
Note that
\begin{align*}
	|v'|^2 = |v-v'|^2 +|v|^2 +2v\cdot (v'-v) &\ge |v-v'|^2 +|v|^2 -2|v-v'||v_\perp|\\
	& \ge |v-v'|^2 +|v|^2 -\frac{2}{\beta}|v-v'||v|\\
	& \ge \left(1-\frac{1}{\beta} \right) (|v-v'|^2+|v|^2),
\end{align*}
and
\begin{align*}
	|v_\shortparallel|^2 =|v|^2 -|v_\perp|^2 \ge \left(1-\frac{1}{\beta^2} \right) |v|^2.
\end{align*}
We have
\begin{align*}
	J \mathbf{1}_{|v|\ge 1 \cap |v_\perp|\le \frac{|v|}{\beta} } &\le \int_{\R^3} \frac{1}{|v-v'|^{\frac{3-\gamma}{2}}} e^{-\frac{|v_\perp|^2}{2}}\frac{\langle v\rangle^\beta}{\left(1+\frac{\beta-1}{\beta}(|v-v'|^2+|v|^2) \right)^{\beta/2}} \left(1+\left(1-\frac{1}{\beta^2}\right)|v|^2 \right)^{\frac{\gamma-1}{4}} dv'\\
	&\le \left(\frac{\beta^2-1}{\beta^2}\right)^{\frac{\gamma-1}{4}}\left(\frac{\beta}{\beta-1}\right)^{\frac{\beta}{2}}\langle v\rangle^{\frac{\gamma-1}{2}}\int_{\R^3} \frac{1}{|v-v'|^{\frac{3-\gamma}{2}}} e^{-\frac{|v_\perp|^2}{2}}\frac{1}{\left(1+\frac{|v-v'|^2}{1+|v|^2} \right)^{\beta/2}} dv'\\
	& \le 2\langle v\rangle^{\frac{\gamma-1}{2}}\int_{\R^3} \frac{1}{|v-v'|^{\frac{3-\gamma}{2}}} e^{-\frac{|v_\perp|^2}{2}}\frac{1}{\left(1+\frac{|v-v'|^2}{1+|v|^2} \right)^{\beta/2}} dv'.
\end{align*}
We make a change of variables
\begin{align*}
	r =|v-v'|,\quad v\cdot (v'-v) = r|v|\cos \theta, \quad dv'=r^2\sin \theta dr d\theta d\phi, \quad |v_\perp|^2 =|v|^2\cos^2 \theta, \quad |v_\shortparallel |^2 = |v|^2 \sin^2\theta
\end{align*}
to obtain
\begin{align*}
	J \mathbf{1}_{|v|\ge 1 \cap |v_\perp|\le \frac{|v|}{\beta}} &\le 4\pi  \langle v\rangle^{\frac{\gamma-1}{2}}\int_0^\infty r^{\frac{\gamma+1}{2}} \frac{1}{\left(1+\frac{r^2}{1+|v|^2} \right)^{\beta/2}} \int_0^\pi e^{-\frac{|v|^2}{2}\cos^2\theta}\sin \theta d\theta dr.
\end{align*}
We make a change of variables
\begin{align*}
	r =\langle  v \rangle z,\quad y=|v|\cos \theta
\end{align*}
to obtain
\begin{align*}
	J \mathbf{1}_{|v|\ge 1 \cap |v_\perp|\le \frac{|v|}{\beta}} 
	& \le 4\pi \langle v\rangle^{\frac{\gamma-1}{2}}\int_0^\infty \langle v\rangle^{\frac{\gamma+1}{2}}z^{\frac{\gamma+1}{2}} \frac{1}{\left(1+z^2 \right)^{\beta/2}} \int_{-|v|}^{|v|} e^{-\frac{y^2}{2}}\frac{1}{|v|}dy\langle v\rangle dz\\
	&\le 4\sqrt{2}\pi \langle v\rangle^{\gamma}\int_0^\infty  \frac{z^{\frac{\gamma+1}{2}}}{\left(1+z^2 \right)^{\beta/2}} \int_{-\infty}^{\infty}e^{-\frac{y^2}{2}}dydz\\
	& = 8\pi^{3/2}\langle v\rangle^{\gamma}\int_0^\infty  \frac{z^{\frac{\gamma+1}{2}}}{\left(1+z^2 \right)^{\beta/2}} dz\\
	& = 4\pi^{3/2}\langle v\rangle^{\gamma}\int_0^\infty  \frac{x^{\frac{\gamma-1}{4}}}{\left(1+x \right)^{\beta/2}} dx.
\end{align*}
Here,
\begin{align} \label{pok22}
	\int_0^\infty  \frac{x^{\frac{\gamma-1}{4}}}{\left(1+x \right)^{\beta/2}} dx \le \int_1^\infty  \frac{x^{\frac{\gamma-1}{4}}}{\left(1+x \right)^{\beta/2}} dx+\int_0^1  \frac{x^{\frac{\gamma-1}{4}}}{\left(1+x \right)^{\beta/2}} dx.
\end{align}
The first term in \eqref{pok22} is bounded by
\begin{align*}
	 \int_1^\infty  \frac{x^{\frac{\gamma-1}{4}}}{\left(1+x \right)^{\beta/2}} dx \le \int_1^\infty \frac{1}{\left(1+x \right)^{\beta/2}} dx = \frac{2}{\beta-2} 2^{-\frac{\beta}{2}+1} \le \frac{2}{\beta-2} \le \frac{2\cdot 4^{\frac{\gamma+3}{4}}}{\beta^{\frac{\gamma+3}{4}}},
\end{align*}
where $\beta> 3$.\\
The last term in \eqref{pok22} is bounded by
\begin{align*}
	\int_0^1  \frac{x^{\frac{\gamma-1}{4}}}{\left(1+x \right)^{\beta/2}} dx &\le \int_0^1  x^{\frac{\gamma-1}{4}}e^{-\frac{\beta}{4}x} dx \le  \frac{4}{\beta}\int_0^\infty \left(\frac{4}{\beta}t \right)^{\frac{\gamma-1}{4}}e^{-t} dt \le \left(\frac{4}{\beta} \right)^{\frac{\gamma+3}{4}}\int_0^\infty t^{\frac{\gamma-1}{4}}e^{-t}dt\\ 
	&\le \frac{4^{\frac{\gamma+3}{4}}}{\beta^{\frac{\gamma+3}{4}}} \left( \int_0^1 t^{\frac{\gamma-1}{4}}dt +\int_1^\infty e^{-t}dt\right) =\frac{4^{\frac{\gamma+3}{4}}}{\beta^{\frac{\gamma+3}{4}}}\left(\frac{4}{\gamma+3}+\frac{1}{e}\right),
\end{align*}
where we have used $(1+x)^{-\frac{\beta}{2}} \le e^{-\frac{\beta}{4}x}$. Thus
\begin{align*}
	J \mathbf{1}_{|v|\ge 1 \cap |v_\perp|\le \frac{|v|}{\beta}} \le 4\pi^{3/2} 4^{\frac{\gamma+3}{4}}\left(2+\frac{4}{\gamma+3}+\frac{1}{e}\right)\frac{1}{\beta^{\frac{\gamma+3}{4}}}\langle v\rangle^{\gamma}.
\end{align*}
Combining all cases, we derive
\begin{align*}
	I_2 &\le  8\sqrt{2}C_b \left(2+\frac{4}{\gamma+3}+\frac{1}{e}\right) \left(\frac{32\pi}{3+\gamma} e^{-3/4} +8\pi \right)\frac{1}{\beta^{\frac{\gamma+3}{4}}}\langle v\rangle^{\gamma}+C_\beta \langle v\rangle^{-1}.
\end{align*}
Hence we complete the proof of the lemma.
\end{proof}

\bigskip

The following two lemmas are used to demonstrate Lemma \ref{popointwiseQ+est}.

\begin{lemma} \cite{Silvestr2016CMP} \label{poCarlemansimple}
	For any non-negative function $F$ (in terms of $v_*$, $v$, $r=|v-v_*|$, $v'$, and/or $v_*'$), we have
	\begin{align*}
		\int_{\R^3} \int_{\S^2} F d\sigma dv_* = 4 \int_{\R^3} \frac{1}{|v-v'|} \int_{\{w : w \cdot (v'-v)=0\}} F \frac{1}{r}dwdv'.
	\end{align*}
	In the right-hand side, we must write $v_* =v'+w$ and $r = \sqrt{|v-v'|^2+|w|^2}$.
\end{lemma}

\bigskip

\begin{lemma} \label{popointwiseQ+est} 
    Suppose $\beta>4$. Then we have
    \begin{align} \label{poQ+estpoint}
        |wQ^+(f, f)(t,y,v)| \leq C_{\beta,\gamma} \langle v\rangle^\gamma  \|wf(t)\|_{L^\infty_{x,v}}\left(\int_{\R^3}\frac{1}{(1+|\eta|)^{\frac{2p}{p-1}(\beta-4)-4}}|wf(t,y,\eta)|^{\frac{2p}{p-1}}d\eta \right)^{\frac{p-1}{2p}},
    \end{align} 
    and
    \begin{align} \label{poQ-estpoint}
    	|wQ^-(f,f)(t,y,v)| 
    \le C_{\beta,\gamma} \langle v\rangle^\gamma \|wf(t)\|_{L^\infty_{x,v}} \left(\int_{\R^3} \frac{1}{(1+|v_*|)^{\frac{2p}{p-1}\beta-\frac{4p+4}{p-1}-\frac{2p}{p-1}\gamma}} |wf(t,y,v_*)|^{\frac{2p}{p-1}}dv_*\right)^{\frac{p-1}{2p}}
    \end{align}
    for some $\frac{4}{3}\le p<\infty$.
\end{lemma}

\begin{proof}
    For the above estimates to hold, the parameters $p$ and $\varepsilon$ must satisfy
    \begin{align*}
        p>1, \quad 0<p\gamma+p\varepsilon \le 2, \quad 1+\frac{p}{p-1} \varepsilon >2, \quad  \frac{2p}{p-1}\varepsilon <3.
    \end{align*}
For example, these conditions are simultaneously satisfied by choosing
$p=\frac{4}{3}$ and $\varepsilon=\frac{5}{16}$.

	First of all, we show \eqref{poQ-estpoint}. We easily compute that
	\begin{align*}
		|\langle v\rangle^\beta Q^-(f,f)(t,y,v)| &\le \int_{\R^3}\int_{\S^2} b(\cos \theta)|v-v_*|^\gamma \frac{1}{\langle v_*\rangle^\beta}|\langle v_*\rangle^\beta f(v_*)||\langle v\rangle^\beta f(v)|d\sigma dv_*\\
		& \le C \langle v\rangle^\gamma \|\langle v\rangle^\beta f\|_{L^\infty_v} \int_{\R^3}	  \frac{1}{\langle v_*\rangle^{\beta-\gamma}}|\langle v_*\rangle^\beta f(v_*)| dv_*\\
		& \le C_\beta \langle v\rangle^\gamma \|\langle v\rangle^\beta f(t)\|_{L^\infty_{x,v}} \left(\int_{\R^3} \frac{1}{(1+|v_*|)^{\frac{2p}{p-1}\beta-\frac{4p+4}{p-1}-\frac{2p}{p-1}\gamma}} |\langle v_*\rangle^\beta f(t,y,v_*)|^{\frac{2p}{p-1}}dv_*\right)^{\frac{p-1}{2p}}
	\end{align*}
	for all $p>1$. We prove \eqref{poQ-estpoint}.
	
	Now, we will show \eqref{poQ+estpoint}.
	We use the fact $\langle v\rangle^{\beta-4} \le C_\beta \langle v'\rangle^{\beta-4}+C_\beta\langle v_*'\rangle^{\beta-4}$ to obtain
	\begin{align*}
		|\langle v\rangle^\beta Q^+(f,f)(t,y,v)| &\le \int_{\R^3}\int_{\S^2} b(\cos \theta)|v-v_*|^\gamma \frac{\langle v\rangle^4}{\langle v_*'\rangle^4\langle v'\rangle^4}\langle v\rangle^{\beta-4}|\langle v_*'\rangle^4f(v_*')||\langle v'\rangle^4f(v')|d\sigma dv_*\\
		& \le C_\beta \int_{\R^3}\int_{\S^2}|v-v_*|^\gamma \frac{\langle v\rangle^4}{\langle v_*'\rangle^4\langle v'\rangle^4}|\langle v_*'\rangle^4 f(v_*')||\langle v'\rangle^\beta f(v')|d\sigma dv_*\\
		& \quad +C_\beta \int_{\R^3}\int_{\S^2}|v-v_*|^\gamma \frac{\langle v\rangle^4}{\langle v_*'\rangle^4\langle v'\rangle^4}|\langle v_*'\rangle^\beta f(v_*')||\langle v'\rangle^4 f(v')|d\sigma dv_*,
	\end{align*}
	where $\beta>4$.
	By rotation, without loss of generality, it suffices to estimate the latter of the two terms. We use the H\"older inequality to obtain
	\begin{align*}
		&\int_{\R^3}\int_{\S^2}|v-v_*|^\gamma \frac{\langle v\rangle^4}{\langle v_*'\rangle^4\langle v'\rangle^4}|\langle v_*'\rangle^\beta f(v_*')||\langle v'\rangle^4 f(v')|d\sigma dv_*\\
		& \le \|\langle v\rangle^\beta f\|_{L^\infty_v} \left(\int_{\R^3}\int_{\S^2}|v-v_*|^{p\gamma+p\varepsilon} \frac{\langle v\rangle^{4p}}{\langle v_*'\rangle^{4p} \langle v'\rangle^{4p}}d\sigma dv_*\right)^{1/p}\\
		&\quad \times  \left( \int_{\R^3}\int_{\S^2}|v-v_*|^{-\frac{p}{p-1} \varepsilon } |\langle v'\rangle^4 f(v')|^{\frac{p}{p-1}}d\sigma dv_*\right)^{\frac{p-1}{p}},
	\end{align*}
	where $p>1$.
	By Lemma \ref{pogaintermestimate}, we deduce
	\begin{align*}
		\left(\int_{\R^3}\int_{\S^2}|v-v_*|^{p\gamma+p\varepsilon} \frac{\langle v\rangle^{4p}}{\langle v_*'\rangle^{4p} \langle v'\rangle^{4p}}d\sigma dv_*\right)^{1/p} \le \langle v\rangle^{\gamma+\varepsilon},
	\end{align*}
	where $0<p\gamma+p\varepsilon \le 2$.\\
	We apply Lemma \ref{poCarlemansimple} to the following:
	\begin{align*}
		&\int_{\R^3}\int_{\S^2}|v-v_*|^{-\frac{p}{p-1} \varepsilon } |\langle v'\rangle^4 f(v')|^{\frac{p}{p-1}}d\sigma dv_*\\
		& \le 4 \int_{\R^3} \frac{|\langle v'\rangle^4 f(v')|^{\frac{p}{p-1}}}{|v-v'|} \int_{\{w: w \cdot (v'-v)=0\}} (\sqrt{|v-v'|^2 + |w|^2})^{-1-\frac{p}{p-1} \varepsilon }dwdv'.
	\end{align*}
	Making a change of variable $w=|v-v'|z$, the above term is bounded by
	\begin{align*}
		&4 \int_{\R^3} |v-v'|^{-\frac{p}{p-1}\varepsilon} |\langle v'\rangle^4 f(v')|^{\frac{p}{p-1}} \int_{\R^2} (\sqrt{1+|z|^2})^{-1-\frac{p}{p-1} \varepsilon }dzdv'\\
		& \le C\int_{\R^3} |v-v'|^{-\frac{p}{p-1}\varepsilon} |\langle v'\rangle^4 f(v')|^{\frac{p}{p-1}}dv',
	\end{align*}
	where $1+\frac{p}{p-1} \varepsilon >2$.
	We use the Cauchy-Schwarz inequality to derive
	\begin{align*}
		&\int_{\R^3}\int_{\S^2}|v-v_*|^{-\frac{p}{p-1} \varepsilon } |\langle v'\rangle^4 f(v')|^{\frac{p}{p-1}}d\sigma dv_* \\
		&\le  C\left(\int_{\R^3} |v-v'|^{-\frac{2p}{p-1}\varepsilon} \langle v'\rangle^{-4}dv'\right)^{1/2} \left(\int_{\R^3} \langle v'\rangle^4 |\langle v'\rangle^4 f(v')|^{\frac{2p}{p-1}}dv'\right)^{1/2}\\
		&\le C \langle v\rangle^{-\frac{p}{p-1}\varepsilon} \left(\int_{\R^3} \langle v'\rangle^{4-\frac{2p}{p-1}(\beta -4)} |\langle v'\rangle^\beta f(v')|^{\frac{2p}{p-1}}dv'\right)^{1/2},
	\end{align*}
	where $\frac{2p}{p-1}\varepsilon <3$. Thus we prove \eqref{poQ+estpoint}.

\end{proof}

\subsection{Nonlinear $L^\infty_{x,v}$ estimate} 
Setting $h=wf=\langle v\rangle^\beta f$, the equation \eqref{FPBE} can be written as
\begin{align}\label{poWpFPBE2}
	\partial_t h + v \cdot \nabla_x h + \nu(v) h=K_wh + wQ\left(f ,f \right).
\end{align}
This subsection is devoted to deriving the nonlinear $L^\infty_{x,v}$ estimate.

\begin{lemma} \label{poLALinftyestimate} 
	Let $\beta>\left(\frac{4C_1}{\nu_0}\right)^{\frac{4}{\gamma+3}}$, where
	\begin{align*}
C_1=8\sqrt{2}C_b \left(2+\frac{4}{\gamma+3}+\frac{1}{e}\right) \left(\frac{32\pi}{3+\gamma} e^{-3/4} +8\pi \right).
	\end{align*}
    Under the a priori assumption, there exists a constant $C_3 \ge 1$, depending on $\beta$, such that
\begin{align}\label{poLALinftyestimate1}
		\|h(t)\|_{L^\infty_{x,v}} &\le 
		C_3\left(1+\bar{M}\right)e^{-\frac{\nu_0}{2}t}\|h_0\|_{L^\infty_{x,v}}\nonumber \\
		&\quad + C_3\left(\frac{1}{N}+\varepsilon C_{\beta,N}+\frac{C_{\beta,N}}{\tilde{N}} \right)\left[\bar{M}+\bar{M}^2+\bar{M}^3\right]\nonumber\\
		&\quad +C_{\beta,N,\varepsilon}\left[H(F_0 \mid \mu )^{1/2} + H(F_0 \mid \mu )\right]\nonumber\\
		& \quad +C_{\beta,N,\varepsilon}\bar{M}^{1+\frac{q+1}{2q}}\left[H(F_0 \mid \mu )^{\frac{q-1}{4q}}+H(F_0 \mid \mu )^{\frac{q-1}{2q}}\right]\nonumber\\
		& \quad +   C_{\beta,N,\tilde{N},\varepsilon}\bar{M}^{2+\frac{q+1}{2q}}\left[H(F_0 \mid \mu )^{\frac{q-1}{4q}}+H(F_0 \mid \mu )^{\frac{q-1}{2q}} \right]
	\end{align}
for all $t\ge 0$, where $\varepsilon>0$ can be arbitrarily small, $N,\tilde{N}>0$ can be arbitrary large, and $q$ is the local H\"older exponent $p$ from Lemma \ref{popointwiseQ+est}.
\end{lemma}

\begin{proof}
	Let $(t,x,v) \in [0, \infty) \times \T^3 \times \R^3$. We apply the Duhamel principle to the weighted Boltzmann equation \eqref{poWpFPBE2}:
	\begin{equation} \label{poLA10}
	\begin{aligned}
			h(t,x,v) &= e^{-\nu(v) t} h_0(x-vt,v)\\
			&\quad + \int_0^t e^{-\nu(v) (t-s)} K_wh(s,x-v(t-s),v) ds\\
			&\quad + \int_0^t e^{-\nu(v) (t-s)} wQ(f,f)(s,x-v(t-s),v) ds\\
		& =: I_1+I_2+I_3.
	\end{aligned}
	\end{equation}
	First of all, we can easily compute that
	\begin{align} \label{poI1}
		 |I_1| \le e^{-\nu_0 t}\|h_0\|_{L^\infty_{x,v}}.
	\end{align}
	Now, let us estimate the term $I_2$ in \eqref{poLA10}:\
	\begin{align*}
		\int_0^{t} e^{-\nu(v) (t-s)}  \int_{\R^3}k_w(v,v')h(s,x-v(t-s),v')dv' ds.
	\end{align*}
	For $|v|\ge N$, we use Lemma \ref{poKESTIMATE} to obtain
	\begin{align} \label{poI2unbounded}
		|I_2 \mathbf{1}_{|v|\ge N}| &\le  \sup_{0\le s \le t}\left[e^{-\frac{\nu_0}{2}(t-s)}\|h(s)\|_{L^\infty_{x,v}}\right]\int_0^t e^{-\frac{\nu(v)}{2}  (t-s)}\left(\frac{C_1}{\beta^{\frac{\gamma+3}{4}}}\langle v\rangle^\gamma+ \frac{C_\beta}{\langle v\rangle} \right)ds  \nonumber\\
		& \le \frac{2C_1}{\nu_0 \beta^{\frac{\gamma+3}{4}}} \sup_{0\le s \le t}\left[e^{-\frac{\nu_0}{2}(t-s)}\|h(s)\|_{L^\infty_{x,v}}\right] +\frac{C_\beta}{N}\bar{M}.
	\end{align}
	where $C_1=8\sqrt{2}C_b \left(2+\frac{4}{\gamma+3}+\frac{1}{e}\right) \left(\frac{32\pi}{3+\gamma} e^{-3/4} +8\pi \right)$.\\
	For $|v| \le N$, we have
	\begin{align*}
		|I_{2}\mathbf{1}_{|v|\le N}| \le \mathbf{1}_{|v|\le N}\int_0^{t} e^{-\nu(v) (t-s)} \int_{\R^3}|k_w(v,v')||h(s,x-v(t-s),v')|dv'ds
	\end{align*}
	and we apply the Duhamel principle to  $|h(s,x-v(t-s),v')|$ in the above integral
	\begin{align*}
		|I_{2}\mathbf{1}_{|v|\le N}| &\le \mathbf{1}_{|v|\le N}\int_0^{t} e^{-\nu(v) (t-s)} \int_{\R^3}|k_w(v,v')|e^{- \nu(v') s}|h_0(x-v(t-s)-v's,v')|dv'ds\\
		&\quad +\mathbf{1}_{|v|\le N}\int_0^{t} e^{-\nu(v) (t-s)}\int_{\R^3}|k_w(v,v')|\\
		& \qquad \times \int_0^{s}e^{-\nu(v') (s-s')}|K_wh(s',x-v(t-s)-v'(s-s'),v')|ds'dv'ds\\
		& \quad +\mathbf{1}_{|v|\le N}\int_0^{t} e^{-\nu(v) (t-s)}\int_{\R^3}|k_w(v,v')|\\
		& \qquad \times\int_0^s e^{-\nu(v') (s-s')}|wQ(f,f)(s',x-v(t-s)-v'(s-s'),v')|ds'dv'ds\\
		& =:I_{21}+I_{22}+I_{23}.
	\end{align*}
	For $I_{21}$, we use Lemma \ref{poKESTIMATE} to obtain
	\begin{align} \label{poI21}
		|I_{21}| \le  C_\beta e^{-\frac{\nu_0}{2}t} \|h_0\|_{L^\infty_{x,v}}\int_0^{t} e^{-\frac{\nu(v)}{2} (t-s)} \langle v\rangle^\gamma ds  \le C_\beta e^{-\frac{\nu_0}{2}t} \|h_0\|_{L^\infty_{x,v}} .
	\end{align} 
	For $I_{22}$, we have
	\begin{align*}
		|I_{22}| &\le \mathbf{1}_{|v|\le N}\int_0^{t} e^{-\nu(v) (t-s)}\int_{\R^3}|k_w(v,v')|\\
		& \qquad \times \int_0^{s}e^{-\nu(v') (s-s')}\int_{\R^3}|k_w(v',v'')||h(s',x-v(t-s)-v'(s-s'),v'')|dv''ds'dv'ds. 
	\end{align*}
	By Lemma \ref{poKESTIMATE}, there exists a bounded function $k_{\beta,N}$ so that
	\begin{align} \label{pokprop}
		k_{\beta,N}(v,v') = k_w(v,v')\mathbf{1}_{\frac{\langle v\rangle}{C_{\beta,N}}\le |v-v'|\le C_{\beta,N}\langle v\rangle}, \quad \int_{\R^3} \left|k_w(v,v')-k_{\beta,N}(v,v') \right|dv' \le \frac{C_{\beta}}{N}\langle v\rangle^\gamma
	\end{align}
	for some large constant $C_{\beta,N}>0$.\\
	We split
	\begin{align*}
		k_w(v,v')k_w(v',v'') &= (k_w(v,v')-k_{\beta,N}(v,v'))k_w(v',v'') +k_{\beta,N}(v,v')(k_w(v',v'')-k_{\beta,N}(v',v'') )\\  &\quad+k_{\beta,N}(v,v')k_{\beta,N}(v',v'')
	\end{align*}
	and 
	\begin{align*}
		|I_{22}| &\le \mathbf{1}_{|v|\le N}\int_0^{t} e^{-\nu(v) (t-s)}\int_{\R^3}|k_w(v,v')-k_{\beta,N}(v,v')|\\
		& \qquad \times \int_0^{s}e^{-\nu(v') (s-s')}\int_{\R^3}|k_w(v',v'')||h(s',x-v(t-s)-v'(s-s'),v'')|dv''ds'dv'ds\\
		& \quad +\mathbf{1}_{|v|\le N}\int_0^{t} e^{-\nu(v) (t-s)}\int_{\R^3}|k_{\beta,N}(v,v')|\\
		& \qquad \times \int_0^{s}e^{-\nu(v') (s-s')}\int_{\R^3}|k_w(v',v'')-k_{\beta,N}(v',v'')||h(s',x-v(t-s)-v'(s-s'),v'')|dv''ds'dv'ds\\
		& \quad +\mathbf{1}_{|v|\le N}\int_0^{t} e^{-\nu(v) (t-s)}\int_{\R^3}|k_{\beta,N}(v,v')|\\
		& \qquad \times \int_0^{s}e^{-\nu(v') (s-s')}\int_{\R^3}|k_{\beta,N}(v',v'')||h(s',x-v(t-s)-v'(s-s'),v'')|dv''ds'dv'ds\\
		&=:I_{22}^1+I_{22}^2+I_{22}^3.
	\end{align*}
	By Lemma \ref{poKESTIMATE} and \eqref{pokprop}, the term $I_{22}^1$ is bounded by
	\begin{align} \label{poI221}
		I_{22}^1 &\le C_\beta \bar{M}\mathbf{1}_{|v|\le N}\int_0^{t} e^{-\nu(v) (t-s)}\int_{\R^3}|k_w(v,v')-k_{\beta,N}(v,v')|\int_0^{s}e^{-\nu(v') (s-s')} \langle v'\rangle^\gamma ds'dv'ds\nonumber\\
		& \le C_\beta\bar{M}\mathbf{1}_{|v|\le N}\int_0^{t} e^{-\nu(v) (t-s)}\frac{1}{N}\langle v\rangle^\gamma ds\nonumber\\
		& \le  \frac{C_\beta}{N}\bar{M}.
	\end{align}
	In a similar way to \eqref{poI221}, we have
	\begin{align} \label{poI222}
		I_{22}^2 \le \frac{C_\beta}{N}\bar{M}.
	\end{align}
	Note that $k_{\beta,N}(v,v')$ and $k_{\beta,N}(v',v'')$ are supported in $\frac{\langle v\rangle}{C_{\beta,N}}\le |v-v'|\le C_{\beta,N}\langle v\rangle$ and $\frac{\langle v'\rangle}{C_{\beta,N}}\le |v'-v''|\le C_{\beta,N}\langle v'\rangle$. Thus, since $|v|\le N$, $k_{\beta,N}(v,v')$ and $k_{\beta,N}(v',v'')$ are supported in $|v|\le N$, $|v'|\le \tilde{C}_{\beta,N}$, and $|v''|\le \tilde{C}_{\beta,N}$ for some constant $\tilde{C}_{\beta,N}>0$.\\
	For $I_{22}^3$, we have
	\begin{align*}
		I_{22}^3 &\le \mathbf{1}_{|v|\le N}\int_0^{t} e^{-\nu(v) (t-s)}\int_{|v'|\le \tilde{C}_{\beta,N}}|k_{\beta,N}(v,v')|\\
		& \qquad \times \int_0^{s}e^{-\nu(v') (s-s')}\int_{|v''|\le \tilde{C}_{\beta,N}}|k_{\beta,N}(v',v'')||h(s',x-v(t-s)-v'(s-s'),v'')|dv''ds'dv'ds\\
		& \le C_{\beta,N}\mathbf{1}_{|v|\le N}\int_0^{t} e^{-\nu(v) (t-s)}\int_{|v'|\le \tilde{C}_{\beta,N}}\\
		& \qquad \times \int_0^{s}e^{-\nu(v') (s-s')}\int_{|v''|\le \tilde{C}_{\beta,N}}|h(s',x-v(t-s)-v'(s-s'),v'')|dv''ds'dv'ds\\
		& = C_{\beta,N}\mathbf{1}_{|v|\le N}\int_0^{t} e^{-\nu(v) (t-s)}\int_{|v'|\le \tilde{C}_{\beta,N}}\\
		& \qquad \times \int_0^{s-\varepsilon}e^{-\nu(v') (s-s')}\int_{|v''|\le \tilde{C}_{\beta,N}}|h(s',x-v(t-s)-v'(s-s'),v'')|dv''ds'dv'ds\\
		& \quad +C_{\beta,N}\mathbf{1}_{|v|\le N}\int_0^{t} e^{-\nu(v) (t-s)}\int_{|v'|\le \tilde{C}_{\beta,N}}\\
		& \qquad \times \int_{s-\varepsilon}^{s}e^{-\nu(v') (s-s')}\int_{|v''|\le \tilde{C}_{\beta,N}}|h(s',x-v(t-s)-v'(s-s'),v'')|dv''ds'dv'ds.
	\end{align*}
	For the first term, it is bounded by
	\begin{align*}
		C_{\beta,N}\mathbf{1}_{|v|\le N}\int_0^{t} \int_0^{s-\varepsilon} e^{-\nu_0  (t-s')}\int_{|v'|\le \tilde{C}_{\beta,N}}\int_{|v''|\le \tilde{C}_{\beta,N}}|h(s',x-v(t-s)-v'(s-s'),v'')|dv''dv'ds'ds
	\end{align*}
	We make a change of variables $v' \mapsto y=x-v(t-s)-v'(s-s')$ with $\left|\frac{dy}{dv'}\right| = (s-s')^3$:
	\begin{align*}
		C_{\beta,N,\varepsilon}\int_{0}^{t}\int_{0}^{s-\varepsilon} e^{-\nu_0  (t-s')}\left(\int_{\T^3}\int_{ |v''|\le \tilde{C}_{\beta,N} }|h(s',y,v'')|dv''dy  \right)ds'ds.
	\end{align*}
	Note that from Lemma \ref{L1L2cont},
	\begin{equation} \label{poRTIRE}
\begin{aligned}
	&\int_{\T^3}\int_{|v''|\le \tilde{C}_{\beta,N} }|h(s,y,v'')|dv'dy\\
	& \le C_{\beta,N}\int_{\T^3}\int_{|v''|\le \tilde{C}_{\beta,N}}|f(s,y,v'')| \mathbf{1}_{\{|f|>\mu\}}dv''dy+C_{\beta,N} \left(\int_{\T^3}\int_{|v''|\le \tilde{C}_{\beta,N}}\frac{1}{4\mu(v'')}|f(s,y,v'')|^2 \mathbf{1}_{\{|f|\le \mu\}} dv''dy\right)^{1/2}\\
	&\le C_{\beta,N} H(F_0 \mid \mu ) + C_{\beta,N}H(F_0 \mid \mu )^{1/2}.
\end{aligned}
\end{equation}
	It follows from \eqref{poRTIRE} that 
	\begin{align}\label{poLA38}
		&C_{\beta,N,\varepsilon}\int_0^{t} \int_0^{s-\varepsilon} e^{-\nu_0  (t-s')}\left[H(F_0 \mid \mu ) + H(F_0 \mid \mu )^{1/2}\right]ds'ds\nonumber\\
		&\le C_{\beta,N,\varepsilon}\left[H(F_0 \mid \mu )^{1/2} + H(F_0 \mid \mu )\right].
	\end{align}
	For the second term, it is bounded by
	\begin{align} \label{poLA37}
		&C_{\beta,N}\mathbf{1}_{|v|\le N}\int_0^{t} e^{-\nu(v) (t-s)}\int_{|v'|\le \tilde{C}_{\beta,N}}\nonumber\\
		& \qquad \times \int_{s-\varepsilon}^{s}e^{-\nu(v') (s-s')}\int_{|v''|\le \tilde{C}_{\beta,N}}|h(s',x-v(t-s)-v'(s-s'),v'')|dv''ds'dv'ds\nonumber\\
		& \le \varepsilon C_{\beta,N} \bar{M} \int_0^{t} e^{-\nu(v) (t-s)}\int_{|v'|\le \tilde{C}_{\beta,N}}dv'ds\nonumber\\
		& \le \varepsilon C_{\beta,N} \bar{M}.
	\end{align}
	\newline
	Gathering \eqref{poLA38} and \eqref{poLA37}, we obtain
	\begin{align} \label{poI223}
		I_{22}^3 \le \varepsilon C_{\beta,N} \bar{M}+C_{\beta,N,\varepsilon}\left[H(F_0 \mid \mu )^{1/2} + H(F_0 \mid \mu )\right].
	\end{align}
	Gathering \eqref{poI221}, \eqref{poI222}, and  \eqref{poI223}, we derive
	\begin{align} \label{poI22}
		|I_{22}| \le  \frac{C_\beta}{N}\bar{M}+\varepsilon C_{\beta,N} \bar{M}+C_{\beta,N,\varepsilon}\left[H(F_0 \mid \mu )^{1/2} + H(F_0 \mid \mu )\right].
	\end{align}
	\newline
	For $I_{23}$, splitting $k_w(v,v') = (k_w(v,v')-k_{\beta,N}(v,v'))+k_{\beta,N}(v,v')$, it is bounded by
	\begin{align*}
		&\mathbf{1}_{|v|\le N}\int_0^{t} e^{-\nu(v) (t-s)}\int_{\R^3}|k_w(v,v')-k_{\beta,N}(v,v')|\\
		& \qquad \times\int_0^s e^{-\nu(v') (s-s')}|wQ(f,f)(s',x-v(t-s)-v'(s-s'),v')|ds'dv'ds\\
		&+	\mathbf{1}_{|v|\le N}\int_0^{t} e^{-\nu(v) (t-s)}\int_{\R^3}|k_{\beta,N}(v,v')|\\
		& \qquad \times\int_0^s e^{-\nu(v') (s-s')}|wQ(f,f)(s',x-v(t-s)-v'(s-s'),v')|ds'dv'ds\\
		& =: I_{23}^1 +I_{23}^2.
	\end{align*}
	By \eqref{pokprop}, note that $k_{\beta,N}(v,v')$ is supported in $\frac{\langle v\rangle}{C_{\beta,N}}\le |v-v'|\le C_{\beta,N}\langle v\rangle$. Thus, since $|v|\le N$, $k_{\beta,N}(v,v')$ is supported in $|v|\le N$ and $|v'|\le \tilde{C}_{\beta,N}$ for some constant $\tilde{C}_{\beta,N}>0$.\\
	For $I_{23}^1$, by \eqref{pokprop} and Lemma \ref{popointwiseQ+est}, we have
	\begin{align} \label{poI231}
		&C_\beta \bar{M}^2 \mathbf{1}_{|v|\le N}\int_0^{t} e^{-\nu(v) (t-s)}\int_{\R^3}|k_w(v,v')-k_{\beta,N}(v,v')|\int_0^s e^{-\nu(v') (s-s')} \langle v'\rangle^\gamma   ds'dv'ds\nonumber\\
		& \le \frac{C_\beta}{N}\bar{M}^2 \int_0^{t} e^{-\nu(v) (t-s)} \langle v\rangle^\gamma ds\nonumber\\
		&\le \frac{C_\beta}{N}\bar{M}^2.
	\end{align}
	For $I_{23}^2$, by Lemma \ref{popointwiseQ+est}, we have
	\begin{align} \label{poLA40}
		&\mathbf{1}_{|v|\le N}\int_0^{t} e^{-\nu(v) (t-s)}\int_{|v'|\le \tilde{C}_{\beta,N}}|k_{\beta,N}(v,v')|\nonumber\\
		& \qquad \times\int_0^s e^{-\nu(v') (s-s')}|wQ(f,f)(s',x-v(t-s)-v'(s-s'),v')|ds'dv'ds\nonumber\\
		& \le C_{\beta,N}\mathbf{1}_{|v|\le N}\int_0^{t} e^{-\nu(v) (t-s)}\int_{|v'|\le \tilde{C}_{\beta,N}}\nonumber\\
		& \qquad \times\int_{s-\varepsilon}^s e^{-\nu(v') (s-s')}|wQ(f,f)(s',x-v(t-s)-v'(s-s'),v')|ds'dv'ds\nonumber\\
		& \quad +C_{\beta,N}\mathbf{1}_{|v|\le N}\int_0^{t} e^{-\nu(v) (t-s)}\int_{|v'|\le \tilde{C}_{\beta,N}}\nonumber\\
		& \qquad \quad \times\int_0^{s-\varepsilon} e^{-\nu(v') (s-s')}|wQ(f,f)(s',x-v(t-s)-v'(s-s'),v')|ds'dv'ds\nonumber\\
		& \le \varepsilon C_{\beta,N} \bar{M}^2 \mathbf{1}_{|v|\le N} \int_0^{t} e^{-\nu(v) (t-s)}\int_{|v'|\le \tilde{C}_{\beta,N}}dv'ds\nonumber\\
		& \quad +C_{\beta,N}\mathbf{1}_{|v|\le N}\int_0^{t} e^{-\nu(v) (t-s)}\int_{|v'|\le \tilde{C}_{\beta,N}}\int_0^{s-\varepsilon} e^{-\nu(v') (s-s')}\langle v'\rangle^\gamma\nonumber\\
		& \qquad \quad \times   \|h(s')\|_{L^\infty_{x,v}}\left(\int_{\R^3}\frac{1}{(1+|\eta|)^{\frac{2p}{p-1}(\beta-4)-4}}|h(s',x-v(t-s)-v'(s-s'),\eta)|^{\frac{2p}{p-1}}d\eta \right)^{\frac{p-1}{2p}}ds'dv'ds\nonumber\\
		& \le \varepsilon C_{\beta,N}\bar{M}^2 \nonumber\\
		& \quad +C_{\beta,N}\bar{M}\mathbf{1}_{|v|\le N}\int_0^{t} e^{-\nu(v) (t-s)}\int_{|v'|\le \tilde{C}_{\beta,N}}\int_0^{s-\varepsilon} e^{-\nu(v') (s-s')}\langle v'\rangle^\gamma \nonumber\\
		& \qquad \quad \times   \left(\int_{\R^3}\frac{1}{(1+|\eta|)^{\frac{2p}{p-1}(\beta-4)-4}}|h(s',x-v(t-s)-v'(s-s'),\eta)|^{\frac{2p}{p-1}}d\eta \right)^{\frac{p-1}{2p}}ds'dv'ds.
	\end{align}
	For the last term, we split the domain for $\eta$ into $|\eta|\le \tilde{N}$ and $|\eta|\ge \tilde{N}$:
	\begin{align} \label{poLA41}
		&\frac{C_{\beta,N}}{\tilde{N}}\bar{M}^2\mathbf{1}_{|v|\le N}\int_0^{t} e^{-\nu(v) (t-s)}\int_{|v'|\le \tilde{C}_{\beta,N}}\int_0^{s-\varepsilon} e^{-\nu(v') (s-s')}\langle v'\rangle^\gamma ds'dv'ds\nonumber\\
		&+C_{\beta,N}\bar{M}\mathbf{1}_{|v|\le N}\int_0^{t} e^{-\nu(v) (t-s)}\int_{|v'|\le \tilde{C}_{\beta,N}}\int_0^{s-\varepsilon} e^{-\nu(v') (s-s')}\langle v'\rangle^\gamma\nonumber\\
		& \qquad \quad \times   \left(\int_{|\eta|\le \tilde{N}}\frac{1}{(1+|\eta|)^{\frac{2p}{p-1}(\beta-4)-4}}|h(s',x-v(t-s)-v'(s-s'),\eta)|^{\frac{2p}{p-1}}d\eta \right)^{\frac{p-1}{2p}}ds'dv'ds\nonumber\\
		&\le \frac{C_{\beta,N}}{\tilde{N}}\bar{M}^2\nonumber\\
		&\quad +C_{\beta,N}\bar{M}\mathbf{1}_{|v|\le N}\int_0^{t}\int_0^{s-\varepsilon} e^{-\nu_0  (t-s')}\nonumber\\
		& \qquad \quad \times   \left(\int_{|v'|\le \tilde{C}_{\beta,N}} \int_{|\eta|\le \tilde{N}}|h(s',x-v(t-s)-v'(s-s'),\eta)|^{\frac{2p}{p-1}}d\eta dv' \right)^{\frac{p-1}{2p}}ds'ds.
	\end{align}
	We make a change of variables $v' \mapsto y=x-v(t-s)-v'(s-s')$ with $\left|\frac{dy}{dv'}\right| = (s-s')^3$:
	\begin{align*}
		C_{\beta,N,\varepsilon}\bar{M}^{1+\frac{p+1}{2p}}\mathbf{1}_{|v|\le N}\int_0^{t}\int_0^{s-\varepsilon} e^{-\nu_0  (t-s')}\left(\int_{\T^3} \int_{|\eta|\le \tilde{N}}|h(s',y,\eta)|d\eta dy \right)^{\frac{p-1}{2p}}ds'ds.
	\end{align*}
	It follows from \eqref{poRTIRE} that 
	\begin{align}\label{poLA45}
		&C_{\beta,N,\varepsilon}\bar{M}^{1+\frac{p+1}{2p}}\mathbf{1}_{|v|\le N}\int_0^{t}\int_0^{s-\varepsilon} e^{-\nu_0  (t-s')}\left[H(F_0 \mid \mu )^{\frac{p-1}{4p}}+H(F_0 \mid \mu )^{\frac{p-1}{2p}}\right]ds'ds\nonumber\\
		& \le C_{\beta,N,\varepsilon}\bar{M}^{1+\frac{p+1}{2p}}\left[H(F_0 \mid \mu )^{\frac{p-1}{4p}}+H(F_0 \mid \mu )^{\frac{p-1}{2p}}\right].
	\end{align}
	Gathering \eqref{poLA40}, \eqref{poLA41}, and \eqref{poLA45}, we derive
	\begin{align} \label{poI232}
		I_{23}^2 \le \varepsilon C_{\beta,N} \bar{M}^2 +\frac{C_{\beta,N}}{\tilde{N}}\bar{M}^2+C_{\beta,N,\varepsilon}\bar{M}^{1+\frac{p+1}{2p}}\left[H(F_0 \mid \mu )^{\frac{p-1}{4p}}+H(F_0 \mid \mu )^{\frac{p-1}{2p}}\right].
	\end{align}
	Gathering \eqref{poI231} and \eqref{poI232}, we obtain
	\begin{align} \label{poI23}
		|I_{23}| \le \frac{C_\beta}{N}\bar{M}^2+\varepsilon C_{\beta,N} \bar{M}^2 +\frac{C_{\beta,N}}{\tilde{N}}\bar{M}^2+C_{\beta,N,\varepsilon}\bar{M}^{1+\frac{p+1}{2p}}\left[H(F_0 \mid \mu )^{\frac{p-1}{4p}}+H(F_0 \mid \mu )^{\frac{p-1}{2p}}\right].
	\end{align}
	\newline
	Inserting \eqref{poI2unbounded}, \eqref{poI21}, \eqref{poI22}, and \eqref{poI23} into $I_{2}$, we derive
	\begin{align} \label{poI2}
		|I_2| &\le \frac{2C_1}{\nu_0\beta^{\frac{\gamma+3}{4}}} \sup_{0\le s \le t}\left[e^{-\frac{\nu_0}{2}(t-s)}\|h(s)\|_{L^\infty_{x,v}}\right]+C_\beta e^{-\frac{\nu_0}{2}t} \|h_0\|_{L^\infty_{x,v}}\nonumber\\
		&\quad +\frac{C_\beta}{N}\bar{M}+\varepsilon C_{\beta,N} \bar{M}+C_{\beta,N,\varepsilon}\left[H(F_0 \mid \mu )^{1/2} + H(F_0 \mid \mu )\right]\nonumber\\
		& \quad +\frac{C_\beta}{N}\bar{M}^2+\varepsilon C_{\beta,N} \bar{M}^2 +\frac{C_{\beta,N}}{\tilde{N}}\bar{M}^2+C_{\beta,N,\varepsilon}\bar{M}^{1+\frac{p+1}{2p}}\left[H(F_0 \mid \mu )^{\frac{p-1}{4p}}+H(F_0 \mid \mu )^{\frac{p-1}{2p}}\right].
	\end{align}
	Lastly, let us estimate $I_{3}$. From Lemma \ref{popointwiseQ+est}, it is bounded by
\begin{align*}
	|I_{3}|& \le C_\beta \bar{M} \int_0^t e^{-\nu(v) (t-s)} \langle v\rangle^\gamma \left(\int_{\R^3}\frac{1}{(1+|\eta|)^{\frac{2p}{p-1}(\beta-4)-4}}|h(s,x-v(t-s),\eta)|^{\frac{2p}{p-1}}d\eta \right)^{\frac{p-1}{2p}} ds\\
	& \le C_\beta \bar{M} \int_0^t e^{-\nu(v) (t-s)} \langle v\rangle^\gamma \left(\int_{|\eta|\ge N }\frac{1}{(1+|\eta|)^{\frac{2p}{p-1}(\beta-4)-4}}|h(s,x-v(t-s),\eta)|^{\frac{2p}{p-1}}d\eta \right)^{\frac{p-1}{2p}} ds\\
	& \quad +C_\beta \bar{M} \int_0^t e^{-\nu(v) (t-s)} \langle v\rangle^\gamma \left(\int_{|\eta|\le N }\frac{1}{(1+|\eta|)^{\frac{2p}{p-1}(\beta-4)-4}}|h(s,x-v(t-s),\eta)|^{\frac{2p}{p-1}}d\eta \right)^{\frac{p-1}{2p}} ds\\
	& \le \frac{C_\beta}{N} \bar{M}^2\\
	& \quad + C_\beta \bar{M} \int_0^t e^{-\nu(v) (t-s)} \langle v\rangle^\gamma \left(\int_{|\eta|\le N }\frac{1}{(1+|\eta|)^{\frac{2p}{p-1}(\beta-4)-4}}|h(s,x-v(t-s),\eta)|^{\frac{2p}{p-1}}d\eta \right)^{\frac{p-1}{2p}} ds.
\end{align*}
	We apply the Duhamel principle to bound the term $|h(s,x-v(t-s),\eta)|$ in the last term of the above
	\begin{align*}
		&\left(\int_{|\eta|\le N }\frac{1}{(1+|\eta|)^{\frac{2p}{p-1}(\beta-4)-4}}|h(s,x-v(t-s),\eta)|^{\frac{2p}{p-1}}d\eta \right)^{\frac{p-1}{2p}}\\
		&\le \left(\int_{|\eta|\le N}\frac{1}{(1+|\eta|)^{\frac{2p}{p-1}(\beta-4)-4}}e^{-\frac{2p}{p-1}\nu(\eta) s}|h_0(x-v(t-s)-\eta s,\eta)|^{\frac{2p}{p-1}} d\eta \right)^{\frac{p-1}{2p}}\\
		& \quad +\left(\int_{|\eta|\le N}\frac{1}{(1+|\eta|)^{\frac{2p}{p-1}(\beta-4)-4}}\left|\int_{0}^{s}e^{-\nu(\eta)(s-s')}K_wh(s',x-v(t-s)-\eta(s-s'),\eta)ds'\right|^{\frac{2p}{p-1}} d\eta \right)^{^{\frac{p-1}{2p}}}\\
		& \quad +\left(\int_{|\eta|\le N}\frac{1}{(1+|\eta|)^{\frac{2p}{p-1}(\beta-4)-4}}\left|\int_{0}^{s}e^{-\nu(\eta)(s-s')}wQ(f,f)(s',x-v(t-s)-\eta(s-s'),\eta)ds'\right|^{\frac{2p}{p-1}} d\eta \right)^{^{\frac{p-1}{2p}}}\\
		&=: \tilde{L}_1+\tilde{L}_2+\tilde{L}_3.
	\end{align*}
	We define
	\begin{align*}
		L_{i} = C_\beta \bar{M} \int_0^t e^{-\nu(v) (t-s)} \langle v\rangle^\gamma \tilde{L}_ids.
	\end{align*}
	For $L_1$, we easily compute 
	\begin{align} \label{poLA54}
		L_1 \le C e^{-\frac{\nu_0}{2}t}\bar{M}\|h_0\|_{L^\infty_{x,v}} .
	\end{align}
	For $L_2$, we first need to estimate $\tilde{L}_2$. We use Lemma \ref{poKESTIMATE} to obtain
	\begin{align*}
		\tilde{L}_2 &\le  \Bigg(\int_{|\eta|\le N}\frac{1}{(1+|\eta|)^{\frac{2p}{p-1}(\beta-4)-4}}\bigg|\int_{0}^{s}e^{-\nu(\eta)(s-s')}\int_{\R^3}|k_w(\eta,\eta')-k_{\beta,N}(\eta,\eta')|\\
		& \qquad \times |h(s',x-v(t-s)-\eta(s-s'),\eta')|d\eta'ds'\bigg|^{\frac{2p}{p-1}} d\eta \Bigg)^{^{\frac{p-1}{2p}}}\\
		& \quad + \Bigg(\int_{|\eta|\le N}\frac{1}{(1+|\eta|)^{\frac{2p}{p-1}(\beta-4)-4}}\bigg|\int_{0}^{s}e^{-\nu(\eta)(s-s')}\int_{\R^3}|k_{\beta,N}(\eta,\eta')|\\
		& \qquad  \quad \times |h(s',x-v(t-s)-\eta(s-s'),\eta')|d\eta' ds'\bigg|^{\frac{2p}{p-1}} d\eta \Bigg)^{^{\frac{p-1}{2p}}}\\
		& \le \frac{C_\beta}{N}\bar{M} \Bigg(\int_{|\eta|\le N}\frac{1}{(1+|\eta|)^{\frac{2p}{p-1}(\beta-4)-4}}\bigg|\int_{0}^{s}e^{-\nu(\eta)(s-s')}\langle \eta\rangle^\gamma ds'\bigg|^{\frac{2p}{p-1}} d\eta \Bigg)^{^{\frac{p-1}{2p}}}\\
		& \quad + C_{\beta,N}\Bigg(\int_{|\eta|\le N}\frac{1}{(1+|\eta|)^{\frac{2p}{p-1}(\beta-4)-4}}\bigg|\left(\int_{s-\varepsilon}^{s}+\int_{0}^{s-\varepsilon}\right)e^{-\nu(\eta)(s-s')}\int_{|\eta'|\le \tilde{C}_{\beta,N}}\\
		& \qquad  \quad \times |h(s',x-v(t-s)-\eta(s-s'),\eta')|d\eta' ds'\bigg|^{\frac{2p}{p-1}} d\eta \Bigg)^{^{\frac{p-1}{2p}}}\\
		& \le \frac{C_\beta}{N}\bar{M}+\varepsilon C_{\beta,N} \bar{M} \\
		& \quad +  C_{\beta,N}\Bigg(\int_{|\eta|\le N}\frac{1}{(1+|\eta|)^{\frac{2p}{p-1}(\beta-4)-4}}\bigg|\int_{0}^{s-\varepsilon}e^{-\nu(\eta)(s-s')}\int_{|\eta'|\le \tilde{C}_{\beta,N}}\\
		& \qquad  \quad \times |h(s',x-v(t-s)-\eta(s-s'),\eta')|d\eta' ds'\bigg|^{\frac{2p}{p-1}} d\eta \Bigg)^{^{\frac{p-1}{2p}}}.
	\end{align*}
	From the H\"older inequality, the last term is bounded by
	\begin{align*}
		&C_{\beta,N}\Bigg(\int_{|\eta|\le N}\frac{1}{(1+|\eta|)^{\frac{2p}{p-1}(\beta-4)-4}}\bigg|\int_{0}^{s-\varepsilon}e^{-\nu_0(s-s')}\int_{|\eta'|\le \tilde{C}_{\beta,N}}|h(s',x-v(t-s)-\eta(s-s'),\eta')|\\
		& \qquad  \quad \times d\eta' ds'\bigg|^{\frac{2p}{p-1}} d\eta \Bigg)^{^{\frac{p-1}{2p}}}\\
		& \le C_{\beta,N} \Bigg(\int_{|\eta|\le N}\frac{1}{(1+|\eta|)^{\frac{2p}{p-1}(\beta-4)-4}}\int_{0}^{s-\varepsilon}e^{-\nu_0(s-s')}\bigg|\int_{|\eta'|\le \tilde{C}_{\beta,N}}|h(s',x-v(t-s)-\eta(s-s'),\eta')|\\
		& \qquad  \quad \times d\eta'\bigg|^{\frac{2p}{p-1}} ds' d\eta \Bigg)^{^{\frac{p-1}{2p}}}\\
		& \le C_{\beta,N}\Bigg(\int_{0}^{s-\varepsilon}e^{-\nu_0(s-s')}\int_{|\eta|\le N}\int_{|\eta'|\le \tilde{C}_{\beta,N}}|h(s',x-v(t-s)-\eta(s-s'),\eta')|^{\frac{2p}{p-1}}d\eta' ds' d\eta \Bigg)^{^{\frac{p-1}{2p}}}\\
		& \le C_{\beta,N,\varepsilon} \bar{M}^{\frac{p+1}{2p}} \left[H(F_0 \mid \mu )^{\frac{p-1}{4p}}+H(F_0 \mid \mu )^{\frac{p-1}{2p}} \right],
	\end{align*}
	where we have used \eqref{poRTIRE} in the last term.
	It follow that
	\begin{align} \label{poLA57}
		L_2 \le \frac{C_\beta}{N}\bar{M}^2+\varepsilon C_{\beta,N} \bar{M}^2 + C_{\beta,N,\varepsilon} \bar{M}^{1+\frac{p+1}{2p}} \left[H(F_0 \mid \mu )^{\frac{p-1}{4p}}+H(F_0 \mid \mu )^{\frac{p-1}{2p}} \right].
	\end{align}
	For $L_3$, we first need to estimate $\tilde{L}_3$. We use Lemma \ref{popointwiseQ+est} to deduce
	\begin{align*}
		\tilde{L}_3
		& \le C_\beta\bar{M}\Bigg(\int_{|\eta|\le N}\frac{1}{(1+|\eta|)^{\frac{2p}{p-1}(\beta-4)-4}}\Bigg|\int_{0}^{s}e^{-\nu(\eta)(s-s')}\langle \eta\rangle^\gamma \Bigg(\bigg(\int_{|\eta'|\ge \tilde{N}}+\int_{|\eta'|\le \tilde{N}} \bigg)\\
		& \quad \times \frac{1}{(1+|\eta'|)^{\frac{2p}{p-1}(\beta-4)-4}}|h(s',x-v(t-s)-\eta(s-s'),\eta')|^{\frac{2p}{p-1}}d\eta' \Bigg)^{\frac{p-1}{2p}} ds'\Bigg|^{\frac{2p}{p-1}} d\eta \Bigg)^{^{\frac{p-1}{2p}}}\\
		& \le \frac{C_{\beta,N}}{\tilde{N}} \bar{M}^2\\
		& \quad + C_\beta\bar{M}\Bigg(\int_{|\eta|\le N}\frac{1}{(1+|\eta|)^{\frac{2p}{p-1}(\beta-4)-4}}\Bigg|\left(\int_{s-\varepsilon}^{s}+\int_0^{s-\varepsilon}\right)e^{-\nu(\eta)(s-s')}\langle \eta\rangle^\gamma\\
		& \qquad \times \Bigg(\int_{|\eta'|\le \tilde{N}}\frac{1}{(1+|\eta'|)^{\frac{2p}{p-1}(\beta-4)-4}}|h(s',x-v(t-s)-\eta(s-s'),\eta')|^{\frac{2p}{p-1}}d\eta' \Bigg)^{\frac{p-1}{2p}} ds'\Bigg|^{\frac{2p}{p-1}} d\eta \Bigg)^{^{\frac{p-1}{2p}}}\\
		& \le \frac{C_{\beta,N}}{\tilde{N}} \bar{M}^2+ \varepsilon C_{\beta,N} \bar{M}^2\\
		& \quad + C_{\beta,N}\bar{M}\Bigg(\int_{|\eta|\le N}\frac{1}{(1+|\eta|)^{\frac{2p}{p-1}(\beta-4)-4}}\Bigg|\int_0^{s-\varepsilon}e^{-\nu_0(s-s')}\\
		& \qquad \quad  \times \Bigg(\int_{|\eta'|\le \tilde{N}}\frac{1}{(1+|\eta'|)^{\frac{2p}{p-1}(\beta-4)-4}}|h(s',x-v(t-s)-\eta(s-s'),\eta')|^{\frac{2p}{p-1}}d\eta' \Bigg)^{\frac{p-1}{2p}} ds'\Bigg|^{\frac{2p}{p-1}} d\eta \Bigg)^{^{\frac{p-1}{2p}}}.
	\end{align*}
	From the H\"older inequality, the last term is controlled by
	\begin{align*}
		&C_{\beta,N}\bar{M}\Bigg(\int_{|\eta|\le N}\frac{1}{(1+|\eta|)^{\frac{2p}{p-1}(\beta-4)-4}}\int_0^{s-\varepsilon}e^{-\nu_0(s-s')}\\
		& \qquad \quad  \times \Bigg(\int_{|\eta'|\le \tilde{N}}\frac{1}{(1+|\eta'|)^{\frac{2p}{p-1}(\beta-4)-4}}|h(s',x-v(t-s)-\eta(s-s'),\eta')|^{\frac{2p}{p-1}}d\eta' \Bigg) ds' d\eta \Bigg)^{^{\frac{p-1}{2p}}}\\
		&\le C_{\beta,N}\bar{M}\Bigg(\int_0^{s-\varepsilon}e^{-\nu_0(s-s')}\int_{|\eta|\le N}\int_{|\eta'|\le \tilde{N}}|h(s',x-v(t-s)-\eta(s-s'),\eta')|^{\frac{2p}{p-1}}d\eta'  d\eta ds' \Bigg)^{^{\frac{p-1}{2p}}}\\
		& \le C_{\beta,N,\tilde{N},\varepsilon}\bar{M}^{1+\frac{p+1}{2p}}\left[H(F_0 \mid \mu )^{\frac{p-1}{4p}}+H(F_0 \mid \mu )^{\frac{p-1}{2p}} \right],
	\end{align*}
	where we have used \eqref{poRTIRE}
	and it follows that
	\begin{align} \label{poLA58}
		L_3 \le \frac{C_{\beta,N}}{\tilde{N}} \bar{M}^3+ \varepsilon C_{\beta,N} \bar{M}^3+C_{\beta,N,\tilde{N},\varepsilon}\bar{M}^{2+\frac{p+1}{2p}}\left[H(F_0 \mid \mu )^{\frac{p-1}{4p}}+H(F_0 \mid \mu )^{\frac{p-1}{2p}} \right].
	\end{align}
	Gathering \eqref{poLA54}, \eqref{poLA57}, and \eqref{poLA58}, we obtain
	\begin{align} \label{poI3}
		|I_3| &\le C e^{-\frac{\nu_0}{2}t}\bar{M}\|h_0\|_{L^\infty_{x,v}}+\frac{C_\beta}{N}\bar{M}^2+\varepsilon C_{\beta,N} \bar{M}^2 \nonumber\\
		& \quad + \frac{C_{\beta,N}}{\tilde{N}} \bar{M}^3+ \varepsilon C_{\beta,N} \bar{M}^3+C_{\beta,N,\tilde{N},\varepsilon}\bar{M}^{2+\frac{p+1}{2p}}\left[H(F_0 \mid \mu )^{\frac{p-1}{4p}}+H(F_0 \mid \mu )^{\frac{p-1}{2p}} \right].
	\end{align}
	Inserting \eqref{poI1}, \eqref{poI2}, and \eqref{poI3} into \eqref{poLA10}, we derive
	\begin{align} \label{poLA59}
		|h(t,x,v)| &\le \left(1+C\bar{M}+C_\beta\right)e^{-\frac{\nu_0}{2}t}\|h_0\|_{L^\infty_{x,v}}+\frac{2C_1}{\nu_0\beta^{\frac{\gamma+3}{4}}} \sup_{0\le s \le t}\left[e^{-\frac{\nu_0}{2}(t-s)}\|h(s)\|_{L^\infty_{x,v}}\right]\nonumber\\
		& \quad + \left(\frac{C_\beta}{N}+\varepsilon C_{\beta,N}+\frac{C_{\beta,N}}{\tilde{N}} \right)\left[\bar{M}+\bar{M}^2+\bar{M}^3\right]\nonumber\\
		&\quad +C_{\beta,N,\varepsilon}\left[H(F_0 \mid \mu )^{1/2} + H(F_0 \mid \mu )\right]\nonumber\\
		& \quad +C_{\beta,N,\varepsilon}\bar{M}^{1+\frac{p+1}{2p}}\left[H(F_0 \mid \mu )^{\frac{p-1}{4p}}+H(F_0 \mid \mu )^{\frac{p-1}{2p}}\right]\nonumber\\
		& \quad +  C_{\beta,N,\tilde{N},\varepsilon}\bar{M}^{2+\frac{p+1}{2p}}\left[H(F_0 \mid \mu )^{\frac{p-1}{4p}}+H(F_0 \mid \mu )^{\frac{p-1}{2p}} \right].
	\end{align}
	Taking the norm $\|\cdot\|_{L^\infty_{x,v}}$ and multiplying both sides by $e^{\frac{\nu_0}{2}t}$, we have
	\begin{align*}
		e^{\frac{\nu_0}{2}t} \|h(t)\|_{L^\infty_{x,v}}
		&\le \left(1+C\bar{M}+C_\beta\right)\|h_0\|_{L^\infty_{x,v}}+\frac{2C_1}{\nu_0\beta^{\frac{\gamma+3}{4}}} \sup_{0\le s \le t}\left[e^{\frac{\nu_0}{2}s}\|h(s)\|_{L^\infty_{x,v}}\right]\nonumber\\
		& \quad + e^{\frac{\nu_0}{2}t} \left(\frac{C_\beta}{N}+\varepsilon C_{\beta,N}+\frac{C_{\beta,N}}{\tilde{N}} \right)\left[\bar{M}+\bar{M}^2+\bar{M}^3\right]\nonumber\\
		&\quad +e^{\frac{\nu_0}{2}t} C_{\beta,N,\varepsilon}\left[H(F_0 \mid \mu )^{1/2} + H(F_0 \mid \mu )\right]\nonumber\\
		& \quad +e^{\frac{\nu_0}{2}t} C_{\beta,N,\varepsilon}\bar{M}^{1+\frac{p+1}{2p}}\left[H(F_0 \mid \mu )^{\frac{p-1}{4p}}+H(F_0 \mid \mu )^{\frac{p-1}{2p}}\right]\nonumber\\
		& \quad +  e^{\frac{\nu_0}{2}t}  C_{\beta,N,\tilde{N},\varepsilon}\bar{M}^{2+\frac{p+1}{2p}}\left[H(F_0 \mid \mu )^{\frac{p-1}{4p}}+H(F_0 \mid \mu )^{\frac{p-1}{2p}} \right],
	\end{align*}
	and taking $\sup_{0\le s \le t}$ and choosing sufficiently large $\beta >0$ so that $\frac{2C_1}{\nu_0\beta^{\frac{\gamma+3}{4}}} < 1/2$,
	\begin{align*}
		\sup_{0\le s \le t}\left[e^{\frac{\nu_0}{2}s}\|h(s)\|_{L^\infty_{x,v}}\right] &\le \left(1+C\bar{M}+C_\beta\right)\|h_0\|_{L^\infty_{x,v}}\nonumber\\
		& \quad + e^{\frac{\nu_0}{2}t} \left(\frac{C_\beta}{N}+\varepsilon C_{\beta,N}+\frac{C_{\beta,N}}{\tilde{N}} \right)\left[\bar{M}+\bar{M}^2+\bar{M}^3\right]\nonumber\\
		&\quad +e^{\frac{\nu_0}{2}t} C_{\beta,N,\varepsilon}\left[H(F_0 \mid \mu )^{1/2} + H(F_0 \mid \mu )\right]\nonumber\\
		& \quad +e^{\frac{\nu_0}{2}t} C_{\beta,N,\varepsilon}\bar{M}^{1+\frac{p+1}{2p}}\left[H(F_0 \mid \mu )^{\frac{p-1}{4p}}+H(F_0 \mid \mu )^{\frac{p-1}{2p}}\right]\nonumber\\
		& \quad +  e^{\frac{\nu_0}{2}t}  C_{\beta,N,\tilde{N},\varepsilon}\bar{M}^{2+\frac{p+1}{2p}}\left[H(F_0 \mid \mu )^{\frac{p-1}{4p}}+H(F_0 \mid \mu )^{\frac{p-1}{2p}} \right].
	\end{align*}
	Hence
	\begin{align*}
		\|h(t)\|_{L^\infty_{x,v}} &\le \left(1+C\bar{M}+C_\beta\right)e^{-\frac{\nu_0}{2}t}\|h_0\|_{L^\infty_{x,v}}\\
		&\quad+  \left(\frac{C_\beta}{N}+\varepsilon C_{\beta,N}+\frac{C_{\beta,N}}{\tilde{N}} \right)\left[\bar{M}+\bar{M}^2+\bar{M}^3\right]\nonumber\\
		&\quad +C_{\beta,N,\varepsilon}\left[H(F_0 \mid \mu )^{1/2} + H(F_0 \mid \mu )\right]\nonumber\\
		& \quad +C_{\beta,N,\varepsilon}\bar{M}^{1+\frac{p+1}{2p}}\left[H(F_0 \mid \mu )^{\frac{p-1}{4p}}+H(F_0 \mid \mu )^{\frac{p-1}{2p}}\right]\nonumber\\
		& \quad +  C_{\beta,N,\tilde{N},\varepsilon}\bar{M}^{2+\frac{p+1}{2p}}\left[H(F_0 \mid \mu )^{\frac{p-1}{4p}}+H(F_0 \mid \mu )^{\frac{p-1}{2p}} \right].
	\end{align*}
\end{proof}

\bigskip

\subsection{Proof of Theorem \ref{Mainresultlargeamp-p}}
\begin{proof}[Proof of Theorem \ref{Mainresultlargeamp-p}]
		Recall that we choose $0<\eta_0\leq1$ from Proposition~\ref{poMainresultsmall} with its exponent $0<p<1$.
        Take
		\begin{align} \label{poLA2}
			T_0 := \frac{2}{\nu_0}\log\left(\frac{4C_3(1+\bar M)M_0}{\eta_0^{1+p}}\right).
		\end{align}
		By \eqref{poLALinftyestimate1}, we get
		\begin{align} \label{poLA3}
			\|h(t)\|_{L^\infty_{x,v}} &\le C_3\left(1+ \bar{M} \right)\|h_0\|_{L^\infty_{x,v}}e^{-\frac{\nu_0}{2} t} +\tilde{E}
		\end{align}
		for all $t \ge 0$, where
		\begin{align*}
			\tilde{E}:&= C_3\left(\frac{1}{N}+\varepsilon C_{\beta,N}+\frac{C_{\beta,N}}{\tilde{N}} \right)\left[\bar{M}+\bar{M}^2+\bar{M}^3\right]\nonumber\\
		&\quad +C_{\beta,N,\varepsilon}\left[H(F_0 \mid \mu )^{1/2} + H(F_0 \mid \mu )\right]\nonumber\\
		& \quad +C_{\beta,N,\varepsilon}\bar{M}^{1+\frac{q+1}{2q}}\left[H(F_0 \mid \mu )^{\frac{q-1}{4q}}+H(F_0 \mid \mu )^{\frac{q-1}{2q}}\right]\nonumber\\
		& \quad +C_{\beta,N,\tilde{N},\varepsilon}\bar{M}^{2+\frac{q+1}{2q}}\left[H(F_0 \mid \mu )^{\frac{q-1}{4q}}+H(F_0 \mid \mu )^{\frac{q-1}{2q}} \right].
       \end{align*}
		We first choose $N>0$ sufficiently large, then choose $\varepsilon >0$ small enough and $\tilde{N}>0$ large enough, and assume $H(F_0 \mid \mu ) \le \varepsilon_0(\eta_0, \bar{M})$ with small enough $\varepsilon_0(\eta_0, \bar{M})$ so that
		\begin{align*}
			\tilde{E} \le  \frac{\eta_0^{1+p}}{4},
		\end{align*}
		and it follows that
		\begin{align} \label{poLA6}
			\|h(t)\|_{L^\infty_{x,v}} \le C_3\left(1+ \bar{M} \right)\|h_0\|_{L^\infty_{x,v}}e^{-\frac{\nu_0}{2} t} +  \frac{\eta_0^{1+p}}{4}
		\end{align}
		for all $t \ge 0$.	\\	
		Assume $H(F_0 \mid \mu ) \le \varepsilon_0$. From \eqref{poLA6}, we get
		\begin{align*}
			\|h(T_0)\|_{L^\infty_{x,v}} \le C_3\left(1+ \bar{M} \right)\|h_0\|_{L^\infty_{x,v}}e^{-\frac{\nu_0}{2} T_0} +  \frac{\eta_0^{1+p}}{4}  \le  \frac{\eta_0^{1+p}}{2}.
		\end{align*}
		It remains to show the exponential decay of the Boltzmann solution $f(t)$ in the weighted $L^\infty_{x,v}$ space. Assume $H(F_0 \mid \mu ) \le \varepsilon_0$. By Theorem \ref{poMainresultsmall}, for all $t \ge T_0$,
		\begin{align} \label{poLA9}
			\|h(t)\|_{L^\infty_{x,v}} \le C_\infty (\|h(T_0)\|_{L^\infty_{x,v}}+\eta_0^{1+p}) e^{-\lambda_\infty(t-T_0)} \le 2 C_\infty \eta_0^{1+p} e^{-\lambda_\infty(t-T_0)}.
		\end{align}
		 Since $\lambda_\infty <\frac{\nu_0}{2}$, it follows from \eqref{poLA6} and \eqref{poLA9} that 
		 \begin{align*}
		 	\|h(t)\|_{L^\infty_{x,v}} \le \max\left\{2,8C_\infty\right\}C_3(1+\bar{M}) M_0e^{-\lambda_\infty t}
		 \end{align*}
		 for all $t \ge 0$.
		 Therefore, we complete the proof of this theorem.
\end{proof}

\newpage

\bibliographystyle{plain}
\bibliography{Refkji.bib}

@Book{Cercignani1994,
  author    = {Cercignani, Carlo and Illner, Reinhard and Pulvirenti, Mario},
  publisher = {Springer Science+Business Media New York},
  title     = {{The Mathematical Theory of Dilute Gases}},
  year      = {1994},
  isbn      = {978-1-4419-8524-8},
  series    = {Applied Mathematical Sciences},
  volume    = {106},
  url       = {http://www.ebook.de/de/product/25194071/carlo_cercignani_reinhard_illner_mario_pulvirenti_mathematical_theory_of_dilute_gases.html},
}

@Article{Ukai1982,
  author    = {Ukai, Seiji and Asano, Kiyoshi},
  journal   = {Publ. Res. Inst. Math. Sci.},
  title     = {{On the Cauchy problem of the Boltzmann equation with a soft potential}},
  year      = {1982},
  number    = {2},
  pages     = {477--519},
  volume    = {18},
  doi       = {10.2977/prims/1195183569},
  fjournal  = {Publications of the Research Institute for Mathematical Sciences},
  publisher = {European Mathematical Publishing House},
}

@Article{Guo2003,
  author    = {Guo, Yan},
  journal   = {Arch. Ration. Mech. Anal.},
  title     = {{Classical Solutions to the Boltzmann Equation for Molecules with an Angular Cutoff}},
  year      = {2003},
  month     = sep,
  number    = {4},
  pages     = {305--353},
  volume    = {169},
  doi       = {10.1007/s00205-003-0262-9},
  fjournal  = {Archive for Rational Mechanics and Analysis},
  publisher = {Springer Science and Business Media {LLC}},
}

@Article{Gualdani2017,
  author    = {Gualdani, Maria Pia and Mischler, St{\'{e}}phane and Mouhot, Cl{\'{e}}ment},
  journal   = {M{\'{e}}moires de la Soci{\'{e}}t{\'{e}} math{\'{e}}matique de France},
  title     = {{Factorization of Non-Symmetric Operators and Exponential H-Theorem}},
  year      = {2017},
  pages     = {1--137},
  volume    = {153},
  doi       = {10.24033/msmf.461},
  publisher = {Societe Mathematique de France},
}

@Article{Kim2011,
  author    = {Kim, Chanwoo},
  journal   = {Comm. Math. Phys.},
  title     = {{Formation and Propagation of Discontinuity for Boltzmann Equation in Non-Convex Domains}},
  year      = {2011},
  month     = nov,
  number    = {3},
  pages     = {641--701},
  volume    = {308},
  doi       = {10.1007/s00220-011-1355-1},
  fjournal  = {Communications in Mathematical Physics},
  publisher = {Springer Science and Business Media {LLC}},
}

@article {Guo2010ARMA,
    AUTHOR = {Guo, Yan},
     TITLE = {Decay and continuity of the {B}oltzmann equation in bounded
              domains},
   JOURNAL = {Arch. Ration. Mech. Anal.},
  FJOURNAL = {Archive for Rational Mechanics and Analysis},
    VOLUME = {197},
      YEAR = {2010},
    NUMBER = {3},
     PAGES = {713--809},
      ISSN = {0003-9527,1432-0673},
   MRCLASS = {35Q20 (35B40 35B65 76N15 76P05 82C40)},
  MRNUMBER = {2679358},
MRREVIEWER = {Mark\ Thompson},
       DOI = {10.1007/s00205-009-0285-y},
       URL = {https://doi.org/10.1007/s00205-009-0285-y},
}

@InCollection{Villani2002,
  author    = {Villani, C{\'{e}}dric},
  booktitle = {Handbook of Mathematical Fluid Dynamics},
  publisher = {Elsevier},
  title     = {{A Review of Mathematical Topics in Collisional Kinetic Theory}},
  year      = {2002},
  pages     = {71--305},
  doi       = {10.1016/s1874-5792(02)80004-0},
}

@book {Carleman1957,
    AUTHOR = {Carleman, T.},
     TITLE = {Probl\`emes math\'ematiques dans la th\'eorie cin\'etique des
              gaz},
    SERIES = {Publ. Sci. Inst. Mittag-Leffler},
    VOLUME = {2},
 PUBLISHER = {Almqvist \& Wiksells Boktryckeri AB, Uppsala},
      YEAR = {1957},
     PAGES = {112},
   MRCLASS = {82.00 (45.00)},
  MRNUMBER = {98477},
MRREVIEWER = {D.\ ter Haar},
}

@Article{Duan2019,
  author    = {Duan, Renjun and Huang, Feimin and Wang, Yong and Zhang, Zhu},
  journal   = {Arch. Ration. Mech. Anal.},
  title     = {{Effects of Soft Interaction and Non-isothermal Boundary Upon Long-Time Dynamics of Rarefied Gas}},
  year      = {2019},
  month     = jul,
  number    = {2},
  pages     = {925--1006},
  volume    = {234},
  doi       = {10.1007/s00205-019-01405-5},
  fjournal  = {Archive for Rational Mechanics and Analysis},
  publisher = {Springer Science and Business Media {LLC}},
}

@Article{Ukai1974,
  author    = {Ukai, Seiji},
  journal   = {Proc. Japan Acad.},
  title     = {{On the existence of global solutions of mixed problem for non-linear Boltzmann equation}},
  year      = {1974},
  month     = jan,
  number    = {3},
  volume    = {50},
  doi       = {10.3792/pja/1195519027},
  publisher = {Project Euclid},
}

@Article{Duan2019b,
  author    = {Duan, Renjun and Wang, Yong},
  journal   = {Adv. Math.},
  title     = {{The Boltzmann equation with large-amplitude initial data in bounded domains}},
  year      = {2019},
  month     = feb,
  pages     = {36--109},
  volume    = {343},
  doi       = {10.1016/j.aim.2018.11.007},
  fjournal  = {Advances in Mathematics},
  publisher = {Elsevier {BV}},
}

@article {DL,
    AUTHOR = {DiPerna, R. J. and Lions, P.-L.},
     TITLE = {On the {C}auchy problem for {B}oltzmann equations: global
              existence and weak stability},
   JOURNAL = {Ann. of Math.},
  FJOURNAL = {Annals of Mathematics. Second Series},
    VOLUME = {130},
      YEAR = {1989},
    NUMBER = {2},
     PAGES = {321--366},
      ISSN = {0003-486X,1939-8980},
   MRCLASS = {82A40 (35Q20 45K05 76P05)},
  MRNUMBER = {1014927},
MRREVIEWER = {Seiji\ Ukai},
       DOI = {10.2307/1971423},
       URL = {https://doi.org/10.2307/1971423},
}

@Article{Imbert2021,
  author        = {Imbert, Cyril and Silvestre, Luis},
  journal       = {J. Amer. Math. Soc.},
  title         = {{Global regularity estimates for the Boltzmann equation without cut-off}},
  year          = {2022},
  archiveprefix = {arXiv},
  doi           = {10.48550/arxiv.1909.12729},
  eprint        = {1909.12729},
  fjournal      = {Journal of the American Mathematical Society},
  primaryclass  = {math.AP},
}

@Article{Shizuta1983,
  author    = {Shizuta, Yasushi},
  journal   = {Comm. Pure Appl. Math.},
  title     = {{On the classical solutions of the boltzmann equation}},
  year      = {1983},
  month     = nov,
  number    = {6},
  pages     = {705--754},
  volume    = {36},
  doi       = {10.1002/cpa.3160360602},
  fjournal  = {Communications on Pure and Applied Mathematics},
  publisher = {Wiley},
}

@Article{Cao2022JFA,
  author    = {Cao, Chuqi},
  journal   = {J. Funct. Anal.},
  title     = {{Cutoff Boltzmann equation with polynomial perturbation near Maxwellian}},
  year      = {2022},
  month     = nov,
  number    = {9},
  pages     = {109641},
  volume    = {283},
  doi       = {10.1016/j.jfa.2022.109641},
  fjournal  = {Journal of Functional Analysis},
  publisher = {Elsevier {BV}},
}

@Article{Briant2016,
  author    = {Briant, Marc and Guo, Yan},
  journal   = {J. Differential Equations},
  title     = {{Asymptotic stability of the Boltzmann equation with Maxwell boundary conditions}},
  year      = {2016},
  month     = dec,
  number    = {12},
  pages     = {7000--7079},
  volume    = {261},
  doi       = {10.1016/j.jde.2016.09.014},
  fjournal  = {Journal of Differential Equations},
  publisher = {Elsevier {BV}},
}

@Article{Guo2015,
  author    = {Guo, Y. and Kim, C. and Tonon, D. and Trescases, A.},
  journal   = {Arch. Ration. Mech. Anal.},
  title     = {{BV-Regularity of the Boltzmann Equation in Non-Convex Domains}},
  year      = {2015},
  month     = nov,
  number    = {3},
  pages     = {1045--1093},
  volume    = {220},
  doi       = {10.1007/s00205-015-0948-9},
  fjournal  = {Archive for Rational Mechanics and Analysis},
  publisher = {Springer Science and Business Media {LLC}},
}

@Article{Arkeryd1972,
  author    = {Arkeryd, Leif},
  journal   = {Arch. Ration. Mech. Anal.},
  title     = {{On the Boltzmann equation, Part I: Existence}},
  year      = {1972},
  month     = jan,
  number    = {1},
  pages     = {1--16},
  volume    = {45},
  doi       = {10.1007/bf00253392},
  fjournal  = {Archive for Rational Mechanics and Analysis},
  publisher = {Springer Science and Business Media {LLC}},
}

@Article{Arkeryd1972a,
  author    = {Arkeryd, Leif},
  journal   = {Arch. Ration. Mech. Anal.},
  title     = {{On the Boltzmann equation part II: The full initial value problem}},
  year      = {1972},
  month     = jan,
  number    = {1},
  pages     = {17--34},
  volume    = {45},
  doi       = {10.1007/bf00253393},
  fjournal  = {Archive for Rational Mechanics and Analysis},
  publisher = {Springer Science and Business Media {LLC}},
}

@Article{Briant2017,
  author    = {Briant, Marc},
  journal   = {Kinet. Relat. Models},
  title     = {{Perturbative theory for the Boltzmann equation in bounded domains with different boundary conditions}},
  year      = {2017},
  number    = {2},
  pages     = {329--371},
  volume    = {10},
  doi       = {10.3934/krm.2017014},
  fjournal  = {Kinetic and Related Models},
  publisher = {American Institute of Mathematical Sciences ({AIMS})},
}

@article {He2024,
    AUTHOR = {He, Ling-Bing and Jiang, Jin-Cheng and Kuo, Hung-Wen and
              Liang, Meng-Hao},
     TITLE = {The {$L^p$} estimate for the gain term of the {B}oltzmann
              collision operator and its application},
   JOURNAL = {Arch. Ration. Mech. Anal.},
  FJOURNAL = {Archive for Rational Mechanics and Analysis},
    VOLUME = {248},
      YEAR = {2024},
    NUMBER = {6},
     PAGES = {Paper No. 112, 55},
      ISSN = {0003-9527,1432-0673},
   MRCLASS = {35Q20 (35D30)},
  MRNUMBER = {4827027},
       DOI = {10.1007/s00205-024-02067-8},
       URL = {https://doi.org/10.1007/s00205-024-02067-8},
}

@article {Guo2010Bounded,
    AUTHOR = {Guo, Yan},
     TITLE = {Bounded solutions for the {B}oltzmann equation},
   JOURNAL = {Quart. Appl. Math.},
  FJOURNAL = {Quarterly of Applied Mathematics},
    VOLUME = {68},
      YEAR = {2010},
    NUMBER = {1},
     PAGES = {143--148},
      ISSN = {0033-569X,1552-4485},
   MRCLASS = {35Q20 (35B45 76P05 82C40)},
  MRNUMBER = {2598886},
MRREVIEWER = {Dehua\ Wang},
       DOI = {10.1090/S0033-569X-09-01180-4},
       URL = {https://doi.org/10.1090/S0033-569X-09-01180-4},
}

@article {Ko2023JSP,
    AUTHOR = {Duan, Renjun and Ko, Gyounghun and Lee, Donghyun},
     TITLE = {The {B}oltzmann equation with a class of large-amplitude
              initial data and specular reflection boundary condition},
   JOURNAL = {J. Stat. Phys.},
  FJOURNAL = {Journal of Statistical Physics},
    VOLUME = {190},
      YEAR = {2023},
    NUMBER = {12},
     PAGES = {Paper No. 189, 46},
      ISSN = {0022-4715,1572-9613},
   MRCLASS = {35Q20 (35L04 35R09 76P05)},
  MRNUMBER = {4670646},
       DOI = {10.1007/s10955-023-03203-6},
       URL = {https://doi.org/10.1007/s10955-023-03203-6},
}

@article {Briant2016DCDS,
    AUTHOR = {Briant, Marc},
     TITLE = {Stability of global equilibrium for the multi-species
              {B}oltzmann equation in {$L^\infty$} settings},
   JOURNAL = {Discrete Contin. Dyn. Syst.},
  FJOURNAL = {Discrete and Continuous Dynamical Systems},
    VOLUME = {36},
      YEAR = {2016},
    NUMBER = {12},
     PAGES = {6669--6688},
      ISSN = {1078-0947,1553-5231},
   MRCLASS = {82C40 (35B35 35Q20 76P05)},
  MRNUMBER = {3567815},
       DOI = {10.3934/dcds.2016090},
       URL = {https://doi.org/10.3934/dcds.2016090},
}

@article {Silvestr2016CMP,
    AUTHOR = {Silvestre, Luis},
     TITLE = {A new regularization mechanism for the {B}oltzmann equation
              without cut-off},
   JOURNAL = {Comm. Math. Phys.},
  FJOURNAL = {Communications in Mathematical Physics},
    VOLUME = {348},
      YEAR = {2016},
    NUMBER = {1},
     PAGES = {69--100},
      ISSN = {0010-3616,1432-0916},
   MRCLASS = {45K05 (35B65 35Q20)},
  MRNUMBER = {3551261},
MRREVIEWER = {Cecilia\ Cavaterra},
       DOI = {10.1007/s00220-016-2757-x},
       URL = {https://doi.org/10.1007/s00220-016-2757-x},
}

@article {Park2026ARMA,
    AUTHOR = {Lee, Donghyun and Park, Sungbin and Yun, Seok-Bae},
     TITLE = {Ill-posedness of the {B}oltzmann-{BGK} model in the
              exponential class},
   JOURNAL = {Arch. Ration. Mech. Anal.},
  FJOURNAL = {Archive for Rational Mechanics and Analysis},
    VOLUME = {250},
      YEAR = {2026},
    NUMBER = {4},
     PAGES = {Paper No. 57, 108},
      ISSN = {0003-9527,1432-0673},
   MRCLASS = {35Q20 (35R25 76P05)},
  MRNUMBER = {5087770},
       DOI = {10.1007/s00205-026-02176-6},
       URL = {https://doi.org/10.1007/s00205-026-02176-6},
}

@Article{Desvillettes,
  author    = {Desvillettes, L. and Villani, C.},
  journal   = {Invent. Math.},
  title     = {{On the trend to global equilibrium for spatially inhomogeneous kinetic systems: The Boltzmann equation}},
  year      = {2005},
  month     = feb,
  number    = {2},
  pages     = {245--316},
  volume    = {159},
  doi       = {10.1007/s00222-004-0389-9},
  fjournal  = {Inventiones Mathematicae},
  publisher = {Springer Science and Business Media {LLC}},
}

@article {Carleman1933Acta.Math.,
    AUTHOR = {Carleman, Torsten},
     TITLE = {Sur la th\'eorie de l'\'equation int\'egrodiff\'erentielle de
              {B}oltzmann},
   JOURNAL = {Acta Math.},
  FJOURNAL = {Acta Mathematica},
    VOLUME = {60},
      YEAR = {1933},
    NUMBER = {1},
     PAGES = {91--146},
      ISSN = {0001-5962,1871-2509},
   MRCLASS = {99-04},
  MRNUMBER = {1555365},
       DOI = {10.1007/BF02398270},
       URL = {https://doi.org/10.1007/BF02398270},
}

@article {Mischer1999AIHP,
    AUTHOR = {Mischler, St\'ephane and Wennberg, Bernt},
     TITLE = {On the spatially homogeneous {B}oltzmann equation},
   JOURNAL = {Ann. Inst. H. Poincar\'e{} C Anal. Non Lin\'eaire},
  FJOURNAL = {Annales de l'Institut Henri Poincar\'e{} C. Analyse Non
              Lin\'eaire},
    VOLUME = {16},
      YEAR = {1999},
    NUMBER = {4},
     PAGES = {467--501},
      ISSN = {0294-1449,1873-1430},
   MRCLASS = {35Q99 (45K05 76P05 82C40)},
  MRNUMBER = {1697562},
MRREVIEWER = {C\'edric\ Villani},
       DOI = {10.1016/S0294-1449(99)80025-0},
       URL = {https://doi.org/10.1016/S0294-1449(99)80025-0},
}

@article {An2026Ann.PDE.,
    AUTHOR = {An, Gayoung and Lee, Donghyun},
     TITLE = {Optimal {$C^{\frac{1}{2}}$} regularity of the {B}oltzmann
              equation in non-convex domains},
   JOURNAL = {Ann. PDE},
  FJOURNAL = {Annals of PDE. Journal Dedicated to the Analysis of Problems
              from Physical Sciences},
    VOLUME = {12},
      YEAR = {2026},
    NUMBER = {1},
     PAGES = {Paper No. 16, 80},
      ISSN = {2524-5317,2199-2576},
   MRCLASS = {35Q20 (35A01 35A02 35F30 82C40)},
  MRNUMBER = {5069632},
       DOI = {10.1007/s40818-026-00243-5},
       URL = {https://doi.org/10.1007/s40818-026-00243-5},
}

@article {Lee2024CPAM,
    AUTHOR = {Kim, Chanwoo and Lee, Donghyun},
     TITLE = {H\"older regularity of the {B}oltzmann equation past an
              obstacle},
   JOURNAL = {Comm. Pure Appl. Math.},
  FJOURNAL = {Communications on Pure and Applied Mathematics},
    VOLUME = {77},
      YEAR = {2024},
    NUMBER = {4},
     PAGES = {2331--2386},
      ISSN = {0010-3640,1097-0312},
   MRCLASS = {76P05},
  MRNUMBER = {4705295},
MRREVIEWER = {Lingda\ Xu},
       DOI = {10.1002/cpa.22167},
       URL = {https://doi.org/10.1002/cpa.22167},
}

@article {Guo2017Invent.Math.,
    AUTHOR = {Guo, Yan and Kim, Chanwoo and Tonon, Daniela and Trescases,
              Ariane},
     TITLE = {Regularity of the {B}oltzmann equation in convex domains},
   JOURNAL = {Invent. Math.},
  FJOURNAL = {Inventiones Mathematicae},
    VOLUME = {207},
      YEAR = {2017},
    NUMBER = {1},
     PAGES = {115--290},
      ISSN = {0020-9910,1432-1297},
   MRCLASS = {35Q20 (35B65)},
  MRNUMBER = {3592757},
MRREVIEWER = {Bertrand\ Lods},
       DOI = {10.1007/s00222-016-0670-8},
       URL = {https://doi.org/10.1007/s00222-016-0670-8},
}

@article {Mouhot2006CMP,
    AUTHOR = {Mouhot, Cl\'ement},
     TITLE = {Rate of convergence to equilibrium for the spatially
              homogeneous {B}oltzmann equation with hard potentials},
   JOURNAL = {Comm. Math. Phys.},
  FJOURNAL = {Communications in Mathematical Physics},
    VOLUME = {261},
      YEAR = {2006},
    NUMBER = {3},
     PAGES = {629--672},
      ISSN = {0010-3616,1432-0916},
   MRCLASS = {82C40 (47N55)},
  MRNUMBER = {2197542},
MRREVIEWER = {Xuguang\ Lu},
       DOI = {10.1007/s00220-005-1455-x},
       URL = {https://doi.org/10.1007/s00220-005-1455-x},
}

@article {Illner1984CMP,
    AUTHOR = {Illner, Reinhard and Shinbrot, Marvin},
     TITLE = {The {B}oltzmann equation: global existence for a rare gas in
              an infinite vacuum},
   JOURNAL = {Comm. Math. Phys.},
  FJOURNAL = {Communications in Mathematical Physics},
    VOLUME = {95},
      YEAR = {1984},
    NUMBER = {2},
     PAGES = {217--226},
      ISSN = {0010-3616,1432-0916},
   MRCLASS = {82A40 (45G10 76P05)},
  MRNUMBER = {760333},
MRREVIEWER = {Carlo\ Cercignani},
       URL = {http://projecteuclid.org/euclid.cmp/1103941523},
}

@article {Nishida1976,
    AUTHOR = {Nishida, Takaaki and Imai, Kazuo},
     TITLE = {Global solutions to the initial value problem for the
              nonlinear {B}oltzmann equation},
   JOURNAL = {Publ. Res. Inst. Math. Sci.},
  FJOURNAL = {Kyoto University. Research Institute for Mathematical
              Sciences. Publications},
    VOLUME = {12},
      YEAR = {1976/77},
    NUMBER = {1},
     PAGES = {229--239},
      ISSN = {0034-5318,1663-4926},
   MRCLASS = {82.45},
  MRNUMBER = {432105},
MRREVIEWER = {Nina\ B.\ Maslova},
       DOI = {10.2977/prims/1195190965},
       URL = {https://doi.org/10.2977/prims/1195190965},
}

@article {Ukai2006,
        AUTHOR = {Ukai, Seiji and Yang, Tong},
     TITLE = {The {B}oltzmann equation in the space {$L^2\cap
              L^\infty_\beta$}: global and time-periodic solutions},
   JOURNAL = {Anal. Appl. (Singap.)},
  FJOURNAL = {Analysis and Applications},
    VOLUME = {4},
      YEAR = {2006},
    NUMBER = {3},
     PAGES = {263--310},
      ISSN = {0219-5305,1793-6861},
   MRCLASS = {35F20 (35B10 35B35 76P05 82C40)},
  MRNUMBER = {2239407},
MRREVIEWER = {Mitsuhiro\ Nakao},
       DOI = {10.1142/S0219530506000784},
       URL = {https://doi.org/10.1142/S0219530506000784},
}

@article {Povzner1962,
    AUTHOR = {Povzner, A. Ja.},
     TITLE = {On the {B}oltzmann equation in the kinetic theory of gases},
   JOURNAL = {Mat. Sb. (N.S.)},
  FJOURNAL = {Matematicheski\u i\ Sbornik. Novaya Seriya},
    VOLUME = {58(100)},
      YEAR = {1962},
     PAGES = {65--86},
      ISSN = {0368-8666},
   MRCLASS = {82.45},
  MRNUMBER = {142362},
MRREVIEWER = {E.\ Bareiss},
}

@article {Elmroth1983,
    AUTHOR = {Elmroth, T.},
     TITLE = {Global boundedness of moments of solutions of the {B}oltzmann
              equation for forces of infinite range},
   JOURNAL = {Arch. Rational Mech. Anal.},
  FJOURNAL = {Archive for Rational Mechanics and Analysis},
    VOLUME = {82},
      YEAR = {1983},
    NUMBER = {1},
     PAGES = {1--12},
      ISSN = {0003-9527},
   MRCLASS = {76P05 (45K05 82A40)},
  MRNUMBER = {684411},
MRREVIEWER = {Reinhard\ Illner},
       DOI = {10.1007/BF00251722},
       URL = {https://doi.org/10.1007/BF00251722},
}

@article {Wennberg1994,
    AUTHOR = {Wennberg, Bernt},
     TITLE = {Regularity in the {B}oltzmann equation and the {R}adon
              transform},
   JOURNAL = {Comm. Partial Differential Equations},
  FJOURNAL = {Communications in Partial Differential Equations},
    VOLUME = {19},
      YEAR = {1994},
    NUMBER = {11-12},
     PAGES = {2057--2074},
      ISSN = {0360-5302,1532-4133},
   MRCLASS = {76P05 (35B65 35Q35 82C40)},
  MRNUMBER = {1301182},
MRREVIEWER = {Giuseppe\ Toscani},
       DOI = {10.1080/03605309408821082},
       URL = {https://doi.org/10.1080/03605309408821082},
}

@article {Desvillettes1990,
    AUTHOR = {Desvillettes, L.},
     TITLE = {Convergence to equilibrium in large time for Boltzmann and B.G.K. equations},
   JOURNAL = {Arch. Rational Mech. Anal.},
  FJOURNAL = {Archive for Rational Mechanics and Analysis},
    VOLUME = {110},
      YEAR = {1990},
    NUMBER = {1},
     PAGES = {73--91},
      ISSN = {1432-0673},
       DOI = {10.1007/BF00375163},
       URL = {https://doi.org/10.1007/BF00375163},
}

@article {Desvillettes1993,
    AUTHOR = {Desvillettes, L.},
     TITLE = {Some applications of the method of moments for the homogeneous
              {B}oltzmann and {K}ac equations},
   JOURNAL = {Arch. Rational Mech. Anal.},
  FJOURNAL = {Archive for Rational Mechanics and Analysis},
    VOLUME = {123},
      YEAR = {1993},
    NUMBER = {4},
     PAGES = {387--404},
      ISSN = {0003-9527},
   MRCLASS = {82C40 (35Q99 76P05)},
  MRNUMBER = {1233644},
MRREVIEWER = {Carlo\ Cercignani},
       DOI = {10.1007/BF00375586},
       URL = {https://doi.org/10.1007/BF00375586},
}

@article {Arkeryd1988ARMA,
    AUTHOR = {Arkeryd, Leif},
     TITLE = {Stability in {$L^1$} for the spatially homogeneous {B}oltzmann
              equation},
   JOURNAL = {Arch. Rational Mech. Anal.},
  FJOURNAL = {Archive for Rational Mechanics and Analysis},
    VOLUME = {103},
      YEAR = {1988},
    NUMBER = {2},
     PAGES = {151--167},
      ISSN = {0003-9527},
   MRCLASS = {82A40 (76P05)},
  MRNUMBER = {946973},
MRREVIEWER = {Carlo\ Cercignani},
       DOI = {10.1007/BF00251506},
       URL = {https://doi.org/10.1007/BF00251506},
}

@article {Wennberg1993,
    AUTHOR = {Wennberg, Bernt},
     TITLE = {Stability and exponential convergence in {$L^p$} for the
              spatially homogeneous {B}oltzmann equation},
   JOURNAL = {Nonlinear Anal.},
  FJOURNAL = {Nonlinear Analysis. Theory, Methods \& Applications. An
              International Multidisciplinary Journal},
    VOLUME = {20},
      YEAR = {1993},
    NUMBER = {8},
     PAGES = {935--964},
      ISSN = {0362-546X,1873-5215},
   MRCLASS = {76P05 (45M10 82C40)},
  MRNUMBER = {1215007},
MRREVIEWER = {Jian\ Ning\ Chen},
       DOI = {10.1016/0362-546X(93)90086-8},
       URL = {https://doi.org/10.1016/0362-546X(93)90086-8},
}

@article {GPV,
    AUTHOR = {Gamba, I. M. and Panferov, V. and Villani, C.},
     TITLE = {Upper {M}axwellian bounds for the spatially homogeneous
              {B}oltzmann equation},
   JOURNAL = {Arch. Ration. Mech. Anal.},
  FJOURNAL = {Archive for Rational Mechanics and Analysis},
    VOLUME = {194},
      YEAR = {2009},
    NUMBER = {1},
     PAGES = {253--282},
      ISSN = {0003-9527,1432-0673},
   MRCLASS = {35Q20 (76P05 82C40)},
  MRNUMBER = {2533928},
MRREVIEWER = {Matthias\ Geissert},
       DOI = {10.1007/s00205-009-0250-9},
       URL = {https://doi.org/10.1007/s00205-009-0250-9},
}

@article {Bobylev,
    AUTHOR = {Bobylev, A. V.},
     TITLE = {Moment inequalities for the {B}oltzmann equation and
              applications to spatially homogeneous problems},
   JOURNAL = {J. Statist. Phys.},
  FJOURNAL = {Journal of Statistical Physics},
    VOLUME = {88},
      YEAR = {1997},
    NUMBER = {5-6},
     PAGES = {1183--1214},
      ISSN = {0022-4715,1572-9613},
   MRCLASS = {82C40 (76P05)},
  MRNUMBER = {1478067},
MRREVIEWER = {Carlo\ Cercignani},
       DOI = {10.1007/BF02732431},
       URL = {https://doi.org/10.1007/BF02732431},
}

@inproceedings {Lu,
    AUTHOR = {Lu, Xuguang},
     TITLE = {A result on uniqueness of mild solutions of {B}oltzmann
              equations},
 BOOKTITLE = {Proceedings of the {F}ourteenth {I}nternational {C}onference
              on {T}ransport {T}heory ({B}eijing, 1995)},
   JOURNAL = {Transport Theory Statist. Phys.},
  FJOURNAL = {Transport Theory and Statistical Physics},
    VOLUME = {26},
      YEAR = {1997},
    NUMBER = {1-2},
     PAGES = {209--220},
      ISSN = {0041-1450,1532-2424},
   MRCLASS = {82C40 (45K05)},
  MRNUMBER = {1457582},
MRREVIEWER = {Carlo\ Cercignani},
       DOI = {10.1080/00411459708221784},
       URL = {https://doi.org/10.1080/00411459708221784},
}

@article {CPDE,
    AUTHOR = {Mouhot, Cl\'ement},
     TITLE = {Quantitative lower bounds for the full {B}oltzmann equation.
              {I}. {P}eriodic boundary conditions},
   JOURNAL = {Comm. Partial Differential Equations},
  FJOURNAL = {Communications in Partial Differential Equations},
    VOLUME = {30},
      YEAR = {2005},
    NUMBER = {4-6},
     PAGES = {881--917},
      ISSN = {0360-5302,1532-4133},
   MRCLASS = {76P05 (35Q72 82B40 82C40 82D05)},
  MRNUMBER = {2153518},
MRREVIEWER = {Laurent\ Desvillettes},
       DOI = {10.1081/PDE-200059299},
       URL = {https://doi.org/10.1081/PDE-200059299},
}

@article {Golse02,
    AUTHOR = {Golse, Fran\c cois and Levermore, C. David},
     TITLE = {Stokes-{F}ourier and acoustic limits for the {B}oltzmann
              equation: convergence proofs},
   JOURNAL = {Comm. Pure Appl. Math.},
  FJOURNAL = {Communications on Pure and Applied Mathematics},
    VOLUME = {55},
      YEAR = {2002},
    NUMBER = {3},
     PAGES = {336--393},
      ISSN = {0010-3640,1097-0312},
   MRCLASS = {35F25 (76P05 82C40)},
  MRNUMBER = {1866367},
MRREVIEWER = {Giuseppe\ Toscani},
       DOI = {10.1002/cpa.3011},
       URL = {https://doi.org/10.1002/cpa.3011},
}

@article {PulvirentiWennberg1997,
    AUTHOR = {Pulvirenti, Ada and Wennberg, Bernt},
     TITLE = {A {M}axwellian lower bound for solutions to the {B}oltzmann
              equation},
   JOURNAL = {Comm. Math. Phys.},
  FJOURNAL = {Communications in Mathematical Physics},
    VOLUME = {183},
      YEAR = {1997},
    NUMBER = {1},
     PAGES = {145--160},
      ISSN = {0010-3616,1432-0916},
   MRCLASS = {82C40 (45K05 76P05)},
  MRNUMBER = {1461954},
MRREVIEWER = {Carlo\ Cercignani},
       DOI = {10.1007/BF02509799},
       URL = {https://doi.org/10.1007/BF02509799},
}

@article {cercig,
    AUTHOR = {Villani, C\'edric},
     TITLE = {Cercignani's conjecture is sometimes true and always almost
              true},
   JOURNAL = {Comm. Math. Phys.},
  FJOURNAL = {Communications in Mathematical Physics},
    VOLUME = {234},
      YEAR = {2003},
    NUMBER = {3},
     PAGES = {455--490},
      ISSN = {0010-3616,1432-0916},
   MRCLASS = {82C40 (76P05)},
  MRNUMBER = {1964379},
MRREVIEWER = {Laurent\ Desvillettes},
       DOI = {10.1007/s00220-002-0777-1},
       URL = {https://doi.org/10.1007/s00220-002-0777-1},
}

@article {DL2,
    AUTHOR = {DiPerna, R. J. and Lions, P.-L.},
     TITLE = {Global solutions of {B}oltzmann's equation and the entropy
              inequality},
   JOURNAL = {Arch. Rational Mech. Anal.},
  FJOURNAL = {Archive for Rational Mechanics and Analysis},
    VOLUME = {114},
      YEAR = {1991},
    NUMBER = {1},
     PAGES = {47--55},
      ISSN = {0003-9527},
   MRCLASS = {76P05 (82C40)},
  MRNUMBER = {1088276},
MRREVIEWER = {Reinhard\ Illner},
       DOI = {10.1007/BF00375684},
       URL = {https://doi.org/10.1007/BF00375684},
}

@article {LuMouhot2012,
    AUTHOR = {Lu, Xuguang and Mouhot, Cl\'ement},
     TITLE = {On measure solutions of the {B}oltzmann equation, part {I}:
              moment production and stability estimates},
   JOURNAL = {J. Differential Equations},
  FJOURNAL = {Journal of Differential Equations},
    VOLUME = {252},
      YEAR = {2012},
    NUMBER = {4},
     PAGES = {3305--3363},
      ISSN = {0022-0396,1090-2732},
   MRCLASS = {35Q20 (35A01 35B35 35R06 82C40)},
  MRNUMBER = {2871802},
MRREVIEWER = {Renjun\ Duan},
       DOI = {10.1016/j.jde.2011.10.021},
       URL = {https://doi.org/10.1016/j.jde.2011.10.021},
}

@Article{DH-entropy,
  author        = {Lee, Donghyun},
  journal       = {Preprint},
  title         = {{Constructive Trend to Equilibrium for Low-Regularity Boltzmann  and BGK equations}},
  year          = {},
  archiveprefix = {arXiv},
  doi           = {},
  eprint        = {},
  primaryclass  = {math.AP},
}

@article {lion1,
    AUTHOR = {Lions, P.-L.},
     TITLE = {Compactness in {B}oltzmann's equation via {F}ourier integral
              operators and applications. {I}, {II}},
   JOURNAL = {J. Math. Kyoto Univ.},
  FJOURNAL = {Journal of Mathematics of Kyoto University},
    VOLUME = {34},
      YEAR = {1994},
    NUMBER = {2},
     PAGES = {391--427, 429--461},
      ISSN = {0023-608X},
   MRCLASS = {35Q35 (35Q99 35S30 76P05 82C40)},
  MRNUMBER = {1284432},
MRREVIEWER = {Andrzej\ Palczewski},
       DOI = {10.1215/kjm/1250519017},
       URL = {https://doi.org/10.1215/kjm/1250519017},
}

@article{He2023,
  title = {On the Cauchy Problem for the Cutoff Boltzmann Equation with Small Initial Data},
  volume = {190},
  issn = {1572-9613},
  url = {http://dx.doi.org/10.1007/s10955-023-03065-y},
  doi = {10.1007/s10955-023-03065-y},
  number = {3},
  journal = {Journal of Statistical Physics},
  publisher = {Springer Science and Business Media LLC},
  author = {He, Ling-Bing and Jiang, Jin-Cheng},
  year = {2023},
  month = Jan,
}

@article{Liu2004b,
  title = {The Green’s function and large‐time behavior of solutions for the one‐dimensional Boltzmann equation},
  volume = {57},
  issn = {1097-0312},
  url = {http://dx.doi.org/10.1002/cpa.20011},
  doi = {10.1002/cpa.20011},
  number = {12},
  journal = {Communications on Pure and Applied Mathematics},
  publisher = {Wiley},
  author = {Liu, Tai‐Ping and Yu, Shih‐Hsien},
  year = {2004},
  month = Feb,
  pages = {1543–1608},
}

@article{MouhotVillani2004,
  author  = {Mouhot, Clément and Villani, Cédric},
  title   = {Regularity theory for the spatially
             homogeneous Boltzmann equation
             with cut-off},
  journal = {Archive for Rational Mechanics
             and Analysis},
  volume  = {173},
  number  = {2},
  pages   = {169--212},
  year    = {2004},
  doi     = {10.1007/s00205-004-0316-7}
}

@article {GLPS,
    AUTHOR = {Golse, Fran\c cois and Lions, Pierre-Louis and Perthame,
              Beno\^it and Sentis, R\'emi},
     TITLE = {Regularity of the moments of the solution of a transport
              equation},
   JOURNAL = {J. Funct. Anal.},
  FJOURNAL = {Journal of Functional Analysis},
    VOLUME = {76},
      YEAR = {1988},
    NUMBER = {1},
     PAGES = {110--125},
      ISSN = {0022-1236},
   MRCLASS = {35Q20 (35B99)},
  MRNUMBER = {923047},
       DOI = {10.1016/0022-1236(88)90051-1},
       URL = {https://doi.org/10.1016/0022-1236(88)90051-1},
}

\end{document}